\documentclass[11pt,letterpaper]{amsart}
\usepackage[margin=3cm]{geometry}
\usepackage{setspace}
\usepackage{amsmath}
\usepackage{mathtools}
\usepackage{amssymb}
\usepackage{amsfonts}
\usepackage{amsthm}
\usepackage{amscd}
\usepackage{todonotes}
\usepackage{bbm}
\usepackage{stmaryrd}
\usepackage{ulem}
\usepackage{mathrsfs}
\usepackage{mathdots}
\usepackage{float}
\usepackage[all]{xy}
\usepackage{wrapfig}
\usepackage{enumitem}

\usepackage{fancyhdr}

\usepackage[pagebackref,colorlinks=true,linkcolor=cyan,citecolor=magenta]{hyperref}

\usepackage[alphabetic,nobysame, initials, backrefs]{amsrefs}
\usepackage{tikz-cd}
\usepackage{tikz}

\tikzset{
  symbol/.style={
    draw=none,
    every to/.append style={
      edge node={node [sloped, allow upside down, auto=false]{$#1$}}}
  }
}

\newcommand{\calm}{\textrm{m}}

\newtheorem*{theorem*}{Theorem}

\newtheorem*{corollary*}
{Corollary}

\usepackage{imakeidx}
\makeindex[columns=2]

\newcommand{\IN}[1]{\index{#1|BH}}

\newcommand{\IC}{\mathcal{I\hskip-2ptC}}
\newcommand{\Loc}{\operatorname{Loc}}
\newcommand{\ceq}{\, := \, }
\newcommand{\tq}{\ \vert\ }

\newcommand{\Z}{\mathbb{Z}}
\newcommand{\Q}{\mathbb{Q}}
\newcommand{\R}{\mathbb{R}}
\newcommand{\BC}{\mathbb{C}}

\newcommand{\iso}{\simeq}

\newcommand{\ol}{\overline}

\newcommand{\SL}{\mathrm{SL}}
\newcommand{\GL}{\mathrm{GL}}
\newcommand{\SO}{\mathrm{SO}}
\newcommand{\Sp}{\mathrm{Sp}}
\newcommand{\St}{\mathrm{St}}
\newcommand{\dualgroup}[1]{\widehat{#1}}
\newcommand{\Lgroup}[1]{{\hskip-2 pt \,^L\hskip-1pt{#1}}}

\newcommand{\reg}{\mathrm{reg}}

\newcommand{\sreg}{\mathrm{sreg}}
\newcommand{\gen}{\mathrm{gen}}
\newcommand{\NEvs}{\mathrm{NEvs}}
\newcommand{\Evs}{\mathrm{Evs}}
\newcommand{\Ev}{\mathrm{Ev}}

\newcommand{\ev}[1]{\mathrm{Ev}_{\hskip-1pt#1}}

\newcommand{\evs}[1]{\mathrm{Evs}_{\hskip-1pt#1}}

\newcommand{\nevs}[1]{\mathrm{NEvs}_{\hskip-1pt#1}}

\newcommand{\nevsmic}[1]{\mathrm{NEvs}^{\hskip-1pt\mic}_{\hskip-1pt#1}}

\newcommand{\geo}{\mathrm{geo}}
\newcommand{\pure}{\mathrm{pure}}
\newcommand{\st}{\mathrm{st}}
\newcommand{\std}{\mathrm{std}}
\newcommand{\pr}{\mathrm{pr}}
\newcommand{\QST}{\mathrm{QST}}

\newcommand{\AZ}{\mathsf{AZ}}

\newcommand{\CTD}{\textbf{(CTD)}}

\newcommand{\CCh}{\textbf{(CC)}}

\newcommand{\TI}{\textbf{(TI)}}

\newcommand{\TW}{\textbf{(TW)}}

\newcommand{\Sta}{\textbf{(St)}}

\newcommand{\Funtheta}{\textbf{(Fun$_{\theta}$)}}
\newcommand{\LC}{\textbf{(LC)}}
\newcommand{\Sha}{\textbf{(Sha)}}

\newcommand{\ECR}{\textbf{(ECR)}}

\newcommand{\Hai}{\textbf{(Hai)}}
\newcommand{\RES}{\textbf{(Res)}}

\newcommand{\mfr}[1]{\mathfrak{#1}}

\newcommand{\RG}{\mathrm{G}}

\newcommand{\Ind}{\mathrm{Ind}}
\newcommand{\Res}{\mathrm{Res}}
\newcommand{\LLC}{\mathrm{LLC}}
\newcommand{\Per}{\mathrm{Per}}

\newcommand{\mic}{\mathrm{mic}}
\newcommand{\codim}{\mathrm{codim}~}

\newcommand{\Trans}{\mathrm{Trans}}

\newcommand{\KS}{\mathrm{KS}}

\newcommand{\ABV}{{\mbox{\raisebox{0pt}{\scalebox{0.5}{$\mathrm{ABV}$}}}}}

\newcommand{\Art}{{\mbox{\raisebox{1pt}{\scalebox{0.6}{$\mathrm{Art}$}}}}}

\newcommand{\etaABV}[2]{{\eta^{\ABV #2}_{#1}}}

\newcommand{\etaArt}[2]{{\eta^{\Art #2}_{#1}}}

\newcommand{\locsys}[2]{\underline{#1}_{#2}}

\newcommand{\stsh}[2]{j_!\locsys{#1}{#2}}

 \newcommand{\Irr}[1]{\mathrm{Irr}({#1})}

\newcommand{\STD}[2]{{\mu_{#1}^{#2}}}
\newcommand{\Rphi}[1]{\mathrm{R}\Phi_{#1}}

\newcommand{\SR}{{M}}
\newcommand{\sr}{{M}}

\newcommand{\Transfer}[2]{{\Trans^{#2}_{#1}}}

\newcommand{\Lift}[2]{{\operatorname{Lift}^{#2}_{#1}}}

\newcommand{\Rep}{\mathrm{Rep}}
\newcommand{\Irrep}{\mathrm{Irrep}}

\newcommand{\half}[1]{\frac{#1}{2}}

\renewcommand{\implies}{\Rightarrow}

\newcommand{\comment}[1]{}

\newcommand{\CC}{\mathbb{C}}

\renewcommand{\O}
{\mathrm{O}}

\newtheorem{theorem}{Theorem}[section]
\newtheorem{corollary}[theorem]{Corollary}
\newtheorem{lemma}[theorem]{Lemma}
\newtheorem{desideratum}[theorem]{Desideratum}
\newtheorem{proposition}[theorem]{Proposition}
\newtheorem{conjecture}[theorem]{Conjecture}
\newtheorem{assu}[theorem]{Assumption}
\newtheorem{ques/conj}[theorem]{Question/Conjecture}

\newtheorem{definition}[theorem]{Definition}
\newtheorem{remark}[theorem]{Remark}

\newtheorem{example}[theorem]{Example}

\newtheorem{assumption}[theorem]{Working Hypotheses}
\newtheorem{hypothesis}[theorem]{Hypothesis}

\newtheorem*{globalcond*}{Global Condition}
\newtheorem*{localcond*}{Local Condition}

\newtheorem*{globalconj*}{Global Conjecture}
\newtheorem*{localconj*}{Local Conjecture}

\newtheorem*{nonzero*}{Conjecture on the non-vanishing of the normalized intertwining operators}
\newtheorem*{holo*}{Conjecture on the holomorphicity of the normalized intertwining operators}

\DeclareMathOperator{\Ad}{Ad}
\DeclareMathOperator{\supp}{supp}

\DeclareMathOperator{\Aut}{Aut}
\DeclareMathOperator{\Hom}{Hom}
\DeclareMathOperator{\Gal}{Gal}

\DeclareMathOperator{\rank}{rank}
\newcommand{\trace}{\operatorname{trace}}

\newcommand{\U}{{\operatorname{U}}}
\renewcommand{\Res}{{\operatorname{Res}}}

\let\oldbullet\bullet
\renewcommand{\bullet}{{\vcenter{\hbox{\tiny$\oldbullet$}}}}

\renewcommand{\1}{\mathbbm{1}}

\numberwithin{equation}{section}

\makeatletter
\let\@wraptoccontribs\wraptoccontribs
\makeatother

\begin{document}

\title[Geometrization of Kottwitz-Shelstad Transfer]{Geometrization of Kottwitz-Shelstad transfer and Vogan's conjecture on local Arthur packets}

\author[Cunningham]{Clifton Cunningham}
\address{Department of Mathematics and Statistics, University of Calgary, Alberta, T2N 1N4,
CANADA}
\email{ccunning@ucalgary.ca}

\author[Hazeltine]{Alexander Hazeltine}
\address{Department of Mathematics\\
Columbia University\\
New York, NY, 10027, USA}
\email{alh2252@columbia.edu}

\author[Liu]{Baiying Liu}
\address{Department of Mathematics\\
Purdue University\\
West Lafayette, IN, 47907, USA}
\email{liu2053@purdue.edu}

\author[Lo]{Chi-Heng Lo}
\address{Department of Mathematics,
National University of Singapore,
119076, Singapore}
\email{ch\_lo@nus.edu.sg}

\author[Ray]{Mishty Ray}
\address{Department of Mathematics and Statistics\\ University of British Columbia\\
Vancouver,
CANADA}
\email{mishtyray@math.ubc.ca}

\author[Xu]{Bin Xu}
\address{Yau Mathematical Sciences Center, Jing Zhai 207, Tsinghua University, 
Haidian District, Beijing 100084}
\email{binxu@tsinghua.edu.cn}

\contrib[with appendices by]{Jörg Schürmann}
\address{Mathematisches Institut\\
University of Münster\\
Einsteinstr. 62, 48149 Münster,
GERMANY}
\email{jschuerm@uni-muenster.de}


\date{\today}


\thanks{C. Cunningham is partially supported by NSERC Discovery Grants RGPIN-2020-05220 and RGPIN-2026-07412. A. Hazeltine is partially supported by an AMS-Simons Travel Grant. B. Liu is partially supported by the NSF Grant DMS-1848058 and the Simons Foundation: Travel Support for Mathematicians. B. Xu is partially supported by the NSFC grant 12371010.
J. Schürmann is funded by the Deutsche Forschungsgemeinschaft
(DFG, German Research Foundation) Project-ID 427320536 - SFB 1442, as well as under
Germany’s Excellence Strategy EXC 2044/2 - 390685587, Mathematics Münster: Dynamics
- Geometry - Structure.}

\begin{abstract}
We give a geometric analogue of Kottwitz-Shelstad transfer for a general class of reductive groups over $p$-adic fields with cyclic component groups of trivial or prime order. We prove that this geometric transfer coincides with Kottwitz-Shelstad transfer, adapted to the language of disconnected groups, assuming the local Langlands correspondence satisfying certain desiderata, including the Kazhdan-Lusztig hypothesis for the ambient group and its endoscopic groups. As an application of this geometrization, we prove Vogan's conjecture on local Arthur packets over {\it p}-adic fields: for quasisplit classical groups and their pure inner forms, all local Arthur packets are ABV-packets, and the distributions attached to these ABV-packets of Arthur type satisfy Arthur's local conjectures on twisted endoscopic character relations. This demonstrates the following general principle: {\it for any reductive group $G$ over a local field, an Arthur packet should be defined as an ABV-packet of Arthur type, equipped with the endoscopic pairing calculated by traces of actions on vanishing cycles or microlocalizations}. 
\end{abstract}

\dedicatory{Dedicated to Jim Arthur, with admiration, on the occasion of his 82nd birthday} 

\maketitle
                     
\setcounter{section}{0}


\section{Introduction}\label{section: introduction}

The Langlands functoriality and the Langlands correspondence are two central problems in the Langlands program. In the seminal work \cite{Art13}, Arthur proved the endoscopic classification of automorphic representations for quasisplit classical groups by comparing trace formulas across different groups. For this, he developed the theory of Arthur packets, which parametrize the local components of discrete automorphic representations and established a local Langlands correspondence. 

In this paper, we introduce a geometric framework for establishing the (local) theory of Langlands functoriality in the case of Kottwitz-Shelstad transfer, assuming the existence of the local Langlands correspondence; see Theorem \ref{theorem: lift = transfer}. 
Our approach is grounded in the theory of ABV-packets, as developed in  
\cites{ABV92,Vog93} and \cite{CFMMX22} and further in 
\cites{Mracek,CFZ2,CFK22,Lo24,CR1,CDFZ25,CR2}. 
These packets are constructed from the arithmetic data of local 
Langlands parameters via geometric methods.

Arthur's conjectures on the fundamental properties of his eponymous packets are stated for disconnected reductive groups in \cite{Art89}*{Conjectures 6.1, 6.2}.
These packets are equipped with stable distributions $\etaArt{\psi}{}$ that satisfy six properties, where $\psi$ is a local Arthur parameter for a reductive group $G$ over a local field of characteristic 0. These properties largely refer to the transfer of $\etaArt{\psi'}{}$ from endoscopic groups $G'$ of $G$, for semisimple $s\in Z_{\dualgroup{G}}(\psi)$.  We let $\Pi_\psi(G)$ denote the local Arthur packet attached to $\psi$.
As discussed in \cite{Art89}*{Remark 2}, these six properties do not uniquely characterize the stable distributions, so  \cite{Art89}*{Conjectures 6.1, 6.2} cannot be taken as a definition of the stable distribution $\etaArt{\psi}{}$ or the packet $\Pi_{\psi}(G)$.

For quasisplit classical groups
 (unitary, odd special orthogonal, symplectic, and even orthogonal groups), 
this ambiguity is addressed in \cites{Art13,Mok15}
using twisted endoscopic transfer to general linear groups.
When combined with other works including \cites{AGIKMS24, Hal26}, the packets introduced in \cites{Art13,Mok15} are known to satisfy all the properties appearing in \cite{Art89}*{Conjectures 6.1, 6.2}. 
Since then, these results have been extended to pure inner forms of these groups by the combined efforts of 
\cites{Art13,KMSW14,Mok15,AG17,MR18,Ish24, Hal26}; currently, these are the only groups for which the properties in \cite{Art89}*{Conjectures 6.1} have been fully established. We remark that the tempered case of Arthur's conjecture has been established for the exceptional group $G_2$ in \cite{Tak26}. 

{As an application of our aforementioned geometrization of Kottwitz-Shelstad transfer, we show that distributions associated with ABV-packets satisfy the analogous twisted endoscopic character relations that the distributions associated with local Arthur packets are predicted to satisfy in \cite{Art89}*{Conjectures 6.1} (Theorem \ref{theorem: Vogan's conjecture}). We apply this result to show that local Arthur packets are ABV-packets for classical groups (Theorem \ref{theorem intro: Vogan's conjecture}). We provide more details on these results, along with highlighting other related results, in the rest of this section.}

\subsection{Vogan's conjecture on Arthur packets}

In order to prove that some candidate packets are Arthur packets, one must associate stable distributions with the putative Arthur packets and then show that they satisfy \cite{Art89}*{Conjectures 6.1, 6.2}. 
The foundational work \cite{Vog93} defines a stable distribution $\etaABV{\phi}{}$ for each $L$-parameter $\phi$ and defines the ABV-packet $\Pi^\ABV_\phi(G)$ as the support of $\etaABV{\phi}{}$. Vogan
implicitly conjectures that if $\phi$ is of Arthur type, then $\etaABV{\phi}{}$ satisfies specific properties with respect to both ordinary and twisted endoscopy, as in \cite{Art89}*{Conjectures 6.1, 6.2}. 
We refer to this as {\it Vogan's conjecture on Arthur packets}\IN{Vogan's conjecture on local Arthur packets}, revisited in \cite{CFMMX22}*{Conjecture 1, \S 8.3}. We emphasize that this conjecture is more than a claim about packets; it is a conjecture about endoscopic transfer, both ordinary and twisted.

In this paper we prove Vogan's conjecture on Arthur packets, assuming the twisted Kazhdan-Lusztig hypothesis for $\GL_N$. 
We note that the twisted Kazhdan-Lusztig hypothesis for $\GL_N$ has been proved for unipotent representations in \cite{CHAI} and the general case will follow from a careful reduction to the unipotent case. 

\begin{theorem}[Vogan's conjecture on local Arthur packets, Theorem~\ref{theorem: Vogan's conjecture}] 
\label{theorem intro: Vogan's conjecture}  
Assume the Kazhdan-Lusztig Hypothesis~\ref{hypothesis: KLH}, for $\GL_N\rtimes{\langle\theta\rangle}$ for outer involution $\theta$ of $\GL_N$ (or $\mathrm{Res}_{E/F}\GL_n\rtimes\langle\theta\rangle$ in the case of unitary groups).
    Let $G$ be any one of the quasisplit groups $\U^E_n$, $\SO_{2n+1}$, $\Sp_{2n}$, or $\O^\alpha_{2n}$ over a $p$-adic field $F$ (see Section \ref{subsection: classical groups} for notation) and let $\delta$ be a pure inner form for $G$.
Let $\psi$ be a local Arthur parameter for $G$ and let $\phi_\psi$ be the Langlands parameter associated to it (see \eqref{eqn: psi to phi psi}).
Then
\[
\Pi_\psi(G_\delta) = \Pi^\ABV_{\phi_\psi}(G_\delta).
\]
Moreover, for all semisimple $s\in Z_{\dualgroup{G}}(\psi)$, the distributions match
\[ \etaArt{\psi,s}{,\delta}= \etaABV{\phi_{\psi},s}{,\delta}.\]
Here the distribution $\etaArt{\psi,s}{,\delta}$ is defined in \cite{Art13} using endoscopic transfer; see Conjecture \ref{conj: Arthur}(2). The distribution $ \etaABV{\phi_{\psi},s}{,\delta}$ is defined using vanishing cycles; in see Section \ref{section: etaABV}. 
See Theorem~\ref{theorem: Vogan's conjecture} for full details. 
\end{theorem}

Note that we have included the disconnected group $\O^\alpha_{2n}$ in the list of classical groups, not $\SO^\alpha_{2n}$. Although the term ``classical group" is not precisely defined in the literature, our rationale is specific: the groups we have listed share the property that there is a lift from $L$-parameters for that group to $L$-parameters for some general linear group that is injective; this property is discussed in Part~\ref{part: IV Lift}.
The failure of this property is largely responsible for the fact that Arthur's work for $\SO^\alpha_{2n}$ treats not irreducible representations of $\SO^\alpha_{2n}$ but orbits of such representations under an outer automorphism $\theta$ with the property that $\O^{\alpha}_{2n} = \SO^\alpha_{2n}\rtimes\langle\theta\rangle$.
Twisted endoscopy from $\SO_{2n}$ to $\GL_{2n}$ is not injective, so it cannot be used to directly define stable distributions for $\SO_{2n}$ attached to local Arthur parameters. 
Instead of working with orbits of representations of $\SO^\alpha_{2n}(F)$, we elect to work with irreducible representations of $\O^\alpha_{2n}(F)$ and pay the price by grappling with disconnected groups. 

Furthermore, from our perspective, the use of disconnected reductive groups is not a bug but a feature. Our geometrization of Kottwitz-Shelstad transfer requires the use of disconnected groups when viewed from the Galois/geometric side of the Langlands correspondence.
Indeed, working with disconnected groups is a leitmotif in this paper, and it provides a convenient framework for Arthur's work. 
For example, what Arthur refers to as the ``third kind of endoscopic datum" for $\SO^\alpha_n$ \cite{Art13}*{p.~14} appears naturally as an endoscopic triple $(G',s,\xi)$ for the disconnected group $\O^\alpha_n$ (see Definition \ref{definition: endoscopic triple}) for a semisimple $s\in \O_n(\CC)$ not lying in the identity component.

\subsection{ABV-packets determine a basis for stable distributions}

To prove Vogan’s conjecture on Arthur packets, we first extend the techniques of \cites{CFMMX22,Vog93} to disconnected groups. Assuming that the component group of $G$ is trivial or cyclic of prime order, and that the local Langlands correspondence holds for $G$, we show that ABV-packets for reductive groups $G$ over $p$-adic fields $F$ are unambiguously defined as the supports of stable distributions $\etaABV{\phi}{}$. Moreover, we prove that these distributions span the space of stable distributions supported on characters of representations, under certain assumptions on the local Langlands correspondence. This establishes \cite{Vog93}*{Conjecture 8.15'}, revisited in \cite{CFMMX22}*{Conjecture 2, Section 8.4}.

\begin{theorem}[Basis for stable distributions, Theorem~\ref{theorem: Conjecture 2 qs}]\label{theorem: ABV distributions form a basis intro}
   Suppose $G$ is a connected reductive group, quasisplit over $F$.
   Assume the local Langlands correspondence for $G$, satisfying desiderata that are specified in Theorem~\ref{theorem: Conjecture 2 qs}. 
   Then the distributions $\etaABV{\phi}{}$ are stable and form a basis for stable distributions in the span of characters of admissible representations with a fixed infinitesimal parameter.
\end{theorem}

\subsection{The weak form of Vogan's conjecture and atomic stability}

We begin by giving a conditional proof of a weak form of Vogan's conjecture that matches packets but not distributions, assuming the following. 

\begin{conjecture}
[Aubert-Zelevinsky duality preserves ABV-packets, Conjecture~\ref{conj AZ}, \cite{CFMMX22}*{Section 10.3.4}]
\label{conjecture intro: AZ}
Let $G$ be a reductive group. 
For any L-parameter $\phi$ of $G$, the Aubert-Zelevinsky dual
$\AZ (\etaABV{\phi}{})$ of $\etaABV{\phi}{}$ is proportional to $\etaABV{\widehat{\phi}}{}$, where ${\hat\phi}$ is the L-parameter obtained from $\phi$ by Pyasetskii duality.
In particular,
\[
\{\widehat{\pi} \ | \ \pi\in\Pi^\ABV_{\phi}(G)\}
=
\Pi^\ABV_{\hat\phi}(G),
\]
where $\widehat{\pi}$ denotes the Aubert-Zelevinsky dual of $\pi$. See Section \ref{sec AZ} for the definitions of these notions of Aubert-Zelevinsky duality.
\end{conjecture}
\noindent
For certain groups $G$ and parameters $\phi$, Conjecture~\ref{conjecture intro: AZ} follows from \cite{AO26}, when combined with results from this paper. 
Using this we prove the weak form of Vogan's conjecture on Arthur packets, stated here.

\begin{theorem}[Weak form of Vogan's conjecture and atomic stability, Theorem~\ref{thm Vogan}]
\label{theorem intro: Vogan's conjecture, atomic}  
Let $G$ be a quasisplit classical group and let $\delta$ be a pure inner form of $G$.
Let $\psi$ be an Arthur parameter for $G$ and let $\phi_\psi$ be the associated Langlands parameter. 
Assume Conjecture~\ref{conjecture intro: AZ} for this parameter $\phi_\psi$.
Then
\[
\Pi_\psi(G) = \Pi^\ABV_{\phi_\psi}(G).
\]
Moreover, the space of stable distributions supported on the set
\begin{align}\label{equation intro: atomic}
    \left\{ \pi \in \Pi(G) \tq \phi_\pi \geq \phi_\psi,\ \phi_{\hat \pi} \geq \phi_{\hat \psi}
\right\}
\end{align}
is one-dimensional, spanned by $\etaABV{\phi_\psi}{}$.
\end{theorem}

We refer to the last assertion of the above theorem as  {\it atomic stability}\IN{atomic stability} of the Arthur packet $\Pi_{\psi}(G)$. 
The main idea in proof of Theorem~\ref{theorem intro: Vogan's conjecture} is to show the similar set for $L$-parameter is a singleton
\[
\left\{ \phi \in \Phi(G) \tq \phi \geq \phi_\psi,\ \widehat{\phi} \geq \phi_{\hat \psi}
\right\}=\{\phi_{\psi}\},
\]
 which already appears in \cite{Lo24} for general linear groups. This idea is not restricted to quasisplit classical groups and their pure inner forms.
However, at the moment, Arthur's conjectures are only known for these groups, so Arthur packets, defined in terms of the distributions satisfying endoscopic character relations, are known only for these groups.

\subsection{Geometrization of endoscopic transfer}

Motivated by Theorem~\ref{theorem intro: Vogan's conjecture, atomic} comparing Arthur packets and ABV-packets, we turn to the heart of  our proof of Vogan's conjecture: the geometrization of Kottwitz-Shelstad transfer. 
For a large class of disconnected reductive groups $G_\delta$ (including pure inner forms of quasisplit classical groups, as above, but also including many other groups), we show that Kottwitz-Shelstad transfer $\Transfer{G',s}{G_\delta} : \mathcal{D}(G')^\st \to \mathcal{D}(G_\delta)$ admits a geometric description, where $G'$ is an endoscopic group for $G$ with $\dualgroup{G'} = Z_{\dualgroup{G}}(s)^0$. 
    Here, $\mathcal{D}(G')^\st$ is the vector space of stable distributions spanned by characters of representations of $G'(F)$.
    See Theorem~\ref{theorem: lift = transfer} (and Part \ref{part: IV Lift} below) for all details on this result and how it depends on the Kazhdan-Lusztig hypothesis for $G$ and $G'$.
This result is considerably more general than required for the proof of Vogan's conjecture on Arthur packets, although the development of this geometrization was motivated by this application.
By applying our geometrization of Kottwitz-Shelstad transfer, we prove that the distributions $\etaABV{\phi}{}$ satisfy Arthur's endoscopic character relations, which is the main step in the proof of Theorem~\ref{theorem intro: Vogan's conjecture}.

\subsection{Generalized Arthur packets}

 Vogan's conjecture supports the following principle:
    \begin{quotation}
    {\it For a reductive group $G$ (not necessarily connected) with an enhanced local Langlands correspondence, an Arthur packet is an ABV-packet $\Pi^\ABV_\phi(G)$ where $\phi$ is a Langlands parameter of Arthur type}.
    \end{quotation}
This paper also shows that ABV-packets coincide with micropackets:
$\Pi^\ABV_{\phi}(G) = \Pi_{\phi}^\mic(G)$
and $\etaABV{\phi}{} =\eta^\mic_{\phi}{}$, 
where $\mic$ refers to analogous constructions defined by microlocalization following \cites{ABV92,Vog93} as opposed to vanishing cycles. This result shows that the terms micropackets and ABV-packets are synonyms for $p$-adic fields, and likewise for the associated stable distributions, in complete generality. This is explained in Appendix~\ref{section: Qmic}.

\subsection*{Relation to other work}
For connected reductive groups over $\mathbb{R}$, micropackets, also called ABV packets, were introduced in \cite{ABV92}. For classical real groups, the fact that Arthur packets are micropackets was established by Adams, Arancibia-Robert, and Mezo in \cite{AAM24}. Their proof uses the theory of twisted endoscopy; however, their work only establishes that the stable distributions attached to local Arthur packets agree with the stable distributions attached to the corresponding ABV packets. In contrast, our result for $p$-adic classical groups (Theorem \ref{theorem intro: Vogan's conjecture}) shows that all of the distributions, stable or otherwise, attached to local Arthur packets agree with the distributions attached to the corresponding ABV packets. Furthermore, \cite{AAM24} provides a geometric description of twisted endoscopic transfer for general linear groups, while the theory we develop treats a more general class of reductive groups, disconnected or otherwise.

For $p$-adic connected reductive groups, Arthur packets have been show to coincide with ABV packets for general linear groups in \cites{CR1, CR2, Lo24, Rid23}.
For tempered parameters of arbitrary groups, this result is established in \cite{CDFZ25}. Further partial progress is given for the exceptional group $G_2$ in \cite{CFZ2}. 
For $p$-adic $F_4$, some micropackets were calculated in \cite{BL26}.
These techniques and results of these works are in the same spirit as ours, but we consider a broader scope. More specifically, the previous aforementioned works do not consider twisted endoscopic transfer, while this it built into our general theory and is a crucial ingredient in our proof of Theorem \ref{theorem intro: Vogan's conjecture}.
 
The above works discuss the situation for ABV-packets. On the other hand, the study of local Arthur packets has made progress for $p$-adic groups outside of the cases of classical groups established in \cites{Art13,KMSW14,Mok15,AG17,MR18,Ish24}. For example, local Arthur packets are constructed for $G_2$ in \cites{GGJ02, GG06, AHRR24} and some twisted endoscopic character identities were studied in \cite{Tak26}. Since $G_2$ is a twisted endoscopic group of $\mathrm{PGSO}_8$, it would be interesting to apply our framework to this situation. There are other cases where local Arthur packets have also been studied but we do not list them here.

Recently, we learned that the Theorem \ref{theorem intro: Vogan's conjecture, atomic} is proved independently by Ciubotaru-Okada in \cite{CO26} in the unipotent setting for general groups. We refer to Remark \ref{rmk CO26} for details.

\subsection*{Acknowledgements}
The authors are grateful to Jeff Adams, Dan Barbasch and David Vogan, as the book \cite{ABV92} and paper \cite{Vog93} were the point of departure for this work. 
Many other people helped with the development of the ideas in this paper, including Andrew Fiori, Ahmed Moussaoui and Qing Zhang, as well as James Steele, Kristaps Balodis and José Cruz. 
We also thank Tasho Kaletha, Freydoon Shahidi, Dihua Jiang, Wee Teck Gan, Taiwang Deng and Nicolas Arancibia for helpful conversations. 

We also would like to thank the following institutes for the warm hospitality where some of the work was done: the Fields Institute for
Research in Mathematical Sciences, the Institute of Advanced Study in Mathematics in Hangzhou, the Yau Mathematical Sciences Center at Tsinghua University, the University of Chicago John W. Boyer Center in Paris, the Institute for Mathematical Sciences and the National University of Singapore, and the Brin Center for Mathematics at the University of Maryland, where these results were announced in September 2026. 

Finally, we thank James Arthur, to whom this paper is dedicated, for a wellspring of ideas that influenced our thinking and our work in many ways.

\subsection*{Use of AI}
No AI was used in the research or the writing of this paper.

\tableofcontents

\section{Overview of the sections of paper}

We now review the structure of this paper,  providing more details on the results and arguments mentioned above.

\subsection{Part~\ref{part: I Langlands}: Langlands classification for disconnected groups}

Let $G^0$ be a connected quasisplit reductive group over a $p$-adic field $F$ equipped with an automorphism $\theta$ and set $G = G^0 \rtimes {\langle{\theta}\rangle}$; groups of this sort are studied in Section~\ref{section: disconnected groups and pure inner forms}. For simplicity, we decide to assume the order of $\theta$ is either $1$ or a prime $\ell$.
Specifically, in Section~\ref{subsection: pure inner forms} we study the algebraic groups $G_\delta$ over $F$ arising from cocycles $Z^1(F,G)$ that appear in the image of the canonical map $Z^1(F,G^0) \to Z^1(F,G)$; we do not assume the automorphism $\theta$ acts non-trivially on these cocycles. 
For pure inner forms $\delta$, given formally in Definition~\ref{definition: pure inner form}, the groups $G_\delta$ need not be quasisplit or connected.
In Sections~\ref{section: classification theorem} and \ref{section: classification for parameters} we find the Langlands classification theorem for the groups $G_\delta$, for  irreducible representations of $G_\delta(F)$ and also for L-parameters; see Theorems~\ref{theorem: standard representations} and \ref{theorem: Langlands classification for L-parameters}.  

We remark that the Langlands classification for representations of quasisplit groups was given in \cite{BJ01}*{Theorem 4.2}. Also, the Langlands classification for $L$-parameters of connected reductive groups is due to \cite{SZ18}*{Section 4.6}. We extend these results to pure inner forms of disconnected groups.

\subsection{Part~\ref{part: II Transfer}: Endoscopic transfer to disconnected groups}

In this part, we rephrase Kottwitz-Shelstad transfer as a linear transformation from stable distributions to invariant distributions on disconnected groups -- namely, groups arising from pure inner forms of quasisplit reductive groups with trivial or prime order cyclic component groups, satisfying additional mild hypotheses. The idea is closely related to Kaletha's treatment for local Langlands correspondence for disconnected groups in \cite{Kal26}.
 
In Section~\ref{section: endoscopic transfer} we adapt twisted transfer as it appears in \cite{KS99} to the disconnected groups $G_\delta$.
Specifically, in Definition~\ref{definition: modified transfer} we introduce transfer
\begin{equation}\label{intro: definition of transfer}
\Transfer{G',s}{G_\delta} : 
\mathcal{D}(G')^\st
\to 
\mathcal{D}(G_\delta)
\end{equation}
from stable distributions $\mathcal{D}(G')^\st$ on $G'(F)$ to invariant distributions on $G_\delta(F)$, where $\delta$ is any pure inner form of $G$ and $(G',s,\xi)$ is an endoscopic triple for $G$.
Note that in this context, $G_\delta$ is generally not connected.
\[
\begin{tikzcd}
&&
K_{\BC}\Pi(G^0_{\delta})_{\theta_{\delta}}
\arrow{rrdd}[swap]{\ell j_!}{\mathrm{Lemma~\ref{lemma: extension by zero}}}
&&
\\ 
&&&&
\\
\arrow{uurr}{{\Transfer{G',s}{G^0_\delta,\theta_\delta}}}[swap]{\eqref{equation: KST}}
 K_{\BC}\Pi(G')^{\st} \arrow{rrrr}{\Transfer{G',s}{G_\delta}}[swap]{\mathrm{Definition~\ref{definition: modified transfer}}}
&&&& K_{\BC}\Pi(G_{\delta})  
\\
\end{tikzcd}
\]
This transfer is defined in such a way that if $\theta_\delta \in \Aut(G_\delta^0)$ generates the component group for $G_\delta$ and if $s$ decomposes as $s = s'\rtimes {\hat\theta}_\delta$, then $\Transfer{G',s}{G_\delta}$ followed by restriction to the component of $\theta_\delta \in G_\delta(F)$ in $G_\delta(F)$ coincides with Kottwitz-Shelstad transfer
\begin{equation}\label{intro: transfer and transfer}
\Transfer{G',s'}{G_\delta^0,\theta_\delta} 
=
\Res^{G_\delta(F)}_{G_\delta^0(F)\theta_\delta} 
\circ 
\Transfer{G',s}{G_\delta},
\end{equation}
taking stable distributions supported on characters of $G'$ to $\theta_\delta$-twisted distributions on 
$G_\delta^0(F)\theta_\delta$. Hereinafter, we view distributions also as virtual representations (see Section \ref{subsection: representations}).

In Section~\ref{section: LLC} we provide comprehensive desiderata for the enhanced local Langlands correspondence $\LLC^\mathfrak{e}(\pi) = (\phi_\pi,\tau_\pi)$ as it applies to the groups $G_\delta$. 
{These desiderata refer, in part, to endoscopic transfer as introduced in Section~\ref{subsection: transfer} and also to (stable) distributions of the form
\begin{equation}\label{intro: proto STD}
\mu_{\phi}:=\sum_{\pi\in \Pi_\phi(G)} 
\rank \tau_\pi \  
\SR(\pi) ,
\end{equation}
where $\SR(\pi)$ is the standard representation of $G(F)$ determined by $\pi \in \Pi_\phi(G)$. 
}
In Definition~\ref{definition: STD} we introduce distributions 
\begin{equation}\label{intro: STD}
\STD{\phi,s}{\delta}
\ceq
\sum_{\pi\in \Pi_\phi(G_\delta)} 
\langle
s,\pi
\rangle_\phi \ 
\SR(\pi)
\end{equation}
for $s\in Z_{\dualgroup{G}}(\phi)$.
When $s=1$ we elide $s$ from the notation thus defining $\STD{\phi}{\delta}$, which exactly recovers
the stable distribution \eqref{intro: proto STD}. 
Likewise, if $\delta$ is the trivial pure inner form of $G$ we suppress it from the notation, which is a convention we follow throughout the paper.

\begin{theorem}[Theorem~\ref{theorem: transfer of STD}]
\label{theorem intro: transfer of STD}
Let $(G',s,\xi)$ be an endoscopic triple for a quasisplit reductive group $G$ and let $\delta$ be a pure inner form for $G$, as in Section~\ref{section: disconnected groups and pure inner forms}. 
Assume a Langlands Correspondence for $G$ satisfying Desiderata~\ref{desideratum: LLC}\ \LC \ and \ECR\  and for $G'$ satisfying Desiderata~\ref{desideratum: LLC}\ \LC \  and \Sta.
Assume the Kazhdan-Lusztig Hypothesis for $G$ and $G'$.
For every Langlands parameter $\phi'$ for $G'$,
\[
\Transfer{G',s}{G_\delta}
(\STD{\phi'}{})
=
e(G_\delta)\ 
\STD{\xi\circ\phi',s}{\delta},
\]
where $\STD{\phi'}{}$ and $\STD{\xi\circ\phi',s}{\delta}$ are instances of the 
distributions above. 
Here $e(G_\delta)$ is Kottwitz sign of the form $G_\delta$ of $G$ as recalled in Definition~\ref{definition: kottwitz sign}.
\end{theorem}

Part~\ref{part: II Transfer} concludes with two results that will be used in Part~\ref{part: IV Lift}.
In Section~\ref{subsection: twisted transfer Bin} we give a calculation showing exactly how transfer acts in the case of $\GL_n$ as a twisted endoscopic group for $\GL_n\times\GL_n$;
this is used in Part~\ref{part: IV Lift} when we prove Theorem~\ref{theorem intro: Vogan's conjecture}.
In Section~\ref{subsection: basis of stable standards Bin} we show that the distributions $\STD{\phi'}{}$ appearing in Theorem~\ref{theorem intro: transfer of STD} give a basis for stable distributions supported by characters of $G'(F)$.
This result will play a role in the proof of our geometrization of endoscopic transfer, presented in Part~\ref{part: IV Lift}.

\subsection{Part~\ref{part: III ABV}: ABV-packets and ABV-distributions} 

In this part of the paper we discuss ABV-packets and ABV-distributions for groups $G_\delta$ that are pure inner forms of quasisplit  reductive groups $G$ as above.
While the work here builds on \cite{CFMMX22}, this treatment makes some important improvements and clarifications, especially drawing on results from Part~\ref{part: VI Schurmann}.

We begin with Definition~\ref{definition: Vogan variety} which, for every infinitesimal parameter $\lambda : W_F \to \Lgroup{G}$, recalls the Vogan variety
\[
V_\lambda
\ceq
\{ x\in \mathrm{Lie}\,\dualgroup{G} \tq \Ad(\lambda(w))(x) = |w|_F\, x, \ \forall w\in W_F 
\},
\]
which comes equipped with the action of the (connected or disconnected) reductive group 
\[
H_\lambda 
\ceq 
Z_{\dualgroup{G}}(\lambda(W_F))
=
\{
g\in \dualgroup{G} \tq \lambda(w) g = g \lambda(w) ,\ \forall w\in W_F
\}
\]
inherited from the adjoint action of $\dualgroup{G}$ on $\mathrm{Lie}\,\dualgroup{G}$.
The rule $\phi \mapsto x_\phi \in V_\lambda$
takes any Langlands parameter $\phi: W'_F \to \Lgroup{G}$
such that
\[
\phi\left(w,\left(
\begin{smallmatrix}|w|_F^{1/2} & 0 \\ 0 & |w|_F^{-1/2}\end{smallmatrix}
\right)\right) 
=\lambda(w)
\]
to the point $x_\phi \in V_\lambda$ defined by
\[
\exp(x_\phi)\rtimes 1 = \phi\left(1,
\begin{pmatrix}1 & 1 \\ 0 & 1\end{pmatrix}\right).
\]
For every such infinitesimal parameter $\lambda$, if we have the enhanced local Langlands correspondence for $G$, then it determines a perfect pairing 
\begin{equation}\label{intro: pairing of GGs}
K_0 \Per_{H_\lambda}(V_\lambda) 
\times 
K_0 \Pi_{\lambda}(G)^\pure 
\to
\mathbb{Z},
\end{equation}
where 
$K_0 \Per_{H_\lambda}(V_\lambda)$ is the Grothendieck group of the abelian category of $H_\lambda$-equivariant peverse sheaves on $V_\lambda$ and where
$\Pi_\lambda(G)^\pure$ is the set of equivalence classes of irreducible representations of $G$ and its pure inner forms with infinitesimal parameter $\lambda$; see Definition~\ref{definition: pairing}. As explained in Definition~\ref{theorem: geometric LLC}, this pairing matches $\pi \in \Pi_\lambda(G)^\pure$ with a simple $H_\lambda$-equivariant perverse sheaf $\IC(\pi)$ on $V_\lambda$.

Now we define ABV-distributions for $G_\delta$ assuming we have an enhanced local Langlands correspondence, which appears here in the rule $\pi \mapsto \IC(\pi)$ attaching simple perverse sheaves to irreducible smooth representations of $G_\delta$.

\begin{definition}[Definitions~\ref{definition: ABV-packet} and \ref{definition: stable etaABV}]
Let $G$ satisfy Definition~\ref{assu: groups we consider} and let $\delta$ be a pure inner form for $G$ as in Section~\ref{subsection: pure inner forms}.
For any Langlands parameter $\phi$ for $G$, the associated ABV-packet for $G_\delta$ is
\begin{equation}\label{intro: definition of ABV-packet}
\Pi^\ABV_\phi(G_\delta)
\ceq
\{
\pi\in \Pi_\lambda(G_\delta) \tq \NEvs_{C_\phi} \IC(\pi) \ne 0
\},
\end{equation}
where $\lambda$ is the infinitesimal parameter for $\phi$, and the associated ABV-distribution is
\begin{equation}\label{intro: definition of ABV stable}
\etaABV{\phi}{,\delta}
:=
\sum_{\pi\in \Pi_\lambda(G_\delta)}
(-1)^{d(\phi)-d(\pi)} 
\
\chi\left(\NEvs_{(x_\phi,y)}\IC(\pi)\right)\ \pi, 
\end{equation}
where:
\begin{enumerate}
    \item 
    $\Lambda_{x_\phi}$ is the vector space appearing in the conormal bundle to $V_{C_\phi}$ at $x_\phi$ and $\Lambda_{x_\phi}^\reg$ is the fibre of $\Lambda^\reg_{C_\phi} \to C_\phi$ above $x_\phi$, both studied in Section~\ref{subsection: normalization} and briefly recalled in Section~\ref{subsection: conormal};
    \item $y\in \Lambda_{x_\phi}^\reg$, any regular covector at ${x_\phi}$;
    \item 
    $d(\phi)\ceq \dim (C_\phi)$ is the dimension of the $H_\lambda$-orbit of $x_\phi$ in $V_\lambda$, as in Definition~\ref{definition: Vogan variety}, and $d(\pi)$ is likewise the dimension of the $H_\lambda$-orbit in $V_\lambda$ determined by the L-parameter of $\pi$ under the local Langlands correspondence; 
    \item 
    $\NEvs_{(x_\phi,y)} = j^*_{(x_\phi,y)} \NEvs_{C_\phi}$ is the stalk of the normalized vanishing cycles functor $\NEvs_{C_\phi} : D^b_{c,H_\lambda}(V_\lambda) \to D^b_{c,H_\lambda}(\Lambda^\reg_{C_\phi})$ that is introduced in Section~\ref{subsection: NEvs} and studied extensively in Section~\ref{section: global vanishing cycles}, which makes significant improvements on the functor $\NEvs_{C_\phi}$ from \cite{CFMMX22};
    \item 
    $\chi\left(\NEvs_{(x_\phi,y)}\IC(\pi)\right)$ is the Euler characteristic of the complex $\NEvs_{(x_\phi,y)}\IC(\pi)$.
    \end{enumerate}
If $\delta$ is trival we suppress it from the notation, thus defining $\Pi_\phi^\ABV(G)$ and $\etaABV{\phi}{}$.
\end{definition}

We prove the following fundamental properties of ABV-packets and ABV-distributions, for disconnected groups.

\begin{theorem}[Fundamental properties of ABV-distributions]
Let $G$ be a reductive group as above.
Suppose we have an enhanced local Langlands correspondence for $G_\delta$.
Let $\phi$ be any Langlands parameter for $G$.
\begin{enumerate}
    \item 
    The distribution $\etaABV{\phi}{,\delta}$ is independent of the choice of $y\in \Lambda_{x_\phi}^\reg$; see Lemma~\ref{lemma: rank and mic}.
    \item 
    The packet $\Pi^\ABV_\phi(G_\delta)$ is the support of $\etaABV{\phi}{,\delta}$ in the basis of characters of irreducible representations of $G_\delta(F)$; {\it cf}, Lemma~\ref{lemma: rank and mic}.
    \item 
    The packet $\Pi^\ABV_\phi(G_\delta)$ is contained in $\Pi_\lambda(G_\delta)$ and contains the L-packet $\Pi_\phi(G_\delta)$; see Theorem \ref{theorem: ABV-LLC compatability}. We refer to the complement of the $L$-packet $\Pi_\phi(G_\delta)$ in the ABV-packet $\Pi^\ABV_\phi(G_\delta)$ as the {\it corona of the $L$-packet}\IN{corona of an $L$-packet}.
    \item 
    If $\pi \in \Pi_\phi^\ABV(G_\delta)$ then $d(\pi) \geq d(\phi)$ with equality if and only if $\pi \in \Pi_\phi(G_\delta)$; see Theorem~\ref{theorem: ABV-LLC compatability}.
    \item 
    If $\pi \in \Pi_\phi(G_\delta)$ then the Euler characteristic of $\NEvs_{(x_\phi,y)}\IC(\pi)$ is $\rank \tau_\pi$, where $\LLC^\mathfrak{e}(\pi) = (\phi,\tau)$ 
    and this coincides with the microlocal Euler characteristic of $\IC(\pi)$; again see Lemma~\ref{lemma: rank and mic}.
    \item
    If $G$ is connected then the distributions $\etaABV{\phi}{}$ form a basis for stable distributions supported on characters of $G(F)$, as $\phi$ ranges over L-parameters for $G$; see Theorem~\ref{theorem: Conjecture 2 qs}. This resolves \cite{Vog93}*{Conjecture 81'}, \cite{CFK22}*{Conjecture 2}.
\end{enumerate}
\end{theorem}

There is a subtle point about ABV-packets that is easy to miss if one is focused only on the case when $\phi$ is of Arthur type; we refer to this as the ``third layer" for the ABV-packet.
To explain this, first observe that when the functor 
\begin{equation}\label{intro: NEvs}
\NEvs_{C_\phi} : D^b_{c,H_\lambda}(V_\lambda) \to D^b_{c,H_\lambda}(\Lambda^\reg_{C_\phi})
\end{equation}
is restricted to equivariant perverse sheaves, its target is the category of $H_\lambda$-equivariant local systems on the regular part of the conormal bundle $\Lambda_{C_\phi}$ in $\Lambda_\lambda$ to $V_{\lambda}$ above $C_\phi$:
\begin{equation}\label{intro: NEvs 1}
\NEvs_{C_\phi} : \Per_{H_\lambda}(V_\lambda) \to \Loc_{H_\lambda}(\Lambda^\reg_{C_\phi}).
\end{equation}
After choosing a base point $y\in \Lambda^\reg_{x_\phi}$, the category of $H_\lambda$-equivariant finite etale covers of $\Lambda^\reg_{C_\phi}$ is a Galois category; the fundamental group of this Galois category is denoted by $A^\reg_{\phi,y}$ and its isomorphism type is denoted by $A^\reg_{\phi}$, giving the equivalence
\begin{equation}
\Loc_{H_\lambda}(\Lambda^\reg_{C_\phi})
\iso \Rep(A_{\phi,y}^\reg).
\end{equation}

On the other hand, the choice of the base point $y\in \Lambda_{x_\phi}^\reg$ also determines a natural restriction of the target of the functor $\NEvs_{C_\phi}$, giving 
\begin{equation}\label{intro: NEvs 2}
{\NEvs_{(x_\phi,y)}}  : \Per_{H_\lambda}(V_\lambda) \to \Loc_{H_\lambda}(H_{\lambda}\cdot(x_\phi,y)).
\end{equation}
Again we see that the choice
of base point $y\in \Lambda_{x_\phi}^\reg$ determines a Galois category and the equivalence
\[
\Loc_{H_\lambda}(H_\lambda\cdot (x_\phi\cdot y)) 
\iso
\Rep(A^\mic_{\phi,y}).
\]
where  
\[
A^\mic_{\phi,y} =
Z_{H_{\lambda_\phi}}(x_\phi,y)/Z_{H_{\lambda_\phi}}(x_\phi,y)^0.
\]
It is the restriction 
\begin{equation}\label{intro: NEvs 3}
{\NEvs_{(x_\phi,y)}}  : \Per_{H_\lambda}(V_\lambda) \to \Rep(A^\mic_{\phi,y})
\end{equation}
that appears in $\etaABV{\phi}{}$, above.

These two groups {$A_{\phi,y}^\mic$} and $A_{\phi,y}^\reg$,  are generally not equal.
The extra layer of information about $\Pi^\ABV_\phi(G_\delta)$ given by $\Pi^\ABV_\phi(G_\delta) \to \Rep(A^\reg_{\phi,y})$, as opposed to $\Pi^\mic_\phi(G_\delta) \to \Rep(A^\mic_{\phi,y})$, is important for understanding how ABV-packets relate to parabolic induction; see \cite{CFK22} for an example showing this. 
These three layers fit together in the following result.

\begin{theorem}[Theorem~\ref{theorem: ABV-LLC compatability}] 
\label{theorem intro: third layer} 
Let $G_\delta$ be as above and assume an enhanced local Langlands correspondence for $G_\delta(F)$.
Let $\phi$ be any Langlands parameter for $G$.
The ABV-packet $\Pi^\ABV_{\phi}(G_\delta)$ comes equipped with 
\begin{equation}\label{intro: fundamental}
\begin{tikzcd}
 \Pi^\ABV_{\phi}(G_\delta) 
 \arrow{rrr}{\NEvs_{C_\phi}}  
 &&& \Rep(A_{\phi,y}^\reg) \arrow{d} & A_{\phi,y}^\reg \\
 \Pi^\mic_\phi(G_\delta) \arrow[equal]{u}
 \arrow{rrr}{\NEvs_{(x_\phi,y)}}
  &&& \Rep(A_{\phi,y}^\mic) & A_{\phi,y}^\mic \arrow{d} \arrow{u} \\
 \Pi_\phi(G_{\delta}) \arrow{u} \arrow{rrr}{\LLC^{\mathfrak{e}}} &&& \Rep(A_\phi) \arrow{u} & A_\phi 
\end{tikzcd}
\end{equation}
where:
    \begin{enumerate}
    \item 
    $y\in \Lambda^\reg_{x_\phi}$;
    \item
    $\Rep(A_\phi) \to \Rep(A_{\phi,y}^\mic)$ is obtained by the pullback along $A_{\phi,y}^\mic \to A_\phi$ which is induced from the inclusion $Z_{H_\lambda}(x_\phi,y) \to Z_{H_\lambda}(x_\phi)$ ;
    and 
    \item $\Rep(A_{\phi,y}^{\reg}) \to \Rep(A_{\phi,y}^\mic)$ is obtained by the pullback along $A_{\phi,y}^\mic \to A_{\phi,y}^\reg$, which corresponds to the functor defined by restriction $\Loc_{H_\lambda}(\Lambda^\reg_{C_\phi}) \to \Loc_{H_\lambda}(H_\lambda\cdot (x_\phi,y))$.
\end{enumerate}
\end{theorem}

 To summarize, $A_{\phi,y}^\reg$ is the ``third layer'' that we referred to earlier, where $A_{\phi,y}^\mic$ is the second layer. It is known that the map $\NEvs_{(x_{\phi,y})}: \Pi_{\phi}^{\mic}(G_{\delta}) \to \Rep(A_{\phi,y}^{\mic})$ is neither injective nor surjective in general. The third layer provides more refined information about the internal structure of $\Pi^\ABV_{\phi}(G_{\delta})$.

Theorem~\ref{theorem: ABV-LLC compatability} can be used to turn $\Pi^\ABV_\phi(G)$ into a multi-set in a way that is consistent with what Arthur has done for his packets. 
This is explained in Remark~\ref{remark: multiset}.

Now we have given the fundamental properties of ABV-packets for reductive groups equipped with an enhanced local Langlands correspondence. 
Our objective in Part~\ref{part: IV Lift} is to show that ABV-distributions satisfy the endoscopic character identities that appear in Arthur's work. 
But before that, and by way of motivation, we end Part~\ref{part: III ABV} with a conditional proof that Arthur packets are ABV-packets for all the groups considered by Arthur (Theorem \ref{theorem intro: Vogan's conjecture, atomic}).

The last part of Theorem~\ref{theorem intro: Vogan's conjecture, atomic} suggests a new characterization of Arthur packets.
Suppose that $G$ has an enhanced local Langlands correspondence and let $\psi$ be a local Arthur parameter of $G$.
\begin{quotation}
{\it An Arthur packet $\Pi_{\psi}(G)$ is the support of a basis element for the one-dimensional space spanned by the characters in the set \eqref{equation intro: atomic}.}
\end{quotation}
Note that this characterization does not require anything other than the local Langlands correspondence $\pi \mapsto \phi_{\pi}$.

\subsection{Part~\ref{part: IV Lift}: The geometrization of Kottwitz-Shelstad transfer} 

Our rephrasing of twisted endoscopic transfer in terms of disconnected groups was motivated by the geometrization presented in this part of the paper.
The main result of Part~\ref{part: IV Lift} concerns the transfer of ABV-distributions under certain conditions, as we now state.

\begin{theorem}[Theorem~\ref{theorem: transfer of ABV-distributions} and Definition~\ref{definition: etaABV}]\label{theorem intro: transfer of ABV-distributions}
Let $G_\delta$ be a pure inner form of a quasisplit group $G$ as in Definition~\ref{assu: groups we consider}. 
Let $(G',s,\xi)$ be an endoscopic triple for $G$ as defined in Definition~\ref{definition: endoscopic triple}.
Assume we have a local Langlands correspondence for $G$ and $G'$ satisfying desiderata specified in Theorem~\ref{theorem: transfer of ABV-distributions}. 
Assume the Kazhdan-Lusztig Hypothesis for $G$ and $G'$.
Let $\phi'$ be a Langlands parameter for $G'$ and set $\phi = \xi\circ\phi'$.
If there is some regular $y'\in \Lambda_{x_{\phi'}}$ such $y\in d\xi(y')$ is regular in $\Lambda_{x_\phi}$, then
\[
\Transfer{G',s}{G_\delta}
(\etaABV{\phi'}{})
=
e(G_\delta)\ 
\etaABV{\phi,s}{},
\]
where $e(G_\delta)$ is Kottwitz-sign of the form $G_\delta$ and where, for any $y\in d\xi(\Lambda_{x_{\phi'}}^\reg) \cap \Lambda_{x_\phi}^\reg$,
\[
\etaABV{\phi,s}{,\delta}
\ceq
\sum_{\pi\in \Pi_\phi^\ABV(G_\delta)}
(-1)^{d(\phi)-d(\pi)}
\trace(s,\NEvs_{(x_\phi,y)}\IC(\pi))\ \pi.
\]
\end{theorem}

The proof of Theorem~\ref{theorem intro: transfer of ABV-distributions} uses the fixed-point formula for vanishing cycles (Theorem~\ref{theorem: Schurmann}) applied to the functor $\Evs_z$ (Theorem~\ref{theorem: FPF}) for $z\in \Lambda_\lambda^\reg$.

The fact that $\etaABV{\phi,s}{,\delta}$ is independent of the choice of $y\in (\Lambda^\reg_{x_\phi})^s$ is perhaps surprising because the coefficients in the expansion refer to trace at $s$ of a representation of $A^\mic_{\phi,y}$ and these groups can vary as $y$ ranges over $(\Lambda^\reg_{x_\phi})^s$.
Indeed, the orbit $H_\lambda\cdot(x_\phi,y)$ and its equivariant fundamental group $A^\mic_{\phi,y}$ and the representation $\NEvs_{(x_\phi,y)}\IC(\pi)$ of $A^\mic_{\phi,y}$ can also change as $y$ varies over $\Lambda_{x_\phi}^\reg$.
One way to address this problem is to use \cite{ABV92}*{Lemma 24.3} which shows the existence of an open $Z_{H_\lambda}(x_\phi)$-subvariety of {\it microlocally generic covectors}\IN{microlocally generic} in $\Lambda^\reg_{x_\phi}$ for which the isomorphism type of the group $A^\mic_{x_\phi,y}$ does not depend on the choice of $y$; the resulting isomorphism type is denoted $A^\mic_\phi$. But the use of \cite{ABV92}*{Lemma 24.3} creates significant technical problems as the characterization of this open subvariety of microlocally generic covectors is almost tautological. 
Instead, we use the third layer in Theorem~\ref{theorem intro: third layer} and the group $A^\reg_{\phi,y}$, whose isomorphism type $A^\reg_{\phi}$ is independent of $y\in A^\reg_{\phi,y}$, to prove that $\etaABV{\phi,s}{,\delta}$ is independent of the choice of $y\in (\Lambda^\reg_{x_\phi})^s$; see Proposition~\ref{proposition: independence of y} for details of this argument.

In order to geometrize Kottwitz-Shelstad transfer, we define a linear transformation
\begin{equation}
\Lift{G',s}{G_\delta}
:
\mathcal{D}(G')^\st \to \mathcal{D}(G_\delta)
\end{equation}
on the space $\mathcal{D}(G')^\st$ of stable distributions on $G'(F)$ in the span of characters of representations. 
{Table~\ref{intro: Lift}} summarizes how this is achieved. 
The square in the middle defines $\Lift{G',s}{G_\delta}$ on the subspace $K_\BC \Pi_{\lambda'}(G')^\st$ of stable distributions determined by an infinitesimal parameter $\lambda' : W_F \to \Lgroup{G'}$, landing in the subspace $K_\BC \Pi_{\lambda}(G_\delta)$ of invariant distributions determined by an infinitesimal parameter $\lambda \ceq \xi\circ\lambda' : W_F \to \Lgroup{G}$.
The triangle on the bottom of Table~\ref{intro: Lift} shows how $\Lift{G',s}{G_\delta}$ is defined on this vector space using geometric operations that all come from the basic fact that the Vogan variety $V_{\lambda'}$ of the fix by the action of $s$ on $V_\lambda$. 
This means that the endoscopic triple $(G',s,\xi)$ for $G$ determines a morphism 
\[
d\xi : V_{\lambda'} = V_{\lambda}^s \to V_{\lambda}
\]
of prehomogeneous vector spaces and therefore induces a functor of derived categories 
\begin{equation}\label{intro: dxi*}
d\xi^*: D^b_{c,H_{\lambda}}(V_\lambda) \to D^b_{c,H_\lambda^s}(V_\lambda^s).
\end{equation}
\begin{table}[ht]
\caption{The definition of $\Lift{G',s}{G_\delta}$, used to geometrize Kottwitz-Shelstad transfer. In this diagram, $G$ is a quasisplit reductive group with cyclic component group, $\delta$ is a pure inner form for $G$, and  $(G',s,\xi')$ is an endoscopic triple for $G$.}

\[
\begin{tikzcd}
K_\CC \Pi_{\lambda'}(G')^\st
\arrow{dd}[swap]{\QST_{G',s}}{\mathrm{Definition~\ref{definition: QST and z-strongly stable}}}
\arrow{rrrr}{{\Lift{G',s}{G_\delta}}}[swap]{\mathrm{Definition~\ref{definition: Lift}}}
&&&&
K_\BC \Pi_\lambda(G_\delta)
\\
&&&&
\\
\arrow{rrrr}{\,^t(T\, \circ\,\calm^*\, \circ\, d\xi^*)}
[swap]{\mathrm{Section~\ref{subsection: KLH}}}
K_\BC \Pi_{\lambda'}(G')^{\pure,\st}
&&&&
K_\BC \Pi_{\lambda}(G)^\pure
\arrow{uu}[swap]{\mathrm{projection}}
\\
&&&&&
\\
K_\CC\Per_{H_{\lambda'}}(V_{\lambda'})  
\arrow[bend right=20]{dd}[swap]{\mathrm{induction}}
&&&& 
K_\CC\Per_{H_{\lambda}}(V_{\lambda})
\arrow{ddll}{d\xi^* \mathrm{(restriction)}}
\arrow{llll}{\mathrm{Definition~\ref{definition: Lift}}}
[swap]{T\,\circ\, \calm^*\,\circ\, d\xi^*}
\\
&&&&
\\
K_\CC\Per_{H^+_{\lambda'}}(V_{\lambda'})
\arrow{uu}[swap]{T}
&&
K_\CC\Per_{H_{\lambda}^{s}}(V_{\lambda}^{s})
\arrow{ll}{\eqref{eq m star}}[swap]{\calm^*}
\arrow[dashed]{lluu}
&&
\end{tikzcd}
\]
\label{intro: Lift}
\end{table}

The subtlety in the definition of $\Lift{G',s}{G_\delta}$ comes from the fact that although $H_{\lambda'}$ is a subgroup of $H_{\lambda}^s$, in general these groups do not coincide. 
To repair this, it is tempting to use the forgetful/change of group functor 
\[
D^b_{c,H_\lambda^s}(V_\lambda^s) 
\to
D^b_{c,H_{\lambda'}}(V_{\lambda'}) 
\]
but that does not provide a geometrization of endoscopic transfer. 
The key to the correct geometrization is the homomorphism $\calm : H^+_{\lambda'}
 \ceq H_{\lambda'} \times S \to H_\lambda^s$
  given by multiplication $\calm(h,s):= hs$, where $S$ is the Zariski closure of the group generated by $s\in \dualgroup{G}$; since $s$ is semisimple, this is a diagonalizable group.
The induced functor $\calm^*: D^b_{c,H_\lambda^s}(V_\lambda^s) \to D^b_{c,H_{\lambda'}^+}(V_{\lambda'})$ and the linear transformation $T: K_0 D^b_{c,H_{\lambda'}^+}(V_{\lambda'}) \to K_0D^b_{c,H_{\lambda'}}(V_{\lambda'})$ appearing in Table~\ref{intro: Lift} are explained in Part~\ref{part: IV Lift}.
This completes the definition of $\Lift{G',s}{G_\delta}$; full details are given in Section~\ref{subsection: Lift}.

\begin{theorem}[Theorem~\ref{theorem: lift = transfer}]\label{theorem intro: lift = transfer}
Let $G$ be a disconnected reductive group satisfying Definition \ref{assu: groups we consider} and $\delta$ be a pure inner form of $G$ (Definition \ref{definition: pure inner form}).
Let $(G',s,\xi)$ be an endoscopic triple for $G$ as in Definition \ref{definition: endoscopic triple}.
Assume a local Langlands correspondence for $G$ and $G'$ satisfying desiderata specified in Theorem~\ref{subsection: basis of stable standards Bin}.
Also assume Hypothesis~\ref{hypothesis: KLH} for $G$ and $G'$.
Then
\[
{
\Lift{G',s}{G_\delta} 
=
\Transfer{G',s}{G_\delta}
}
\]  
on the space $\mathcal{D}(G')^\st$
of stable distributions spanned by characters of representations of $G'(F)$.
Equivalently, under the pairing of Grothendieck groups determined by the enhanced Langlands correspondence,  
\[
\langle 
\Transfer{G',s}{G_\delta}(\mu) ,
\mathcal{F}
\rangle
=
\langle 
\QST_{G'}(\mu), {T\circ \calm^*\circ d\xi^*}
\mathcal{F}
\rangle
\]
where
    $\mu$ is any stable distribution on $G'(F)$ spanned by characters of representations of $G'(F)$ with infinitesimal parameter $\lambda'$, 
    and $\mathcal{F} \in D^b_{c,H_\lambda}(V_{\xi\circ\lambda'})$.
Full details are given in Theorem~\ref{theorem: lift = transfer}.
\end{theorem}

Another ingredient in the proof of Theorem~\ref{theorem intro: transfer of ABV-distributions} is a very general fixed-point formula for vanishing cycles, Theorem~\ref{theorem: Schurmann}, from which we deduce the following result.

\begin{theorem}[Theorem~\ref{theorem: FPF}]\label{theorem intro: FPF}
Let $G$ is a disconnected reductive group of the sort considered in Section~\ref{subsection: assumptions on groups}.
Let $(G',s,\xi)$ be an endoscopic triple for $G$ as in Section~\ref{subsection: endoscopic triple}.
Let $\phi'$ be a Langlands parameter for $G'$ and set $\phi = \xi \circ \phi'$.
Then
\[
(-1)^{d(\phi)}\trace(s,\NEvs_{(x_{\phi},y)}\mathcal{F})
=
(-1)^{d(\phi')} \trace(s,\NEvs_{(x_{\phi'},y)} d\xi^* \mathcal{F}),
\]
for all $y\in \Lambda_{x_{\phi'}}^\reg \cap \Lambda_{x_\phi}^\reg$ and for all $\mathcal{F}\in D^b_{c,H_\lambda}(V_\lambda)$.
\end{theorem}

\subsection{Part~\ref{part: V Arthur}: Application to Arthur packets}

In this part of the paper we prove Vogan's conjecture on Arthur packets, Theorem~\ref{theorem: Vogan's conjecture}.
We recalled this result above, as Theorem~\ref{theorem intro: Vogan's conjecture}, but here we express it in a different form.

\begin{theorem}[Theorem~\ref{theorem: Vogan's conjecture}] Assume the Kazhdan-Lusztig Hypothesis~\ref{hypothesis: KLH}, for $\GL_N\rtimes{\langle\theta\rangle}$ for outer involution $\theta$ of $\GL_N$ (or $\mathrm{Res}_{E/F}\GL_n\rtimes\langle\theta\rangle$ in the case of unitary groups).
Let $G$ be the quasisplit classical group $\U^E_n$, $\SO_{2n+1}$, $\Sp_{2n}$, or $\O^\alpha_{2n}$; let $\delta$ be a pure inner form for $G$.
    For every local Arthur parameter $\psi$ for $G$ and semisimple $s\in Z_{\dualgroup{G}}(\psi)$,
\[
\etaArt{\psi,s}{,\delta}
=
\etaABV{\phi_\psi,s}{,\delta},
\]
where $\phi_\psi$ is the Langlands parameter determined by the local Arthur parameter $\psi$.
Consequently,
\[
\Pi_\psi(G_\delta)
=
\Pi_{\phi_\psi}^\ABV(G_\delta).
\]
\end{theorem}

The proof of Theorem~\ref{theorem: Vogan's conjecture} occupies Part~\ref{part: V Arthur}. 
We begin by recalling Arthur's theory of packets and associated stable distributions for quasisplit classical groups, adapting is as we go to the disconnected reductive groups $\O^{\alpha}_{2n}$ and $ \Res_{E/F}\GL_N\rtimes{\langle\theta\rangle}$, for trivial or quadratic extensions $E/F$. 
For example, for any local Arthur parameter $\psi$ for $\GL_N$ fixed by $\widehat{\theta}$, we have 
\begin{equation}
    \Pi_{\psi^+}(\GL_N^+)
    =
    \left\{
    \pi_\psi^+, \pi_\psi^-
    \right\}  
\end{equation}
where $\pi_\psi^+$ is the Whittaker-normalised extension of the Arthur type representation $\pi_\psi$ of $\GL_N(F)$ to $\GL_N^+(F)$ and where 
\begin{equation}
\Ind^{\GL_N^+}_{\GL_N} \pi_\psi = \pi_\psi^+\oplus \pi_\psi^-.
\end{equation}
Our framework for Arthur packets of disconnected groups gives
\[
\etaArt{\psi^+}{}
=
\pi_\psi^+ + \pi_\psi^-
\qquad\text{and}\qquad
\etaArt{\psi^+,{\hat\theta}}{}
=
\pi_\psi^+ - \pi_\psi^-
.
\]

In Sections~\ref{subsection: classical groups} and \ref{subsection: G'} we review the theory of endoscopy for quasisplit classical groups and $\Res_{E/F}\GL_N\rtimes{\langle\theta\rangle}$,
as summarized in Table~\ref{table intro: endo summary}.
 Note that Arthur's ``third kind" of endoscopy fits uniformly into this framework: $(\Sp_{2n_1}\times\Sp_{2n_2},{\hat\theta},\xi)$ is an endoscopic triple for the disconnected group $\O^{\alpha}_{2n} = \SO^\alpha_{2n}\rtimes\langle\theta\rangle$ 
for a semisimple ${\hat\theta}\in \dualgroup{\O^\alpha_{2n}} = 
\O_{2n}(\BC)$ such that ${\hat\theta}\not\in \SO_{2n}(\BC)$.
Theorem~\ref{theorem intro: Arthur} recalls Arthur's endoscopic character relations for quasi-split classical groups $G$ and the ambient groups $G^+$ in Table~\ref{table intro: endo summary}.

\begin{table}[!ht]
\caption{Endoscopy for quasisplit classical groups.}
\scalebox{1.0}{
\begin{tabular}{|c|c|} 
 \hline
Endoscopic group $G^0$ & Ambient group $G^+$ \\
\hline
\hline
$\U^E_n$ & $\Res_{E/F}\GL_n\rtimes{\langle\theta\rangle}$ 
\\
$\U^{E}_{n_1}\times\U^{E}_{n_2}$ & $\U^{E}_{n_1+n_2}$ 
\\
\hline
$\SO_{2n+1}$ & $\GL_{2n}\rtimes{\langle\theta\rangle}$ 
\\
$\SO_{2n_1+1} \times \SO_{2n_2+1}$ & $\SO_{2(n_1+n_2)+1}$ 
\\
\hline
$\Sp_{2n}$ & $\GL_{2n+1}\rtimes{\langle\theta\rangle}$ 
\\
$\Sp_{2n_1} \times \SO^{\alpha}_{2n_2}$ & $\Sp_{2(n_1+n_2)}$ 
\\
\hline
$\SO^\alpha_{2n}$ & $\GL_{2n}\rtimes{\langle\theta\rangle}$ 
\\
$\SO^{\alpha}_{2n_1} \times \SO^{\alpha}_{2n_2}$ & $\O^{\alpha}_{2(n_1+n_2)}$ 
\\
$\Sp_{2n_1} \times \Sp_{2n_2}$ & 
$\O^\alpha_{2(n_1+n_2)}$
\\
\hline
\end{tabular}
}
\label{table intro: endo summary}
\end{table}

\begin{theorem}[\cites{Art13,KMSW14,Mok15,AG17,MR18,Ish24}]
\label{theorem intro: Arthur}
Let $G$ be the quasisplit classical group $\U^E_{n}$, $\SO_{2n+1}$, $\Sp_{2n}$ or $\O^\alpha_{2n}$. For any local Arthur parameter $\psi$ of $G$, there exists distributions $\etaArt{\psi,s}{}$ for semisimple $s \in Z_{\dualgroup{G}}(\psi)$, that satisfies Arthur's endoscopic character relations; see Conjecture \ref{conj: Arthur}. 
In particular, if $G^+$ is the ambient group which realizes $(G^0,s,\xi)$ as a endoscopic triple for $G^+$ as summarized in Table~\ref{table intro: endo summary}, then
\[
\Transfer{G^0,s}{G^+}\,
\frac{1}{[G:G^0]} \Res^{G}_{G^0}
(\etaArt{\psi}{}
)
=
\etaArt{\xi\circ\psi,s}{}.
\]
These relations define the function $\Pi^\Art_\psi(G) \to \Rep(A_\psi)$, compatible with the enhanced local Langlands correspondence as in the following commutative diagram.
\[
\begin{tikzcd}
 \Pi^\Art_{\psi}(G) 
 \arrow{rrr}{\pi\mapsto \langle \ , \pi \rangle_\psi^\Art} 
 &&& \Rep(A_\psi) 
 & A_\psi \arrow{dd} \\
  &&&
  & 
  \\
 \Pi_{\phi_\psi}(G) \arrow{uu} \arrow{rrr}{\LLC^{\mathfrak{e}}} &&& \Rep(A_{\phi_\psi}) \arrow{uu} & A_{\phi_\psi}
\end{tikzcd}
\]
\end{theorem}

Up to this point, Part~\ref{part: V Arthur} has re-cast the existing theory of Arthur packets for quasisplit classical groups in the framework of disconnected groups, as developed in Part~\ref{part: I Langlands}. 
Now, we specialize the work done in Part~\ref{part: IV Lift} to the case of local Arthur parameters for these groups. 
We use the Langlands correspondence for the quasisplit classical groups $G$ to define ABV-packets $\Pi^\ABV_{\phi_\psi}(G)$.
For comparison with Arthur's theory, we define ABV-packets for $\O_{2n}^\alpha$, without a local Langlands correspondence for $\SO_{2n}^\alpha$.
In the proof of Vogan's conjecture on Arthur packets, Theorem~\ref{theorem: Vogan's conjecture}, we show that the ABV-distributions $\etaABV{\phi_\psi}{}$ satisfy the same endoscopic character relations as the distributions $\etaArt{\psi}{}$ satisfy, always remembering that $G$ need not be connected, in our framework.

The following commuting diagram
\begin{equation}\label{intro: fundamental A-ABV}
\begin{tikzcd}
\Pi^\ABV_{\phi_\psi}(G) 
 \arrow{rrrr}{\NEvs_{C_\psi}} 
 &&&& \Rep(A_{\phi_\psi}^\reg)
 \arrow{d}
& A_{\phi_\psi}^\reg & \\
\Pi_\psi^\Art(G) \arrow[equal]{u} 
\arrow{rrrr}{\pi \mapsto \langle\ , \pi \rangle_\psi}
&&&& \Rep(A_{\psi}) 
& A^\mic_{\phi_\psi,y_\psi}  \arrow[equal]{r} \arrow{u}\arrow{d} & A_{\psi} \\
\Pi_{\phi_{\psi}}(G) 
\arrow{u} 
\arrow{rrrr}{\LLC^{\mathfrak{e}}} 
 &&&& \Rep(A_{\phi_\psi}) \arrow{u} & A_{\phi_\psi} & 
\end{tikzcd}
\end{equation}
shows how Theorem~\ref{theorem intro: Arthur} relates to Theorem~\ref{theorem: ABV-LLC compatability}. 
The Arthur parameter $\psi$ determines a base point $z_\psi = (x_{\phi_\psi}, y_\psi) \in \Lambda^\reg_{C_{\phi_\psi}}$ that is used crucially in the equality above.

For $\SO^\alpha_{2n}$, there is another way to relate ABV-packets and Arthur packets, which comes from the following more general result.

\begin{theorem}[Theorem~\ref{theorem: tranfer from the identity}]
Let $G$ be a reductive group of the sort considered in Section~\ref{subsection: assumptions on groups} and $\delta$ any pure inner form of $G$.
    Assume the local Langlands correspondence for $G_\delta$ and $G_\delta^0$. 
For any Langlands parameter $\phi$ for $G$,
\[
\etaABV{\phi}{,\delta}
=
\Ind_{G_\delta^0}^{G_\delta} 
(
\etaABV{\phi^0}{}
)
\]
\end{theorem}

Together with other work in this paper, we show that applying this to $G = \O^\alpha_{2n}$ gives the following result.

\begin{theorem}
Assume the local Langlands correspondence for $\SO^\alpha_{2n}$.
Then
\[
\etaArt{\psi}{}
=
\Ind_{\SO^\alpha_{2n}}^{\O^\alpha_{2n}} 
(
\etaABV{\phi_{\psi^0}}{}
)
\]
for all local Arthur parameters $\psi$ for $\O^\alpha_{2n}$.
where $\psi^0$ is the associated local Arthur parameter for $\SO^\alpha_{2n}$.
In particular, 
\[
\Pi^\Art_\psi(\O^\alpha_{2n})
=
\Ind_{\SO^\alpha_{2n}}^{\O^\alpha_{2n}}\ (\Pi^\ABV_{\phi_{\psi^0}}(\SO^\alpha_{2n})).
\]
\end{theorem}

\subsection{Part~\ref{part: VI Schurmann}: Appendices on geometry and functors, by Jörg Schürmann}

This part of the paper makes important improvements to the theory and practice of vanishing cycles as they appear in \cite{CFMMX22}.
Some of these improvements play a crucial role in the proof of Theorem~\ref{theorem: transfer of ABV-distributions}, which is one of the most fundamental results of this paper. 

Appendix~\ref{section: equivariant vanishing cycles} proves a very general fixed-point formula for (local) vanishing cycles which is used in Part~\ref{part: IV Lift} to study the transfer of ABV-distributions. It also establishes an inverse Kashiwara-Dubson formula in this language, which is used in our proof that ABV-distributions give a basis for stable distributions.

Appendix~\ref{section: global vanishing cycles} introduces the functor
\begin{equation}
\Ev_{C_\phi}
:
D^b_{c,H_\lambda}(V_\lambda)
\to
D^b_{c,H_\lambda}(\Lambda_{C^*_\phi})
\end{equation}
which is an improved version of a functor by the same name introduced in \cite{CFMMX22}, from which $\NEvs_{C_\phi}: \Per_{H_\lambda}(V_\lambda) \to \Loc_{H_\lambda}(\Lambda^\reg_{C_\phi})$ is defined as
explained in Part~\ref{part: VI Schurmann}.
The improvements made here regarding $\Ev_{C_\phi}$ are significant. 
Especially, the target of this functor is the equivariant derived category of sheaves on  $\Lambda_{C^*_\phi}$, not only on the open subvariety $\Lambda^\reg_{C^*_\phi} = \Lambda^\reg_{C_\phi}$. 
We prove
\begin{equation}\label{intro: EFPF}
\trace(s,\Ev_{(x_{\phi},y)}\mathcal{F})
=
 \trace(s,\Ev_{(x_{\phi'},y)} d\xi^* \mathcal{F}),
\end{equation}
for all $y\in \Lambda_{x_\phi}$ and  for all $\mathcal{F}\in D^b_{c,H_\lambda}(V_\lambda)$,
where $\Ev_{(x_\phi,y)}$ is the stalk, in the sense above, of the functor $\Ev_{C_\phi}$.
This is a consequence of the following more general result.

\begin{theorem}[Theorem~\ref{theorem: Schurmann}]
Let $M$ be a smooth complex algebraic variety equipped with the action of a diagonalizable group $G$. 
For $\sigma\in G$ (automatically semisimple), let $f: M\to \BC$ be $\sigma$-invariant morphism of complex varieties: $f =f\circ \sigma$.
Then, for every $\mathcal{F}\in D^b_{c,G}(M;k)$,
\[
\trace\left(\sigma, (\Rphi{f}[-1]{ \mathcal{F}})_x \right)
=
\trace\left(\sigma, (\Rphi{f\vert_{M^\sigma}}[-1]{ \mathcal{F}\vert_{M^\sigma}})_x \right),
\]
for all $x\in M^\sigma$.
\end{theorem}

Appendix~\ref{section: Qmic} 
studies the functor $Q^\mic$ that was introduced in \cite{ABV92}*{Theorem~24.8} and gives a complete proof of that result. 
This functor was used in \cite{ABV92} to define micropackets $\Pi_\phi^\mic(G)$ for real groups $G$ and \cite{Vog93} suggests the same approach to define micropackets for $p$-adic groups. 
Here we also prove that these two approaches produce the same packets: {\it ABV-packets, defined by vanishing cycles, are precisely the same as micropackets, defined by microlocalization}.

\part{Langlands classification for disconnected groups}\label{part: I Langlands}

In this part, we recall the theory of twisted endoscopy for disconnected groups under certain conditions. In Section~\ref{section: endoclassical} we verify these assumptions for classical groups and certain other groups occurring in the theory of twisted endoscopy.

\section{Disconnected groups and their pure inner forms}\label{section: disconnected groups and pure inner forms}

In this section, we discuss disconnected reductive groups, their pure inner forms, their standard parabolic subgroups, and related notions.

\subsection{Some quasisplit disconnected groups}\label{subsection: assumptions on groups}

Let $F$ be a nonarchimedean local field of characteristic zero, $\ol{F}$ be an algebraic closure of $F$, and set $\Gamma_F=\mathrm{Gal}(\ol{F}/F)$, the absolute Galois group of $F$. We fix an non-trivial additive character of $F$ throughout the paper.
In this paper, we consider quasisplit disconnected groups $G$, in the sense of \cite{Kal26}*{Section 3.1}, and furthermore, we require that the component group is either trivial or cyclic of prime order.

\begin{definition}\label{assu: groups we consider}
Let $G^0$ be a connected quasisplit reductive group over \(F\) and let \(\theta\) be an automorphism of \(G^0\) defined over \(F\). This automorphism is either trivial or has a prime order \(\ell\).
Let \((B,T,\{X_{\alpha}\})\) be an \(F\)-pinning of \(G^0\) that is preserved under the automorphism \(\theta\).
From this data, we define \[G := G^0 \rtimes \langle \theta \rangle.\] Every group \(G\) considered in this paper arises in this way. Consequently, whenever we write \(G\), we have implicitly fixed the data
\[
G^0, \theta, (B,T,\{X_{\alpha}\})
\]
where we view \(\theta\) as a chosen generator of \(\pi_0(G)\) that fixes the \(F\)-pinning \((B,T,\{X_{\alpha}\})\). The $F$-pinning \((B,T,\{X_{\alpha}\})\) and the additive character of $F$ (chosen and fixed, above) determine a Whittaker datum $\mathfrak{w}$ of the quasisplit group $G^0$.
\end{definition}

We will often encounter cyclic groups of prime order $\ell$, e.g., $\langle\theta\rangle$. It is convenient to fix some notation now. Let $x$ be a generator of a cyclic group of prime order $\ell$ fixed in the context. We fix a primitive $\ell$-th root of unity $\zeta$ and consider the character $\vartheta$ of $\langle x\rangle$ defined by $\vartheta(x):=\zeta.$ Note that the set of irreducible representations of $\langle x \rangle$ is $\langle\vartheta\rangle$, which is also a cyclic group of order $\ell.$
Since $\theta$ is defined over $F$,  $\Gamma_F$ acts trivially on $\langle\theta\rangle$.

\begin{remark}\label{remark: trivial a}
The pair $(G^0,\theta)$ is an instance of the data $(G^0,\theta,\mathbf{a})$ appearing in \cite{KS99}*{Section 2.1} for trivial $\mathbf{a} \in H^1(W_F,Z(\dualgroup{G^0}))$. 
The conditions we have given above relate to properties explained in \cite{KS99}*{Section 1}. 
\end{remark}

We also consider pure inner forms of $G$, in a sense that will be explained in Section~\ref{subsection: pure inner forms}.

\subsection{Pure inner forms}\label{subsection: pure inner forms}

In this section, we define the notion of pure inner forms for the quasisplit groups $G$ satisfying Definition~\ref{assu: groups we consider}; our approach is informed by \cite{Kal26}, but not exactly the same.

Recall that a pure inner form for the connected group $G^0$ is a cocycle $\delta^0\in Z^1(F,G^0)$ and it determines a form of $G^0$, denoted by $G^0_{\delta^0}$, defined by the image of $\delta^0$ under $Z^1(F,G^0) \to Z^1(F,\Aut(G^0))$. More precisely, we have $G^0_{\delta^0}(\overline{F})=G^0(\overline{F})$ and the Galois action is given by $\sigma \mapsto \sigma_{\delta^0}:=\textrm{Int}(\delta^0(\sigma)) \circ \sigma.$ We will call this the $\delta^0$-twisted Galois action.

\begin{remark}
Pure inner forms for quasisplit connected reductive groups $G^0$ were introduced in \cite{Vog93}*{Definition~2.6}, where they were called pure rational forms of $G^0\otimes_{F}{\bar F}$.
\end{remark}

Similarly, every $\delta \in Z^1(F,G)$ determines a form of $G$, denoted by $G_\delta$, defined by the image of $\delta$ under $Z^1(F,G) \to Z^1(F,\Aut(G))$. 
It is important to note that its identity component, $G_\delta^0$, may not be an inner form to $G^0$ (see Example~\ref{example: inner form}). This is  one motivator for Definition~\ref{definition: pure inner form} to restrict to a subset of $Z^1(F,G)$.

\begin{example}
\label{example: inner form}
Let $G = \mathrm{O}_{2n}$ such that $G^0 = \SO_{2n}$ is split (see \ref{subsection: classical groups} for notation). Since $H^1(F, G)$ parametrizes the isomorphism classes of nondegenerate quadratic spaces over $F$ of dimension $2n$ (\cite{KMRT98}*{Equation (29.28)}), there exists a $\delta \in Z^{1}(F, G)$ such that $G_{\delta} \cong \mathrm{O}(V)$ for some quadratic space $V$, where ${\rm disc}\,V \neq (-1)^n$ mod $(F^{\times})^{2}$. Then its identity component $G^0_{\delta} = \SO(V)$ is an outer form of $G^0$.
\end{example}

\begin{definition}\label{definition: pure inner form}
    Suppose that $G$ satisfies Definition~\ref{assu: groups we consider}.
    A {\it pure inner form} of $G$\IN{$\delta$, pure inner form of $G$} is a cocycle $\delta \in Z^1(F,G)$ that lies in the image of the injection $Z^1(F,G^0) \to Z^1(F,G)$. Note that, by this definition, a pure inner form $\delta$ of $G$ uniquely determines a pure inner form of $G^0$, which we denote by $\delta^0$.
\end{definition}

\begin{lemma}
If $\delta$ is a pure inner form of $G$ and determines a pure inner form $\delta_0$ of $G^0$, then $G_\delta^0 = G^0_{\delta^0}$ as algebraic groups over $F$. 
\end{lemma}
\begin{proof}
    Note that $G_\delta^0(\overline{F})= G^0(\overline{F})= G^0_{\delta^0}(\overline{F})$. The $\delta$-twisted action on $G_\delta^0(\overline{F})$ is the same as the $\delta^0$-twisted action on $G^0_{\delta^0}(\overline{F})$ under this identification. Thus, $G_\delta^0 = G^0_{\delta^0}$ as schemes over $F$. 
\end{proof}

\begin{remark}
    It is an interesting problem to characterize those forms of $G^0$ that are produced as $G_{\delta}^0$ for arbitrary $\delta\in Z^1(F,G)$. We will not consider this here as it will not be needed in our main results.
\end{remark}

In the case at hand, the automorphism $\theta\in \Aut_F(G^0)$ (automorphism in the category of $F$-group schemes) acts on $Z^1(F,G^0)$ and on $H^1(F,G^0)$.
Our definition takes this into account, as we now explain.
While $Z^1(F,G^0)\to Z^1(F,G)$ is injective, it is certainly not true that $H^1(F,G^0)\to H^1(F,G)$ is injective.
The following lemma addresses this fact. 
Let $H^1(F,G^0)_{\langle\theta\rangle}$ denote the coinvariants of $H^1(F,G^0)$ under the action of $\langle\theta\rangle$; thus, 
\[
H^1(F,G^0)_{\langle\theta\rangle}:= H^1(F,G^0)/\sim,
\]
where the equivalence relation is defined by $[\delta_1] \sim [\delta_2]$ if and only if there exists a $\theta^k \in \langle \theta \rangle$ such that $[\delta_2]= \theta^k [\delta_1]$.

\begin{lemma}\label{lem: coinvariants}
    For $H^1(F,G^0)_{\langle\theta\rangle}$ as defined above, 
    there is a canonical bijection
    \[
H^1(F,G^0)_{\langle\theta\rangle} \rightarrow \mathrm{image}(H^1(F,G^0)\to H^1(F,G)).
    \]
    We define $H^1(F,G)_{ \theta }:=\mathrm{image}(H^1(F,G^0)\to H^1(F,G))$; \IN{$H^1(F,G)_{ \theta }$, classes of pure inner forms} this gives the isomorphism 
    \[
H^1(F,G^0)_{\langle\theta\rangle}\iso H^1(F,G)_{\theta}.
    \]
\end{lemma}

\begin{proof}
We show that $[\delta_1], [\delta_2] \in H^1(F,G^0)$ have the same image in $H^1(F,G)$ if and only if $[\delta_1]= \theta^k [\delta_2]$ for some $k$. Indeed, if they have the same image, then there exists a $g \rtimes \theta^k \in G(\overline{F})$ such that 
\[ \delta_1(\sigma){ \rtimes 1}= (g \rtimes \theta^k)^{-1} (\delta_2(\sigma){ \rtimes 1)} \sigma(g \rtimes \theta^k)\]
for all $\sigma \in \Gamma_F$. Recall that $\sigma$ acts trivially on $1 \rtimes \theta$. Thus, we have 
\[ \delta_1(\sigma){ \rtimes 1}= (1 \rtimes \theta^k)^{-1} (g\rtimes 1)^{-1} (\delta_2(\sigma){ \rtimes 1)} \sigma (g\rtimes 1)  (1\rtimes \theta^k)= \theta^{{-}k} ( \delta_3 (\sigma)){ \rtimes 1},\]
where $\delta_3(\sigma):= g^{-1} \delta_2(\sigma) \sigma(g)$. Note that $g \in G^0$; hence $[\delta_2]=[\delta_3]$ in $H^1(F,G^0)$.
We conclude that $[\delta_1]= \theta^{{-}k} [\delta_3]= \theta^{{-}k} [\delta_2]$ in $H^1(F,G^0)$.
\end{proof}

\begin{remark}    
    The action of ${\langle\theta\rangle}$ on $H^1(F,G^0)$ may indeed be non-trivial. For example, consider the split group $G^0 =\mathrm{PSO}_{2n}$ with $\theta$ denoting an order 2 outer automorphism. Then, assuming $n\geq 2$, there exists $[{\delta^0}]\in H^1(F,\mathrm{PSO}_{2n})$ such that $\theta([{\delta^0}])  \ne [{\delta^0}]$. More specifically, consider the unique nontrivial element $[\gamma^0]\in H^1(F,\SO_{2n})$, where $\SO_{2n}$ is split over $F$. Via the map $Z^1(F,\SO_{2n})\rightarrow Z^1(F,\mathrm{PSO}_{2n})$, we consider the induced element $[\overline{\gamma^0}]\in H^1(F,\mathrm{PSO}_{2n}).$ Then $\theta$ fixes $[\overline{\gamma^0}]$ as $[\gamma^0]$ must be fixed by an outer automorphism of $\SO_{2n}$ (since it would fix the trivial element and $|H^1(F,\SO_{2n})|=2$). However, $|H^1(F,\mathrm{PSO}_{2n})|=4$ and the remaining nontrivial elements of $H^1(F,\mathrm{PSO}_{2n})$ are exchanged by $\theta$. Indeed,
    this phenomenon is discussed in \cite{Art13}*{Section 9.4, p. 552} (in the context of inner forms of $\SO_{2n}$ which correspond to pure inner forms of $\mathrm{PSO}_{2n}$).
\end{remark}

Although it will not appear again until we study endoscopic transfer, 
here we define Kottwitz sign $e(G_\delta)$ of the algebraic group $G_\delta$ over $F$ as Kottwitz sign for the identity component $e(G_{\delta}^0)$, as in \cite{Kal26}*{Remark 4.13.2}. 

\begin{definition}\label{definition: kottwitz sign}
   Let $G$ be as in Definition~\ref{assu: groups we consider} and let $\delta$ be a pure inner form of $G$ as defined above in this section. Let $r(G_\delta^0)$ denote the $F$-rank of the derived group of $G_\delta^0$. Then Kottwitz sign of $G_\delta$ is  $e(G_{\delta}):= e(G_\delta^0):=(-1)^{r(G^0)-r(G_\delta^0)}$.\IN{$e(G_{\delta})$, Kottwitz sign} 
\end{definition}

\subsection{Rational points for pure inner forms}\label{subsection: rational points for pure inner forms}
Next we discuss the rational points of the pure inner forms that we introduced in the previous section. 
Let $\delta$ be a pure inner form of $G$ as in Definition \ref{definition: pure inner form} and determine a pure inner form $\delta^0$ of $G^0$. According to 
\cite{Kal26}*{\S 4.1}, we have a short exact sequence of the $F$-points
 \begin{equation}\label{eqn: exact seq pure inner forms first form}
1 \to G^0_{\delta^0}(F) \to G_{\delta}(F) \to \langle \theta \rangle^{[\delta^0]} \to 1,     
 \end{equation}
where $\langle \theta \rangle^{[\delta^0]}$ is the subgroup of $\langle \theta \rangle $ that fixes the class $[\delta^0] \in H^1(F,G^0)$. By our condition that $\langle \theta \rangle $ is of prime order $\ell$, if $\theta[\delta^0] \neq [\delta^0]$, then $G_{\delta}(F)$ is isomorphic to $G^{0}_{\delta^0}(F)$. We give a canonical isomorphism between $G^0_{\delta^0} $ and $G^0_{\theta(\delta^0)}$ as algebraic groups over $F$ in the following lemma.

\begin{lemma}\label{lem: F-point iso for non-fixed delta}
    The conjugation by $1 \rtimes \theta$ gives a canonical isomorphism between $G^0_{\delta^0} \cong G^0_{\theta(\delta^0)}$ defined over $F$.
\end{lemma}
\begin{proof}
Clearly, the conjugation gives an isomorphism on $\overline{F}$-points. We check that it intertwines the $\delta$-twisted action and $\theta(\delta)$-twisted action of $\Gamma_F$ on $G^0(\overline{F})$. Let $g\in G^0_{\delta^0} (\overline{F})$ and $\sigma\in\Gamma_F$. Recall that $(1\rtimes\theta)(g\rtimes 1)(1\rtimes\theta^{-1})=\theta(g)\rtimes 1.$ Therefore, 
\begin{align*}
(1\rtimes\theta)(\sigma_{\delta^0}(g)\rtimes 1)(1\rtimes\theta^{-1})&= \theta(\sigma_{\delta^0}(g))\rtimes 1 \\
&=\theta(\delta^0(\sigma)(\sigma(g))\delta^0(\sigma)^{-1})\rtimes 1 \\
&=\theta(\delta^0(\sigma))\theta(\sigma(g))\theta(\delta^0(\sigma)^{-1}) \rtimes 1 \\
&=\theta(\delta^0(\sigma))\sigma(\theta(g))\theta(\delta^0(\sigma)^{-1}) \rtimes 1 \\
&=(\mathrm{Int}(\theta(\delta^0(\sigma)))\circ\sigma)(\theta(g)) \rtimes 1 \\
&=\sigma_{\theta(\delta^0)}(\theta(g)) \rtimes 1.
\end{align*}
We conclude that the conjugation by $1 \rtimes \theta$ gives an isomorphism between $G^0_{\delta^0} \cong G^0_{\theta(\delta^0)}$ defined over $F$.
\end{proof}

Now we consider the case that $\theta[\delta^0] = [\delta^0]$, where we have $G_{\delta}(F)/ G^0_{\delta^0}(F) \cong \langle \theta \rangle$. We will fix a lift $\theta_{\delta} \in G_\delta(F)$ of $\theta$ as follows. By definition, there exists a $g_\delta\in G^0(\overline{F})$ such that $\theta(\delta^0(\sigma))=g_\delta^{-1}\delta^0(\sigma)\sigma(g_\delta)$ for any $\sigma\in\Gamma_F.$ We fix such $g_\delta$ and set $\theta'_{\delta} = g_{\delta} \rtimes \theta \in G_{\delta}(\overline{F})$ and
$\theta'_{\delta^0}:=\mathrm{Int}(g_\delta)\circ\theta$, an automorphism of $G^0_{\delta^0}$.

\begin{lemma}\label{lemma: automorphism preserves Borel}
   In the above setting, $\theta'_{\delta}\in G_{\delta}(F)$ and the automorphism $\theta'_{\delta^0}$ of $G^0_{\delta^0}$ is defined over $F$. 
\end{lemma}

\begin{proof} 
We show that $\sigma_{\delta}(g_{\delta} \rtimes \theta)= g_{\delta} \rtimes \theta$, where recall that $\sigma_\delta=\textrm{Int}(\delta(\sigma)) \circ \sigma$. Indeed,
\begin{align*}
    \sigma_{\delta}(g_{\delta} \rtimes \theta)&= \delta(\sigma) \sigma(g_{\delta} \rtimes \theta) \delta(\sigma)^{-1}\\
    &= \delta(\sigma) \sigma(g_{\delta} \rtimes 1) (1 \rtimes \theta) \delta(\sigma)^{-1} (1 \rtimes \theta)^{-1} (1 \rtimes \theta)\\
    &=\delta(\sigma) \sigma(g_{\delta} \rtimes 1) \theta( \delta(\sigma)^{-1} ) (1 \rtimes \theta)\\
    &=\delta(\sigma) \sigma(g_{\delta} \rtimes 1) \sigma(g_{\delta}\rtimes 1)^{-1} \delta(\sigma)^{-1} (g_{\delta}\rtimes 1) (1 \rtimes \theta)\\
    &=g_{\delta}\rtimes \theta.
\end{align*}
Here we again use the condition that $\theta$ is defined over $F$. This completes the proof of the lemma.
\end{proof}

Next, we take a minimal $F$-parabolic subgroup $P_{\delta^0,\min}^0$ of $G_{\delta^0}^0$ with a Levi factor $M^0_{\delta^0, \min}$. Since $\theta_{\delta^0}(P_{\delta^0,\min}^0)$ is also a minimal $F$-parabolic subgroup, there exists a $b_{\delta^0} \in G_{\delta^0}^0(F)$ such that 
\[
b_{\delta^0} (\theta'_{\delta^0}(P_{\delta^0,\min}^0) )b_{\delta^0}^{-1}=P_{\delta^0,\min}^0 \quad \text{and} \quad b_{\delta^0} (\theta'_{\delta^0}(M_{\delta^0,\min}^0) )b_{\delta^0}^{-1}=M_{\delta^0,\min}^0.
\]

\begin{definition}\label{def: theta delta}
In the above setting, let $\theta_{\delta} := b_{\delta^0} \theta'_{\delta} \in G_{\delta}(F)$,  and $\theta_{\delta^0}:=\mathrm{Int}(b_{\delta^0})\circ\theta'_{\delta^0}$, an $F$-automorphism of $G^0_{\delta^0}$ fixing the minimal $F$-parabolic subgroup $P^0_{\delta^0, \min}$ and Levi factor $M^0_{\delta^0, \min}$. \IN{$\theta_{\delta}$, $\theta_{\delta^0}$, $F$-automorphism of $G_\delta^0$} 
\end{definition}

Recall that $\theta$ is assumed to be of finite order $\ell$. This may not be true for $\theta_{\delta^0}$; however, it is true for its image in $\mathrm{Out}(G^0_{\delta^0}).$ Indeed, we have
\[
\theta_{\delta}^{\ell} = (b_{\delta^0}g_{\delta} \rtimes \theta)^{\ell}= \left(\prod_{i=0}^{\ell-1}\theta^i(b_{\delta^0}g_\delta) \right)\rtimes 1\in G^0(\overline{F}) \cap G_{\delta}(F) = G_{\delta}^0(F). 
\]
Thus, $\theta_{\delta^0}^\ell=\mathrm{Int}(\prod_{i=0}^{\ell-1}\theta^i(b_{\delta^0}g_\delta))$ is an inner automorphism of $G^0_{\delta^0}$. Since $\theta_{\delta}^{\ell}$ normalizes $P_{\delta^0,\min}^0$ and $M_{\delta^0,\min}^0$, then $\theta_{\delta}^{\ell} \in M^0_{\delta^0, \min}(F)$.

\subsection{Parabolic subgroups}\label{subsection: parabolic subgroups}

In this section, we review the parametrization of parabolic subgroups for pure inner forms $G_\delta$, and introduce the notation of standard parabolic subgroups following \cite{BJ01}*{\S 3}.

First, we review the case for the connected group $G_{\delta^0}^0$.
Recall that we have fixed an $F$-pinning $(B,T,\{X_{\alpha}\})$ of $G^0$ preserved under the action of $\theta$. We let $\Delta$ denote the corresponding set of simple roots, which admits a commuting action of $\Gamma_F$ and $\langle \theta \rangle$.
Following \cite{Bor79}*{Section 3}, fix a minimal parabolic $F$-subgroup $P_{\delta^0,\min}^0$ of $G^0_{\delta^0}$ that is fixed by $\theta_{\delta^0}$, as in Definition \ref{def: theta delta}, and let $\Delta_{\delta^0} \subseteq\Delta$ denote the set of simple roots corresponding to its Levi subgroup.
Let $p(G^0)$ denote the set of conjugacy classes of parabolic subgroups of $G^0$ over ${\bar F}$ and let $p(G^0)_F$ denote the subset of conjugacy classes which are defined over $F$, {\it i.e.}, for any  $P\in[P]\in p(G^0)_F$, we have that ${}^\gamma P\in [P]$ for any $\gamma\in\Gamma_F$. 
We let $p(G^0_{\delta^0}/F)$ denote the set of conjugacy classes of parabolic subgroups of $G^0_{\delta^0}(F)$. 
We have an inclusion $p(G^0_{\delta^0}/F)\hookrightarrow p(G^0)_F$. 
The set $p(G^0)_F$ corresponds to $\Gamma_F$-stable subsets of $\Delta,$ while the set $p(G^0_{\delta^0}/F)$ corresponds to $\Gamma_F$-stable subsets of $\Delta$ which contain $\Delta_{\delta^0}$.
Given a $\Gamma_F$-stable set $\Delta_{\delta^0}\subseteq {S}\subseteq\Delta,$ let $P_{S}^0$ denote the corresponding parabolic $F$-subgroup of $G^0_{\delta^0}$ containing $P_{\delta^0,\min}^0= M_{\delta^0,\min}^0 N_{\delta^0,\min}^0$.  We let $P_{S}^0=M_{S}^0 N_{S}^0$ be the Levi decomposition with Levi subgroup $M_{S}^0 \supseteq M_{\delta^0,\min}^0$. We call such $P_{S}^0$ and $M_S^0$ a standard parabolic subgroup and a standard Levi subgroup of $G^0_{\delta^0}$.

\begin{lemma}
\label{lemma: injection of pure inner forms of Levi}
For any Levi subgroup $M^0$ of $G^0$, the natural map
\[
H^{1}(F, M^0) \rightarrow H^{1}(F, G^0)
\]
is injective.
\end{lemma}

\begin{proof}
Using Kottwitz isomorphism \cite{Kot84}*{Proposition 6.4} we form the commutative diagram
\[
\xymatrix{
H^{1}(F, M^0) \ar[r]^{\simeq \quad \quad \quad} \ar[d] & \Hom(\pi_0(Z(\dualgroup{M^0})^{\Gamma_F}), \BC^{\times }) \ar[d] \\
H^{1}(F, G^0) \ar[r]^{\simeq \quad \quad \quad} & \Hom(\pi_0(Z(\dualgroup{G^0})^{\Gamma_F}), \BC^{\times }),
}
\]
where the right vertical arrow is induced by the natural inclusion
\(
Z(\dualgroup{G^0}) \rightarrow Z(\dualgroup{M^0}).
\)
Then it suffices to show that this induces a surjection
\[
\pi_0(Z(\dualgroup{G^0})^{\Gamma_{F}}) \rightarrow \pi_0(Z(\dualgroup{M^0})^{\Gamma_F}),
\]
i.e., $Z(\widehat{M^0})^{\Gamma_{F}} = (Z(\widehat{M^0})^{\Gamma_{F}})^0 
Z(\widehat{G^0})^{\Gamma_{F}}$.

Consider the $\Gamma_{F}$-equivariant exact sequence (here $\Delta^{M^0}$ is the simple roots contained in $M^0$)
\[
1 \rightarrow Z(\widehat{G^0}) \rightarrow Z(\widehat{M^0}) \rightarrow \prod_{\alpha \in \Delta \backslash \Delta^{M^0}} \BC^{\times}_{\alpha} \rightarrow 1
\]
 where $\BC^{\times}_{\alpha} = \BC^{\times}$ and for any $\alpha \in \Delta$, we take
\[
Z(\widehat{M^0}) \rightarrow \BC^{\times}_{\alpha}
\]
to be the restriction of the simple root of $(\widehat{G^0}, \widehat{T^0})$ corresponding to coroot $\alpha^{\vee}$. By taking $\Gamma_{F}$-invariants, we get
\[
1 \rightarrow Z(\widehat{G^0})^{\Gamma_{F}} \rightarrow Z(\widehat{M^0})^{\Gamma_{F}} \rightarrow \left(\prod_{\alpha \in \Delta \backslash \Delta^{M^0}} \BC^{\times}_{\alpha}\right)^{\Gamma_{F}}. 
\]
Since the $\Gamma_{F}$-action only permutes $\BC^{\times}_{\alpha}$, then $(\prod_{\alpha \in \Delta \backslash \Delta^{M^0}} \BC^{\times}_{\alpha})^{\Gamma_{F}}$ is connected. Hence the map \[
Z(\widehat{M^0})^{\Gamma_{F}} \rightarrow \left(\prod_{\alpha \in \Delta \backslash \Delta^{M^0}} \BC^{\times}_{\alpha}\right)^{\Gamma_{F}}
\] is surjective.
It follows that $Z(\widehat{M^0})^{\Gamma_{F}}/Z(\widehat{G^0})^{\Gamma_{F}}$ is connected. Then 
\[
(Z(\widehat{M^0})^{\Gamma_{F}})^0 \rightarrow Z(\widehat{M^0})^{\Gamma_{F}}/Z(\widehat{G^0})^{\Gamma_{F}}
\]
is surjective, which implies that $Z(\dualgroup{M^0})^{\Gamma_{F}} = (Z(\dualgroup{M^0})^{\Gamma_{F}})^0 Z(\dualgroup{G^0})^{\Gamma_{F}}$.
\end{proof}

Let $M^0(\delta^0)$ be the standard Levi subgroup of $G^0$ determined by $\Delta_{\delta^0}$. Since $M^0_{\delta^0, \min}$ is fixed by $\theta_{\delta^0}$, then $M^0(\delta^0)$ is fixed by $\theta$.

\begin{lemma}
\label{lemma: pure inner forms of Levi}
There exists a unique $[\delta^0_{M^0(\delta^0)}] \in H^1(F, M^0(\delta^0))$  such that its image in $H^1(F, G^0)$ is  $[\delta^0]$. 
\end{lemma}

\begin{proof}
Let $S$ be the maximal central split torus of $M^0_{\delta^0, \min}$ and $T^0_{\delta^0}$ a maximal torus of $G^0_{\delta^0}$ containing $S$. By definition, $G^0_{\delta^0}(\overline{F}) = G^0(\overline{F})$. Since $G^0$ is quasisplit, there exists $g \in G^0(\overline{F})$ such that $T^0(\delta^0) := gT^0_{\delta^0}g^{-1}$ is a maximal torus of $G^0$ over $F$ and the isomorphism 
\[
T^0_{\delta^0} \rightarrow T^0(\delta^0)
\]
induced by the conjugation is defined over $F$ (cf. \cite{KS99}*{Lemma 3.3B}). 

Let $\delta^0_1(\sigma) := g\delta^0(\sigma)\sigma(g^{-1})$. Then $\delta^0_1$ is cohomologous to $\delta^0$ and it induces an isomorphism $G^0_{\delta^0} \cong G^0_{\delta^0_1}$ over $F$, which sends $S$ (resp. $T^0_{\delta^0}$) to $S_1$ (resp. $T^0_{\delta^0_1}$) in $G^0_{\delta^0_1}$. Under $G^0_{\delta^0_1}(\overline{F}) = G^0(\overline{F})$, $T^0_{\delta^0_1} = T^0(\delta^0)$ as $F$-tori. It follows that $\delta^0_1(\sigma) \in T^0(\delta^0)(\overline{F})$, i.e., $\delta^0_1$ lies in the image of 
\[
Z^1(F, T^0(\delta^0)) \rightarrow Z^1(F, G^0).
\]
Let $M^0_{\delta^0_1, \min} = Z_{G^0_{\delta^0_1}}(S_1)$. By viewing $S_1$ as a subtorus of $T^0(\delta^0)$, we also define $M^0(\delta^0_1) = Z_{G^0}(S_1)$. Then $M^0_{\delta^0_1, \min}(\overline{F}) = M^0(\delta^0_1)(\overline{F})$. Since $\delta^0_1$ lies in the image of 
\[
Z^1(F, M^0(\delta^0_1)) \rightarrow Z^1(F, G^0),
\]
let us denote its preimage in $Z^1(F, M^0(\delta^0_1))$ by $\delta^0_{M^0(\delta^0_1)}$. In this way, $M^0_{\delta^0_1, \min}$ is associated with the pure inner form $\delta^0_{M^0(\delta^0_1)}$ of $M^0(\delta^0_1)$. Finally, we can conjugate $M^0(\delta^0_1)$ to $M^0(\delta^0)$ and $\delta^0_{M^0(\delta^0_1)}$ to $\delta^0_{M^0(\delta^0)}$. This proves the existence, and the uniqueness follows from Lemma~\ref{lemma: injection of pure inner forms of Levi}. 
\end{proof}

\begin{remark}
\label{remark: pure inner forms of Levi}
We can associate $M^0_{\delta^0, \min}$ with a pure inner form $\delta^0_{M^0(\delta^0)}$ of $M^0(\delta^0)$. 
\end{remark}

Next, we move to the disconnected setting. Henceforth, we assume that $\pi_0(G_{\delta})$ is a cyclic group of prime order $\ell$. 
Fix a Borel subgroup $B^0_{\delta^0}$ and a maximal torus $T^0_{\delta^0}$ of $G^0$ over $\overline{F}$ such that
$
P_{\delta^0,\min}^0 \supseteq B^0_{\delta^0} \supseteq T^0_{\delta^0}.$
Let $\Delta(B^0_{\delta^0},T^0_{\delta^0})$ denote the set of simple roots associated to the pair $(B^0_{\delta^0},T^0_{\delta^0})$.
For $\tau \in \Gamma_F$ or $\tau \in \langle \theta_{\delta^0} \rangle$, the action of $\tau$ on $\Delta(B^0_{\delta^0},T^0_{\delta^0})$ is defined as follows. There exists an element $g_{\tau} \in G_{\delta^0}^0(\overline{F})$ such that
$
g_{\tau}\bigl(\tau(B^0_{\delta^0}),\tau(T^0_{\delta^0})\bigr)g_{\tau}^{-1}
=
(B^0_{\delta^0},T^0_{\delta^0}),
$ so $\widetilde{\tau}:= \operatorname{Int}(g_{\tau}) \circ \tau$ preserves $(B^0_{\delta^0},T^0_{\delta^0})$. Since $\widetilde{\tau}$ preserves $T^0_{\delta^0}$, it acts on the roots. Since it preserves $b_{\delta^0}$, it permutes the corresponding positive roots, hence permuting simple roots $\Delta(B^0_{\delta^0},T^0_{\delta^0})$. We define the action of  $\tau$ on $\Delta(B^0_{\delta^0},T^0_{\delta^0})$ to be this permutation given by $\widetilde{\tau}$. 
Now choose an element in $G^0(\overline{F})$ conjugating $(B^0_{\delta^0},T^0_{\delta^0})$ to $(B,T)$ and identify the simple roots. Then the $\delta$-twisted $\Gamma_F$-action and the $\theta_{\delta^0}$-action on $\Delta(B^0_{\delta^0},T^0_{\delta^0})$ coincide with the usual $\Gamma_F$-action and $\theta$-action on $\Delta$. We shall make this identification in the following discussion.

By the choice of $\theta_{\delta^0}$ and $P_{\delta^0,\min}^0$ in Definition \ref{def: theta delta}, for any $\Gamma_F$-invariant subset $\Delta_{\delta^0} \subseteq S \subseteq \Delta$, we have 
\[  \theta_{\delta^0}(P_{S}^0) = P_{\theta(S)}^0.\]
Following \cite{BJ01}, we define 
\[C({S})=\{c\in\langle \theta\rangle \ | \ c{S}={S}\}.\]
For any subgroup $D$ of $C(S)$, we define a parabolic subgroup and a Levi subgroup by
\[
P_{{S},D}:=\begin{cases}
P_S^0 & \text{ if }  D=\{1\}, \\
\langle P_{S}^0, \theta_{\delta}\rangle & \text{ otherwise,}
\end{cases} \ \ \ \ \ \ \ \ \ 
M_{{S},D}:=\begin{cases}
M_S^0 & \text{ if }  D=\{1\}, \\
\langle M_{S}^0, \theta_{\delta} \rangle & \text{ otherwise.}
\end{cases}
\]
Note that $P_{S,D}=M_{S,D} N_{S}^0$. A subgroup $P$ of $G_\delta$ is said to be a {\it parabolic subgroup} of $G_\delta$\IN{parabolic subgroup of $G_\delta$} if it is conjugate (by $G_\delta$) to some $P_{{S},D}$. 
The group $\langle\theta\rangle$ acts on the collection of pairs $\{({S},D)\}$. Following \cite{BJ01}*{\S3}, we fix a total order on $ \Delta$ from now on, which induces a lexicographic order on the subsets of $\Delta$. Then, for each orbit of pairs $({S},D)$, we fix the representative that is maximal with respect to the lexicographic order.
We say a parabolic subgroup (resp. Levi subgroup) $P_{{S},D}$ (resp. $M_{{S}, D}$) is {\it standard}\IN{standard parabolic subgroup} if $({S},D)$ is the representative we chose for its orbit. By Lemma~\ref{lemma: pure inner forms of Levi}, we can associate $M_{S, D}$ with a pure inner form of the corresponding standard Levi subgroups of $G$ determined by $S, D$.

\begin{example}
Here we present standard parabolic subgroups for the split group $G^0=\SO_4$ with outer conjugation $
    \theta$ given by $c\in \mathrm{O}_4$, an element with $\det(c)=-1.$
    Then $G=\SO_4\rtimes\langle \theta\rangle\simeq \mathrm{O}_4$. 
    Label the set of simple roots by $\Delta=\{\alpha_1,\alpha_2\}$. Since $G^0$ is split, any subset of $\Delta$ is $\Gamma_F$-stable. However, the $\langle \theta\rangle$-orbit of $\{\alpha_1\}$ is $\{\{\alpha_1\},\{\alpha_2\}\}$. Thus, while $\SO_4$ has 4 standard parabolic subgroups $\mathrm{O}_4$ has 5 standard parabolic subgroups according to the following table (where we choose an order on $\Delta$ via $\alpha_1>\alpha_2$).
\begin{center}
    \begin{tabular}{|c|c|}
    \hline
  $({S}, D)$  & $\mathrm{Standard \ Levi \ subgroup}$  \\
    \hline $(\Delta, \langle \theta \rangle)$    &  $\mathrm{O}_4\simeq \SO_4\rtimes \langle \theta \rangle$ \\
    \hline $(\Delta, \{1\})$    & $\SO_4$   \\
    \hline $(\{\alpha_1\}, \{1\})$    & $\GL_2$ \\
    \hline $(\emptyset,\langle \theta \rangle)$    & $\GL_1\times\GL_1\rtimes\langle \theta \rangle$ \\
    \hline $(\emptyset, \{1\})$    & $\GL_1\times\GL_1$ \\
    \hline
    \end{tabular}
    \end{center}
The nontrivial cohomology class in $H^1(F,G^0)$ has a $\theta$-fixed representative $\delta^0 \in Z^1(F,G^0)$. Let $\delta$ be the corresponding pure inner form of $G$. Then $G_\delta= G^0_{\delta^0} \rtimes \langle \theta \rangle$. The minimal parabolic $F$-subgroup of $G^0_{\delta^0}$ is $G^0_{\delta^0}$. Thus $\Delta_{\delta^0}=\Delta$. We have the following standard parabolic subgroups.
\begin{center}
    \begin{tabular}{|c|c|}
    \hline
  $({S}, D)$  & $\mathrm{Standard \ Levi \ subgroup}$  \\
    \hline $(\Delta, \langle \theta \rangle)$    &  $G_\delta\simeq G^0_{\delta^0}\rtimes \langle \theta \rangle$ \\
    \hline $(\Delta, \{1\})$    & $G^0_{\delta^0}$   \\
    \hline 
    \end{tabular}
    \end{center}
\end{example}

\begin{example} 
Here we present standard parabolic subgroups for the non-split, quasisplit group $G^0=\SO_4^\alpha$, where $\alpha$ denotes a choice of the nontrivial class in $F^\times/(F^\times)^2$. Again, $
    \theta$ is the outer conjugation given by $c\in \mathrm{O}_4^\alpha$, an element with $\det(c)=-1.$
    We still have $G=\SO_4^\alpha\rtimes\langle \theta\rangle\simeq \mathrm{O}_4^\alpha$. Label the set of simple roots by $\Delta=\{\alpha_1,\alpha_2\}$. Unlike the previous example, $G^0$ is nonsplit and the only $\Gamma_F$-stable subsets of $\Delta$ are $\emptyset$ and $\Delta$. The standard parabolic subgroups $\mathrm{O}_4^\alpha$ are given in the following table.
\begin{center}
    \begin{tabular}{|c|c|}
    \hline
  $({S}, D)$  & $\mathrm{Standard \ Levi \ subgroup}$  \\
    \hline $(\Delta, \langle \theta \rangle)$    &  $\mathrm{O}_4^\alpha\simeq \SO_4^\alpha\rtimes \langle \theta \rangle$ \\
    \hline $(\Delta, \{1\})$    & $\SO_4^\alpha$   \\
    \hline $(\emptyset,\langle \theta \rangle)$    & $\GL_1\times\mathrm{O}_2^\alpha\simeq\GL_1\times\SO_2^\alpha\rtimes\langle \theta \rangle$ \\
    \hline $(\emptyset, \{1\})$    & $\GL_1\times\SO_2^\alpha$ \\
    \hline
    \end{tabular}
    \end{center}
The case for $G_\delta$, where $[\delta^0]\in H^1(F,G^0)$ denotes the nontrivial cohomology class is similar.
\end{example}

\section{Standard representations and 
the Langlands classification theorem}\label{section: classification theorem}

In this section, we state the Langlands classification theorem for the reductive groups $G_\delta$ that appeared in Section~\ref{subsection: pure inner forms}.
This was done in \cite{BJ01} for trivial pure inner form $\delta$.

\subsection{\texorpdfstring{Notations for representations}{}}\label{subsection: restriction induction}
Let $P_{S,D}$ be a standard parabolic subgroup of $G_{\delta}$ with the standard Levi subgroup $M_{S,D}$. We let $\Pi(M_{S,D})$ (resp. $\Pi^{\textrm{fl}}(M_{S,D})$) denote the set of isomorphism classes of smooth irreducible (resp. smooth finite length) representations of $M_{S,D}(F)$. In this section, we introduce operations on these representations and describe relations between $\Pi(M_{S,D})$ and $\Pi(M_{S,1})$.
These will be needed in the following sections in order to state and prove the classification theorem.

Recall that $\theta_{\delta^0}(P_{S,1})= P_{\theta(S),1}$ and $\theta_{\delta^0}(M_{S,1})= M_{\theta(S),1}$. For $\pi \in \Pi(M_{S,1})$ and $\theta_{\delta^0}^k \in \langle \theta_{\delta^0} \rangle$, we
let ${}^{\theta_{\delta^0}^k}\pi \in \Pi(M_{\theta^k(S),1})$ be the representation acting on the same underlying space as $\pi$ by
\[ {}^{\theta_{\delta^0}^k}\pi(g):= \pi( \theta_{\delta^0}^k(g)).\]
Since $\theta_{\delta}^{\ell} \in M_{S,1}(F)$, then ${}^{\theta_{\delta^0}^{\ell}}\pi \simeq \pi$. 

Next, suppose $\theta(S)=S$. 
In this case, let 
\[ \Ind_{M_{S,1}}^{M_{S,\langle \theta\rangle}}: \Pi^{\textrm{fl}}(M_{S,1}) \to \Pi^{\textrm{fl}}(M_{S,\langle \theta \rangle}),\ \ \ \ \Res_{M_{S,1}}^{M_{S,\langle \theta\rangle}}: \Pi^{\textrm{fl}}(M_{S,\langle \theta \rangle}) \to \Pi^{\textrm{fl}}(M_{S,1})\]
denote the usual induction and restriction. Namely, for $\pi \in \Pi^{\textrm{fl}}(M_{S,1}) $, let $V$ denote the underlying space of $\pi$. Then the underlying space of $\Ind_{M_{S,1}}^{M_{S,\langle \theta\rangle}}\pi$ is the space of functions
\[ W:=\{f:M_{S,\langle \theta\rangle}(F) \to V\ | \ f(m^0 m)= \pi(m^0) f(m) \text{ for all }m^0 \in M_{S,1}(F) \},\]
and $M_{S,\langle \theta \rangle}$ acts by right translation. It is clear that each function $f$ in $W$ is determined by the collection of values $\{f( \theta_{\delta}^k )\}_{k=0,\ldots, \ell-1}$. Also, 
for $0 \leq i\leq \ell-1$, the right translation action of $M_{S,1}(F)$ leaves the subspace 
\[ W_i:=\{f \in W\ | \ f(\theta_{\delta}^k)=0 \text{ for all }k \in \{0,\ldots,\ell-1\}\setminus \{i\}  \} \]
invariant, and the linear transformation $f \mapsto f(\theta_{\delta}^i)$ gives an isomorphism between the  $M_{S,1}(F)$-action on $W_i$ and ${}^{\theta_{\delta^0}^i}\pi$. Therefore, we have
\begin{align}\label{eq finite index restriction}
   \Res_{M_{S,1}}^{M_{S,\langle \theta \rangle}} \circ \Ind_{M_{S,1}}^{M_{S,\langle \theta \rangle}}\pi \simeq \bigoplus_{i=0}^{\ell-1} {}^{\theta_{\delta^0}^{i}}\pi.
\end{align}

Now we discuss when the induction is irreducible.

\begin{lemma}\label{lem finite index induction restriction}
   Suppose that $\theta(S)=S$.  Set\IN{$\Pi^{\theta_{\delta^0}}(M_{S,1})$}
\begin{align}
\label{eq Pi theta}
    \Pi^{\theta_{\delta^0}}(M_{S,1})
    := 
    \{
    \pi \in \Pi(M_{S,1})\ | \ {}^{\theta_{\delta^0}}\pi \simeq \pi
    \}. 
\end{align}
Then the following holds.
\begin{enumerate}
    \item We have $\Ind_{M_{S,1}}^{M_{S,\langle \theta \rangle}}\pi $ is reducible if and only if $\pi \in \Pi^{\theta_{\delta^0}}(M_{S,1})$.
    \item For $\pi \in \Pi^{\theta_{\delta^0}}(M_{S,1})$, take any irreducible summand $\pi_0$ of $\Ind_{M_{S,1}}^{M_{S,\langle \theta \rangle}}\pi$. Fix any non-trivial character $\vartheta: M_{S,\langle \theta \rangle}(F) \to M_{S,\langle \theta \rangle}(F)/M_{S,1}(F) \to \BC^{\times}$ and define $\pi_i:= \pi_0 \otimes \vartheta^i$ for $0 \leqslant i \leqslant \ell-1$. Then we have 
\begin{align}\label{eq finite index induction}
\Ind_{M_{S,1}}^{M_{S,\langle \theta \rangle}}\pi \simeq \bigoplus_{i=0}^{\ell-1} \pi_i.\end{align}
\item For $\pi \in \Pi(M_{S,1})$ and $m \in M_{S,\langle \theta \rangle}(F) \setminus M_{S,1}(F)$, we have 
\[ \mathrm{trace}( m, \Ind_{M_{S,1}}^{M_{S,\langle \theta \rangle}}\pi )=0.\]
\end{enumerate}

\end{lemma}
\begin{proof}
    Suppose that $\pi $ is not in $ \Pi^{\theta_{\delta^0}}(M_{S,1})$. Let $W' \subseteq W$ be a nonzero $M_{S,\langle \theta \rangle}(F)$-invariant subspace. We show that $W'=W$. Indeed, by \eqref{eq finite index restriction}, we may find a subspace $ W'' \subseteq W'$ on which the $M_{S,1}(F)$-action is isomorphic to ${}^{\theta_{\delta^0}^{i}}\pi$ for some $i$. Let $W_k'':= \theta_{\delta}^k \cdot W'' \subseteq W'$ . Then the $M_{S,1}(F)$-action on $W_k''$ is isomorphic to ${}^{\theta_{\delta^0}^{i+k}}\pi$. Thus, from \eqref{eq finite index restriction}, we see that $W'=W$.

    Conversely, suppose that $\pi \in \Pi^{\theta_{\delta^0}}(M_{S,1})$. We fix an isomorphism $T$ between  ${}^{\theta_{\delta^0}}\pi$ and $\pi$. Namely, $T$ is a nonzero linear transformation from $V$ to $V$ such that for all $v \in V$ and $m \in M_{S,1}(F)$, we have
\begin{align}
\label{eq: intertwining operator}
T\big(\pi(\theta_{\delta^0}(m)) v\big) = \pi(m) T(v)
\end{align}
and
\begin{align}
\label{eq: normalization of intertwining operator}
T^{-\ell} = \pi(\theta_{\delta}^{\ell}).
\end{align}
We consider the following subspaces of $W$
\[
W_i':= \{f \in W\ | \ f(\theta_{\delta}^k)=\vartheta(\theta_{\delta}^k)^i T^{-k}(f(1)) \ \text{ for } 0 \leqslant k \leqslant \ell - 1  \} \quad (0 \leqslant i \leqslant \ell - 1).
\]
One can easily check that all $W'_i$ are nonzero and they are linearly independent. By \eqref{eq: intertwining operator} and \eqref{eq: normalization of intertwining operator}, $W'_i$ is invariant under the action of $M_{S,\langle \theta \rangle}(F)$, and we denote this action by $\pi_i$. Moreover, the map
\[
W'_{i} \rightarrow V, \quad f \mapsto f(1)
\]
induces an isomorphism $\pi_i|_{M_{S,1}} \cong \pi$. Hence, $\pi_i$ is irreducible. This completes the proof of Part (1). From the description of $W_i'$, we have $\pi_{i} \simeq \pi_0 \otimes \vartheta^i$. By \eqref{eq finite index restriction}, $W = \bigoplus_{0 \leqslant i \leqslant \ell - 1} W'_i$. This proves Part (2).

For Part (3), we can let $m = m^0 \theta^k_\delta$ for $m^0 \in M_{S,1}(F)$ and $0 < k \leqslant \ell - 1$. Then it sends $W_i$ to $W_{i+k}$ by right translation, where we have extended the definition of $W_i$ to all $i \in \mathbb{Z}$ by periodicity $W_{i+\ell}=W_i$. Since $W_i \neq W_{i+k}$, then the right translation by $m$ must have zero trace. This completes the proof of Part (3) and concludes the proof of the lemma.
\end{proof}

Finally, since any $\sigma \in \Pi(M_{S,\langle \theta \rangle})$ appears in $ \Ind_{M_{S,1}}^{M_{S,\langle \theta \rangle}} \circ \Res_{M_{S,1}}^{M_{S,\langle \theta \rangle}} \sigma$, we have 
\begin{align}\label{eq parametrization of rep of disconnected group}
    \Pi(M_{S,\langle \theta \rangle})= \left(\bigsqcup_{ \pi \in \Pi^{\theta_{\delta^0}}(M_{S,1})} \{\pi_0 ,\ldots, \pi_{\ell-1}\}\right) \bigsqcup \left(\bigsqcup_{\pi \in \Pi(M_{S,1})\setminus \Pi^{\theta_{\delta^0}}(M_{S,1})} \{\Ind_{M_{S,1}}^{M_{S,\langle \theta \rangle}}\pi\}\right).
\end{align}

\subsection{Standard triples for representations}\label{subsection: standard representation}

In this subsection we state the Langlands classification theorem for irreducible representations of $G_\delta(F)$.
For trivial $[\delta^0]\in H^1(F,G^0)$, which is to say, for the group $G(F)$, this result follows from \cite{BJ01}. Thus, for the rest of the section we assume that $G_{\delta}(F)/G^0_{\delta}(F)$ is non-trivial.
We begin by recalling the ingredients needed to define standard triples for representations of $G_\delta(F)$.

Consider a standard parabolic subgroup $P_{{S},D}$ of $G_\delta$ with standard Levi subgroup $M_{{S},D}$. 
Consider $P^0_{S}=G^0_\delta\cap P_{{S},D}$ with Levi decomposition $P^0_{S}=M_{{S}}^0N_{{S}}^0.$
Let $A_{{S}}$ be the split component of $Z(M_{{S}}^0)$ and $X(M_{{S}}^0)_F$ be the group of $F$-rational characters of $M_{{S}}^0$. 
Further consider the real Lie algebra and its dual
\[
\mathfrak{a}_{{S}}=\Hom(X(M_{{S}}^0)_F,\R),\ \ \mathfrak{a}^*_{{S}}=X(M_{{S}}^0)_F\otimes\R.
\]
We have an action of $\langle\theta\rangle$ on $\mathfrak{a}^*_{{S}}$ which is induced from the action on $\Delta.$
Let $H_{M_S^0}: M_{S}^0 \to \mathfrak{a}_S$ be the homomorphism satisfying
\[ q^{\langle \chi, H_{M_S^0}(m) \rangle}= |\chi(m)|.\]
Then for $\nu \in \mathfrak{a}_S^{\ast}$,  we define the character $\exp \nu$ of the connected Levi subgroup $M_S^0$ by
\[ \exp \nu (m):= q^{\langle \nu, H_{M_S^0}(m) \rangle}.\]
Note that the image of $\exp \nu$ is contained in $\R_{>0}$. When $\theta(\nu)=\nu$, by abuse of notation, we also use $\exp \nu$ to denote its unique extension to $M_S$ whose value remains positive (any other extension differs by a twist by $\vartheta^j$).

We identify $\mathfrak{a}_{S}^*$ with a subspace of $\mathfrak{a}_{\Delta_{\delta^0}}^*$, where $\Delta_{\delta^0} \subseteq \Delta$ corresponds to the fixed minimal parabolic subgroup of $G_\delta^0$.
Consider the standard inner product $\langle-,-\rangle:\mathfrak{a}_{\Delta_{\delta^0}}^*\times\mathfrak{a}_{\Delta_{\delta^0}}^*\rightarrow\mathbb{R}.$
We define another inner product $\langle-,-\rangle_\delta:\mathfrak{a}_{\Delta_{\delta^0}}^*\times\mathfrak{a}_{\Delta_{\delta^0}}^*\rightarrow\mathbb{R}$ via
\[
\langle x,y\rangle_\delta=\sum_{\rho\in \langle \theta\rangle}\langle \rho \cdot x, \rho \cdot y\rangle. 
\]
The inner product $\langle-,-\rangle_\delta$ is invariant under $\langle \theta\rangle$. By restriction, we obtain a $C({S})$-invariant inner product $\langle-,-\rangle_\delta:\mathfrak{a}_{S}^*\times \mathfrak{a}_{S}^*\rightarrow\mathbb{R}$. 

Now, consider the following chamber inside $ \mathfrak{a}_{{S}}^*$:
\[
\mathfrak{a}_{{S},+}^*:=\{x\in\mathfrak{a}_{S}^* \ | \ \langle x, \alpha\rangle_\delta>0,\ \forall \alpha \in \Delta_{S}\}.
\]
Unlike the connected case, there may exist $\nu_1, \nu_2 \in \mathfrak{a}_{{S},+}^*$ that are permuted under $\langle \theta\rangle$. 
In this case, we recall that we have fixed a total order on $\Delta$, which induces a lexicographic order $\geq$ on $ \mathfrak{a}_{\Delta_{\delta^0}}^{\ast}$ (see the last paragraph of \cite{BJ01}*{\S 3}). We consider the subchamber
\[
\mathfrak{a}_{{S},+}^*(C({S})):=\{x\in\mathfrak{a}_{{S},+}^* \ | \ x\leq\rho\cdot x, \ \forall \rho\in C({S})\}.
\]

\begin{definition}\label{definition: standard triple for representation}
A {\it standard triple} for $G_\delta$\IN{$(P,\pi_t,\nu)$, standard triple for $G_\delta$} consists of $(P,\pi_t,\nu)$ where:
    \begin{enumerate}
        \item[(0)] 
        $P=MN$ is a standard parabolic subgroup of $G_\delta$ and $P^0=M^0_{S} N_S^0$;
        \item[(1)] $\nu\in\mathfrak{a}_{{S},+}^*(C(S))$;
        \item[(2)]  $M=M_{{S},C({S},\nu)}$, where $C({S},\nu):=\{\rho\in C({S}) \ | \ \rho \cdot\nu=\nu\}$;
        \item[(3)]
        $\pi_t$ is an irreducible {\it tempered representation} of $M(F)$, which is to say, all (equivalently, any one of) the irreducible components of $\Res^{M}_{M^0}(\pi)$ are tempered.    \end{enumerate}
A {\it standard representation}\IN{$\SR(\pi)$, $\SR(P,\pi_t,\nu)$, standard representation} of $G_\delta(F)$ is a representation of the form
\[
\SR(P,\pi_t,\nu):=\Ind_{P}^{G_\delta}(\exp\nu\otimes\pi_t)
\]
for some standard triple $(P,\pi_t,\nu)$ for $G_\delta$. Here $\Ind_{P}^{G_\delta}$ is normalized parabolic induction.
\end{definition}

Now we state the Langlands classification theorem for representations of $G_{\delta}(F)$.

\begin{theorem}\label{theorem: standard representations}
For the group $G_\delta(F)$:
    \begin{enumerate}
        \item 
        every standard representation $\SR(P,\pi_t,\nu)$ has a unique irreducible quotient;
        \item 
        any irreducible representation $\pi$ of $G_\delta(F)$ occurs as the quotient of some standard representation $\SR(P,\pi_t,\nu)$ of $G_\delta(F)$.
    \end{enumerate}
    In particular, there is a bijection between $\Pi(G_{\delta})$ and the set of standard triples for $G_\delta$.
    If $\pi$ is the unique quotient of $M(P,\pi_t,\nu)$, then we also write $\SR(\pi):=M(P,\pi_t,\nu)$.
\end{theorem}

When $\delta$ is trivial, this is exactly \cite{BJ01}*{Theorem 4.2}. In fact, remarks in \cite{BJ01}*{p. 242} suggest that the quasisplit assumption in \cite{BJ01} is mainly for convenience: there exists a representative of $\pi_0(G)$ that fixes a Borel subgroup of $G^0$. 
For a pure inner form $\delta$, we have fixed a representative $\theta_{\delta}$ that fixes a minimal $F$-parabolic subgroup $P_{\delta,\min}$ in Definition \ref{def: theta delta}. Then, Theorem \ref{theorem: standard representations} follows from the same argument in \cite{BJ01}, which we omit. 
Here we only recall the construction of the standard triple and standard representation for $\pi \in \Pi(G_{\delta})$ from that of an irreducible subquotient $\pi^0$ of $\Res_{G_{\delta}^0}^{G_{\delta}}(\pi)$.
To that end, let $\pi^0$ be an irreducible representation of $G^0_\delta(F)$ with standard triple $( P_{{S},1},  \sigma_t^0,  \nu )$. In particular, $\pi^0$ is the unique quotient of the standard representation 
\[
\SR(\pi^0)=\Ind_{P_{{S},1}}^{G^0_\delta} (\exp \nu \otimes \sigma_{t}^0).
\] 
Since $(\theta_{\delta} (P_{{S},1})= P_{\theta({S}),1}, {}^{\theta_{\delta}}\sigma_t^0, \theta(\nu))$ is also a standard triple of $G_{\delta}^0$, the exactness of conjugation by $\theta_{\delta}$ implies that it must be the standard triple for ${}^{\theta_{\delta}}\pi^0$. Therefore, uniqueness of the standard triple for the connected group $G_{\delta}^0$ implies that $\pi^0 \simeq {}^{\theta_{\delta}}{\pi}^0$ if and only if 
\[
( P_{{S},1}, \sigma_t^0,  \nu  )=( P_{\theta({S}),1} , {}^{\theta_{\delta}}\sigma_t^0, \theta ( \nu)).
\]

\begin{definition}\label{def: induction of standards to disconnected}
    Suppose $\pi^0 \in \Pi(G_{\delta}^0)$ corresponds to the standard triple $( P_{{S},1},  \sigma_t^0,  \nu )$. We assign a standard triple for each irreducible summand of $\Ind_{G^0_\delta}^{G_\delta}(\pi^0)$ case by case as follows.
\begin{enumerate}
 \item[Case 1.] 
 Suppose that $\pi^0\not\simeq {}^{\theta_{\delta}}{\pi^0}$and consider the irreducible representation $\pi:= \Ind_{G^0_{\delta}}^{G_{\delta}}(\pi^0)$. The standard representation $\SR(\pi)$ for $\pi$ from Theorem~\ref{theorem: standard representations} is given by
    \[ 
    \Ind_{G^0_\delta}^{G_\delta} ( \SR(\pi^0))= \Ind_{P_{{S},1}}^{G_\delta} (\exp \nu \otimes \sigma_t^0)=\Ind_{P_{\theta^k({S}),1}}^{G_\delta} (\exp \theta^k(\nu) \otimes {}^{\theta_{\delta}^k}\sigma_t^{0}),
    \]
    for any $0 \leq k \leq \ell-1$. 
    In this case, the standard triple is given case by case as follows.
    \begin{itemize}
        \item [(i)]
        If the set $ \{\theta^k({S}) \tq 0 \leq k \leq \ell-1\}$ is not a singleton, then there exists a unique $k'$ such that $P_{\theta^{k'}({S}),1}$ is standard. Then the standard triple for $\pi$ is 
        \[
        (P_{\theta^{k'}({S}),1},  {}^{\theta_{\delta}^{k'}}\sigma_t^0, \theta^{k'}(\nu)).
        \]
        \item [(ii)]  
        If $\theta ({S})={S}$ but $\theta(\nu) \neq \nu$, then we take $ k' $ such that $\theta^{k'}(\nu)$ is the unique minimal element of $\{ \theta^k(\nu) \tq 0 \leq k \leq \ell-1\}$ under the lexicographical order we fixed. Then the standard triple for $\pi$ is 
        \[
        (P_{\theta^{k'}({S}),1},  {}^{\theta_{\delta}^{k'}}\sigma_t^0, \theta^{k'}(\nu)).
        \]
        \item [(iii)] 
        If $\theta ({S})={S}$ and $\theta (\nu)=\nu$, then we must have ${}^{\theta_{\delta}}\sigma_t^0 \not\simeq \sigma_t^0$ since $\pi^0 \not\simeq {}^{\theta_{\delta}}{\pi}^0$. Then $\Ind_{M_{{S},1}}^{M_{{S}, \langle \theta \rangle}}(\sigma_t^0)$ is irreducible (and tempered) and the standard triple for $\pi$ is given by 
        \[
        (P_{{S},\langle \theta \rangle},  \Ind_{M_{{S},1}}^{M_{{S}, \langle \theta \rangle}}(\sigma_t^0) , \nu).
        \]
    \end{itemize}

\item [Case 2.] 
Suppose that $\pi^0 \simeq {}^{\theta_{\delta}}{\pi}^0$ and write $\Ind_{G^0_\delta}^{G_\delta}(\pi^0)=\bigoplus_{i=0}^{\ell -1} \pi_i$ where $\pi_{i}= \pi_j\otimes \vartheta^{i-j} $. Write $\Ind_{M_{{S},1}}^{M_{{S}, \langle \theta \rangle}}(\sigma_t^0)= \bigoplus_{i=0}^{\ell -1} \sigma_{i}$ similarly so 
    $\sigma_i = \sigma_j \otimes \vartheta^{i-j}$.
Re-labeling if necessary, we assume that $\pi_0$ is an irreducible subquotient of $\Ind_{P}^{G_\delta}(\sigma_0\otimes \vartheta^i)$. Then the standard representation of $\pi_i$ is $\SR(\pi_i)= \Ind_{P_{{S},\langle \theta \rangle}}^{G_\delta} ( \sigma_{i})$ and the standard triple for $\pi_i$ is 
    \[
    (P_{{S}, \langle \theta \rangle}, \sigma_i, \nu).
    \]
\end{enumerate}
\end{definition}

\begin{remark}\label{remark: BJ01}
The assumption that $\langle \theta\rangle$ is cyclic of prime order, from Definition~\ref{assu: groups we consider}, simplifies Definition~\ref{def: induction of standards to disconnected}, considerably. 
    In the notation of \cite{BJ01}*{\S4}, the subgroup $C_1$ is always trivial, which avoids the issue of choosing representative $d(\Phi,\nu_1)$ for each class $d \in C_2/C_1$ in $C_2$ in  \cite{BJ01}*{Proposition 4.5}.  
\end{remark}

\begin{remark}
    In Case 1 in Definition~\ref{def: induction of standards to disconnected}, $\pi\simeq \Ind_{G^0_{\delta}(F)}^{G_{\delta}(F)}({}^{\theta_{\delta}^k}\pi^0)$ for any $0 \leq k \leq \ell-1$. However, the assignment of the standard triple is independent of $k$.
\end{remark}

\subsection{Whittaker normalization}\label{sec Whittaker normalization}

In \eqref{eq parametrization of rep of disconnected group}, we have classified the representations of $G_\delta(F)$ in terms of those of $G^0_\delta(F)$.
In this section we choose a base point $\pi_0$ in the set $\{\pi_i \tq 0 \leq i \leq \ell -1\}$ for $\pi \in \Pi^{\theta_{\delta}}(G_{\delta}^0)$. This is the main difficulty to extending the local Langlands correspondence from $G^0_\delta$ to $G_\delta$. 

However, in the special case when the group $G_\delta$ is quasisplit and the standard representation $\SR(\pi^0)$ is generic, there is a canonical choice of $\pi^+ \in \{\pi_i \tq 0 \leq i \leq \ell-1 \}$ with respect to the fixed Whittaker datum $\mathfrak{w}$, known as the Whittaker normalization. We review the definition of the Whittaker normalization, following \cite{Art13}*{\S 2.2}. Recall that $\mathfrak{w}$ is the Whittaker datum given by the $\theta$-invariant $F$-pinning and the additive character of $F$ that we fixed in Section~\ref{subsection: assumptions on groups}.

Now let $\pi^0 \in \Pi^{\theta}(G^0)$. As in \eqref{eq finite index induction} and \eqref{eq finite index restriction}, we have
\[
\Ind_{G^0}^{G}(\pi^0)= \bigoplus_{i=0}^{\ell -1} \pi_i, \ \ \ \ \ \ \ \Res_{G^0}^{G}(\pi_i) = \pi^0.
\]
Thus, the representations $\{\pi_i\}_{0 \leq i \leq \ell-1}$ may be regarded as $\ell$ distinct extensions of $\pi^0$ from $G^0(F)$ to $G(F)$, or, equivalently, as $\ell$ different $\theta$-actions on the underlying space $V_{\pi^0}$ of $\pi^0$.
Let $T: V_{\pi^0} \to V_{\pi^0}$ be a nonzero intertwining operator realizing the isomorphism $\,^{\theta}\pi^0 \simeq \pi^0$. That is, for all $v \in V_{\pi^0}$ and $g \in G^0(F)$,
\[
T\big(\pi^0(\theta(g)) v\big) = \pi^0(g) T(v).
\]
The choice of $T$ is unique up to scaling by Schur's lemma, and there are $\ell$ choices for $T$ such that $T^{\ell} = \operatorname{id}$ on $V_{\pi^0}$ (if $T$ is one such, then so is $\zeta_\ell^{i} T$ for $1 \leq i \leq \ell-1$).

To define $\pi^+$ from $\pi$ when $\SR(\pi^0)$ is generic, let us begin by considering the case in which $\pi^0$ itself is generic. In this situation, there exists a nonzero Whittaker functional $\omega: V_{\pi^0} \rightarrow \BC$ with respect to the Whittaker datum $\mathfrak{w}$, unique up to scalar. The Whittaker datum $\mathfrak{w}$ is $\theta$-invariant, since so is the $F$-pinning. As a consequence, $\omega \circ T$ is again a Whittaker functional on $V_{\pi^0}$. Consequently, $\omega \circ T = c\, \omega$ for some $c \in \BC^\times$. We fix the normalization of $T$ so that
\begin{align}\label{eq Whittaker normalization intertwining}
    \omega \circ T = \omega.
\end{align}
Then we define {\it $\pi^+$}\IN{$\pi^+$} to be the representation of $G$ whose underlying space is
\[ \{f \in \Ind_{G^0}^{G}(\pi^0) \ | \ f( (1 \rtimes \theta)^k)= T^{-k}(f(1)) \ \text{ for all }k  \},\]
which is irreducible from the proof of Lemma \ref{lem finite index induction restriction}.

Now consider a general irreducible representation $\pi^0$ of $G^0(F)$ such that the standard representation $\SR(\pi^0)$ is generic and $\pi^0 \simeq {}^{\theta}\pi^0$. Using the notation from Definition \ref{def: induction of standards to disconnected}, we have
\[
\SR(\pi^0) = \Ind_{P_{{S},1}}^{G^0}(\sigma^0)
\qquad \text{and} \qquad
\SR(\pi_i) = \Ind_{P_{{S},\langle \theta\rangle}}^{G}(\sigma_i).
\]
 Let $\mathfrak{w}_M$ denote the Whittaker datum for $M^0(F) = M_{{S},1}(F)$ determined by the fixed $F$-pinning and the chosen additive character of $F$. Since $\SR(\pi^0)$ is generic with respect to $\mathfrak{w}$, it follows that $\sigma^0$ is generic with respect to $\mathfrak{w}_M$. Write 
\[\Ind_{M_{{S},1}}^{M_{{S}, \langle \theta \rangle}}(\sigma^0)= \bigoplus_{i=0}^{\ell -1} \sigma_{i}.\]
From the preceding case, we obtain a canonical choice $\sigma^+ \in \{\sigma_i \tq 0 \leq i \leq \ell-1\}$. We then define $\pi^+ \in \{\pi_i \tq 0 \leq i \leq \ell-1\}$ to be the unique irreducible quotient of $\Ind_{P_{{S},\langle \theta\rangle}}^{G}(\sigma^+)$.

Finally, if ${}^{\theta}{\pi^0} \not\simeq \pi^0$ then $\Ind_{G^0}^{G} (\pi^0)$ is irreducible, so we simply set $\pi^+ = \Ind_{G^0}^{G} (\pi^0)$ in this case.

\begin{remark}
Let $\Pi^{\mfr{w}-\textrm{Mgen}}(G^0)$ denote the set of irreducible representations of $G^0$ whose standard representations are generic with respect to $\mfr{w}$, and let $\Pi^{\mfr{w}-\textrm{Mgen}}(G)$ be the set of irreducible representations of $G$ that appear as irreducible subquotients of $\Ind_{G^0}^{G}(\pi^0)$ for some $\pi^0 \in \Pi^{\mfr{w}-\textrm{Mgen}}(G^0)$. 
Using Whittaker normalization, the local Langlands correspondence for $G^0$ extends naturally to assign enhanced $L$-parameters to the representations in $\Pi^{\mfr{w}-\textrm{Mgen}}(G)$. See $\mathbf{(Sha)}$ and $\mathbf{(LC)}$ in Desideratum~\ref{desideratum: LLC} below.

In general, $\Pi^{\mfr{w}-\textrm{Mgen}}(G)$ is only a proper subset of $\Pi(G)$. However, when $G^0 = \GL_n(F)$, we have $\Pi^{\mfr{w}-\textrm{Mgen}}(G^0) = \Pi(G^0)$ and $\Pi^{\mfr{w}-\textrm{Mgen}}(G) = \Pi(G)$. Thus, for $G = \GL_n(F) \rtimes \langle \theta \rangle$, the local Langlands correspondence is known.\end{remark}

\section{Langlands classification for parameters}\label{section: classification for parameters}

\subsection{L-groups, Langlands parameters and their enhancements}
\label{subsection: parameters}

In this section, we define $L$-groups,
Langlands parameters and enhanced Langlands parameters for the groups $G_\delta$ appearing in Section~\ref{subsection: pure inner forms}. The main result of Proposition~\ref{proposition: enhancement by purity} partition these according to pure inner forms.

First, we define $L$-groups for pure inner forms of group $G$ as in Section \ref{subsection: assumptions on groups}. Recall that we have fixed a generator $\theta$ of $\pi_0(G)$ and an $F$-pinning $(B,T,\{X_{\alpha}\})$ of $G^0$ invariant under $\theta$-action. This allows us to define an automorphism $\hat{\theta}$ on $\dualgroup{G^0}$ dual to $\theta$. (See \cite{Kal26}*{Section 4.5} for example.) Let $W_F$ be the Weil group of $F$ and $W'_F := W_F \times \SL_2(\BC)$ the Weil-Deligne group.

\begin{definition}\label{definition: L-group}
Suppose that $G$ satisfies Definition~\ref{assu: groups we consider}.
    The {\it dual group} for $G$\IN{$\dualgroup{G}$, dual group for $G$} is $\dualgroup{G} := \dualgroup{G^0}\rtimes\langle{\hat\theta}\rangle$.
    The {\it $L$-group} of $G$\IN{$\Lgroup{G}$, $L$-group of $G$} is defined as $\Lgroup{G} := \dualgroup{G}\rtimes W_F$, where $W_F$ acts on $\dualgroup{G}$ through the action of $W_F$ on $\dualgroup{G^0}$.
\end{definition}

\begin{remark}\label{remark: dual group and L group}\ 
    \begin{enumerate}
    \item Let $M$ be a standard Levi subgroup of $G$, as discussed in Section \ref{subsection: parabolic subgroups}. Then there is a natural embedding between the $L$-groups ${}^LM\to \Lgroup{G}$.
        \item Recall from condition (2) in Section \ref{subsection: assumptions on groups} that $\Gamma_F$ commutes with $\theta$. Hence, $W_F$ acts trivially on $\langle{\hat\theta}\rangle$, so we may write 
$
\Lgroup{G}  = \Lgroup{G^0}\rtimes\langle{\hat\theta}\rangle.
$
\item As a semi-direct product of linear algebraic groups, $\dualgroup{G}$ is also a linear algebraic group.
    \end{enumerate}
\end{remark}

Next, we define Langlands parameters, $L$-parameters, and several related notions.

\begin{definition}\label{definition: parameters}\ 
\begin{enumerate}
    \item  A {\it Langlands parameter for $G$}\IN{$\phi$, Langlands parameter for $G$} is an admissible homomorphism $\phi: W_F' \to \Lgroup{G}$ such that $\phi(w') = \phi^0(w')\rtimes 1_{\langle{\hat\theta}\rangle}$ for a Langlands parameter $\phi^0$ for $G^0$; note that $\phi^0$ uniquely determines $\phi$, and vise versa. 
  Thus, the following diagram commutes
\[
\begin{tikzcd}
W'_F \arrow{rr}{\phi} \arrow{dr}[swap]{\phi^0} && \Lgroup{G}\\
& \Lgroup{G^0} \arrow{ur} &
\end{tikzcd}
\]
where $\Lgroup{G^0} \to \Lgroup{G}$ is the natural embedding from Definition~\ref{definition: L-group}.
An {\it $L$-parameter for $G$}\IN{$[\phi]$, $L$-parameter for $G$} is the $\widehat{G}$-conjugacy class of a Langlands parameter; we will use the notation $[\phi]$ for the $L$-parameter defined by the Langlands parameter $\phi$.
We write $\Phi(G)$ for the set of $L$-parameters for $G$.

\item
We call a Langlands parameter an {\it infinitesimal parameter}\IN{$\lambda$, infinitesimal parameter for $G$} if it is trivial on $\SL_2(\BC)$. We regard this as a homomorphism from $W_F$ to $\Lgroup{G}$. We will use the letter $\lambda$ for an infinitesimal parameter. For a Langlands parameter $\phi$, we associate an infinitesimal parameter $\lambda_{\phi}$\IN{$\lambda_{\phi}$, infinitesimal parameter of a Langlands parameter $\phi$} by 
\[ \lambda_{\phi}(w):= \phi\left(w, \begin{pmatrix}
    |w|^{\half{1}} & \\ &|w|^{\half{-1}}
\end{pmatrix}\right).\]
Fixing an infinitesimal parameter $\lambda$, we define the set of $L$-parameters \IN{$\Phi_{\lambda}(G)$, the set of $L$-parameters with infinitesimal parameter $\lambda$} 
\[\Phi_{\lambda}(G):= \{[\phi] \in \Phi(G)\ | \ \lambda_{\phi}=\lambda\}. \] 
We will typically write $\phi$ for both a Langlands parameter and its equivalence class when it is clear in context.
\item For every Langlands parameter $\phi$ for $G$ we set
\begin{equation}\label{eq: disconnected component group langlands parameter}
A_\phi := Z_{\dualgroup{G}}(\phi)/Z_{\dualgroup{G}}(\phi)^0 = \pi_0(Z_{\dualgroup{G}}(\phi))
\end{equation}
 and call it the {\it component group for $\phi$} (or for $[\phi]$)\IN{$A_\phi$, component group for $\phi$}. 

 \item Let $\widehat{A_\phi}$ be the set of isomorphism classes of irreducible representations of $A_\phi$.\IN{$\widehat{A_\phi}$, isomorphism classes of irreducible representations of $A_\phi$}
An {\it enhanced $L$-parameter for $G$}\IN{$([\phi],\tau)$, $(\phi,\tau)$, enhanced $L$-parameter} is a pair $([\phi],\tau)$, where $[\phi]$ is an $L$-parameter for $G$ and $\tau\in \widehat{A_\phi}$.
We will often use the simplified notation $(\phi,\tau)$ for an enhanced $L$-parameter.
We write $\Phi^\mathfrak{e}(G)$ for the set of enhanced $L$-parameters for $G$. 
\end{enumerate}
\end{definition}

    By Definition~\ref{definition: parameters}, there is a natural bijection $\phi \mapsto \phi^0$ between Langlands parameters for $G$ and Langlands parameters for $G^0$; however, this does {\it not} define a bijection between $L$-parameters. Rather, $\phi \mapsto \phi^0$ defines the bijection
   \begin{equation}\label{equation: equivalence disconnected g}
        \Phi(G) 
        \cong
    \Phi(G^0)/\langle{\hat\theta}\rangle,
    \end{equation} 
  where ${\hat\theta}$ acts on Langlands parameters for $G^0$ by
${}^{\hat\theta}\phi^0 := {\hat\theta}\circ\phi^0$. Observe that ${\hat\theta}$ acts on Langlands parameters for $G$ by conjugation:
$
({}^{\hat\theta}\phi)(w') := (1\rtimes{\hat\theta})\phi(w')(1\rtimes{\hat\theta})^{-1}.
$
These actions are compatible in the following sense:
\[
({}^{\hat\theta}\phi)(w') = ({}^{\hat\theta}\phi^0)(w') \rtimes 1, 
\qquad \forall w'\in W'_F.
\]
Note that $\phi^0$ and ${}^{\hat\theta}\phi^0$ 
determine Langlands parameters for $G$ in the same $\widehat{G}$-conjugacy class.

Definition \ref{definition: parameters}(3) includes the connected case: for $\phi^0 \in \Phi(G^0)$, its component group is $A_{\phi^0} =  \pi_0(Z_{\dualgroup{G^0}}(\phi^0))$. To describe the relation between $A_{\phi}$ and $A_{\phi^0}$, we give a lemma on the relation between component groups of centralizers in $\dualgroup{G^0}$ and $\dualgroup{G}$.

\begin{lemma}\label{lem centralizer}
    Let $\Omega$ be a non-empty subset of $\dualgroup{G^0}$. We have the following short exact sequence for the centralizer subgroups
    \begin{align*}
        1 \to Z_{\dualgroup{G^0}}(\Omega) \to Z_{\dualgroup{G}}(\Omega) \to \mathrm{im}(Z_{\dualgroup{G}}(\Omega) \to \pi_0(\widehat{G}) ) \to 1,
    \end{align*}
    and a short exact sequence of their component groups
      \begin{align*}
        1 \to \pi_0(Z_{\dualgroup{G^0}}(\Omega)) \to \pi_0(Z_{\dualgroup{G}}(\Omega)) \to \mathrm{im}(Z_{\dualgroup{G}}(\Omega) \to \pi_0(\widehat{G}) ) \to 1.
    \end{align*}
\end{lemma}
\begin{proof}
    Since $Z_{\dualgroup{G^0}}(\Omega)= Z_{\dualgroup{G}}(\Omega) \cap \dualgroup{G^0}$, the first short exact sequence holds. The second short exact sequence follows from the observation that $Z_{\dualgroup{G^0}}(\Omega)^0= Z_{\dualgroup{G}}(\Omega)^0$ and the second isomorphism theorem.
\end{proof}

Now we apply Lemma~\ref{lem centralizer} in the case $\Omega= \textrm{im}(\phi)$.
\begin{corollary}\label{cor: component group  l pardisc}
For any Langlands parameter $\phi$ for $G$ there is a short exact sequence 
\[
1 \to A_{\phi^0} \to A_{\phi} \to \langle {{\hat\theta}} \rangle^{[\phi^0]} \to 1,
\]
where $\phi^0$ is the Langlands parameter for $G^0$ determined by $\phi$. 
Here, $[\phi^0]$ is the $L$-parameter for $G^0$, i.e., the $\dualgroup{G^0}$-equivalence class of $\phi^0$, and $\langle {\dualgroup{\theta}} \rangle^{[\phi^0]}$ is the stabilizer of $[\phi^0]$ in $\langle {{\hat\theta}} \rangle$. 
\end{corollary}
\begin{proof}
By Lemma \ref{lem centralizer}, it suffices to show that $\textrm{im}(Z_{\dualgroup{G}}(\phi) \to \pi_0(\widehat{G}) )$ is non-trivial if and only if so is $ \langle {{\hat\theta}} \rangle^{[\phi^0]}$, which is also equivalent to ${}^{{\hat\theta}}\phi^0 \sim \phi^0$ as Langlands parameters of $G^0$.

Suppose first that ${}^{{\hat\theta}}\phi^0 \sim \phi^0$. Then there exists a $g \in \dualgroup{G^0}$ such that $ g ({}^{\hat{\theta}}\phi^0 )g^{-1}= \phi^0$. Equivalently, $g \rtimes \hat{\theta} \in Z_{\dualgroup{G}}(\phi)$. Hence, we have 
\[\textrm{im}(Z_{\dualgroup{G}}(\phi) \to \pi_0(\widehat{G})) \simeq \Z/\ell\Z \simeq \langle {{\hat\theta}} \rangle^{[\phi^0]}.\]  Conversely, if $\textrm{im}(Z_{\dualgroup{G}}(\phi) \to \pi_0(\widehat{G}))$ is not trivial, then there exists a $g\rtimes \hat{\theta}^k \in Z_{\dualgroup{G}}(\phi)$ where $0<k <\ell$. This implies that ${}^{\hat{\theta}^k} \phi^0 \sim \phi^0 $. Hence, ${}^{\hat{\theta}} \phi^0 \sim \phi^0 $ as $\ell$ is a prime. This completes the proof of the corollary.
\end{proof}

Let 
\begin{align*}
    H^1(F,G^0) &\to \Hom(\pi_0(Z(\dualgroup{G^0})^{\Gamma_F}), \BC^{\times }),\\
    [\delta^0]& \mapsto \chi_{[\delta^0]}
\end{align*}
denote Kottwitz isomorphism for $G$. In the following lemma, we show that this bijection respects the $\theta$-action on $H^1(F,G^0)$ and $\hat{\theta}$-action on $ \Hom(\pi_0(Z(\dualgroup{G^0})^{\Gamma_F}), \BC^{\times })$.

\begin{lemma}\label{lemma: Kottwitz isomorphism}
     Suppose that $G$ satisfies Definition~\ref{assu: groups we consider}. Then, for any $[\delta^0] \in H^1(F,G^0)$, we have 
     \[ \chi_{\theta[\delta^0]}= {}^{\hat{\theta}} \chi_{[\delta^0]}.\]
\end{lemma}

\begin{proof}
    Let $\delta^0\in Z^1(F,G^0).$ Then $\theta[\delta^0]:=[\theta\circ\delta^0].$ 
    On the other hand, we have that $\hat{\theta}$ is an automorphism of $\dualgroup{G^0}$ which commutes with the action of $\Gamma_F.$ Thus, it restricts to an automorphism of $Z(\dualgroup{G^0})^{\Gamma_F}$ which then induces its action on $\pi_0(Z(\dualgroup{G^0})^{\Gamma_F}).$ We have that ${}^{\hat{\theta}} \chi_{[\delta^0]}:= \chi_{[\delta^0]}\circ \hat{\theta}.$ Let $\kappa_{G^0}: H^1(F,G^0) \to \Hom(\pi_0(Z(\dualgroup{G^0})^{\Gamma_F}),\BC^\times)$ denote Kottwitz isomorphism. We must show that $\kappa_{G^0}([\theta\circ\delta^0])=\chi_{[\delta^0]}\circ \hat{\theta}.$

    Consider the connected reductive group $G_\delta^0$ defined over $F$, which we obtain via the $\delta^0$-twisted Galois action. Then there exists $F$-maximal tori $S'$ of $G_\delta^0$ (with respect to the $\delta^0$-twisted action) and $S$ of $G^0$ (with respect to the non-twisted action) such that there exists an isomorphism $\phi:S'\rightarrow S$ defined over $F$ and given by $\phi(t)=g^{-1}tg$ for some $g\in G^0(\overline{F})$. We show that there exists $\delta_1^0\in Z^1(F,S)$ which is cohomologous with $\delta^0$ in $H^1(F,G^0)$ (via the map $H^1(F,S)\rightarrow H^1(F,G)$ induced from the inclusion).

    Indeed, since $\phi$ is defined over $F$, we have $\phi(\sigma_{\delta^0}(t))=\sigma(\phi(t)),$ for any $t\in S'(\overline{F})$ and $\sigma\in \Gamma_F.$ 
    Set $y:=\phi(t)\in S(\overline{F}).$
    From the definitions, we obtain
    \[
g^{-1}\delta^0(\sigma)\sigma(t)\delta^0(\sigma)^{-1}g=\sigma(y)
    \]
    and hence
    \[
    g^{-1}\delta^0(\sigma)\sigma(t)\delta^0(\sigma)^{-1}g=\sigma(y).
    \]
    Using that $t=gyg^{-1}$ and
    setting $\delta_1^0(\sigma):=g^{-1}\delta^0(\sigma)\sigma(g)$, we rewrite the above equation as
    \[
    \delta_1^0(\sigma)\sigma(y) \delta_1^0(\sigma)^{-1} = \sigma(y).
    \]
    Since $\phi$ is an isomorphism, varying $y\in S(\overline{F})$ shows that $\delta_1^0(\sigma)\in Z_{G^0}(S)=S$. Furthermore, the definition implies that $[\delta_1^0]=[\delta^0]$ in $H^1(F,G).$ 

    Functoriality for Kottwitz isomorphism implies that we have the following commutative diagram (more generally, for any $F$-torus $S$ of $G^0$)
    \[\begin{tikzcd}
        H^1(F,S) \ar[r, "\kappa_S"] \ar[d, "i"] & \Hom(\pi_0(Z(\dualgroup{S})^{\Gamma_F}),\BC^\times) \ar[d, "r"] \\
        H^1(F,G^0) \ar[r, "\kappa_{G^0}"]  & \Hom(\pi_0(Z(\dualgroup{G^0})^{\Gamma_F}),\BC^\times),
    \end{tikzcd}\]
    where $i$ denotes the homomorphism induced from the inclusion $S\rightarrow G^0$ and $r$ denotes the homomorphism given by restriction.

    Although $\theta$ preserves a maximal torus, it may fail to preserve $S.$ Still, we have that $S'':=\theta(S)$ is a maximal $F$-torus in $G^0$. Furthermore, we have $\delta_2^0:=\theta\circ \delta_1^0\in Z^1(F,S'')$ and $[\delta_2]=[\theta\circ \delta_1^0]=\theta[\delta_1^0]=\theta[\delta^0]$ in $H^1(F,G^0).$
     Next, we invoke functoriality of Kottwitz isomorphism (which for tori is the Tate-Nakayama duality). Specifically, we have a homomorphism of tori (obtain by restriction) $\theta:S\rightarrow S''$ which induces homomorphisms $\theta_*:H^1(F,S)\rightarrow H^1(F,S'')$ and $\hat{\theta}:\dualgroup{S''}\rightarrow\dualgroup{S}$ (obtained by restricting $\hat{\theta}$ from $\dualgroup{G^0}$). Functoriality for the duality implies that 
     \[
     \kappa_{S''}(\theta[\delta_1^0])= \chi_{[\delta_2^0]}=\chi_{[\delta_1^0]}\circ\hat{\theta}.
     \]
     Finally, since $[\delta_2]=[\theta\circ \delta_1^0]=\theta[\delta_1^0]=\theta[\delta^0]$ in $H^1(F,G^0)$, using the above commutative diagram shows that $\kappa_{G^0}([\theta\circ\delta^0])=\chi_{[\delta^0]}\circ \hat{\theta}.$ This completes the proof of the lemma.
\end{proof}

\begin{proposition}
    \label{proposition: enhancement by purity}
   Suppose that $G$ satisfies Definition~\ref{assu: groups we consider}. Let $\phi$ be a Langlands parameter for $G$. 
There is a canonical function
\[
\widehat{A_\phi} \to \mathrm{im}(H^1(F,G^0)\to H^1(F,G))
\]
defined by pullback along $Z(\dualgroup{G^0})^{\Gamma_F} \to Z_{\dualgroup{G}}(\phi)$ together with Kottwitz isomorphism.
This defines the partition
\[
\widehat{A_\phi}
=
\coprod_{[\delta]\in H^1(F,G^0)_{\langle\theta \rangle}}\widehat{A_\phi}[\delta],
\]
into the fibres of $\widehat{A_\phi} \to \mathrm{im}(H^1(F,G^0)\to H^1(F,G))$.
\end{proposition}

\begin{proof}

The injections
\[ Z(\dualgroup{G^0})^{\Gamma_F} \to Z_{\dualgroup{G^0}}(\phi) \to Z_{\dualgroup{G}}(\phi) \]
induce homomorphisms (not necessarily injective)
\[ \pi_0(Z(\dualgroup{G^0})^{\Gamma_F})  \to A_{\phi^0} \to  A_{\phi}. \]
Let $\tau \in \widehat{A_{\phi}}$. Let $\chi$ be an irreducible constituent of $\tau|_{\pi_0(Z(\widehat{G^0})^{\Gamma})}$, which is a character. We let $[\delta_{\chi}] \in H^1(F,G^0)$ be its image under Kottwitz isomorphism, and let $[\delta_{\chi}]_{\langle \theta \rangle}$ be the image of $[\delta_{\chi}]$ in the coinvariants $H^1(F,G^0)_{\langle \theta \rangle}$. We show that $[\delta_{\chi}]_{\langle \theta\rangle}$ is independent of the choice of the irreducible constituent $\chi$. By Lemma \ref{lemma: Kottwitz isomorphism}, it is equivalent to show that if $\chi_1, \chi_2$ are two irreducible constituents of $\tau|_{\pi_0(Z(\widehat{G^0})^{\Gamma})}$, then $\chi_1= {}^{\hat{\theta}^k}\chi_2$ for some $k$.

First, we pull back $\tau$ along $ A_{\phi^0} \to A_{\phi}$. Fix an irreducible subquotient $\tau^0$ of $\tau|_{A_{\phi^0}}$.
By Corollary \ref{cor: component group  l pardisc}, either $A_{\phi^0}=A_{\phi}$ or $A_{\phi^0}$ is a subgroup of $A_{\phi}$ of index $\ell$. Suppose that $\tau|_{A_{\phi^0}}$ is not irreducible; then by the same discussion in Section \ref{subsection: restriction induction}, there exists a $g \rtimes \hat{\theta} \in Z_{\widehat{G}}(\phi)$ such that
\[ \tau|_{A_{\phi^0}}= \bigoplus_{i=0}^{\ell-1} {}^{(\overline{g \rtimes \hat{\theta}})^i}\tau^0.\]
Here $\overline{g \rtimes \hat{\theta}}$ is the image of $g \rtimes \hat{\theta}$ in $A_{\phi}$.

Next, we pull back along $\pi_0(Z(\dualgroup{G^0})^{\Gamma}) \to A_{\phi^0}$.
Since $Z(\dualgroup{G^0})^{\Gamma_F}$ is contained in the center of $Z_{\dualgroup{G^0}}(\phi)$, the image of $\pi_0(Z(\dualgroup{G^0})^{\Gamma_F})$ is contained in the center of $A_{\phi^0}$. Thus, we have that $\tau^0|_{\pi_0(Z(\dualgroup{G^0})^{\Gamma_F})}$ is an isotypic direct sum of the pullback of the central character of $\tau^0$. Also, since $g \rtimes 1 \in \dualgroup{G^0}$, we have 
\begin{align*}
  (  {}^{(\overline{g \rtimes \hat{\theta}})^i}\tau^0)|_{\pi_0(Z(\dualgroup{G^0})^{\Gamma})}= {}^{\hat{\theta}^i}(\tau^0|_{\pi_0(Z(\dualgroup{G^0})^{\Gamma})})
\end{align*}
This completes the proof of the proposition.
\end{proof}

With Proposition \ref{proposition: enhancement by purity}, we define $L$-parameters for pure inner form $G_{\delta}$.

\begin{definition}\label{def L-par relevance}
    Let $\delta$ be a pure inner form of $G$. We define \IN{$\Phi(G_{\delta})$, $L$-parameters of $G_{\delta}$} \IN{$\Phi^{\mfr{e}}(G,\delta)$, enhanced $L$-parameters of $G_{\delta}$}
    \[ \Phi(G_{\delta}):= \{\phi \in \Phi(G) \ | \ \widehat{A_{\phi}}[\delta] \neq \emptyset \},\ \ \ \ \Phi^{\mfr{e}}(G,\delta):= \{(\phi,\varepsilon) \in \Phi^{\mfr{e}}(G) \ | \ \varepsilon \in \widehat{A_{\phi}}[\delta] \}.\]
\end{definition}
Note that $\Phi^{\mfr{e}}(G,\delta)$ depends on both $G_{\delta}$ and the class $[\delta] \in H^1(F,G)$, while $\Phi(G_{\delta})$ only depends on the algebraic group $G_{\delta}$ by the following remark.

\begin{remark}\label{remark: relevant parameters}
An $L$-parameter $\phi^0$ of $G^0$ is called $G_{\delta}^0$-relevant if every parabolic subgroup ${}^L P$ that $\phi^0$ factors through is conjugate to an $L$-group of some standard $F$-parabolic subgroup of $G_{\delta}$. See, e.g., \cite{Bor79} for more details. It is known, for instance by \cite{Kal16}*{Lemma 5.7}, that a Langlands parameter \(\phi^0\) for \(G^0\) is \(G^0_{\delta}\)-relevant if and only if \(\widehat{A_{\phi^0}}[\delta]\) is nonempty. From the proof of Proposition \ref{proposition: enhancement by purity}, we see that \(\widehat{A_{\phi}}[\delta]\) is nonempty if and only if \(\widehat{A_{\phi^0}}[\delta^0]\) is nonempty. Therefore, we may define that \(\phi \in \Phi(G)\) is \(G_{\delta}\)-relevant if and only if \(\phi^0\) is \(G_{\delta}^0\)-relevant. Then \(\Phi(G_{\delta})\) is precisely the subset of \(\Phi(G)\) consisting of those parameters that are \(G_{\delta}\)-relevant.
\end{remark}

If $\theta[\delta^0] \neq [\delta^0]$, then $G_{\delta}(F)=G_{\delta^0}^0(F)$. From the viewpoint of the enhanced local Langlands correspondence (see Section \ref{section: LLC}), it is natural to expect that $\widehat{A_{\phi^0}}[\delta^0]$ is in bijection with $ \widehat{A_{\phi}}[\delta]$ for any $\phi^0 \in \Phi(G^0)$. We prove in the following corollary that induction gives the desired bijection.
\begin{corollary}\label{cor component group theta[delta]neq [delta]}
    Let $\delta$ be a pure inner form of $G$. Suppose that $\theta[\delta^0] \neq [\delta^0]$. Then for all $\phi^0 \in \Phi(G^0)$, the induction $\Ind_{A_{\phi^0}}^{A_{\phi}}$ gives a bijection $\widehat{A_{\phi^0}}[\delta^0]$ between $ \widehat{A_{\phi}}[\delta]$. Moreover, for any $\tau \in \widehat{A}_{\phi}[\delta]$ and $s \in A_{\phi} \setminus A_{\phi^0}$, we have $\textrm{trace}(s, \tau)=0$.
\end{corollary}
\begin{proof}
    By Corollary \ref{cor: component group  l pardisc}, we see that $A_{\phi^0}$ is always a subgroup of $A_{\phi}$, so the induction makes sense. If $[A_{\phi}:A_{\phi^0}]=1$, then the statement trivially holds. If $[A_{\phi}:A_{\phi^0}]=\ell$ and $\tau^0 \in \widehat{A_{\phi^0}}[\delta^0]$, then ${}^{\hat{\theta}^i}\tau^0 \in \widehat{A_{\phi^0}}[\theta^i(\delta^0)]$ by Lemma \ref{lemma: Kottwitz isomorphism}. In particular, ${}^{\hat{\theta}}\tau^0 \neq \tau^0$, and hence $ \Ind_{A_{\phi^0}}^{A_{\phi}} \tau^0=: \tau$ is irreducible and lies in $\widehat{A_{\phi}}[\delta]$. By the finite group version for \eqref{eq finite index restriction}, we have 
    \[ \tau|_{A_{\phi^0}}= \bigoplus_{i=0}^{\ell-1} {}^{\hat{\theta}^i}\tau^0. \]
    This shows the injectivity.   
    To show the surjectivity, suppose that $\tau \in \widehat{A_{\phi}}[\delta]$. Then there exists a subquotient $\tau^0$ of $ \tau|_{A_{\phi^0}}$ such that $\tau^0 \in \widehat{A_{\phi^0}}[ \theta^k (\delta^0)]$ for some $k$. Then, ${}^{\hat{\theta}^{-k}}\tau^0 \in \widehat{A_{\phi^0}}[\delta^0]$ and $\tau= \Ind_{A_{\phi^0}}^{A_{\phi}} {}^{\hat{\theta}^{-k}}\tau^0$. $\Ind_{A_{\phi^0}}^{A_{\phi}}$ gives a bijection $\widehat{A_{\phi^0}}[\delta^0]$ between $ \widehat{A_{\phi}}[\delta]$.

    The second part of the statement is non-empty only when $[A_{\phi}: A_{\phi^0}]= \ell$. The same argument in Lemma \ref{lem finite index induction restriction}(3) implies that $\textrm{trace}(s, \tau)=0$.
    This completes the proof of the corollary.
\end{proof}

\subsection{Standard triples for parameters}\label{subsection: standard parameters}

In this section, we prove the Langlands classification theorem for $L$-parameters for disconnected groups.
The main result is Theorem~\ref{theorem: Langlands classification for L-parameters}, which establishes a bijection between $L$-parameters and standard triples for $L$-parameters for $G$.

\begin{definition}\label{definition: standard parameters}
    A {\it standard triple for  $\Phi(G)$}\IN{$(P,[\phi_t],\nu)$, standard triple for  $\Phi(G)$} consists of $(P,[\phi_t],\nu)$ where:
    \begin{enumerate}
        \item[(0)] $P=MN$ is a standard parabolic subgroup of $G$ with $P^0=M^0_{S} N_S^0$ for ${S}\subseteq \Delta$;
       \item[(1)] $\nu\in\mathfrak{a}_{{S},+}^*(C(S))$;
        \item[(2)]  $M=M_{{S},C({S},\nu)}$;
        \item[(3)]
        $[\phi_t]$ is a {\it tempered}
        $L$-parameter for $M$, which is to say, the image of the restriction of $\phi_t$ to $W_F$ is bounded. 
    \end{enumerate}
With reference to Section~\ref{subsection: parabolic subgroups}, this means $P = P_{{S},D}$ for $D = C({S},\nu)$.
\end{definition}

We may identify $\nu$ as an unramified cocharacter $\nu: W_F \to Z(\widehat{M_S^0})^{\Gamma}$ (see \cite{SZ18}*{Corollary 4.3}). We let $\phi_t \otimes \nu$ be the corresponding Langlands parameter that sends $(w,x)$ to $\phi_t(w,x)\nu(w)$. 

\begin{theorem}\label{theorem: Langlands classification for L-parameters}
    There is a bijection, called the Langlands classification, between $\Phi(G)$ and standard triples of $L$-parameters for $G$:
    \[
    [\phi]
    \leftrightarrow
    (P,[\phi_t],\nu),
    \]
given by the relation $[\phi]= [ ({}^LM_{S,C(S,\nu)} \hookrightarrow \Lgroup{G}) \circ (\phi_t \otimes \nu)]$. 
\end{theorem}

\begin{proof}
Let $\phi$ be a Langlands parameter for $G$.
 If $\theta$ is trivial, then $G_\delta$ is connected and Theorem~\ref{theorem: Langlands classification for L-parameters} is proved in \cite{SZ18}*{Section 4.6}.
    Consequently, we now assume that $\theta$ is nontrivial.
Consider the Langlands parameter $\phi^0$ for $G^0$ determined by $\phi$, using Definition~\ref{definition: parameters}.
Using \cite{SZ18}*{Section 4.6}, let $(P^0, [\phi^0_t], \nu)$ be the standard triple for $[\phi^0] \in \Phi(G^0)$. Note that we already have 
\[[\phi^0]= [ ({}^LM_S^0 \hookrightarrow \Lgroup{G}^0) \circ (\phi_t^0 \otimes \nu)]\]
as $L$-parameter of $G^0$. However, $(P^0, [\phi^0_t], \nu)$ may not be a standard triple according to Definition \ref{definition: standard parameters}. We are going to modify the triple $(P^0, [\phi^0_t], \nu)$ case by case. Observe that
    $({}^{\theta^i} P^0= P_{\theta^i(S),1}, [{}^{{\hat\theta}^i}\phi_t], {}^{\theta^i}\nu)$ is again a standard triple of $G^0$, so it must be the standard triple for $[{}^{{\hat\theta}^i}{\phi^0}]\in \Phi(G^0)$. Also,
    \[ [\phi]= [ ({}^LM_{S,1} \hookrightarrow \Lgroup{G}) \circ (\phi_t^0 \otimes \nu)]=  [ ({}^LM_{\theta^i(S),1} \hookrightarrow \Lgroup{G}) \circ ({}^{{\hat\theta}^i}\phi_t^0 \otimes ({}^{{\hat\theta}^i}\nu))] \]
    as $L$-parameters of $G$ for all $i=0,\ldots, \ell-1$.

    First, suppose that $\theta(S)\neq S$. Then there is exactly one $j\in\{0,\dots,\ell-1\}$ such that $P_{\theta^j({S}),1}$ is a standard parabolic subgroup of $G$. Note that $C({\theta^j(S)},\nu)$ is trivial. Thus, $( P_{\theta^j(S),1}, [{}^{{\hat\theta}^j}\phi_t^0], {}^{\theta^j}\nu) $ is a standard triple of $G$, and we assign it to the Langlands parameter $\phi$ in this case.

    Next, suppose that $\theta(S)=S$ and $^\theta \nu \neq \nu$. Similarly, there exists a unique $j\in\{0,1,\dots,\ell-1\}$ such that ${}^{\theta^j}\nu\geq {}^{\theta^i}\nu$ for any $i\in\{0,1,\dots,\ell-1\}$ under our fixed lexicographic order (see Section \ref{subsection: standard representation}). Namely,  $\{{}^{\theta^i}\nu \ | \ i=0,1,\dots,l-1\} \cap \mathfrak{a}_{S,+}^{\ast}(C(S))= \{{}^{\theta^j}\nu\}$. Note that $C({S},\nu)$ is trivial and so $P_{\theta^j(S),1}=P_{S,1}=P_{{S},C({S},\nu)}$. Again, the triple $( P_{\theta^j(S),1}, [{}^{{\hat\theta}^j}\phi_t^0], {}^{\theta^j}\nu) $ is a standard triple of $G$, and we assign it to the Langlands parameter $\phi$ in this case

    Finally, suppose that both $\theta(S)=S$ and ${}^\theta\nu= \nu$. Then, $C(S,\nu)= \langle \theta \rangle$. Let $[\phi_t]$ denote the $\dualgroup{G}$-conjugacy class of $\phi_t^0$. Then, $(P_{S,\langle \theta \rangle}, [\phi_t], \nu)$ is a standard triple, and we assign it to the Langlands parameter $\phi$ in this case.

  We have now assigned a standard triple of $\Phi(G)$ to every Langlands parameter $\phi$ of $G$. It is not hard to check that this assignment only depends on the $\dualgroup{G}$-conjugacy class of $\phi$, which we omit. Thus, we have obtained a well-defined function from $\Phi(G)$ to the set of standard triples. Then, it is straightforward from the above construction that this function is the inverse of the function 
  \[ (P_{S,C(s,\nu)}, [\phi_t],\nu) \mapsto [ ({}^LM_{S,C(S,\nu)} \hookrightarrow \Lgroup{G}) \circ (\phi_t \otimes \nu)], \]
 which we also omit. This completes the proof of the theorem.  
\end{proof}

\begin{remark}\label{remark: LC for L-parameters of pure inner forms}
    If $[\phi] \in \Phi(G_{\delta}) \subseteq \Phi(G)$ and $[\phi] \leftrightarrow (P_{S,D}, [\phi_t], \nu)$ under the Langlands classification, then $S \supseteq \Delta_{\delta^0}$ since $\phi$ factors through ${}^LP_{S,D}$ (see Remark \ref{remark: relevant parameters}). 
\end{remark}

\begin{lemma}\label{lemma: LC for enhanced parameters}
         If  $(P,[\phi_t],\nu)$ is the standard triple for $[\phi]\in \Phi(G)$ under Theorem \ref{theorem: Langlands classification for L-parameters}, then the map $ A_{\phi_t} \to A_\phi$ induced from the $L$-embedding is an isomorphism.
\end{lemma}

 \begin{proof}
Consider the Langlands parameter $\phi^0$ for $G^0$ determined by $\phi$ via Definition~\ref{definition: parameters}. Let $(P_{\phi^0}, [\phi_t^0], \nu_{\phi^0})$ be the standard triple for $[\phi^0]\in \Phi(G^0)$. If $P_{\phi}$ is connected, then this is proved in \cite{SZ18}*{Section 7.2}. Therefore, we assume that $P=P_{S,\langle \theta \rangle}$ is not connected, so $\theta(S)=S$ and $\theta(\nu)=\nu$. Let $\dualgroup{M}$ denote the Levi of the dual group of the Levi $M_{S,\langle\theta\rangle}$ of $P$. Consider the following commutative diagram. 
\[\begin{tikzcd}
    1 \ar[r]& A_{\phi^0\ar[r]}&\ar[r] A_{\phi} & \textrm{im}(Z_{\dualgroup{G}}(\phi^0) \to \pi_0(\dualgroup{G}))\ar[r] & 1\\
    1 \ar[r]& A_{\phi^0_t}\ar[r]\ar[u]& A_{\phi_t}\ar[r] \ar[u]& \textrm{im}(Z_{\dualgroup{M}}(\phi_t^0) \to \pi_0(\dualgroup{M}) \ar[r] \ar[u]& 1
\end{tikzcd}
\]
The left vertical arrow is an isomorphism by the connected case. We show that the right vertical arrow is always an isomorphism as well. Hence, $A_{\phi} \simeq A_{\phi_t}$. 

From the proof of Theorem \ref{theorem: Langlands classification for L-parameters}, we may assume that $(P_{S,1}, [\phi_t^0], \nu)$ is the standard triple of $\phi^0$, and the following equality of Langlands parameters holds (not up to $\dualgroup{G^0}$-conjugacy) 
\[ \phi^0= ({}^L M_{S,1} \hookrightarrow \Lgroup{G}^0) \circ (\phi_t^0 \otimes \nu).\]
We recall from the construction of standard triple of the connected case (\cite{SZ18}*{Proposition 6.1}) that $\nu$ corresponds to the hyperbolic part $z_{\phi^0}$ of $\phi^0(\textrm{Frob})$ (\cite{SZ18}*{Proposition 5.3}) and ${}^LM_{S,1}= Z_{\Lgroup{G}}(z_{\phi^0})$. Suppose that there exists an $m \rtimes \hat{\theta} \in Z_{\dualgroup{M}}(\phi_t^0)$. Then, since $\theta(\nu)=\nu$  and $m \in \dualgroup{M_{S,1}} \subseteq Z_{\Lgroup{G}}(z_{\phi^0})$, we see that ${}^{m \rtimes \theta} \nu= \nu$. Therefore, we have ${}^{m \rtimes \hat{\theta}} \phi^0= \phi^0$, and the right vertical arrow is an isomorphism between cyclic groups of order $\ell$. Conversely, suppose that there exists an $g \rtimes \hat{\theta} \in Z_{\dualgroup{G}}(\phi^0)$. Then $g\rtimes \hat{\theta}$ centralizes the hyperbolic part of $\phi^0$. Therefore, we have ${}^{g \rtimes \hat{\theta}} \nu= \nu$, and hence ${}^g \nu =\nu$ since $^{\hat{\theta}} \nu=\nu $. This implies that $g \in Z_{\Lgroup{G}}(z_{\phi^0})={}^LM_{S,1}$. Thus, in fact $g \rtimes \hat{\theta} \in Z_{\dualgroup{M}}(\phi_t^0)$. We again conclude that the right vertical arrow is an isomorphism between cyclic groups of order $\ell$. This completes the proof of the lemma.
\end{proof}

\begin{definition}\label{definition: standard triple for enhanced parameters}
    The {\it standard triple} for an enhanced $L$-parameter $(\phi,\tau)$ is the triple 
    \[
    (P_\phi,([\phi_t],\tau_t),\nu_\phi), 
    \]
    where $(P_\phi,[\phi_t],\nu_\phi)$ is the standard triple for $[\phi]\in \Phi(G)$ under Theorem \ref{theorem: Langlands classification for L-parameters} and $\tau= \xi \circ \tau_t \in \widehat{A}_{\phi_t}$ where $\xi: A_{\phi_t} \simeq A_{\phi}$ is the isomorphism induced from the $L$-embedding. \IN{$(P_\phi,([\phi_t],\tau_t),\nu_\phi)$, standard triple for an enhanced $L$-parameter}
\end{definition}

\part{Endoscopic transfer to disconnected groups}\label{part: II Transfer}

In this section, we define the endoscopic transfer from stable distributions of endoscopic groups $G'(F)$ to invariant distributions of disconnected group $G_\delta(F)$, where $\delta$ is a pure inner form of $G$ (Definition \ref{definition: pure inner form}).

\section{Endoscopic transfer}\label{section: endoscopic transfer}

We re-cast Kottwitz-Shelstad transfer in the language of certain disconnected groups. 
The main result is Proposition~\ref{proposition: induction commutes with transfer} which shows that endoscopic transfer commutes with induction, building on the well known result for Kottwitz-Shelstad transfer.

\subsection{Endoscopic triples}\label{subsection: endoscopic triple}
This section presents the class of twisted endoscopic data that we consider in this paper. 

\begin{definition}\label{definition: endoscopic triple}
Suppose that $G$ satisfies Definition~\ref{assu: groups we consider}.
An {\it endoscopic triple} $(G',s,\xi)$ for $G$\IN{$(G',s,\xi)$, endoscopic triple for $G$} consists of
\begin{enumerate}
    \item
    a connected quasisplit group $G'$ defined over $F$;
    \item 
    a semisimple element $s$ in $\dualgroup{G}$ (recall that $\dualgroup{G}$ is a linear algebraic group by Remark \ref{remark: dual group and L group}(3)) which is $\hat{\theta}$-semisimple if $s\in\dualgroup{G^0}$;
    \item 
    an $L$-embedding  $\xi:\Lgroup{G'}\rightarrow\Lgroup{G}$ such that $\xi$ restricted to $\dualgroup{G'}$ is an isomorphism to $Z_{\dualgroup{G^0}}({s})^0$, the identity component of the ${s}$-fixed points of $\dualgroup{G^0}$, and the image of $\xi$ lies in $Z_{\Lgroup{G}}(s)$.
\end{enumerate}
Any algebraic group $G'$ arising in this way is called an endoscopic group for $G$.
\end{definition}

\begin{remark}
    Note that the above definition considers only endoscopic groups attached to a quasisplit group $G$, but these also allow for endoscopic transfer from $G'$ to pure inner forms $G_\delta$, and 
        we allow $s \in \widehat{G^0}$, the identity component of $\dualgroup{G}$ but also $s\in \widehat{G}\setminus \widehat{G^0}$.
    This will be explained when we modify Kottwitz–Shelstad transfer in Definition~\ref{definition: modified transfer};
    see also Remark~\ref{rmk transfer disconnected}.
\end{remark}

\begin{remark}
    Suppose that $s=g_s\rtimes\hat{\theta}^{n_s}$ is semisimple and $n_s\neq 0$. Then $g_s$ is $\hat{\theta}^{n_s}$-semisimple in the sense of \cite{KS99}*{Definition 2.1.3} (and since the order of $\theta$ is prime,  $G^0\rtimes\langle\theta^{n_s}\rangle$ satisfies Definition~\ref{assu: groups we consider} and is equal to $G$). Indeed, if $s$ is semisimple then the restriction of $\mathrm{Int}(g_s)\circ\hat{\theta}^{n_s}$ to the derived group $(\dualgroup{G^0})_{\mathrm{der}}$ is semisimple, but this is exactly the definition of being $\hat{\theta}^{n_s}$-semisimple. 
\end{remark}

We provide some examples of endoscopic triples in Table~\ref{table: endosclassical}.
In Section~\ref{section: endoclassical} we explain how our results apply to some of these groups.
\begin{table}
\begin{center}
\caption{Examples of endoscopic triples $(G',s,\xi)$ for quasisplit $G = G^0\rtimes\langle{\theta}\rangle$, as in Definition~\ref{definition: endoscopic triple}.
}\label{table: endosclassical}
\scalebox{0.9}{
\begin{tabular}{ |l|l|l|l|l| l | }
 \hline
$G'$ & $ \Lgroup{G'} = \dualgroup{G'}\rtimes W_F $ & $\Lgroup{G} = (\Lgroup{G^0}\rtimes W_F)\rtimes\langle{\hat\theta}\rangle$ & $G = G^0\rtimes\langle{\theta}\rangle$ & $\hat{\theta}$ & $\ell$ \\ 
 \hline\hline
$\U_{n}$ & $\GL_{n} \rtimes W_F$ & $(\GL_{n}^{\times 2} \rtimes W_F)\rtimes\langle\hat{\theta}\rangle$ & $\Res_{E/F} \GL_{n}\rtimes\langle{\theta}\rangle$ & \S \ref{subsubsection: A}  & 2\\
\hline
$\SO_{2n+1}$ & $\Sp_{2n} \times W_F$ & $(\GL_{2n} \times W_F)\rtimes\langle\hat{\theta}\rangle$ & $\GL_{2n}\rtimes\langle\theta\rangle $ & \S \ref{subsubsection: B} & 2 \\ 
 $\mathrm{GSpin}_{2n+1}$ & $\mathrm{GSp}_{2n} \times W_F $ & $ (\GL_{2n}\times\GL_1 \times W_F)\rtimes\langle\hat{\theta}\rangle$ & $(\GL_{2n}\times\GL_1) \rtimes\langle{\theta}\rangle$ & \cite{Moe11b}*{\S2.1} & 2\\
 \hline
 $\Sp_{2n}$ & $\SO_{2n+1} \times W_F$ & $(\GL_{2n+1} \times W_F)\rtimes\langle\hat{\theta}\rangle$ & $\GL_{2n+1}\rtimes\langle{\theta}\rangle$ & \S \ref{subsubsection: C} & 2 \\ 
 \hline
  $\SO_{2n}^\alpha$ & $\SO_{2n} \rtimes W_F$ & $(\GL_{2n} \times W_F)\rtimes\langle\hat{\theta}\rangle$ & $\GL_{2n}\rtimes\langle{\theta}\rangle$ & 
\S\ref{subsubsection: D} & 2\\
    $\mathrm{GSpin}^\alpha_{2n}$ & $\mathrm{GO}_{2n} \rtimes W_F$ & $(\GL_{2n}\times\GL_1 \times W_F)\rtimes\langle\hat{\theta}\rangle$ & $( \GL_{2n}\times\GL_1)\rtimes\langle{\theta}\rangle$ & \cite{Moe11b}*{\S2.1} & 2 \\
\hline
  $\mathrm{G}_{2}$ & $\mathrm{G}_{2} \times W_F$ & $(\mathrm{Spin}_{8} \times W_F)\rtimes\langle\hat{\theta}\rangle$ & $\mathrm{PGSO}_{8}\rtimes\langle{\theta}\rangle$ & \cite{Gan25} & 3\\
  \hline
    $\mathrm{F}_{4}$ & $\mathrm{F}_4 \times W_F$ & $(\mathrm{E}_6  \times W_F)\rtimes\langle\hat{\theta}\rangle$ & $\mathrm{E}_{6}\rtimes\langle{\theta}\rangle$ & \cite{Hel01}*{p. 514} & 2\\
 \hline
\end{tabular}
}
\end{center}
\end{table}

\begin{remark}\label{remark: endoscopic KS}
Definition~\ref{definition: endoscopic triple} is closely related to \cite{KS99}*{Section 2.1}, as we now explain.
First, it is straightforward that our definition above implies the following commutative diagram
\[
\begin{tikzcd}
    \Lgroup{G'} \arrow{rr}{\xi} \arrow{dr}[swap]{\xi^0} && \Lgroup{G} \\
    & \Lgroup{G^0}\arrow{ur} & 
\end{tikzcd}
\] 
where $\Lgroup{G^0}\to \Lgroup{G}$ is the obvious inclusion and where $\xi^0$ is part of an endoscopic datum $(G',\Lgroup{G'},s,\xi^0)$ for $G^0$ as in \cite{KS99}.
Write $s=s'\rtimes\hat{\theta}^k$.
The data $(G,\theta,\textbf{a})$ of Kottwitz and Shelstad is given in our setting by $(G_\delta^0,\theta_\delta,1)$, so in particular, we assume that $\textbf{a}$ is the trivial class in $H^1(F,Z(\dualgroup{G^0}))$;
see Remark~\ref{remark: trivial a}.
Also, their endoscopic datum $(H,\mathcal{H},s,\xi_H)$ corresponds to $(G',\Lgroup{G'},s',\xi^0)$ in Definition~\ref{definition: endoscopic triple}; thus, in our setting, we have that $\mathcal{H} = \Lgroup{G'}$ and so this information is already encoded in $G'$. This restriction allows us to avoid the case of $z$-extensions. 
Finally, note that ${s'\rtimes{\hat\theta}^k}$ lies in $\dualgroup{G}$ but not in $\dualgroup{G^0}$ and, for $g\in \dualgroup{G^0}$, 
\[
{}^{{s'\rtimes{\hat\theta}^k}}(g\rtimes{\hat\theta}^0) = g\rtimes{\hat\theta}^0
\qquad\textrm{iff}\qquad
(\operatorname{Int}(s')\circ{\hat\theta}^k)(g) = g.
\]
\end{remark}

\subsection{Representations and invariant distributions}
\label{subsection: representations}

In this subsection, we discuss invariant distributions and the Grothendieck groups of representations for representations of the groups $G_\delta$ appearing in Section~\ref{subsection: pure inner forms}. 

Let $\delta\in Z^1(F,G)$ be a pure inner form of $G$. 
Recall that $\Pi(G_\delta)$ denotes the set of equivalence classes of irreducible admissible complex representations of $G_\delta(F)$ and let $K_{\BC}\Pi(G_\delta)$ denote the corresponding Grothendieck group. 
Consider the Hecke algebra $C_c^\infty(G_\delta(F))$. 
For $f\in C_c^\infty(G_\delta(F)),$ we fix a choice of Haar measure and consider the Harish-Chandra distribution character $\Theta(\pi)$ for the representation $\pi$ defined by 
\[
\Theta(\pi)(f):=\mathrm{trace}(\pi(f)).
\]
We remark that the foundational results on these distributions in the connected case are primarily due to Harish-Chandra and have been extended to disconnected groups by Clozel in \cite{Clo87}.
The distribution $\Theta(\pi)$ is an {\it invariant distribution}, {\it i.e.}, it is invariant under conjugation by any $g\in G_\delta(F).$
It is well-known that for $\pi_1,\pi_2\in\Pi(G_\delta),$ if $\pi_1\not\simeq\pi_2,$ then $\Theta({\pi_1})$ and $\Theta({\pi_2})$ are linearly independent. 
Moreover, if $\Theta(\pi_1)=\Theta(\pi_2)$ then $\pi_1\simeq\pi_2$. 
Consequently, we identify the span of 
$
\{
\Theta(\pi) \ | \ \pi\in\Pi(G_\delta)
\}
$ 
with $K_{\BC}\Pi(G_\delta)$ via the linear transformation $\Theta_{\pi}\mapsto[\pi]$. It is known that $K_{\BC}\Pi(G_{\delta})$ is dense in $\widehat{I}(G_{\delta})$, the space of invariant distributions on $G_{\delta}$; see \cite{Kaz86a}*{Theorem 1}, which does not assume $G_\delta$ is connected.

\subsection{Twisted invariant distributions}\label{sec twisted invariant distribution}

In this subsection, we recall the definition of twisted invariant distributions of $G_{\delta}^0(F)$ with respect to $\theta_{\delta}^k$ in Definition \ref{def: theta delta}, and define a certain section from twisted invariant distributions of $G_{\delta}^0(F)$ to invariant distributions of $G_{\delta}(F)$ in Lemma \ref{lemma: extension by zero}. Throughout the section, we assume that $\theta[\delta^0]=[\delta^0]$; otherwise the discussion is trivial.

Recall from Lemma \ref{lemma: automorphism preserves Borel} and Definition \ref{def: theta delta} that we have fixed $b_{\delta^0}g_{\delta}\rtimes \theta \in G_{\delta}(F)$, and the corresponding automorphism $\theta_{\delta}:= \textrm{Int}(b_{\delta^0}g_{\delta}) \circ \theta$ on $G^0_\delta(F)$.

\begin{definition}\label{definition: twisted invariant distributions}
A distribution of $G^0_{\delta}(F)$ is called {\it twisted invariant} with respect to $\theta_{\delta}^k$ if it is invariant under the $\theta_{\delta}^k$-conjugation: 
\[ \Ad_{\theta_{\delta}}(g)x:= gx\theta_{\delta}^k(g^{-1}), \forall g, x \in G^0_{\delta}(F).  \]
Let $\widehat{I}(G^0_{\delta})_{\theta_{\delta}}$ denote the space of twisted invariant distributions.
Note that if $\theta_{\delta}$ is trivial or $k=0$, then $\widehat{I}(G^0_{\delta})_{\theta_{\delta}}= \widehat{I}(G^0_{\delta})$, the space of invariant distributions of $G^0_{\delta}(F)$. \IN{$\widehat{I}(G^0_{\delta})_{\theta_{\delta}}$, the space of twisted invariant distributions}
\end{definition}

We provide an important source of twisted invariant distributions. Suppose that $\pi$ is a representation of $G_{\delta}(F)$. 
Consider the distribution $\Theta(\pi)_{\theta_{\delta}^k}$ of $G^0_\delta(F)$, called the {\it twisted character} of $\pi$ \IN{$\Theta(\pi)_{\theta^k}$, twisted character of $\pi$} with respect to $\theta_{\delta}^k$, that sends $h \in \mathcal{C}_{c}^{\infty}(G^0_{\delta}(F))$ to the trace of the finite rank operator 
\[ 
\pi(h)_{\theta_{\delta}^k}:= \int_{ G^0_\delta(F)} h(g) \pi( (g \rtimes 1)  (b_{\delta^0}g_{\delta} \rtimes  \theta)^k) dg.  
\]
Indeed, $\Theta(\pi)_{\theta_{\delta}^k}$ is the restriction of $\Theta(\pi)$ (see Section \ref{subsection: representations}) to $\mathcal{C}_{c}^{\infty}(G^0_{\delta}(F)\cdot \theta_{\delta}^k),$ composing with the identification $\mathcal{C}_{c}^{\infty}(G^0_{\delta}(F)\cdot \theta_{\delta}^k)\simeq \mathcal{C}_{c}^{\infty}(G^0_{\delta}(F))$ given by translation with $\theta_{\delta}^k$.
It can be checked by a straightforward computation that $\Theta(\pi)_{\theta_{\delta}^k}$ is a twisted invariant distribution with respect to $\theta_{\delta}^k$ of $G^0_\delta(F)$. Thus, we obtain a map 
$\Theta (\cdot)_{\theta_{\delta}^k}: K_{\BC} \Pi (G_\delta) \to \widehat{I}(G^0_\delta)_{\theta_{\delta}^k}$.
We let $K_{\BC}\Pi(G^0_{\delta})_{\theta_{\delta}^k}$\IN{$K_{\BC}\Pi(G^0_{\delta})_{\theta_{\delta}^k}$} denote the image of this linear transformation, thus defining
\begin{equation}\label{equation: Theta theta}
\Theta (\cdot)_{\theta_{\delta}^k}: K_{\BC} \Pi (G_\delta) \to K_{\BC}\Pi(G^0_{\delta})_{\theta_{\delta}^k}.
\end{equation}
Note that if $\theta_{\delta}$ is trivial then  $\Theta(\cdot)_{\theta_{\delta}}$ is simply the identity map. 

In the rest of this section, we assume that $\theta_{\delta}$ is non-trivial and construct a section $j_{!,k}$ of $\Theta(\cdot)_{\theta_{\delta}^k}$ given by {\it extension by zero} explicitly, for $0 \leq k \leq \ell-1$. We will use the notations introduced in Section \ref{subsection: restriction induction}. We recall that we have fixed a non-trivial character $\vartheta$ of $\langle \theta \rangle$ and $\zeta:= \vartheta(\theta)$, which is a primitive $\ell$-th root of unity.

First, we consider the case that $k=0$. By definition, we have 
\[ \Theta(\Pi)_{\theta_{\delta}^0}= \Theta( \Res_{G_{\delta}^0}^{G_{\delta}}\Pi). \]
Thus, $K_{\BC}\Pi(G_{\delta}^0)_{\theta_{\delta}^0}$ is simply the image of $\Res_{G_{\delta}^0}^{G_{\delta}}$ in $K_{\BC}\Pi(G_{\delta}^0)$. By the discussion in Section \ref{subsection: restriction induction}, it is not hard to see that 
\[ j_{!,0}:= \frac{1}{\ell} \Ind_{G_{\delta}^0}^{G_{\delta}}\]
gives a section of $\Theta(\cdot)_{\theta_{\delta}^0}$.

Next, we consider the case that $0 < k \leq \ell-1$. We first define this section on a collection of one-dimensional subspaces of $K_{\BC}\Pi(G^0_{\delta})_{\theta_\delta^k}$, which we introduce now. For each $\pi \in \Pi^{\theta_{\delta}}(G_{\delta}^0)$, we take any irreducible subquotient $\pi_0$ of $\Ind_{G^0_{\delta}}^{G_{\delta}} \pi$, and let $K_{\BC}(\pi)_{\theta_\delta^k}$ be the one-dimensional vector space spanned by $\Theta(\pi_0)_{\theta_\delta^k}$. This does not depend on the choice of $\pi_0$ since $\Theta(\pi_i)_{\theta_\delta^k}= \zeta^{k(i-l)}\Theta(\pi_l)_{\theta_\delta^k}$ where $\pi_i={}^{\theta_{\delta}^i}\pi_0$ is any other irreducible subquotient of $\Ind_{G^0_{\delta}}^{G_{\delta}} \pi$.  
Next, we define the following elements in $K_{\BC}\Pi(G_{\delta})$ for $0 \leq a \leq \ell-1$
\begin{align}\label{eq extension by zero}
    s^{a}(\pi_0) := \sum_{l=0}^{\ell -1} \frac{\zeta^{ -la } \pi_l}{\ell}.
\end{align}
Again since $ \Theta(\pi_l)_{\theta_\delta^{a}}= \zeta^{la}\cdot \Theta(\pi_0)_{\theta_\delta^{a}}$, we see that 
\begin{align}\label{eq twisted distribution section}
    \Theta( s^{a}(\pi_0))_{\theta_\delta^{b}}=\begin{cases}
        \Theta(\pi_0)_{\theta_\delta^{a}} & \text{ if }a=b,\\
        0 & \text{ otherwise.}
    \end{cases}
\end{align} 
Next we define $j_{!,k}$, the linear transformation from $K_{\BC}(\pi)_{\theta_\delta^k}$ to $K_{\BC}\Pi(G_{\delta})$,  by
\[ j_{!,k}( \Theta(\pi_0)_{\theta_\delta^k}):= s^k( \pi_0). \]
This definition is independent of the choice of $\pi_0$. Indeed, if we let $\sigma_0= \pi_i$ and $\sigma_l= \sigma_0 \otimes \vartheta^{l} $, then 
\[ j_{!,k}( \Theta(\sigma_0)_{\theta_\delta^k})= j_{!,k}( \zeta^{ik} \Theta(\pi_0)_{\theta_\delta^k})=\zeta^{ik} \cdot s^k(\pi_0)= \sum_{l=0}^{\ell -1}  \frac{\zeta^{-(l-i)k}\pi_l}{\ell}=\sum_{l=0}^{\ell -1}  \frac{\zeta^{-(l-i)k}\sigma_{l-i}}{\ell}= s^k(\sigma_0).  \]
 In the following lemma, we check that gluing the $j_{!,k}$ on this collection of one-dimensional subspaces gives a section of $\Theta(\cdot)_{\theta_{\delta}^k}$.

\begin{lemma}\label{lemma: extension by zero}
In the above setting with $0 <k \leq \ell-1$, the following holds. 
    \begin{enumerate}
    \item For any representation $\pi$ of $G^0_{\delta}(F)$, we have $\Theta(\Ind_{G^0_{\delta}}^{G_{\delta}}(\pi))_{\theta_\delta^{k}}=0$ for any $0<k\leq \ell-1.$
    
    \item Define 
        
\[
K_{\BC}\Pi(G_{\delta})^{\theta_\delta^k}:= \sum_{\pi \in \Pi^{\theta_{\delta^0}}(G^0_{\delta})} j_{!,k}(K_{\BC} (\pi)_{\theta_\delta^k}).
\]
        Then the restriction of $\Theta(\cdot)_{\theta_\delta^k}$ gives a bijection from $ K_{\BC}\Pi(G_{\delta})^{\theta_\delta^k}$ onto $K_{\BC}\Pi(G_{\delta}^0)_{\theta_\delta^k}$.
        
        \item  
        We have $K_{\BC}\Pi(G^0_{\delta})_{\theta_\delta^k}= \bigoplus_{\pi \in \Pi^{\theta_{\delta^0}}(G^0_{\delta})} K_{\BC}(\pi)_{\theta_\delta^k}$, and the linear transformation 
        \begin{align*}
            j_{!,k}: K_{\BC}\Pi(G^0_{\delta})_{\theta_\delta^k} &\to K_{\BC}\Pi(G_{\delta}),\\
            \Theta(\pi_0)_{\theta_\delta^k} &\mapsto s^k(\pi_0)
            \end{align*}  \IN{{$j_{"!,k}$}, extension by zero}
        gives a well-defined section of $\Theta(\cdot)_{\theta_\delta^k}$. 
    \end{enumerate}
\end{lemma}

\begin{proof}
Part (1) follows from the same argument as Lemma \ref{lem finite index induction restriction}(3), which we omit. For Part (2), recall from Section \ref{subsection: restriction induction} that $K_{\BC}\Pi(G_{\delta})$ has a basis 
\[
    \{
    \pi_0,\ldots,\pi_{\ell-1} \ \vert\ \pi\in \Pi^{\theta_{\delta^0}}(G^0_{\delta})
    \} 
\ \bigcup \ 
    \{
    \Ind_{G_\delta^0(F)}^{G_\delta(F)} \pi \ \vert\ \pi \in \Pi(G_\delta^0) \setminus \Pi^{\theta_{\delta^0}}(G^0_{\delta})
    \}.
\]
By Part (1) and  $\Theta(\pi_l)_{\theta_\delta^k}= \zeta^{lk}\cdot \Theta(\pi_0)_{\theta_\delta^k}$ for $\pi \in K_{\BC}\Pi^{\theta_{\delta^0}}(G^0_{\delta})$, we see that $K_{\BC}\Pi(G^0_{\delta})_{\theta_\delta^k}$ is generated by $\{\Theta(\pi_0)_{\theta_\delta^k}\ \vert\ \pi \in \Pi^{\theta_{\delta^0}}(G^0_{\delta})\}$. This proves the surjectivity. For the injectivity, let $ K'$ be the subspace of $K_{\BC}\Pi(G_{\delta})$ spanned by 
\[
    \{
    s^{l}(\pi_0)
        \ \vert\ 
    \pi \in \Pi^{\theta_{\delta^0}}(G^0_{\delta}),\ l \in \{0,\ldots, \ell-1\}\setminus \{k\}
    \} 
\cup 
    \{
    \Ind_{G^0_{\delta}}^{G_{\delta}} \pi
        \ \vert\ 
    \pi \in \Pi(G^0_{\delta}) \setminus \Pi^{\theta_{\delta^0}}(G^0_{\delta})
    \}.
\]
We claim that
\begin{enumerate}
    \item [(i)] 
    $ K_{\BC}\Pi(G_\delta)= K_{\BC}\Pi(G_\delta)^{\theta_\delta^k} \oplus K'$, and
    \item [(ii)] 
    $K'$ is exactly the kernel of $\Theta(\cdot)_{\theta_\delta^k}$.
\end{enumerate}
These claims imply the desired injectivity.
 Now we prove Claims (i) and (ii). 
 
 For Claim (i), observe that for $\pi \in \Pi^{\theta_{\delta^0}}(G_\delta^0)$, any representation in  $\{\pi_0,\ldots, \pi_{\ell-1}\}$ can be expressed as a linear combination (in the Grothendieck group) of elements in $\{s^l(\pi_0)\}_{0 \leq l \leq \ell-1}$. Thus, $ K_{\BC}\Pi(G_{\delta})= K_{\BC}\Pi(G_{\delta})^{\theta_\delta^k}+ K'
 $ and it remains to show that $K_{\BC}\Pi(G_{\delta})^{\theta_\delta} \cap K'= 0$. Suppose that $\Pi \in K_{\BC}\Pi(G_{\delta})^{\theta_\delta^k} \cap K'$. Then since $\Pi \in K_{\BC}\Pi(G_{\delta})^{\theta_\delta^k}$, Equation~\eqref{eq twisted distribution section} implies that $\Theta(\Pi)_{\theta_\delta^l}=0$ for $l\in \{0,\ldots, \ell-1\}\setminus \{k\}$. Since $\Pi \in K'$, Part (1) and Equation~\eqref{eq twisted distribution section} implies  that $\Theta(\Pi)_{\theta_\delta^k}=0$ as well. Thus, $\Theta(\Pi)=0$ as an invariant distribution of $G_{\delta}(F)$; hence $\Pi=0$. 
 
For Claim (ii), we have already shown that $K'$ is contained in the kernel of $\Theta(\cdot)_{\theta_\delta^k}$ in the previous paragraph. Conversely, suppose that $\Pi$ lies in the kernel and write
\begin{align*}
    & \Theta(\Pi)_{\theta_\delta^l} \\
    & = \begin{cases}
        \mathop{\sum}\limits_{\pi \in \Pi^{\theta_{\delta^0}}(G^0_{\delta})} \frac{1}{\ell} a_{\pi, 0} \Theta(\pi_0)_{\theta_\delta^0} + \mathop{\sum}\limits_{\sigma \in \Pi(G^0_{\delta})\setminus \Pi^{\theta_{\delta^0}}(G^0_{\delta})}  {a_{\sigma,0}} \Theta(\Ind_{G^0_{\delta}(F)}^{G_{\delta}(F)} \sigma)_{\theta_\delta^0} & \text{ if }l=0,\\
        \mathop{\sum}\limits_{\pi \in \Pi^{\theta_{\delta^0}}(G^0_{\delta})} a_{\pi, l} \Theta(\pi_0)_{\theta_\delta^l} & \text{ if }l \in \{1,\ldots, \ell-1\} \setminus\{k\},
    \end{cases}
\end{align*}
where $a_{\pi,l} \neq 0$ for only finitely many $\pi$. Then consider the following element in $K'$ (note that  $s^0(\pi_0)= \frac{1}{\ell}\Ind_{G^0_{\delta}}^{G_{\delta}} \pi$ for $\pi \in \Pi^{\theta_{\delta^0}}(G^0_{\delta})$)
\begin{align*}
   \Pi':= \sum_{\pi\in \Pi(G^0_\delta)} a_{\pi,0} \cdot \Ind_{G_\delta^0}^{G_\delta} \pi + \sum_{l \in \{1,\ldots,\ell-1\}\setminus \{k\}} \sum_{\pi \in \Pi^{\theta_{\delta^0}}(G^0_\delta)} a_{\pi, l} \cdot s^{l}(\pi_0).  
\end{align*}
Then $\Theta(\Pi- \Pi')_{\theta_\delta^l}=0$ for any $0 \leq l \leq \ell-1$ by Part (1) and Equation~\eqref{eq twisted distribution section}. Therefore, we have $\Pi=\Pi' \in K'$ which completes the verification of Claim (ii).

Part (3) is an immediate consequence of Equation~\eqref{eq twisted distribution section} and Part (2). 
Note that if $\theta_\delta$ is trivial, then $j_{!}$ is just the identity map.
This completes the proof of the lemma.
\end{proof}

\subsection{Stable invariant distributions and endoscopic transfer}\label{subsection: transfer}

In this subsection, we recall the notation of endoscopic transfer from stable distributions of $G'(F)$ to twisted invariant distributions of $G_\delta(F)$; here $(G',s,\xi)$ is an endoscopic triple for $G$, as in Definition~\ref{definition: endoscopic triple}. We modify this transfer by composing it with $j_{!,k}$, defined in the previous section, and then scaling it with $\ell$ in Definition \ref{definition: modified transfer}. The main result in this subsection is Proposition~\ref{proposition: induction commutes with transfer} which shows that endoscopic transfer, as we define it, commutes with parabolic induction.

Let $\Theta$ be an invariant distribution on $C_c^\infty(G'(F))$. 
We say that $\Theta$ is a {\it stable invariant distribution} if $\Theta$ lies in the closure of the span of stable orbital integrals (see \cite{Art89}*{\S3}). 
Recall from Section~\ref{subsection: representations} that $\widehat{I}(G')$ denotes the space of invariant distributions of $G'(F)$ and we identify $K_{\BC}\Pi(G')$ as a subspace of $\widehat{I}(G')$. 
Let $\widehat{I}(G')^\st$
\IN{$\widehat{I}(G')^\st$, stable invariant distributions} 
 denote the subspace of stable invariant distributions of $G'(F)$ and let $K_{\BC}\Pi(G')^\st$ denote the intersection of $K_{\BC}\Pi(G')$ and $\widehat{I}(G')^\st$. 

Let $(G', s =s'\rtimes \hat{\theta}^k, \xi)$ be an endoscopic triple of $G$ with $0 \leq k \leq \ell-1$, and let $\delta$ be a pure inner form of $G$. Suppose that $\theta^{k}[\delta^0]=[\delta^0]$, so that $\theta_{\delta}^k$ defines an automorphism on $G_{\delta}^0$. From \cite{KS99}, Kottwitz-Shelstad defined the transfer attached to the endoscopic triple $(G', s , \xi)$ of $G$, which we call the {\it Kottwitz-Shelstad transfer}. This is a linear transformation
\begin{equation}\label{equation: KST}
\Transfer{G',s'}{G^0_{\delta},\theta_\delta^k}: 
\widehat{I}(G')^\st 
\to 
\widehat{I}(G^0_{\delta})_{\theta_\delta^k}.
\end{equation}
We remark that the endoscopic triples we consider in Definition \ref{definition: endoscopic triple} are not the most general setting from \cite{KS99}.

\begin{remark}\label{rmk Kottwitz-Shelstad transfer}\ 
    \begin{enumerate}
   \item It is known (\cite{Art96}*{Theorems 6.2, (6.2), Lemma 5.2})
that the restriction of $\Transfer{G',s'}{G^0_{\delta},\theta_\delta^k}$ to $K_{\BC}\Pi(G')^{\st}$ has image in $K_{\BC}\Pi(G_{\delta}^0)_{\theta_{\delta}^k}$. We shall mainly consider this restriction in the discussion.
   \item If $k=0$, then $(G',s,\xi)$ is also an endoscopic triple for the connected group $G^0$. Then the transfer attached to the triple $(G',s,\xi)$ of $G$ and of $G^0$ coincides. This case is also known as Langlands-Shelstad transfer. We shall sometimes refer to the transfer for the case that $0< k \leq \ell-1$ as twisted transfer. 
   
        \item The definition of the transfer depends on a choice of Whittaker datum for $G$, which we have fixed from Definition~\ref{assu: groups we consider}.

    \end{enumerate}
\end{remark}

Now we modify the transfer in \eqref{equation: KST} by composing with $j_{!,k}$ together and then scaled by $\ell$.

\begin{definition}\label{definition: modified transfer}
Suppose that $G$ satisfies Definition~\ref{assu: groups we consider}. Let $\delta$ be a pure inner form of $G$. For each endoscopic triple $(G', s=s'\rtimes \hat{\theta}^k, \xi)$, the
{\it endoscopic transfer}  $\Transfer{G',s}{G_\delta}$\IN{$\Transfer{G',s}{G_\delta}$, endoscopic transfer} from stable distributions of $G'(F)$ to invariant distributions of $G_\delta(F)$ is the linear transformation defined as follows: If $\theta[\delta^0] \neq [\delta^0]$, then set $\Transfer{G',s}{G_\delta} \equiv 0$. Otherwise, we define it to be the composition of $\Transfer{G',s'}{G^0_{\delta},\theta_\delta^k}$ with the section $j_{!,k}$ defined in Section \ref{sec twisted invariant distribution}, scaled by $\ell$. Thus, the following diagram commutes. 
\[
\begin{tikzcd}
&&
K_{\BC}\Pi(G^0_{\delta})_{\theta_{\delta}^k}
\arrow{rrdd}[swap]{\ell j_{!,k}}{\mathrm{Lemma~\ref{lemma: extension by zero}}}
\arrow{rr}{\iso}
&&
K_{\BC}\Pi(G_{\delta}) |_{G_{\delta}^0 \cdot \theta_{\delta}^k} 
\\
&&&&
\\
\arrow{uurr}{{\Transfer{G',s'}{G^0_\delta,\theta_\delta^k}}}[swap]{\eqref{equation: KST}}
 K_{\BC}\Pi(G')^{\st} \arrow{rrrr}{\Transfer{G',s}{G_\delta}}
&&&& K_{\BC}\Pi(G_{\delta}) 
\arrow{uu}[swap]{{\frac{1}{\ell}} \cdot \mathrm{restriction}} 
\end{tikzcd}
\]
\end{definition}

\begin{remark}\label{rmk transfer disconnected}\ 
\begin{enumerate}
    \item When $\theta[\delta^0] \neq [\delta^0]$, $G_{\delta}(F)$ is connected, so it does not come with an automorphism of the identity component. This is not in the setting of Kottwitz-Shelstad transfer. Therefore, we set $\Transfer{G',s}{G_\delta}  \equiv 0$, which is also suggested by the endoscopic character relations for tempered packets (see Remark \ref{rmk LLC(1)}.
    \item When $k=0$, $(G',s,\xi)$ is also an endoscopic triple for the connected group $G_{\delta}^0$. Then,
    \[ \Transfer{G',s}{G_\delta}= \Ind_{G_{\delta}^0}^{G_{\delta}}\circ \Transfer{G',s}{G_\delta^0}\]
    by Remark \ref{rmk Kottwitz-Shelstad transfer}(2) and $j_{!,0}= \frac{1}{\ell}\Ind_{G_{\delta}^0}^{G_{\delta}}$ as defined in the previous section.
\end{enumerate}

\end{remark}

We do not recall the explicit definition of $\Transfer{G',s}{G^0_{\delta}}$ since we do not use it in this paper (except in as special case in Proposition \ref{prop twisted transfer from Siegel GL to disconencted GL} below). Instead, 
 we will use the fundamental property of endoscopic transfer that it commutes with parabolic induction, which we now explain.  
 Let $(G',s,\xi)$ be an endoscopic triple for $G$ as defined in Section~\ref{subsection: endoscopic triple}, and fix a pure inner form $\delta \in Z^1(F,G)$. Let $P'$ be a standard parabolic subgroup of $G'$ with the standard Levi subgroup $M'$. Then $M'$ is an endoscopic group of a standard Levi subgroup $M$ of $G$. More explicitly, we write $\dualgroup{M'}=Z_{\dualgroup{G'}}(\dualgroup{T'})$ for some torus $\dualgroup{T'}$ of $\dualgroup{G'}$. Let $M$ is the standard Levi of $G$ with dual group $\dualgroup{M}:= Z_{\dualgroup{G}}(\dualgroup{T'})$. Then $\dualgroup{M'}=Z_{\dualgroup{M}}(s)^0$; hence $M'$ is an endoscopic group of $M$. Let $P$ be the standard parabolic subgroup of $G$ that contains $M$.

\begin{proposition}\label{proposition: induction commutes with transfer}
In the above setting, the following diagram commutes.
\[
\begin{tikzcd}
    K_{\BC}\Pi(G')^{\st} 
        \arrow{rr}{\Transfer{G',s}{G_\delta}} 
    && K_{\BC}\Pi(G_\delta) 
         \\
    K_{\BC}\Pi(M')^{\st} 
        \arrow{rr}{\Transfer{M',s}{M}} 
        \arrow{u}{\Ind_{P'}^{G'}}
    && K_{\BC}\Pi(M)^{\st}
        \arrow{u}[swap]{\Ind_{P}^{G_\delta}}
\end{tikzcd}
\]
\end{proposition}

\begin{proof}
Write $s=s'\rtimes \hat{\theta}^k$. 
If $\theta^k[\delta^0]\neq [\delta^0]$, then both transfers are zero and the diagram commutes. Thus, we assume that $\theta^k[\delta^0]= [\delta^0]$.
Note that if $G_{\delta}$ is disconnected, then so is $M$ since $s  \in  Z_{\dualgroup{G^0}}(\dualgroup{T'})$. Consider the following diagram.
\[
\begin{tikzcd}
    K_{\BC}\Pi(G')^{\st} 
        \arrow{rr}{\Transfer{G',s'}{G_\delta^0,\theta_{\delta}^k}} 
    && K_{\BC}\Pi(G_\delta)_{\theta_{\delta}^k}  &&K_{\BC}\Pi(G_\delta) \arrow{ll}[swap]{\Theta(\cdot)_{\theta_{\delta}^k}}
         \\
    K_{\BC}\Pi(M')^{\st} 
        \arrow{rr}{\Transfer{M',s'}{M^0, \theta_{\delta}^k}} 
        \arrow{u}{\Ind_{P'}^{G'}}
    && K_{\BC}\Pi(M)^{\st}_{\theta_{\delta}^k} 
        \arrow{u}[swap]{(\Ind_{P}^{G_\delta})_{\theta_{\delta}^k}} &&\arrow{u}[swap]{\Ind_{P}^{G_\delta}} \arrow{ll}[swap]{\Theta(\cdot)_{\theta_{\delta}}}K_{\BC}\Pi(M)^{\st} 
\end{tikzcd}
\]
Here $(\Ind_{P}^{G_\delta})_{\theta_{\delta}^k}$  is the unique linear transformation such that the right square commutes. It is well-known that the left square commutes (see    
{\cite{KS99}, \cite{Hir04}*{Lemma 6.1}, or \cite{MW16}*{Section I.3}} for example). Therefore, it suffices to show that the following diagram commutes.
\[
\begin{tikzcd}
    K_{\BC}\Pi(G_\delta)_{\theta_{\delta}^k} \arrow{rr}{j_{!,k}} &&K_{\BC}\Pi(G_\delta) 
         \\
    K_{\BC}\Pi(M)^{\st}_{\theta_{\delta}^k} 
        \arrow{u}[swap]{(\Ind_{P}^{G_\delta})_{\theta_{\delta}^k}} \arrow{rr}{j_{!,k}} &&K_{\BC}\Pi(M)^{\st} \arrow{u}[swap]{\Ind_{P}^{G_\delta}} 
\end{tikzcd}
\]
As we have shown that $j_{!,k}$ is a section of $\Theta(\cdot)_{\theta_{\delta}^k}$, it suffices to show that $\Ind_{P}^{G_{\delta}}$ preserves the image of $j_{!,k}$. 
From the proof of Lemma \ref{lemma: extension by zero} and the definition of $j_{!,0}$, we know that the image of $j_{!,k}$ consists of elements  
\[ s^k(\pi_0)= \sum_{j=0}^{\ell-1} \frac{\zeta^{-jk} \pi_j}{\ell},\]
where $\pi_0$ is any virtual representation of $G_{\delta}$ or $M$, $\pi_j:= \pi_0 \otimes \vartheta^j$ where $\vartheta$ is a fixed non-trivial character of $\pi_0(G_{\delta})$ or $\pi_0(M)$ with $\zeta= \vartheta(\theta_{\delta})$. Since $\Ind_{P}^{G_{\delta}}( \pi_0 \otimes \vartheta^j)\cong \Ind_{P}^{G_{\delta}}(\pi_0) \otimes \vartheta^j$, we have
\[\Ind_{P}^{G_{\delta}} (s^k(\pi_0))= \sum_{j=0}^{\ell-1} \frac{\zeta^{-jk} \Ind_{P}^{G_{\delta}}(\pi_j)}{\ell}=\sum_{j=0}^{\ell-1} \frac{\zeta^{-jk} \Ind_{P}^{G_{\delta}}(\pi)_j}{\ell}= s^k(\Ind_{P}^{G_{\delta}}(\pi_0)), \]
which again lies in the image of $j_{!,k}$. 
\end{proof}

\section{The local Langlands correspondence}\label{section: LLC}

In this section, we recall the conjectural {\it local Langlands correspondence} for {quasisplit} disconnected groups and their pure inner forms, based largely on \cite{Kal26}*{Conjectures 4.2, 4.7}. 
For any pure inner form $\delta \in Z^1(F,M)$ of a Levi $M$ of $G$, the correspondence matches irreducible representations of $M_\delta(F)$, up to equivalence, with enhanced $L$-parameters $\Phi^{\mfr{e}}(M,\delta)$ for $M_{\delta}$ (see Definition \ref{def L-par relevance}).
We formulate the desired properties of this correspondence in Desiderata~\ref{desideratum: LLC}.
Some of these Desiderata will serve as hypotheses for our results in Section~\ref{section: endoclassical}, where we verify that they hold for all classical groups.

\begin{definition}\label{definition: LLC}
    Suppose that $G$ satisfies Definition~\ref{assu: groups we consider}. 
Recall that this means we have fixed a $\theta$-invariant $F$-pinning of $G^0$. Together with the fixed nontrivial additive character of $F$, this defines a Whittaker datum $\mfr{w}_{M^0}$ for every standard Levi subgroup $M^0$ of $G^0$. 
 A {\it local Langlands correspondence of $G$} \IN{$\LLC^{\mfr{e}}_{G,{\delta}}$, local Langlands correspondence} is a bijection,
\[
    \begin{aligned}
\LLC^{\mfr{e}}_{G,{\delta}}: \Pi(G_{\delta}) &\to \Phi^{\mfr{e}}(G,\delta),\\
    \pi& \mapsto (\phi_{\pi}, \tau_{\pi}),\\
    \pi(\phi,\tau)& \mapsfrom (\phi,\tau).
    \end{aligned}
\]
We define the {\it $L$-packet}\IN{$\Pi_{\phi}(G_{\delta})$, $L$-packet} by
    \[ 
    \Pi_{\phi}(G_{\delta}):= \{\pi \in \Pi(G_{\delta})\ | \ \phi_{\pi}= \phi \}= \{\pi(\phi,\tau)\ | \ \tau \in \widehat{A_{\phi}}[\delta]\}.
    \] 
\end{definition}

Although the following desiderata are well known in the literature, we have decided to collect them all here for the convenience of the reader.
There is some overlap with desiderata  \cite{Kal26}*{Conjectures 4.2, 4.7, 7.1, A.1}.
There are also redundancies among the desiderata below, reflecting when and how we use them. Some of these desiderata concerns the relation of the local Langlands correspondence of $G$ and a related group $H$ (e.g. Levi subgroups of endoscopic groups of $G$), and we assume the existence of local Langlands correspondence for these groups.

\begin{desideratum}[Various desiderata for local Langlands correspondence]\label{desideratum: LLC}
Suppose that $G= G^0 \rtimes \langle \theta \rangle$ satisfies Definition~\ref{assu: groups we consider}.

\begin{enumerate}
 \item[\CCh](Central characters) For $\phi^0 \in \Phi(G^0)$, the elements of $\Pi_{\phi}(G^0)$ have the same central character.

\item [\TW] (Compatibility with twists by automorphisms and quasi-characters) Let $\theta$ be an automorphism of $G^0$ fixing the $F$-pinning $(B, T, \{X_{\alpha}\}_{\alpha \in \Delta})$ and let $\hat{\theta}$ be the dual automorphism. Let $\mathbf{a}$ be an element in $Z^1(W_{F}, Z(\widehat{G^0}))$ associated with a quasi-character $\omega$ of $G^0(F)$ by the Langlands isomorphism
$
H^{1}(W, Z(\widehat{G^0})) \rightarrow \Hom(G^0(F), \BC^{\times}).
$
Recall that for any $\phi^0 \in \Phi(G^0)$, we denote its $\widehat{\theta}$-conjugate by $\,^{\widehat{\theta}} \phi^0$ and its product with $\mathbf{a}$ by $\phi^0 \otimes \mathbf{a}$. Then
\[
\Pi_{\,^{\widehat{\theta}}\phi^0}(G^0) = \,^\theta\Pi_{\phi^0}(G^0) \quad \text{ and } \quad  \Pi_{\phi^0 \otimes \mathbf{a}}(G^0) = \Pi_{\phi^0}(G^0) \otimes \omega.
\]

 \item [\Funtheta](Functoriality with respect to $\theta$-action)
If $M$ is a connected Levi subgroup of $G$ then 
\[
    \LLC^{\mfr{e}}_{M}(\pi)= (\phi_{\pi}, \tau_{\pi})
    \qquad\implies\qquad
    \LLC^{\mfr{e}}_{{}^\theta M}( {}^{\theta_{\delta}}{\pi})= ({}^{{\hat\theta}}\phi_{\pi}, {}^{{\hat\theta}}\tau_{\pi}).
\]

\item [\CTD] (Classification of tempered and discrete series representations)
Let $\Pi_{\mathrm{temp}}(G^0)$ (resp. $\Pi_{2}(G^0)$) be the set of isomorphism classes of irreducible tempered (resp. discrete series) representations of $G^0(F)$. Denote the set of tempered (resp. discrete) Langlands parameters by $\Phi_{\mathrm{bdd}}(G^0)$ (resp. $\Phi_2(G^0)$). Then there there are disjoint decompositions into $L$-packets:
\[
\Pi_{\mathrm{temp}}(G^0) = \bigsqcup_{\phi^0 \in \Phi_{\mathrm{bdd}}(G^0)} \Pi_{\phi^0}(G^0)
\qquad
and
\qquad
\Pi_{2}(G^0) = \bigsqcup_{\phi^0 \in \Phi_{2}(G^0)} \Pi_{\phi}(G^0).
\]

\item [\TI] (Compatibility with tempered parabolic induction) Suppose $\phi^0 \in \Phi_{\mathrm{bdd}}(G^0)$ factors through $\phi_{M} \in \Phi_{\mathrm{bdd}}(M)$ for a Levi subgroup $M$ of $G^0$. Let $P$ be any parabolic subgroup of $G^0$ with Levi factor $M$. Then $\Pi_{\phi}(G^0)$ consists of irreducible constituents of 
\[
\{{\rm Ind}^{G^0}_{P}  \pi_{M} \ | \ \pi_{M} \in \Pi_{\phi_{M}}(M)\}.
\]

\item [\Hai] (Haines conjecture)
Let $M$ be a Levi subgroup of $G$ and let $\delta$ be a pure inner form of $M$. We regard $\delta$ as a pure inner form of $G$ by considering its image under the injection $Z^1(F,M) \to Z^1(F,G)$; hence $M_{\delta}$ is isomorphic to a Levi subgroup of $G_{\delta}$. For every $\sigma \in \Pi(M_{\delta})$ and every irreducible subquotient $\pi$ of $\Ind_{M_{\delta}}^{G_{\delta}} (\sigma)$, we have 
\[ \lambda_{\phi_{\pi}}= ({}^LM \to \Lgroup{G}) \circ \lambda_{\phi_{\sigma}}.\]
 \item [\LC] (Compatibility with Langlands classification for parameters) 
    Let $\delta$ be a pure inner form of $G$ and let $\pi \in \Pi(G_{\delta})$.
    Suppose that
    \[
    \LLC^{\mfr{e}}_{G,{\delta}}(\pi)= (\phi, \tau).
    \]
    If $(P_{\pi}, \pi_t, \nu_{\pi})$ is the standard triple for $\pi$ as in 
    Theorem \ref{theorem: standard representations}, and if
     $(P_{\phi}, ([ \phi_t], \tau_t), \nu_{\phi})$ is the standard triple for $(\phi, \tau)$ as in Lemma \ref{lemma: LC for enhanced parameters},
    then $P_\pi= P_\phi$,  $\nu_\pi = \nu_\phi$,
and    \[
    \LLC^{\mfr{e}}_{M,{\delta_M}}( \pi_t)=(\phi_t,\tau_t),
    \]
    where $M_{\pi}$ is the standard Levi subgroup of $P_\pi=P_\phi$ with an associated pure inner form $\delta_M$ of the corresponding standard Levi subgroup $M$ of $G$, and
    $\tau_t \in \widehat{A_{\phi_t}}[\delta_M]$.
    
    \item [\Sta] 
    (Stability of tempered packets) 
   For a Levi subgroup $M$ of a quasisplit $G$ and for every tempered Langlands parameter $\phi$ of $M$, the linear combination
    \[  
    \sum_{
    \tau' \in \widehat{A_{\phi}}[1]
    } 
    \dim(\tau)\ \pi(\phi,\tau)
    \]
    is a stable distribution on $M(F)$, which is independent of the choice of Whittaker datum.
    \item [\ECR] 
    (Endoscopic character relations)
    Let  $(G',s,\xi)$ be an endoscopic triple for $G$ as in Definition~\ref{definition: endoscopic triple}, and $\phi'$ be a tempered Langlands parameter for $G'$. Let $\phi:= \xi \circ \phi'$ be the corresponding Langlands parameter of $G$. Then for every pure inner form $\delta $ of $G$, we have
    \[
    \Transfer{G',s}{G_{\delta}} 
    \left( 
\sum_{
    \tau' \in \widehat{A_{\phi'}}[1]
    } 
    \dim(\tau')\ \pi(\phi',\tau')
\right)
    =  e(M_\delta)
    \sum_{\tau \in \widehat{A_{\phi}}[\delta]} \trace(s ,  \tau)  \pi(\phi, \tau),
    \]
    where we identify the semisimple element $s\in \dualgroup{G}$ with its image under the quotient $Z_{\dualgroup{G}}(\phi) \to Z_{\dualgroup{G}}(\phi)/Z_{\dualgroup{G}}(\phi)^\circ =: A_\phi$ and where $\Transfer{G',s}{G_\delta}$ was defined in Definition \ref{definition: modified transfer}.
    
    \item [\RES] Suppose that $M$ is a disconnected Levi subgroup of $G$. Let  $\phi^0$ be a Langlands parameter of $M^0$ and let $\phi$ be the corresponding Langlands parameter of $M$. Let $\delta^0\in Z^1(F,M^0)$ and let $\delta \in Z^1(F,M)$ be the corresponding pure inner form of $M$. Then
    \begin{itemize}
        \item for all $\tau^0 \in \widehat{A}_{\phi^0}[\delta^0]$,
        \[ \Ind_{M^0_{\delta^0}}^{M_{\delta}} \pi(\phi^0, \tau^0)= \bigoplus_{\tau \in \widehat{A}_{\phi}} m(\tau, \Ind_{A_{\phi^0}}^{{A_{\phi}}}\tau^0) \pi(\phi,\tau),\]
      where $m(\tau, \Ind_{A_{\phi^0}}^{{A_{\phi}}}\tau^0):= \dim(\Hom(\tau, \Ind_{A_{\phi^0}}^{{A_{\phi}}}\tau^0) ).$
        \item For all $\tau \in \widehat{A}_{\phi}[\delta]$,
\[  
        \Res_{M^0_{\delta^0}}^{M_{\delta}} \pi(\phi, \tau)= \bigoplus_{\phi^0_1} \bigoplus_{\tau^0_1 \in \widehat{A}_{\phi^0_1}} m(\tau_{1}^0, \Res_{A_{\phi^0_1}}^{{A_{\phi}}}\tau) \pi(\phi^0_1,\tau^0_1),
        \]
        where $\phi^0_1$ ranges over Langlands parameters of $M^0$ that maps to $\phi$, and 
        \[
        m(\tau_{1}^0, \Res_{A_{\phi^0_1}}^{{A_{\phi}}}\tau):= \dim(\Hom ( \tau^0_1, \Res_{A_{\phi^0_1}}^{{A_{\phi}}}\tau)).
        \]
    \end{itemize}

 \item [\Sha]
    (Shahidi conjecture) 
    Let  $\phi^0$ be a tempered Langlands parameter of $G^0$ and let $\phi$ be the corresponding Langlands parameter of $G$. Then any irreducible subquotient of $\Res_{G^0}^G \pi(\phi,\1)$ is a $\mfr{w}_{G^0}$-generic representation of $G^0(F)$, where $\1$ is the trivial representation of $A_\phi$. Moreover, if $A_{\phi} \neq A_{\phi^0}$, then $\pi(\phi,\1)$ is the extension of $\Res_{G^0}^G \pi(\phi,\1)$ using the Whittaker normalization. See Section~\ref{sec Whittaker normalization} for the definition of the Whittaker normalization.
\end{enumerate}
\end{desideratum}

\begin{remark}\label{rmk LLC(1)}
When the Levi subgroup $M$ is disconnected and $\theta[\delta^0] \neq [\delta^0]$, the right hand side of the displayed equation in \ECR\ is zero, by Corollary \ref{cor component group theta[delta]neq [delta]}. This is one reason why we define $ \Transfer{M',s}{M_{\delta}} \equiv 0$ for this case in Definition \ref{definition: modified transfer}. On the other hand, in this case we have $G_{\delta}(F) \simeq G^0_{\delta^0}(F)$, and \RES\  reduces the enhanced local Langlands correspondence of $G_{\delta}$ to that of the connected Levi subgroup $G^0_{\delta^0}$. 
\end{remark}

\begin{remark}\label{rmk LLC(2)}
    Let $\mathcal{H}_G$ be the minimal set of isomorphism classes of quasisplit groups that contains $G$ and is closed under taking endoscopic groups and Levi subgroups. We may state our Desiderata for all pair $(H,\delta)$ where $H \in \mathcal{H}_G$ and $\delta$ is a pure inner form of $H$, especially for \Sta \ and \ECR.
\end{remark}

\section{Transfer of stable standard distributions}\label{section: stable standards}

\subsection{Standard distributions}\label{subsection: standard distributions}

Throughout this section, $G$ satisfies Definition~\ref{assu: groups we consider}, and $\delta$ is a pure inner form of $G$ (Definition \ref{definition: pure inner form}).
Recall that in Definition~\ref{definition: LLC} we fixed a Whittaker normalized, local Langlands correspondence.
Thus, each $\pi \in \Pi_\phi(G_\delta)$ corresponds, uniquely, to an irreducible representation $\tau_\pi \in \widehat{A}_\phi[\delta]$.
For semisimple $s\in Z_{\dualgroup{G}}(\phi)$, we set
\IN{$\langle \cdot,\pi \rangle_\phi$, character of representation of component group of $L$-parameter}
\begin{equation}\label{equation: pairing phi}
    \langle s,\pi \rangle_\phi \coloneqq \trace(s,\tau_\pi), \qquad
    \LLC^\mathfrak{e}(\pi) = (\phi_\pi,\tau_\pi),
\end{equation}    
    where we again write $s$ for its image under the quotient $Z_{\dualgroup{G}}(\phi) \to Z_{\dualgroup{G}}(\phi)/Z_{\dualgroup{G}}(\phi)^\circ =: A_\phi$.

\begin{definition}
\label{definition: STD}
With $G$, $\delta$ and $\phi$ as above, and for $s\in Z_{\dualgroup{G}}(\phi)$, define $\STD{\phi,s}{\delta} \in K_\CC \Pi(G_{\delta})$ by
\IN{$\STD{\phi,s}{\delta}, \STD{\phi,s}{\pure}, \STD{\phi}{}, \STD{\phi}{\delta}$, standard distributions} 
\begin{equation*}
\STD{\phi,s}{\delta}
= 
\sum_{\pi\in \Pi_\phi(G_\delta)}
 \langle s,\pi \rangle_\phi\ \SR(\pi).
\end{equation*}
We also introduce the notation
\begin{equation*}
     \STD{\phi,s}{\pure} 
     \coloneqq
     \mathop{\sum}\limits_{\pi\in \Pi_\phi(G)^\pure} e(\pi) \langle s,\pi \rangle_\phi\ \SR(\pi) \in K_\CC \Pi_\phi(G)^\pure, 
\end{equation*}
where \IN{$e(\pi)$, Kottwitz sign of the group carrying $\pi$} $e(\pi)\ceq e(G_\delta)$ for $\pi \in \Pi(G_\delta)$; see Definition \ref{definition: kottwitz sign}. 
When $s=1$ or $[\delta]$ is trivial, we will generally omit it from the notation, thus defining $\STD{\phi}{\delta}$ and, in particular, 
\[
\STD{\phi}{}
\ceq
\sum_{\pi\in \Pi_\phi(G)} \rank(\tau_\pi) \SR(\pi).
\]
\end{definition}

Using the local Langlands correspondence for $G$, Definition~\ref{definition: STD} may equivalently be written in the form
\[
\mu_{\phi,s}^{\delta}
= 
\sum_{\sigma \in \widehat{A_\phi}[\delta]} \trace(s,\sigma)\ \SR(\phi, \sigma)
\]
and
\[
\mu_{\phi,s}^{\pure}
= 
\sum_{\sigma \in \widehat{A_\phi}} 
e(\sigma) \trace(s,\sigma)\ 
\SR(\phi, \sigma)
\]
where $e(\sigma)$ \IN{$e(\sigma)$, Kottwitz sign of $\sigma\in \widehat{A_\phi}$}
is Kottwitz sign of the image of $\sigma$ under the homomorphism of Proposition~\ref{proposition: enhancement by purity}, so that 
\[
e(\pi) = e(\sigma)
\qquad \text{for} \qquad \LLC^\mathfrak{e}(\pi) = (\phi,\sigma).
\]

\begin{remark}
When $\phi$ is tempered, $\STD{\phi}{\delta}$ is the distribution in Desiderata \ref{desideratum: LLC}\Sta. 
\end{remark}

\subsection{Transfer of stable standard distributions}\label{subsection: transfer of STD}

\begin{theorem}\label{theorem: transfer of STD}
Suppose that $G$ satisfies Definition~\ref{assu: groups we consider} and let $(G',s,\xi)$ be an endoscopic triple as defined in Section~\ref{subsection: endoscopic triple}.
Desiderata~\ref{desideratum: LLC} \ECR,\  \LC\ and \Sta\  imply the following.
    For any $\delta \in Z^1(F,G^0)$ and any Langlands parameter $\phi'$ for $G'$,
    \begin{equation*}
        \Transfer{G',s}{G_\delta}(\STD{\phi'}{}) = e(G_\delta)\ 
        \STD{\xi\circ\phi',{{s}}}{\delta}.
    \end{equation*}
\end{theorem}

\begin{proof}
When $\phi$ is tempered, this is a direct consequence of \ECR \ since tempered representations coincide with their standard modules. If $\phi$ is non-tempered, then the claim follows from the tempered case, \LC, and the compatibility of endoscopic transfer and parabolic induction (Proposition~\ref{proposition: induction commutes with transfer}).
\end{proof}

\section{Some results on endoscopic transfer}\label{section: Bin} 

We conclude Part~\ref{part: II Transfer} by proving the following two results which require a more in-depth discussion of endoscopic transfer.
First, we compute twisted endoscopic transfer in a specific setting.

\begin{proposition}\label{prop twisted transfer from Siegel GL to disconencted GL}
Let $G^0 =\GL_m \times  \GL_m$ and let $\theta$ be the involution of $G^0$ given by
\[ \theta(g_1, g_2) = ({}^t g_2^{-1}, {}^t g_1^{-1}).\]
It determines an endoscopic triple $(\GL_m, \hat{\theta}, \xi)$ for $G = (\GL_m\times\GL_m)\rtimes\langle\theta\rangle$.
For any irreducible representation $\sigma$ of $\GL_m$, we have 
\[ \Transfer{\GL_m, \hat{\theta}}{G} (\sigma) = (\sigma \boxtimes \sigma^{\vee})^+,\]
    where $(\sigma \boxtimes \sigma^{\vee})^+$ is the Whittaker normalization introduced in Section \ref{sec Whittaker normalization}.
\end{proposition}

Second, we show that stable standard distributions form a basis for the space of stable distributions supported on  $K_{\BC}\Pi_{\lambda}(G)$.

\begin{theorem}\label{thm: stable standards form a basis - quasisplit}
Let $G$ be a connected quasisplit group equipped with a local Langlands correspondence satisfying Desiderata~\ref{desideratum: LLC}\ \TW, \Sta, \LC, \CTD, \CCh, \ECR, \TI\ and \Hai.
The set of stable distributions $\{\mu_{\phi} \tq \phi \in \Phi_{\lambda}(G)\}$ forms a basis of stable distributions supported on $K_{\BC}\Pi_{\lambda}(G)$.
\end{theorem}

\subsection{A particular twisted endoscopic transfer}
\label{subsection: twisted transfer Bin}

Let $F$ be a $p$-adic field and $\overline{F}$ be its algebraic closure. Let $H = GL_m$ over $F$ and $G^0 = H \times H$. Let $\tau$ be the involution of $H$ given by transpose inverse. Let $\theta$ be the involution of $G^0$ defined by
\[
\theta(x, y) = (\tau(y), \tau(x)).
\]
Let $G = G^0 \rtimes \langle \theta \rangle$. Let $\hat{\tau}$ be the dual automorphism on $\widehat{H}$. Then the embedding
\[
\widehat{H} \rightarrow \widehat{H} \times \widehat{H} = \widehat{G^0}, \quad h \mapsto (h, \hat{\tau}(h))
\]
extends to an L-embedding 
\[
\xi: {}^LH \rightarrow {}^LG,
\]
such that $(H, \hat{\theta}, \xi)$ forms an endoscopic triple of $G$. 

For any irreducible smooth representation $(\pi, V)$ of $H(F)$, we denote by $(\pi^{\tau}, V)$ the representation of $H(F)$ defined by 
\[
\pi^{\tau}(y)(v) : = \pi(\tau(y))(v), \quad y \in H(F), v \in V.
\] 
Let $\Pi := \pi \boxtimes \pi^{\tau}$ be the irreducible representation of $G^0(F)$ on $V \otimes V$. It is $\theta$-invariant with a natural intertwining operator
\[
A_{\theta}: V \otimes V \rightarrow V \otimes V, \quad v \otimes w \mapsto w \otimes v.
\]
This gives an extension of $\Pi$ to $G(F)$, denoted by $\Pi^{+}$. The Harish-Chandra distribution character of $\pi$ is 
\[
\Theta^{H}(\pi)(f) = {\rm Tr} \, \pi(f), \quad f \in \mathcal{C}^{\infty}_{c}(H(F)).
\]
It can be represented by a locally integrable function $\Theta^{H}(\pi)$ on $H(F)$, which is smooth over the regular semisimple elements of $H(F)$. We take $\theta$-twisted Harish-Chandra distribution character of $\Pi$ to be
\[
\Theta^{G}(\Pi^{+})_{\theta}(f) = {\rm Tr} \, (\Pi(f) \circ A_{\theta}), \quad f \in \mathcal{C}^{\infty}_{c}(G^0(F)).
\]
It can be represented by a locally integrable function $\Theta^{G}(\Pi^{+})_{\theta}$ on $G(F)$, which  is smooth over the $\theta$-regular $\theta$-semisimple elements of $G^0(F)$. 
We would like to show 
\begin{align}
\label{eq: twisted endoscopic transfer GL}
{\rm Trans}_{H, \hat{\theta}}^{G^0, \theta} \, \Theta^{H}(\pi) = \Theta^{G}(\Pi^{+})_{\theta}.
\end{align}
Note Proposition~\ref{prop twisted transfer from Siegel GL to disconencted GL} follows from this.

Define the norm map by
\[
{\rm Nm}: G^0 \rightarrow H, \quad (x, y) \mapsto x \tau(y).
\]
It is surjective, and the fibre over any $x \in H(\overline{F})$ is $\{(x\tau(y), y^{-1}) \, | \, y \in H(\overline{F})\} \cong H(\overline{F})$. For any $\delta = (x, y) \in G^0(F)$, let $\gamma = {\rm Nm}(\delta)$. Define
\[
H_{\gamma} := \{h \in H(\overline{F}) \, | \, h^{-1} \gamma h = \gamma \}
\]
and
\[
(G^{0})^{\theta}_{\delta} := \{g \in G^0(\overline{F}) \, | \, g^{-1}\delta \theta(g) = \delta \}
\]
as algebraic groups over $F$. Then there is an isomorphism of algebraic groups
\[
H_{\gamma} \xrightarrow{\simeq} (G^0)^{\theta}_{\delta}, \quad h \mapsto  (h, y\tau(h)y^{-1}).
\]
Hence, $\delta$ is $\theta$-regular $\theta$-semisimple if and only if $\gamma$ is regular semisimple. At last, \eqref{eq: twisted endoscopic transfer GL} can be restated as in the following Proposition.

\begin{proposition}
\label{prop: character relation p-adic}
For any $\theta$-regular $\theta$-semisimple elements $g \in G^0(F)$,  
\[
\Theta^{H}(\pi)({\rm Nm}(g)) = \Theta^{G}(\Pi^{+})_{\theta}(g).
\]
\end{proposition}

We will prove Proposition~\ref{prop: character relation p-adic} in the rest of this section. Define
\[
\mathcal{T}: \mathcal{C}^{\infty}_{c}(G^0(F)) \rightarrow \mathcal{C}^{\infty}_{c}(H(F))
\]
by
\[
(\mathcal{T}f)(x) = \int_{H(F)}f(x\tau(y), y^{-1})dy.
\]
The following lemma shows that $\mathcal{T}$ can be extended to a map $L^{1}(G^0(F)) \rightarrow L^{1}(H(F))$.

\begin{lemma}
\label{lemma: change variable p-adic}
For any $f \in L^{1}(G^0(F))$, the integral defining
\(
(\mathcal{T}f)(x)
\)
is absolutely convergent for almost all $x$ and $\mathcal{T}f \in L^{1}(H(F))$. Moreover,
\[
\int_{G^0(F)} f(g) dg = \int_{H(F)} (\mathcal{T}f)(h)dh.
\]
\end{lemma}

\begin{proof}
The map
\[
\varphi: H(F) \times H(F) \rightarrow G^0(F) = H(F) \times H(F), \quad (x, y) \mapsto (x\tau(y), y^{-1}),
\]
is a change of variables, which preserves the measure. Since $f \in L^{1}(G^0(F))$, then so is $f \circ \varphi$. It follows that
\[
\int_{G(F)} f(g)dg = \int_{H(F) \times H(F)} (f \circ \varphi)(x, y)dxdy =  \int_{H(F)} \Big(\int_{H(F)} (f \circ \varphi)(x, y)dy \Big)dx,
\]
where 
\[
(\mathcal{T}f)(x) = \int_{H(F)} (f \circ \varphi)(x, y)dy
\]
is absolutely convergent for almost all $x$ and defines a function in $L^{1}(H(F))$. This proves the lemma.
\end{proof}

\begin{lemma}
\label{lemma: distribution relation p-adic}
For any $f \in \mathcal{C}^{\infty}_{c}(G^0(F))$, 
\[
\Theta^{H}(\pi)(\mathcal{T}f) = \Theta^{G}(\Pi^{+})_{\theta}(f).
\]
\end{lemma}

\begin{proof}
For any $f \in \mathcal{C}^{\infty}_{c}(G^0(F))$, there exists an open compact subgroup $K$ of $H(F)$ such that the image of $\Pi(f)$ lies in $V^{K} \otimes V^{K}$ and the image of $\pi((\mathcal{T}f))$ lies in $V^{K}$. Since $\mathcal{C}^{\infty}_{c}(G^0(F)) \cong \mathcal{C}^{\infty}_{c}(H(F)) \otimes \mathcal{C}^{\infty}_{c}(H(F))$, by making $K$ smaller we can further assume that for any $x \in H(F)$ and $v \in V$,
\[
\int_{H(F)} f(x, y)\pi^{\tau}(y)v \, dy
\]
lies in $V^{K}$. Let $(\tilde{\pi}, \widetilde{V})$ be the smooth dual of $\pi$. Since $\pi$ is admissible, then so is $\tilde{\pi}$. Hence, both $V^{K}$ and $\widetilde{V}^{K}$ are finite dimensional. We choose a basis $\{e_i\}_{i \in I}$ of $V^{K}$ and let $\{\tilde{e}_i\}_{i \in I}$ be the dual basis in $\widetilde{V}^{K}$ under the restriction of the natural pairing
$
\langle \, ,  \, \rangle : V \times \widetilde{V} \rightarrow \BC.
$
Then
\begin{align*}
\Theta^{G}(\Pi^{+})_{\theta}(f) & = {\rm Tr} \, (\Pi(f) \circ A_{\theta}) = \sum_{i, j} \langle \Pi(f) \circ A_{\theta} (e_i \otimes e_j), \tilde{e}_i \otimes \tilde{e}_j \rangle \\
& = \sum_{i, j} \langle \Pi(f) (e_j \otimes e_i), \tilde{e}_i \otimes \tilde{e}_j \rangle \\
& = \sum_{i, j} \int_{H(F) \times H(F)} f(x, y) \langle \pi(x)e_{j}, \tilde{e}_i \rangle \cdot \langle \pi^{\tau}(y) e_{i}, \tilde{e}_j \rangle dxdy \\
& = \sum_{i, j} \int_{H(F)}  \langle e_{j}, \tilde{\pi}(x^{-1})\tilde{e}_i \rangle \Big( \int_{H(F)}f(x,y)\langle \pi^{\tau}(y) e_{i}, \tilde{e}_j \rangle dy \Big)dx
\end{align*}
By our assumption, 
\[
w_i(x) := \int_{H(F)} f(x,y)\pi^{\tau}(y)e_i \in V^{K} \quad (i \in I).
\]
Hence,
\[
w_i(x) = \sum_{j} \langle w_i(x), \tilde{e}_j \rangle e_j = \sum_{j} \Big(\int_{H(F)} f(x,y) \langle \pi^{\tau}(y)e_i, \tilde{e}_j \rangle dy\Big) e_j.
\]
Substitute this into the previous expansion, we get
\begin{align*}
\Theta^{G}(\Pi^{+})_{\theta}(f) & = \sum_{i} \int_{H(F)} \langle w_i(x), \tilde{\pi}(x^{-1})\tilde{e}_i \rangle dx \\
& = \sum_{i} \int_{H(F)} \Big( \int_{H(F)} f(x, y) \langle \pi^{\tau}(y) e_{i}, \tilde{\pi}(x^{-1})\tilde{e}_i \rangle dy \Big) dx  \\
& = \sum_{i} \int_{H(F)} \Big( \int_{H(F)} f(x, y) \langle \pi(x\tau(y)) e_{i}, \tilde{e}_i \rangle dy \Big) dx  \\
& = \sum_{i} \int_{H(F) \times H(F)} f(x, y) \langle \pi(x\tau(y)) e_{i}, \tilde{e}_i \rangle dxdy \\
&= \sum_{i} \int_{G^{0}(F)} f(g) \langle \pi({\rm Nm}(g)) e_{i}, \tilde{e}_i \rangle dg \\
& = \sum_{i} \int_{H(F)} (\mathcal{T}f)(x) \langle \pi(x)e_{i}, \tilde{e}_i \rangle dx \\
& = {\rm Tr} \, \pi(\mathcal{T}f) \\
& = \Theta^{H}(\pi)(\mathcal{T}f) .
\end{align*}
The first equality of the last line follows from Lemma~\ref{lemma: change variable p-adic} and the fact that $\mathcal{T}$ sends the function to
\(
f(g) \langle \pi({\rm Nm}(g)) e_{i}, \tilde{e}_i \rangle
\)
to $(\mathcal{T}f)(x) \langle \pi(x)e_{i}, \tilde{e}_i \rangle$. 
\end{proof}

Now we are ready to give the proof of Proposition~\ref{prop: character relation p-adic}.

\begin{proof}
By definition,
\[
\Theta^{H}(\pi)(f) = \int_{H(F)} f(h) \Theta^{H}(\pi)(h)dh, \quad f \in \mathcal{C}^{\infty}_{c}(H(F))
\]
and
\[
\Theta^{G}(\Pi^{+})_{\theta}(f) = \int_{G^0(F)} f(g) \Theta^{G}(\Pi^{+})_{\theta}(g) dg, \quad f \in \mathcal{C}^{\infty}_{c}(G^0(F)).
\]
By Lemma~\ref{lemma: change variable p-adic}, for any $f \in \mathcal{C}^{\infty}_{c}(G^0(F))$,
\[
\Theta^{H}(\pi)(\mathcal{T}f)  = \int_{H(F)} (\mathcal{T}f)(h) \Theta^{H}(\pi)(h)dh = \int_{G^0(F)} f(g) \Theta^{H}(\pi)({\rm Nm}(g)) dg,
\]
where $\mathcal{T}$ sends 
\[
f(g)\Theta^{H}_{\pi}({\rm Nm}(g)) \in L^{1}(G^0(F))
\]
to 
\[
(\mathcal{T}f)(h) \Theta^{H}_{\pi}(h) \in L^{1}(H(F)).
\]
By Lemma~\ref{lemma: distribution relation p-adic}, we get
\[
\int_{G^0(F)} f(g) \Theta^{G}(\Pi^{+})_{\theta}(g) dg = \int_{G^0{F}} f(g) \Theta^{H}(\pi)({\rm Nm}(g)) dg, \quad f \in \mathcal{C}^{\infty}_{c}(G^0(F)).
\]
Hence, 
\[
\Theta^{G}(\Pi^{+})_{\theta}(g) = \Theta^{H}(\pi)({\rm Nm}(g))
\]
for all $\theta$-regular $\theta$-semisimple elements of $g \in G^0(F)$.
\end{proof}

\subsection{A basis for stable distributions}\label{subsection: basis of stable standards Bin}

Let $F$ be a $p$-adic field and $G$ be a quasisplit connected reductive algebraic group over $F$. We will not consider disconnected groups and their pure inner forms in this subsection. We fix an $F$-pinning $(B, T, \{X_{\alpha}\})$ of $G$ and a nontrivial additive character of $F$ throughout this subsection. This gives rise to a Whittaker datum $\mathfrak{m}$ of $G(F)$. The goal is to prove Theorem~\ref{thm: stable standards form a basis - quasisplit} based on some standard hypotheses on the local Langlands correspondence. More precisely,  we assume Desiderata~\ref{desideratum: LLC}\  \Sta, \LC, \CTD, \CCh, \ECR, \TI\ and \Hai.

Since we only consider quasisplit groups in this subsection, we would like to simplify some notations. For $\phi \in \Phi(G)$, let $\mathcal{S}_{\phi} := Z_{\widehat{G}}(\phi)/Z_{\widehat{G}}(\phi)^{0}Z(\widehat{G})^{\Gamma_{F}}$. There is an exact sequence
\[
 \pi_0(Z(\widehat{G})^{\Gamma_F}) \rightarrow A_{\phi} \rightarrow \mathcal{S}_{\phi} \rightarrow 1.
\]
By our definition, we can identify $\widehat{\mathcal{S}_{\phi}} = \widehat{A_{\phi}}[1]$.

We will denote $\mathcal{C}^{\infty}_c(G(F))$ by $\mathcal{H}(G)$ as in \cite{Art96}. Let $I(G)$ (resp. $SI(G)$) be the space of normalized (resp. stable) orbital integrals of $\mathcal{H}(G)$ over the strongly regular semisimple elements of $G(F)$. By definition, there are natural surjections 
\[
\mathcal{H}(G) \rightarrow I(G) \xrightarrow{\mathcal{T}^G} SI(G).
\]
Denote the set of conjugacy classes (resp. stable conjugacy classes) of strongly regular semisimple elements of $G(F)$ by $\Gamma(G)$ (resp. $\Sigma(G)$). Then $I(G)$ (resp. $SI(G)$) can be viewed as a subspace of functions over $\Gamma(G)$ (resp. $\Sigma(G)$). Note that there is a natural surjection
\[
\Gamma(G) \rightarrow \Sigma(G).
\]
For any $a_{G} \in SI(G)$, we denote its pullback to $\Gamma(G)$ again by $a_{G}$. For $\gamma \in \Gamma(G)$ and the corresponding $\sigma \in \Sigma(G)$, $a_{G}(\gamma) := a_{G}(\sigma)$. Let $n(\sigma)$ be the number of conjugacy classes in $\sigma$. Then 
\[
\gamma \mapsto n(\sigma)^{-1}a_{G}(\gamma), \quad \gamma \in \Gamma(G)
\]
defines an element in $I(G)$. In this way, we get a section of $\mathcal{T}^{G}$. Let $I_{cusp}(G)$ (resp. $SI_{cusp}(G)$) be the subspace of functions supported on the elliptic elements $\Gamma_{ell}(G)$ (resp. $\Sigma_{ell}(G)$). From the construction of the section of $\mathcal{T}^{G}$, we also get a surjection
\[
\mathcal{T}^{G}_{cusp}: I_{cusp}(G) \rightarrow SI_{cusp}(G)
\]
by restriction. Let $\mathcal{H}_{cusp}(G)$ be the preimage of $I_{cusp}(G)$ in $\mathcal{H}(G)$.

Let 
\[
T(G) = \{(L, \pi, r) \, | \, L \in \mathcal{L}^{G}, \pi \in \Pi_2(L), r \in R^{G}_{\pi}\} / W^{G},
\]
where $\mathcal{L}^{G}$ is the set of Levi subgroups of $G$ containing $T$ and $R^{G}_{\pi}$ is the representation theoretic $R$-group. There is a decomposition
\[
T(G) = \bigsqcup_{\{M\}} \Big( T_{ell}(M) / W(M) \Big).
\]
Here $M$ ranges over $G(F)$-conjugacy classes of Levi subgroups of $G$, $W(M)$ is the relative Weyl group and
\[
T_{ell}(M) = \{(L, \pi, r) \, | \, L \in \mathcal{L}^{M}, \pi \in \Pi_2(L), r \in R^{M}_{\pi, reg}\} / W^{M},
\]
where $\mathcal{L}^{M}$ is the set of Levi subgroups of $M$ containing $T$, and 
\[
R^{M}_{\pi, reg} = \{r \in R^{M}_{\pi} \, | \, \mathfrak{a}^{r}_L = \mathfrak{a}_M \}.
\] 

For a character $\chi$ of a central subgroup $Z$ of $G(F)$, we can define spaces of $\chi^{-1}$-equivariant functions
\[
\mathcal{H}(G, \chi), \, \mathcal{H}_{cusp}(G, \chi), \, I(G, \chi), \, I_{cusp}(G, \chi), \, SI(G, \chi), \, SI_{cusp}(G, \chi), \]
analogously and also the subsets 
\[
T(G, \chi), \, T_{ell}(G, \chi)
\] 
by requiring the restriction of the central character of $\pi$ to $Z$ is equal to $\chi$. Now we let $\chi$ be a character of the center $Z_{G}(F)$. For $\tau = (L, \pi, r) \in T(G, \chi)$, one can associates a distribution character
\[
\Theta(\tau)(f) := {\rm Tr} (R_{P}(r, \pi) \circ I_{P}(\pi, f)), \quad f \in \mathcal{H}(G, \chi),
\]
where $P$ is any Levi subgroup of $G$ with Levi factor $L$, $I_{P}(\pi, f)$ is the action of $f$ on the parabolically induced representation $\Ind^{G}_{P}(\pi)$ and $R_{P}(r, \pi)$ is certain self-intertwining operator on this induced representation, which is uniquely determined up to some scalar. It is representable by some locally integrable function, still denoted by $\Theta(\tau)$. Let
\[
I(\tau, \gamma) := |D_{G}(\gamma)|^{1/2} \Theta(\tau)(\gamma), \quad \gamma \in \Gamma(G),
\]
where $D_{G}(\gamma)$ is the Weyl discriminant of $\gamma$.

For $\phi \in \Phi_2(G, \chi)$, let 
\[
\Theta(\phi) := \sum_{\rho \in \widehat{\mathcal{S}_{\phi}}} {\rm dim}(\rho) \, \Theta(\pi(\phi, \rho))
\]
and
\[
S(\phi, \gamma) := |\mathcal{S}_{\phi}|^{-1}|D_{G}(\gamma)|^{1/2} \Theta_{\phi}(\gamma), \quad \gamma \in \Gamma(G).
\]
Since $\Theta_{\phi}$ is stable by Desideratum~\ref{desideratum: LLC} \Sta, then $S(\phi, \gamma)$ only depends on the stable conjugacy class $\sigma$ of $\gamma$, and we will write $S(\phi, \sigma) := S(\phi, \gamma)$. Let $a_{G}(\tau, \gamma) := I(\tau, \gamma)$ for $\gamma \in \Gamma_{ell}(G)$. By \cite{Art93}*{Section 6}, we can normalize $R_{P}(r, \pi)$ for $\Theta(\tau)$ so that $\{a_{G}(\tau, \gamma)\tq \tau \in T_{ell}(G, \chi)\}$ forms an orthonormal basis of $I_{cusp}(G, \chi)$. Their images under the $\chi^{-1}$-equivariant analogue of $\mathcal{T}^{G}_{cusp}$
\[
\mathcal{T}^{G}_{cusp, \chi}: I_{cusp}(G, \chi) \rightarrow SI_{cusp}(G, \chi)
\]
span $SI_{cusp}(G, \chi)$. By definition,
\[
S(\phi, \gamma) = |\mathcal{S}_{\phi}|^{-1} \sum_{\rho \in \widehat{\mathcal{S}_{\phi}}} {\rm dim}(\rho) \, a_{G}(\pi(\phi, \rho), \gamma), \quad \gamma \in \Gamma_{ell}(G).
\]
By Desideratum~\ref{desideratum: LLC} ($\mathbf{CTD}$) and ($\mathbf{CC}$), it lies in $I_{cusp}(G, \chi)$. Its image in $SI_{cusp}(G, \chi)$ is
\[
a_{G}(\phi, \sigma) := n(\sigma) S(\phi, \sigma), \quad \sigma \in \Sigma_{ell}(G).
\]

\begin{lemma}
$\{a_{G}(\phi, \sigma)\tq \phi \in \Phi_2(G, \chi)\}$ is orthonormal.
\end{lemma}

\begin{proof}
By definition (cf. \cite[]{Art96}),
\begin{align*}
\langle a_{G}(\phi_1, \sigma), a_{G}(\phi_2, \sigma) \rangle_{\Sigma_{ell}(G)} & = \int_{\Sigma_{ell}(G)} n(\sigma)^{-1} a_{G}(\phi_1, \sigma) \overline{a_{G}(\phi_2, \sigma)} d\sigma \\
& = \int_{\Gamma_{ell}(G)} n(\sigma)^{-2} a_{G}(\phi_1, \gamma) \overline{a_{G}(\phi_2, \gamma)} d\gamma \\
& = \int_{\Gamma_{ell}(G)} S(\phi_1, \gamma) \overline{S(\phi_2, \gamma)} d\gamma.
 \end{align*} 
Then it follows from the fact that $\{a_{G}(\tau, \gamma)\}_{\tau \in T_{ell}(G, \chi)}$ is orthonormal and 
\[
|\mathcal{S}_{\phi}|^2 = \sum_{\rho \in \widehat{\mathcal{S}_{\phi}}} ({\rm dim}(\rho))^2. \qedhere
\]
\end{proof}

\begin{lemma}
For $\phi \in \Phi_2(G, \chi)$ and $\tau \in T_{ell}(G, \chi)$,
\begin{align*}
\langle a_{G}(\phi, \sigma), \mathcal{T}^{G}_{cusp, \chi}(a_{G}(\tau, \gamma)) \rangle_{\Sigma_{ell}(G)} = \begin{cases} {\rm dim}(\rho)|\mathcal{S}_{\phi}|^{-1} \quad & \text{ if } \tau = \pi(\phi, \rho) \text{ for some } \rho\\
0 \quad & \text{ otherwise. }
\end{cases} 
\end{align*}
\end{lemma}

\begin{proof}
It follows from
\begin{align*}
\langle a_{G}(\phi, \sigma), \mathcal{T}^{G}_{cusp, \chi}(a_{G}(\tau, \gamma)) \rangle_{\Sigma_{ell}(G)} & = \langle n(\sigma)S(\phi, \sigma), \mathcal{T}^{G}_{cusp, \chi}(a_{G}(\tau, \gamma) \rangle_{\Sigma_{ell}(G)} \\
&= \langle S(\phi, \gamma), a_{G}(\tau, \gamma) \rangle_{\Gamma_{ell}(G)}.
\end{align*}
and the fact that $\{a_{G}(\tau, \gamma)\vert \tau \in T_{ell}(G, \chi)\}$ forms an orthonormal basis of $I_{cusp}(G, \chi)$.
\end{proof}

\begin{proposition}
\label{prop: elliptic stable basis}
$\{a_{G}(\phi, \sigma) \vert\phi \in \Phi_2(G, \chi)\}$ is an orthonormal basis of $SI_{cusp}(G, \chi)$.
\end{proposition}

\begin{proof}
It remains to show that they span $SI_{cusp}(G, \chi)$. If not, then we can find $0 \neq b_{G} \in SI_{cusp}(G, \chi)$ orthogonal to $a_{G}(\phi, \sigma)$ for all $\phi \in \Phi_2(G, \chi)$. The endoscopic transfer gives an isomorphism
\[
I_{cusp}(G, \chi) \xrightarrow{\sim} I^{\mathcal{E}}_{cusp}(G, \chi) := \bigoplus_{G' \in \mathcal{E}_{ell}(G)} SI_{cusp}(G', \chi')^{{\rm Out}_G(G')},
\] 
$\chi'$ is a character of some central subgroup $Z'$ of $G'(F)$ determined by $\chi$ and ${\rm Out}_G(G')$ is the outer automorphism group of the endoscopic datum represented by $G'$ (see \cite{Art96}*{Proposition 3.5}). As a result, we can choose $f \in \mathcal{H}_{cusp}(G, \chi)$ such that the image of the endoscopic transfer $f^{G'}$ in $SI_{cusp}(G', \chi')$ is zero for $G \neq G'$ and the image of $f$ in $SI_{cusp}(G, \chi)$ is $b_{G}$. Since $b_G \neq 0$, the image of $f$ in $I_{cusp}(G, \chi)$ is nonzero, hence there exists $\tau \in T_{ell}(G, \chi)$ such that $\Theta^{G}(\tau)(f) \neq 0$. By the inversion formula of the endoscopic character relations (cf. Desideratum~\ref{desideratum: LLC} \ECR), we can write 
\begin{align}
\label{eq: elliptic stable basis}
\Theta^{G}(\tau)(f) = \sum_{G' \in \mathcal{E}_{ell}(G)} \sum_{\phi' \in \Phi_{bdd}(G')} c(G', \phi') \Theta^{G'}(\phi')(f^{G'})
\end{align}
for some constants $c(G', \phi') \in \BC$. For $G \neq G'$, $\Theta^{G'}(\phi')(f^{G'}) = 0$ . By Desideratum~\ref{desideratum: LLC} ($\mathbf{TI}$), we also know that $\Theta^G(\phi) =0$ for any $\phi \in \Phi_{bdd}(G) \backslash \Phi_2(G)$. At last, for $\phi \in \Phi_2(G, \chi)$, 
\[
\Theta^{G}(\phi)(f) = |\mathcal{S}_{\phi}| \, \langle b_G(\sigma), a_{G}(\phi, \sigma)\rangle_{\Sigma_{ell}(G)} = 0
\]
by our assumption. Hence, the right hand side of \eqref{eq: elliptic stable basis} is zero, which leads to a contradiction.
\end{proof}

Let $\lambda$ be an infinitesimal parameter of $G$ and $V_{\lambda}$ be the corresponding Vogan variety. Let $\Phi_{\lambda}(G)$ be the set of Langlands parameters with infinitesimal parameter $\lambda$ and $\Pi_{\lambda}(G)$ be the set of isomorphism classes of irreducible admissible representations of $G(F)$ with infinitesimal parameter $\lambda$. Let $K_{\BC}\Pi_{\lambda}(G)$ be the subspace of $K_{\BC}\Pi(G)$ spanned by $\Pi_{\lambda}(G)$. For $\phi \in \Phi(G)$ and $\rho \in \widehat{\mathcal{S}_{\phi}}$, let $\SR(\phi, \rho)$ be the standard representation that contains $\pi(\phi, \rho)$ as the Langlands quotient. The standard stable distribution associated with $\phi$ has been defined to be
\[
\Theta^{G}(\phi) = \sum_{\rho \in \widehat{\mathcal{S}_{\phi}}} {\rm dim}(\rho) \, \Theta^{G}(\SR(\phi, \rho)).
\]
Now are ready to give the proof of Theorem~\ref{thm: stable standards form a basis - quasisplit}.

\begin{proof}[Proof of Theorem~\ref{thm: stable standards form a basis - quasisplit}]
Since parabolic induction preserves infinitesimal parameters by Desideratum~\ref{desideratum: LLC} \Hai, $\SR(\phi, \rho)$ belongs to $K_{\BC}\Pi_{\lambda}(G)$ and they form a basis  by \cite{BorelWallach:2000}*{XI, Lemma 2.13}. The linear independence of $\{\Theta^{G}(\phi)\tq \phi \in \Phi_{\lambda}(G)\}$ follows from this by the linear independence of irreducible characters. Next we shall show that they span the space of stable distributions supported on $K_{\BC}\Pi_{\lambda}(G)$. Suppose 
\[
S^{G} := \sum_{\pi \in \Pi_{\lambda}(G)} c_{\pi} \Theta^{G}(\pi).
\] 
is a nonzero stable distribution. For any $\pi \in \Pi_{\lambda}(G)$, we have
\[
\Theta^{G}(\pi) = \sum_{\phi \in \Phi_{\lambda}(G)} \, \sum_{\rho \in \widehat{\mathcal{S}_{\phi}}} m^{\pi}_{(\phi, \rho)} \, \Theta^{G}(\SR(\phi, \rho)).
\]
Then
\begin{align*}
S^{G} & = \sum_{\pi \in \Pi_{\lambda}(G)} c_{\pi} \Big( \sum_{\phi \in \Phi_{\lambda}(G)} \, \sum_{\rho \in \widehat{\mathcal{S}_{\phi}}} m^{\pi}_{(\phi, \rho)} \, \Theta^{G}(\SR(\phi, \rho))\Big) \\
& = \sum_{\phi, \rho, \pi} c_{\pi} \, m^{\pi}_{(\phi, \rho)} \Theta^{G}(\SR(\phi, \rho)).
\end{align*}
By the inversion formula of endoscopic character relations (cf. Desideratum~\ref{desideratum: LLC} \ECR),
\[
\Theta^{G}(\SR(\phi, \rho))(f) = \sum_{(G', \phi')} m^{(\phi, \rho)}_{\phi'} \Theta^{G'}(\phi')(f^{G'}), \quad f \in \mathcal{H}(G)
\]
where $f^{G'}$ is the endoscopic transfer of $f$, $\phi$ factors through $\phi' \in \Phi(G')$ for some elliptic endoscopic group $G'$ of $G$ and $\phi'$ does not factor through any common Levi subgroup of ${}^LG'$ and ${}^LG$ unless $\phi' = \phi$, in which case $m^{(\phi, {\rho})}_{\phi} = |\mathcal{S}_{\phi}|^{-1}$. Then we can further write
\[
S^{G}(f) = \sum_{G'} \sum_{\phi'} \Big( \sum_{\rho, \pi} m^{(\phi, \rho)}_{\phi'} m^{\pi}_{(\phi, \rho)} c_{\pi} \Big) \Theta^{G'}(\phi')(f^{G'}).
\]
We separate the stable distributions on $G(F)$ from those on $G'(F)$ ($G' \neq G$) as follows
\begin{align}
\label{eq: stable}
&S^{G}(f) - \sum_{\phi \in \Phi_{\lambda}(G)} |\mathcal{S}_{\phi}|^{-1} \Big( \sum_{\rho, \pi} m^{\pi}_{(\phi, \rho)} c_{\pi} \Big) \Theta^{G}(\phi)(f) \\
&= \sum_{G' \neq G} \sum_{\phi'} \Big( \sum_{\rho, \pi} m^{(\phi, \rho)}_{\phi'} m^{\pi}_{(\phi, \rho)} c_{\pi} \Big) \Theta^{G'}(\phi')(f').
\end{align}
Suppose the left hand side defines a nonzero distribution, then we would like to make a contradiction by plugging some $f \in \mathcal{H}(G)$ such that the left hand is nonzero and the right hand side is zero. The existence of such test functions follows from Arthur's spectral characterization of the image of $\mathcal{H}(G)$ under the endoscopic transfer (cf. \cite{Art96}), which will be recalled afterwards. As a result,
\[
S^{G} = \sum_{\phi \in \Phi_{\lambda}(G)} |\mathcal{S}_{\phi}|^{-1} \Big( \sum_{\rho, \pi} m^{\pi}_{(\phi, \rho)} c_{\pi} \Big) \Theta^{G}(\phi).
\]
This finishes the proof.
\end{proof}

We would like to review Arthur's spectral characterization of the image of $\mathcal{H}(G)$ under the endoscopic transfer (cf. \cite{Art96}). Define $I^{\mathcal{E}}(G)$ to be the subspace of 
\[
\bigoplus_{G' \in \mathcal{E}_{ell}(G)} SI(G')^{{\rm Out}_{G}(G')},
\]
consisting of $a^{\mathcal{E}}_{G} = (a_{G'})_{G'}$ such that $a_{G'_1}, a_{G'_2}$ have the same image in $SI(M')$ under the descent map if $M'$ is a common Levi subgroup of $G'_1, G'_2 \in \mathcal{E}_{ell}(G)$. The Paley-Wiener space on $T(G)$ can be written as
\[
PW(T(G)) = \bigoplus_{\{M\}} PW(T_{ell}(M))^{W(M)}.
\]
Let
\[
T^{\mathcal{E}}(G) = \bigsqcup_{\{M\}} \Big( T^{\mathcal{E}}_{ell}(M) / W(M) \Big)
\]
where
\[
T^{\mathcal{E}}_{ell}(M) = \{ (M' , \phi'_M) \, | \, M' \in \mathcal{E}_{ell}(M), \phi'_M \in \Phi_{2}(M')/{\rm Out}_{M}(M')\}.
\]
In \cite{Art96}, the local Langlands correspondence is not assumed and $\Phi_2(G, \chi)$ is defined simply as an index set for some chosen basis of $SI_{cusp}(G, \chi)$, which is invariant under twists by outer automorphisms and unramified characters fixing the central character $\chi$. For $\phi \in \Phi_2(G, \chi)$, one can define a distribution $\Theta^{G}(\phi)$ on $SI_{cusp}(G, \chi)$ by integration against the associated basis element over $\Sigma_{ell}(G)$. Then  
\[
\Theta^{G}(\phi)(f) = \sum_{\tau \in T_{ell}(G, \chi)} m(\phi, \tau) \Theta^{G}(\tau)(f), \quad f \in \mathcal{H}_{cusp}(G),
\]
where the right hand side is a finite sum. One of the main results of \cite{Art96} is that the right hand side defines a stable distribution on $\mathcal{H}(G)$. By Proposition~\ref{prop: elliptic stable basis} and Desideratum~\ref{desideratum: LLC} ($\mathbf{Tw}$), we can choose this basis to be $\{a_{G}(\phi, \sigma)\}_{\phi \in \Phi_2(G, \chi)}$ under the local Langlands correspondence. Similarly, the stable distributions with $\phi'_M \in \Phi_{2}(M')$ in \cite{Art96} can be chosen to be those given by the local Langlands correspondence.

The Paley-Wiener space on $T^{\mathcal{E}}(G)$ is  written as
\[
PW(T^{\mathcal{E}}(G)) = \bigoplus_{\{M\}} PW(T^{\mathcal{E}}_{ell}(M))^{W(M)}.
\]
For $G' \in \mathcal{E}_{ell}(G)$, 
\[
\Phi_{bdd}(G') = \bigsqcup_{\{M'\}} \Big( \Phi_2(M') / W(M') \Big).
\]
Similarly, the Paley-Wiener space on $\Phi_{bdd}(G')$ is written as
\[
PW(\Phi_{bdd}(G')) = \bigoplus_{\{M'\}} PW(\Phi_2(M'))^{W(M')}.
\]
Since a Levi subgroup $M'$ of $G'$ can be viewed as an elliptic endoscopic group of some Levi subgroup $M$ of $G$, then one can show that $PW(T^{\mathcal{E}}(G))$ is equal to a subspace of
\[
\bigoplus_{G' \in \mathcal{E}_{ell}(G)} PW(\Phi_{bdd}(G'))^{{\rm Out}_{G}(G')},
\]
consisting of $b^{\mathcal{E}}_{G} := (b_{G'})_{G'}$ such that $b_{G'_1}, b_{G'_2}$ are equal over $\Phi_{2}(M')$ if $M'$ is a common Levi subgroup of $G'_1, G'_2 \in \mathcal{E}_{ell}(G)$. In \cite{Art96}, it is proved that there is a commutative diagram
\begin{align}
\label{diag: transfer}
\xymatrix{
I(G) \ar[r] \ar[d] & I^{\mathcal{E}}(G) \ar[d] \\
PW(T(G)) \ar[r] & PW(T^{\mathcal{E}}(G)).
}
\end{align}
Here the top and bottom arrows are geometric and spectral endoscopic transfers respectively, and it is also proved implicitly in \cite{Art96} that they are isomorphisms. The left vertical map is the Fourier transform and it is an isomorphism by the trace Paley-Wiener theorem. The right vertical map is given by the stable analogue and it is also an isomorphism, which  follows from this commutative diagram. In particular, $SI(G) \cong PW(\Phi_{bdd}(G))$. 

\begin{remark}
Since we have assumed the local Langlands correspondence, which includes the endoscopic character relations (cf. Desideratum~\ref{desideratum: LLC} \ECR), the commutativity of this diagram becomes a direct consequence of that. In practice, one often needs to establish this diagram first (cf. \cite{Art13}). In particular, it follows from this diagram and the definition of spectral endoscopic transfer map
\[
PW(T(G)) \rightarrow PW(T^{\mathcal{E}}(G))
\]
defined in \cite{Art96} that the endoscopic transfer sends $K_\BC\Pi(G')^{st}$ to $K_{\BC} \Pi(G)$.
\end{remark}

At last, we would like to show the existence of certain test functions used in the proof Theorem~\ref{thm: stable standards form a basis - quasisplit} as a consequence of the above discussion. For any $b_{G} \in PW(\Phi_{bdd}(G))$, we can find $f \in \mathcal{H}(G)$ such that the image of $f$ in $PW(\Phi_{bdd}(G))$ is $b_{G}$ and $\Theta^{G'}(\phi')(f') = 0$ for all $\phi' \in T^{\mathcal{E}}(G) \backslash \Phi_{bdd}(G)$ and endoscopic transfer $f^{G'}$ of $f$. By analytic continuation, we know that $\Theta^{G'}(\phi')(f^{G'})$ on the right hand side of \eqref{eq: stable} all vanish for such choice of $f$.

\part{ABV-packets and distributions for disconnected groups}\label{part: III ABV}

Here we extend the theory of ABV-packets and associated distributions from \cites{ABV92, Vog93, CFMMX22} to include the disconnected groups discussed in Part~\ref{part: I Langlands}.

\section{Geometric perspective on the local Langlands correspondence}\label{section: Vogan's LLC}

\subsection{Geometric perspective on the local Langlands correspondence}\label{subsection: Vogan's LLC}

In this section we explain how the local Langlands correspondence in Section~\ref{section: LLC} for the groups $G_\delta$, as $\delta$ ranges over pure inner forms for $G$, assemble together such that irreducible representations of these groups correspond to simple equivariant perverse sheaves on a common space. 
This is a generalization of \cite{Vog93} and \cite{CFMMX22}*{Chapter 4} from connected groups to disconnected groups.

\begin{definition}\label{definition: pure packet}
Suppose that $G$ satisfies Definition~\ref{assu: groups we consider} and that there exists a local Langlands correspondence for $G$ and its pure inner forms as in Definition~\ref{definition: LLC}.
Recall the definition of $H^1(F,G)_{\theta}$ from Section~\ref{subsection: pure inner forms} and set
\IN{$\Pi(G)^\pure$, set of representations of pure inner forms}
\[
    \Pi(G)^\pure 
    \ceq
    \coprod_{[\delta]\in H^1(F,G)_{\theta}}
    \Pi(G_\delta).
\]
Define \IN{$\LLC^{\mfr{e}}$, pure/Vogan Langlands correspondence}
\[ \LLC^{\mfr{e}}: \Pi(G)^{\pure} \to \Phi^{\mfr{e}}(G)  \]
by $\LLC^{\mfr{e}}(\pi):= \LLC^{\mfr{e}}_{G,\delta}(\pi)$ for $\pi \in \Pi(G_{\delta})$.
Let $\phi$ be a Langlands parameter. The {\it pure $L$-packet}\IN{$\Pi_\phi(G)^\pure$, pure $L$-packet} for $\phi$ is 
\[
    \Pi_\phi(G)^\pure 
    \ceq
    \coprod_{[\delta]\in H^1(F,G)_{\theta}}
    \Pi_{\phi}(G_\delta).
\]
For use in Theorem~\ref{theorem: geometric LLC} below, for each infinitesimal parameter $\lambda$ of $G$, we also set 
\[
\Pi_\lambda(G)^\pure 
\ceq
\coprod_{\phi \in \Phi_{\lambda}(G)} \Pi_{\phi}(G)^\pure, \qquad \Pi_\lambda(G_{\delta}) 
\ceq
\Pi_{\lambda}(G)^{\pure} \cap \Pi(G_{\delta}).
\]
\IN{$\Pi_\lambda(G_{\delta})$, set of representations with fixed infinitesimal parameter}
\end{definition}

\begin{lemma}
The local Langlands correspondence determines a bijection 
$
\Pi_\lambda(G)^\pure 
\leftrightarrow \Phi_\lambda^\mathfrak{e}(G)$, where 
\[
\Phi_\lambda^\mathfrak{e}(G) \ceq \{ (\phi,\tau) \in \Phi^\mathfrak{e}(G) \tq \lambda_\phi \equiv \lambda,\ \tau\in \widehat{A_\phi}\}.
\]
\end{lemma}

\begin{proof}
    For every $\delta \in Z^1(F,G)$ in the image of $Z^1(F,G^0)\to Z^1(F,G)$,
    the local Langlands correspondence determines a bijection 
    $
\Pi_\lambda(G_\delta) 
\leftrightarrow \Phi^\mathfrak{e}_\lambda(G, \delta)$, where 
\IN{$\Phi^\mathfrak{e}_\lambda(G,\delta)$, enhanced $L$-parameters of $G_\delta$ with infinite parameter $\lambda$}
\[
\Phi^\mathfrak{e}_\lambda(G,\delta) \ceq \{ (\phi,\tau) \in \Phi^\mathfrak{e}(G,\delta) \tq \lambda_\phi \equiv \lambda,\ \tau\in \widehat{A_\phi}[\delta]\}.
\]
The proof now follows from Proposition~\ref{proposition: enhancement by purity}: 
\[
\widehat{A_\phi} = \coprod_{[\delta]\in H^1(F,G)_{\theta}} \widehat{A_\phi}[\delta]. \qedhere
\] 
\end{proof}

\begin{definition}\label{definition: Vogan variety}
Let $G$ be a group as in Definition~\ref{assu: groups we consider}.
Let $\lambda : W_F \to \Lgroup{G}$ be an infinitesimal parameter $\lambda$ of $G$ as in Definition \ref{definition: parameters}(2). 
    \begin{enumerate}
        \item  The {\it Vogan variety}\IN{$V_\lambda$, Vogan variety} for $\lambda$ is 
\[
    V_\lambda \ceq \{x \in \dualgroup{\mathfrak{g}} 
    ~ \vert ~ 
    \Ad(\lambda(w))x=|w|x, ~\forall w\in W_F \},
\]
where $\dualgroup{\mathfrak{g}}:= \mathrm{Lie}(\widehat{G})$.
The Vogan variety $V_\lambda$ comes equipped with the conjugation/adjoint action of the {\it stabilizer group}\IN{$H_\lambda$, stabilizer group in $\dualgroup{G}$ of $\lambda$} 
$H_\lambda\ceq Z_{\widehat{G}}(\lambda)$.
\item Let $C$ be an $H_{\lambda}$-orbit of $V_{\lambda}$. 
For any base point $x \in C$, the {\it $H_{\lambda}$-equivariant fundamental group for $C$}\IN{$A_C$, equivariant fundamental group for $C\subseteq V_\lambda$} is 
\begin{equation}\label{eqn: equivariant fundamental group}
    A_{C,x}:=\pi_0(Z_{H_{\lambda}}(x))
    := Z_{H_{\lambda}}(x)/Z_{H_{\lambda}}(x)^0.
\end{equation}
Since the isomorphism class of $A_{C,x}$ is independent of the choice of $x$, we often write $A_{C}$ for $A_{C,x}$ if the choice of base point $x$ is implicit or unimportant.

\item Let $C_1, C_2$ be two $H_{\lambda}$-orbits of $V_{\lambda}$. The {\it closure ordering} on orbits is defined by
\[ 
C_1 \geq C_2 \Longleftrightarrow \overline{C_1} \supseteq C_2.
\]

\item We will use the function \IN{$x_\phi \in V_\lambda$}
$
\phi \mapsto x_\phi \in V_{\lambda_\phi}
$
from Langlands parameters with infinitesimal parameter $\lambda$ to $V_\lambda$ defined by 
\[
 \phi(1,e) = (\exp(x_\phi),1) \in \Lgroup{G},
\qquad
\text{for}
\quad
e := 
\begin{pmatrix}
    1 & 1 \\ 0 & 1 
\end{pmatrix}.
\] 
We let $C_{\phi}$ 
\IN{$C_\phi$, $H_{\lambda}$-orbit of $V_{\lambda}$ that contains $x_{\phi}$}
denote the $H_{\lambda}$-orbit of $V_{\lambda}$ that contains $x_{\phi}$ and set 
\IN{$d(\phi) = \dim(C_\phi)$, dimension of $C_\phi$}
$d(\phi) \ceq \dim C_{\phi}$.
    \end{enumerate}
\end{definition}

\begin{remark}\label{remark: Vogan varieties for general G}
Definition~\ref{definition: Vogan variety} can be made for reductive groups more generally than we have done above, but the definition requires the L-group $\Lgroup{G} = \dualgroup{G}\rtimes W_F$ which is not well documented without conditions on the component group of $G$. It is for this reason that we have studied $\dualgroup{G} = \dualgroup{G^0}\rtimes \langle {\hat\theta} \rangle$ in Part~\ref{part: I Langlands}.
\end{remark}

The following properties are proved in \cite{CFMMX22}*{Propositions 5.6, 4.2, Lemma 4.4} for connected $G$. 
Here we generalize them to the disconnected setting.

\begin{lemma}\label{lemma: bijection parameters and orbits disconnected}
    Fix an infinitesimal parameter $\lambda$ of $G$.
     \begin{enumerate}
   \item 
   The complex variety $V_\lambda$ is affine and stratified into finitely many orbits under the $H_\lambda$-action. 
   \item The assignment $\phi \mapsto C_\phi$ as in Definition \ref{definition: Vogan variety}(4) determines a bijection between the set 
        $\Phi_{\lambda}(G)$ and the set of $H_{\lambda}$-orbits in $V_{\lambda}$. 
    \item  For any $L$-parameter $\phi \in \Phi_{\lambda}(G)$, we have an isomorphism between the $H_{\lambda}$-equivariant component group for $C_{\phi}$ and the component group $A_{\phi}=\pi_0(Z_{\dualgroup{G}}(\phi))$ for the $L$-parameter $\phi$:
    $$A_{C_{\phi}}\simeq A_{\phi}.$$ 
\item $V_{\lambda^0} = V_{\lambda}$ and $H_{\lambda^0}$ is a finite-index subgroup of $H_{\lambda}$ and $A_{\phi^0}$ is a finite-index subgroup of $A_{\phi}$.
    \end{enumerate}
\end{lemma}

\begin{proof}
Part (1) follows from the same proof as \cite{CFMMX22}*{Proposition 5.6}, which we omit. For Parts (2) and (3), we compare with the connected case. Let $\lambda^{0}$ denote the corresponding infinitesimal parameter of the identity component $G^0$. Note that $V_{\lambda}= V_{\lambda^0}$ by definition. By Lemma \ref{lem centralizer}, we have a short exact sequence
\[ 1 \to H_{\lambda^0} \to H_{\lambda} \to \textrm{im}(H_{\lambda}  \to \pi_0(\dualgroup{G})) \to 1.\]

Suppose that ${}^{\hat{\theta}}\lambda^0 \not\sim \lambda^0$. Then $\textrm{im}(H_{\lambda}  \to \pi_0(\dualgroup{G}))$ is trivial and so
$H_{\lambda}=H_{\lambda^0}$. We obtain 
\[ \Phi_{\lambda}(G)= \left(\bigsqcup_{i=0}^{\ell-1} \Phi_{{}^{\hat{\theta}^i} \lambda^0} (G^0)\right)/ \langle \hat{\theta} \rangle. \]
Thus, the natural function $\Phi_{\lambda^0}(G) \to \Phi_{\lambda}(G)$ is a bijection. This proves Part (2) in this case. Part (3) follows immediately from Corollary \ref{cor: component group  l pardisc} since 
\[ A_{C_{\phi}}= A_{C_{\phi^0}} \simeq A_{\phi^0} \simeq A_{\phi}.\]

Suppose that ${}^{\hat{\theta}}\lambda^0 \sim \lambda^0$. Take $g \in \dualgroup{G^0}$ such that $ g {}^{\hat{\theta}} \lambda^0 g^{-1}=  \lambda^0$. Then $g \rtimes \hat{\theta} \in H_{\lambda}$ and we have
\[H_{\lambda} \backslash V_{\lambda}= \langle g \rtimes \hat{\theta}\rangle \backslash (H_{\lambda^0} \backslash V_{\lambda^0}) \simeq \Phi_{\lambda^0}(G^0)/ \langle  g \rtimes \hat{\theta}\rangle = \Phi_{\lambda}(G). \]
Part (3) follows from Corollary \ref{cor: component group  l pardisc}  and an analogous short exact sequence 
\[ 1 \to \pi_0(Z_{\dualgroup{G^0}}(\phi^0)) \to \pi_0(Z_{\dualgroup{G}}(\phi)) \to \textrm{im}(\pi_0(Z_{\dualgroup{G}}(\phi)) \to \pi_0(\dualgroup{G})) \to 1,\]
whose proof is similar to that in Corollary \ref{cor: component group  l pardisc}, which we omit.

Part (4) also follows from Corollary~\ref{cor: component group  l pardisc}.
\end{proof}

\begin{remark}
All parts of Lemma~\ref{lemma: bijection parameters and orbits disconnected} may be proved without reduction to the analogous statements for $G^0$. 
We use this reduction here because it gives an explicit description of the $A_\phi$ in the case at hand.
\end{remark}

The choice of base point $x_\phi$ for $C_\phi$ and Lemma \ref{lemma: bijection parameters and orbits disconnected} determine an equivalence of categories
\begin{equation}\label{equation: base point}
\Loc_{H_\lambda}(C_\phi) \simeq \Rep(A_{C_\phi})
\simeq
\Rep(A_\phi),
\end{equation}
where $\Loc_{H_\lambda}(C_\phi)$ denotes that monoidal abelian category of $H_\lambda$-equivariant (complex) local systems on $C_\phi$.
Thus, isomorphism classes of simple equivariant local systems on $C_\phi$ are parametrized by $\widehat{A_\phi}$ as defined in Definition \ref{definition: parameters}(4).
For $\tau \in \Rep(A_\phi)$, we write 
\[
\locsys{\tau}{C_\phi} \in \Loc_{H_{\lambda}}(C_\phi)
\]
\IN{$\locsys{\tau}{C_\phi}$, local system on $C_\phi$}
for the local system matching under \eqref{equation: base point}; when it causes no confusion, we will use the same notation when matching $\widehat{A_\phi}$ with  isomorphism classes of simple objects in $\Loc_{H_{\lambda}}(C_\phi)$.

Let $\Per_{H_\lambda}(V_\lambda)$ denote the category of $H_\lambda$-equivariant perverse sheaves on $V_\lambda$, which appears as the heart of a $t$-structure on the $H_\lambda$-equivariant bounded derived category $D^b_{c,H_\lambda}(V_\lambda)$ of constructible sheaves on $V_\lambda$. 
For more general exposition on perverse sheaves and derived categories for use in representation theory, see, for example, \cite{Ach21}*{\S 3.1, 3.3, 6.4}. Using \eqref{equation: base point}, we see that simple objects in $\Per_{H_\lambda}(V_\lambda)$ are all of the form $\IC(\locsys{\tau}{C_\phi})$ for simple $\tau \in \dualgroup{A_\phi}$ with $\lambda_\phi \equiv \lambda$, where $\IC(\locsys{\tau}{C_\phi})$ denotes 
the corresponding intersection cohomology complex:
\begin{equation}\label{eq: simple perverse sheaves classified}
\Per_{H_\lambda}(V_\lambda)_{\text{/iso}}^{\text{simple}}=\{\IC(\locsys{\tau}{C_\phi}) \tq [\phi]\in \Phi_{\lambda}(G),\  \tau \in \widehat{A_\phi} \}.
\end{equation}
In particular, this establishes a bijection of sets
$\Per_{H_{\lambda}}(V_{\lambda})^{\mathrm{simple}}_{/\mathrm{iso}} \leftrightarrow \Phi^\mathfrak{e}_\lambda(G)$,

\begin{theorem}\label{theorem: geometric LLC}
Suppose that $G$ satisfies Definition~\ref{assu: groups we consider}.
Fix an infinitesimal parameter $\lambda$ for $G$ as in Definition \ref{definition: parameters}(2).
The local Langlands correspondence (Definition~\ref{definition: LLC}) determines a bijection
\IN{$\IC(\pi)$, perverse sheaf corresponding to $\pi$ under $\LLC^\mathfrak{e}$}
\[
\begin{array}{rcl}
\IC : 
\Pi_{\lambda}(G)^\pure &\rightarrow& \Per_{H_{\lambda}}(V_{\lambda})_{/\text{iso}}^{\text{simple}}\\
\pi &\mapsto& \IC(\pi)
\end{array}
\]
by the condition
\[
\LLC^\mathfrak{e}(\pi) = (\phi,\tau)
\qquad\iff\qquad
\IC(\pi) = \IC(\locsys{\tau}{C_\phi}).
\]
If $\phi$ is tempered then $C_\phi$ is the open $H_{\lambda_\phi}$-orbit in $V_{\lambda_\phi}$, which is unique, and $\pi(\phi,\1)$ is generic with respect to the fixed Whittaker datum; here, $\1$ is the trivial character of $A_\phi$.
\end{theorem}

\begin{proof}
This is a direct consequence of combining Section~\ref{section: LLC} with the classification of equivalence classes of simple perverse sheaves in \eqref{eq: simple perverse sheaves classified}. 
\end{proof}

We recall Vogan's perspective on the local Langlands correspondence for the connected reductive algebraic group $G^0$, established in \cite{Vog93}. As explained in \cite{CFMMX22}*{Proposition 4.2}, there is a bijection $[\phi^0] \mapsto C_{\phi^0}$ between equivalence classes of Langlands parameters of $G^0$ that share the infinitesimal parameter $\lambda^0$, and the $H_{\lambda^0}$-orbits $C_{\phi}$ in $V_{\lambda^0}$. 
Let $\Pi_{\lambda^0}(G^0)^\pure$ denote the subset of representations in $\Pi(G^0)^\pure$ whose Langlands parameters share the infinitesimal parameter $\lambda^0$. 
Let $\pi(\phi^0,\tau)$ be the representation satisfying $\LLC^{\mathfrak{e}}_{G^0}(\pi(\phi^0,\tau))=(\phi^0, \tau)$, where $\tau \in \Irr{A_{\phi^0}}$, where $\LLC_{G^0}^{\mathfrak{e}}$ is described in \ref{desideratum: LLC}. 
Now $A_{C_{\phi^0}}=A_{\phi^0}$ by
\cite{CFMMX22}*{Lemma 4.4}.
Vogan's geometric perspective on the Local Langlands correspondence [for $G^0$] can be summarized as the bijection
\[
\begin{array}{rcl}
 \Pi_{\lambda^0}(G^0)^\pure 
    &\to& \Per_{H_{\lambda^0}}(V_{\lambda^0})_{\text{/iso}}^{\text{simple}} \\
    \pi(\phi, \tau) 
    &\mapsto&
    \IC(\locsys{\tau}{C_{\phi}}), \\
    \pi &\mapsto& \IC(\pi).
\end{array}
\]

\begin{remark}\label{Remark: Fixing where the tempered rep goes wrt whittaker normalization}
    Desideratum~\ref{desideratum: LLC}\Sha\ gives us the following. Let $[\phi^0]\in \Phi_{\lambda^0}(G^0)$ denote the $L$-parameter for $G^0$. If $\phi^0$ is tempered, then 
    \[
    \mathcal{IC}(\pi(\phi,\1))= \mathcal{IC}(\pi(\phi^0,\1)^+)=\IC(\locsys{\1}{C_\phi}),
    \]
    where $\pi(\phi^0,\1)^+$ is specified in Section~\ref{section: LLC} and $\1$ is the trivial character of $A_\phi$ and where $\locsys{\1}{C_\phi}$ is the constant sheaf of $C_\phi$.
\end{remark}

\subsection{Pairing of Grothendieck groups}\label{subsection: pairing}

In this section we write down a pairing between the Grothendieck groups $K_0\Pi_\lambda(G)^\pure$ and $K_0\Per_{H_\lambda}(V_\lambda)$ (see below) determined by the enhanced local Langlands correspondence in terms of a basis of simple objects. 
We explain the relationship between the pairing we use and the pairing introduced in \cite{Vog93}*{Section (8.11')(a)} and \cite{ABV92}*{Theorem 1.24} (for the connected case) in Remark \ref{remark: equivalence of pairings}.
Recall that we have defined Kottwitz sign $e(G_{\delta})$ for $G_{\delta}$ in Definition \ref{definition: kottwitz sign}, and for $\pi \in \Pi(G)^\pure$, $e(\pi)= e(G_{\delta})$ if $\pi \in \Pi(G_{\delta})$. 

\begin{definition}\label{definition: pairing}
Let $K_0\Pi_\lambda(G)^\pure$ be the $\Z$-span of irreducible representations in $\Pi_\lambda(G)^\pure$ as defined in Section~\ref{subsection: Vogan's LLC}; let $K_0\Per_{H_\lambda}(V_\lambda)$ be Grothendieck group of $\Per_{H_\lambda}(V_\lambda)$. Using Theorem \ref{theorem: geometric LLC}, we define the following pairing between groups:
\begin{equation}\label{equation: pairing}
    \begin{array}{rcl}
        K_0\Pi_\lambda(G)^\pure \times K_0\Per_{H_\lambda}(V_\lambda)
            &\to& \mathbb{Z}\\
        \langle \pi, \IC(\locsys{\tau}{C_\phi}) \rangle 
            &\mapsto& 
            \begin{cases}
                e(\pi)\ (-1)^{d(\pi)} & \text{ if }\LLC^\mathfrak{e}(\pi) = (\phi,\tau), \\
                0 & \text{else},
            \end{cases}
    \end{array}
\end{equation}
where $d(\pi) \ceq \dim (C_{\phi_\pi})$.
\IN{$d(\pi) = d(\phi_\pi)$, dimension of $C_{\phi_\pi}$}
\IN{$\langle \pi, \IC(\locsys{\tau}{C_\phi}) \rangle$, pairing between groups $K_0\Pi_\lambda(G)^\pure \times K_0\Per_{H_\lambda}(V_\lambda)$}
\end{definition}

We will use this pairing in Section~\ref{section: etaABV}.

\begin{remark}\label{remark: equivalence of pairings}
    The pairing in \eqref{equation: pairing} is different than the pairing defined in \cite{Vog93}:
    \[
    \begin{array}{rcl}
        K_0\Pi_\lambda(G)^\pure \times K_0\Per_{\dualgroup{G}}(X_\lambda)
            &\to& \mathbb{Z}\\
        \langle \pi, \IC(\locsys{\tau}{C_\phi}) \rangle 
            &\mapsto& 
            \begin{cases}
                e(\pi)(-1)^{d_X(\pi)} & \text{ if }\LLC^\mathfrak{e}(\pi) = (\phi,\tau), \\
                0 & \text{else},
            \end{cases}
    \end{array}
\]
where $d_X(\pi)$ denotes the dimension of the $\dualgroup{G}$-orbit of $\phi_\pi$ in $X_\lambda := (\dualgroup{G}\times V_\lambda)/H_\lambda$, and $\lambda$ is the infinitesimal parameter of $\pi$.
This pairing, and that of Definition~\ref{definition: pairing}, are equivalent by
induction equivalence of equivariant derived categories as in \cite{Vog93}*{Prop 8.8}, which is precisely \cite{ABV92}*{Prop 7.14}. This fits with the induction equivalence of equivariant perverse sheaves only up to a shift fitting with the differences in signs above.
\end{remark}

We end this section by introducing the notation for projecting to forms.

\begin{definition}\label{def project to form}
Let $\delta$ be a pure inner form of $G$. We write
\[\operatorname{proj}_{\delta}: K_{0}\Pi(G)^\pure \to K_{0}\Pi(G_{\delta}) \hookrightarrow K_{0}\Pi(G)^\pure\]
for projection composed with inclusion. We use the same notation 
\IN{$\operatorname{proj}_{\delta}$, projection}
$\operatorname{proj}_{\delta}$ for its restriction to $K_{0}\Pi_{\lambda}(G)^\pure$.
\end{definition}

We now describe the transpose of $\operatorname{proj}_{\delta}$  via the pairing in Definition \ref{definition: pairing}.

\begin{lemma}\label{lem: transpose of projection IC}
   Let ${}^t \operatorname{proj}_{\delta}: K_{0}\Per_{H_{\lambda}}(V_{\lambda}) \to K_{0}\Per_{H_{\lambda}}(V_{\lambda})$ denote the transpose of $\operatorname{proj}_{\delta}$ via the pairing in Definition \ref{definition: pairing}. Then ${}^t \operatorname{proj}_{\delta}$ is the projection onto the direct summand of $K_{0}\Per_{H_{\lambda}}(V_{\lambda}) $ generated by
   $ \{\IC(\locsys{\tau}{C_{\phi}} ) \tq \phi\in \Phi_{\lambda}(G), \ \tau \in \widehat{A_{\phi}}[\delta]\}$.
\end{lemma}
\begin{proof}
    This is immediate from the requirement that if $\pi \in \Pi(G_{\delta})$ then $\LLC(\pi) \in \Phi^{\mfr{e}}(G,\delta)$.
\end{proof}

\subsection{Kazhdan-Lusztig hypothesis}\label{subsection: KLH}

In this section, we explain the Kazhdan-Lusztig hypothesis for $p$-adic groups and use it to evaluate the pairing in Definition \ref{definition: pairing} on a basis of standard objects. 

Let $D^b_{c,H_\lambda}(V_\lambda)$ be the triangulated category of $H_\lambda$-equivariant complexes on $V_\lambda$, constructible with respect to the stratification of $V_\lambda$ into $H_\lambda$-orbits; see the appendices for more on this category.
The triangulated category $D^b_{H_\lambda}(V_\lambda)$ 
 has another heart consisting of constructible sheaves on $V_\lambda$; see, for example,
\cite{Ach21}*{Section 6.4, Exercise 6.4.3}.
For every local system $\locsys{\tau}{C}$ supported on $C$, we obtain a constructible sheaf $\stsh{\tau}{C}$ on $V_\lambda$ via extension by zero.
We refer to such $\stsh{\tau}{C}$ as \textit{standard sheaves}. The sheaves $\stsh{\tau}{C}$, as $\tau \in \widehat{A_C}$ and $C$ varies over $H_\lambda$-orbits in $V_\lambda$, form another basis of $K_0\Per_{H_\lambda}(V_\lambda)$. In the same vein, $K_0\Pi_{\lambda}(G)^\pure$ has another basis, consisting of standard representations $\SR(\pi)$.

\begin{hypothesis}[Kazhdan-Lusztig Hypothesis for $G$]\label{hypothesis: KLH}
    Suppose that $G$ satisfies Definition~\ref{assu: groups we consider} and that there exists a pure local Langlands correspondence for $G$ satisfying Desiderata \ref{desideratum: LLC}\textbf{(Hai)}. We further assume that for any infinitesimal parameter $\lambda$ for $G$, 
  the pairing defined in \eqref{equation: pairing} is equivalent to the pairing defined by
\begin{equation*}\label{equation: KLH}
    \begin{array}{rcl}
        K_0\Pi_\lambda(G)^\pure \times K_0\Per_{H_\lambda}(V_\lambda)
            &\to& \mathbb{Z}\\
        \langle \sr(\pi), \stsh{\tau}{C}\rangle 
            &\mapsto& 
            \begin{cases}
                e(\pi) & \text{ if }\LLC^{\mathfrak{e}}(\pi) = (\phi,\tau), \\
                0 & \text{else},
            \end{cases}
    \end{array}
\end{equation*}
where $\stsh{\tau}{C}$ denotes a standard sheaf in $D^b_{H_\lambda}(V_\lambda)$ on the $H_\lambda$-orbit $C_\phi$ and where $e(\pi)=e(G_{\delta}^0)$ is Kottwitz sign where $\pi \in \Pi(G_{\delta})$. 
\end{hypothesis}

\begin{remark}\label{remark: KLH status}
    For $G=\GL_N$,  Hypothesis~\ref{hypothesis: KLH} is a theorem by \cite{bal25}.
    For unipotent representations of quasisplit classical groups and their pure inner forms, this is a theorem by \cite{Sol25}*{Theorem 5.4}.
        For unipotent representations of $G(F) = \GL_N(F)\rtimes{\langle\theta\rangle}$, Hypothesis~\ref{hypothesis: KLH} is proved in \cite{CHAI}.
\end{remark}

We give a simple illustration of Hypothesis~\ref{hypothesis: KLH}.
From Section~\ref{subsection: KDF}, recall the definition of the linear functional $\chi^\textrm{loc}_C$ on the Grothendieck group $K_0 D^b_{H_\lambda}(V_\lambda)$.

\begin{lemma}\label{lemma: mu and chiloc}
Suppose that $G$ satisfies Definition~\ref{assu: groups we consider}. 
Assume Hypothesis~\ref{hypothesis: KLH} for $G$. 
Let $\phi$ be a Langlands parameter for $G$ and let $\lambda$ be its infinitesimal parameter.
Recall the definition $\STD{\phi}{\pure}$ from Definition~\ref{definition: STD}.
Then, for all $\mathcal{F} \in D^b_{c,H_\lambda}(V_\lambda)$,
\[
\langle \STD{\phi}{\pure}, \mathcal{F} \rangle
=
\chi^\textrm{loc}_{C_\phi}(\mathcal{F}),
\]
 where we have used the pairing of Definition~\ref{definition: pairing}.
\end{lemma} 

\begin{proof}
Without loss of generality, we may take  $\mathcal{F}$ from the basis for $K_0 D^b_{H_\lambda}(V_\lambda)$ adapted to standard representations using Hypothesis~\ref{hypothesis: KLH}, so suppose now that $\mathcal{F} =  j_!\locsys{\sigma}{C_{\phi_0}}$ for $\sigma \in \widehat{A_{\phi_0}}$. 
Then $\chi^{\textrm{loc}}_{C_\phi}(\mathcal{F}) = \rank(\sigma)$ if $[\phi_0] = [\phi]$ and $\chi^{\textrm{loc}}_{C_\phi}(\mathcal{F}) = 0$ otherwise. Similarly, if $[\phi_0] \ne [\phi]$ then
\[
\begin{array}{rcl}
\langle
\STD{\phi}{\pure}, \mathcal{F}
\rangle
&=&
\sum_{\tau\in \widehat{A_\phi}} e(\tau) \rank(\tau) 
	\langle \SR(\phi,\tau), \stsh{\sigma}{C_{\phi_0}} \rangle  = 0,
\end{array}
\]
while, if $[\phi_0] = [\phi]$ then
\[
\begin{array}{rcl}
\langle
\STD{\phi}{\pure}, \mathcal{F}
\rangle
&=&
\sum_{\tau\in \widehat{A_\phi}} e(\tau) \rank(\tau) 
	\langle \SR(\phi,\tau), \stsh{\sigma}{C_{\phi_0}} \rangle \\
&=& 	e(\sigma) e(\sigma) \rank(\sigma) = \rank(\sigma). 
\end{array}
\]
Thus,  $\langle
\STD{\phi}{\pure}, \mathcal{F}
\rangle 
= \chi^{\textrm{loc}}_{C_\phi}(\mathcal{F})$ in all cases. 
\end{proof}

Finally, we compute the transpose of $\operatorname{proj}_{\delta}$ defined in Definition \ref{def project to form} using the formulation of the pairing in Kazhdan-Lusztig hypothesis.

\begin{lemma}\label{lem: transpose of projection standard}
    Suppose that $G$ satisfies Definition~\ref{assu: groups we consider} and assume Hypothesis~\ref{hypothesis: KLH} for $G$. Let 
    \IN{${}^t \operatorname{proj}_{\delta}$, transpose of projection}
    ${}^t \operatorname{proj}_{\delta}: K_{0}\Per_{H_{\lambda}}(V_{\lambda}) \to K_{0}\Per_{H_{\lambda}}(V_{\lambda})$ denote the transpose of $\operatorname{proj}_{\delta}$ via the pairing in Definition \ref{definition: pairing}. Then ${}^t \operatorname{proj}_{\delta}$ is the projection onto the direct summand of $K_{0}\Per_{H_{\lambda}}(V_{\lambda}) $ generated by
   $ 
   \{
   \stsh{\tau}{C_{\phi}} 
   \tq 
   \phi\in \Phi_{\lambda}(G), \tau \in \widehat{A_{\phi}}[\delta] 
   \} $.
\end{lemma}
\begin{proof}
     For $\pi \in \Pi(G)^\pure$, we have $\pi \in \Pi(G_{\delta})$ if and only if $\SR(\pi) \in K_0\Pi(G_{\delta})$. Thus, we have 
\[ \operatorname{proj}_{\delta} (\SR(\pi))= \begin{cases}
     \SR(\pi) & \text{ if }\pi \in \Pi(G_{\delta}),\\
     0 & \text{ otherwise.}
\end{cases}\]
Then the result follows from the formulation of the pairing in Hypothesis \ref{hypothesis: KLH}.
\end{proof}

\begin{remark}
Although it is not needed in this paper, it follows from \cite{CS2} that the abelian category $\Per_{H_\lambda}(V_\lambda)$ decomposes into subcategories indexed by classes of pure inner forms of $G$:
 \[
    \Per_{H_\lambda}(V_\lambda)
    =
    \mathop{\bigoplus}_{[\delta]\in H^1(F,G)_{\theta}}\Per_{H_\lambda}^{[\delta]}(V_\lambda)
    \]
    where the direct sum is taken over 
    $
    H^1(F,G^0)_{\langle\theta\rangle} 
    \cong
    \mathrm{im}(H^1(F,G^0)\to H^1(F,G)),
    $
    using Lemma~\ref{lem: coinvariants}.
    The following diagram summarizes the correspondence:
\[
\begin{tikzcd}
    \Pi_\lambda(G)^\pure 
    \arrow{rr}{\pi \mapsto \IC(\pi)} 
    \arrow[equal]{d}
    &&
    \Per_{H_{\lambda}}(V_{\lambda})_{/\text{iso}}^{\text{simple}}
    \arrow[equal]{d}\\
    \mathop{\coprod}\limits_{[\delta]\in H^1(F,G)_{\theta}} \Pi_\lambda(G_\delta)^\pure
    &&
    \mathop{\coprod}\limits_{[\delta]\in H^1(F,G)_{\theta}}
    \Per^{[\delta]}_{H_{\lambda}}(V_{\lambda})_{/\text{iso}}^{\text{simple}} \\
    \Pi_\lambda(G_\delta) 
    \arrow{rr}{\pi\mapsto \IC(\pi)}[swap]{\textrm{Section~\ref{section: LLC}}} 
    &&
    \Per^{[\delta]}_{H_{\lambda}}(V_{\lambda})_{/\text{iso}}^{\text{simple}}
\end{tikzcd}
\]
\end{remark}

\section{ABV-packets}
\label{section: ABV-packets}

In this section we define ABV-packets for the algebraic groups $G_\delta$ {with $G$} satisfying Definition~\ref{assu: groups we consider}. We note that this was done for connected reductive groups $G_\delta^0$ in \cite{CFMMX22}, building on \cite{Vog93}.

\subsection{Conormal variety and conormal bundles}\label{subsection: conormal}

In this section, we collect some basic facts about the conormal variety $\Lambda_\lambda$ to the Vogan variety $V_\lambda$, concerning: bundles, components, orbit duality, and regular covectors.

Let $\lambda$ be any infinitesimal parameter for $G$.
The {\it conormal variety to $V_\lambda$}
\IN{$\Lambda_\lambda$, conormal variety to $V_\lambda$} is 
\begin{equation}\label{definiition: conormal}
\Lambda_\lambda 
\ceq
\{ (x,y)\in V_\lambda\times V_{\lambda}^* \tq [x,y]=0 \},
\end{equation}
where $[x,y]$ is the Lie bracket in $\dualgroup{\mathfrak{g}}$.
We also define the {\it conormal bundle}\IN{$\Lambda_{C_\phi}$, conormal bundle above $C_\phi$} above $C_\phi$ by
\[
\Lambda_{C_\phi} 
\ceq
\{ (x,y)\in \Lambda_{\lambda} \tq x\in C_\phi \}
=
\{ (x,y)\in C_\phi \times V_{\lambda}^* \tq [x,y] =0 \}.
\] 
Observe that
\[
\Lambda_\lambda 
= 
\coprod_{[\phi]\in \Phi_\lambda(G)}
\Lambda_{C_\phi}
=
\bigcup_{[\phi]\in \Phi_\lambda(G)}
\overline{\Lambda_{C_\phi}},
\] where, $\overline{\Lambda_{C_\phi}}$ is the Zariski closure.
This is explained in Subsection \ref{subsection: normalization} of the Appendix.

We remark that the subvarieties $\overline{\Lambda_{C_\phi}}\subseteq \Lambda_\lambda$ need not be irreducible components of $\Lambda_\lambda$; indeed, the $H_\lambda$-orbits $C_\phi \subseteq V_\lambda$ need not be connected.
However, in the category of varieties with $H_\lambda$-group actions, every irreducible component of $\Lambda_\lambda$ is indeed of the form  $\overline{\Lambda_{C_\phi}}$, as $\phi$ ranges over $\Phi_\lambda(G)$; see \cite{ABV92}*{Lemma 19.2(b)}. 
From this observation, the usual definition of Pyasetskii duality of $H_\lambda$-orbits in $V_\lambda$ extends to the present case, since, by the same argument, for every $H_\lambda$-orbit $D\subseteq V_\lambda^*$, the closure of bundle $\Lambda^*_{D} \ceq (D\times V_\lambda) \cap \Lambda^*_\lambda$ is likewise an irreducible  component of $\Lambda^*_\lambda \ceq \{ (y,x)\in V_\lambda^*\times V_\lambda \tq [y,x]=0 \}$ as a variety with $H_\lambda$-action. Using the isomorphism $\Lambda^*_\lambda \to \Lambda_\lambda$ given by $V^*_\lambda \times V_\lambda \to V_\lambda\times V^*_\lambda$, we can pair $C_\phi$ uniquely with some $H_\lambda$-orbit $C_\phi^*$ in $V^*_\lambda$. This defines Pyasteskii duality for orbits that need not be connected.
This is explained in greater generality in Subsection \ref{subsection: PDGV} from the Appendix.

It will be useful to observe that, just as an orbit $C_\phi$ is an induction space $C_\phi \iso H_{\lambda}\times_{Z_{H_\lambda}(x_\phi)}\{x_\phi\}$, each conormal bundle is an induction space by the isomorphism
\begin{equation}\label{equation: induction space}
H_\lambda \times_{Z_{H_\lambda}(x_\phi)} \Lambda_{x_\phi} \mathop{\longrightarrow}\limits^{\iso} \Lambda_{C_\phi}  
\end{equation}
given by $[h,y]_{Z_{H_\lambda}(x)} \mapsto h\cdot(x,y)$, where,
\IN{$\Lambda_{x_\phi}$, covectors at $x_\phi \in V_\lambda$}
\[
\Lambda_{x_\phi} \ceq \{ y\in V_\lambda^* \tq [x_\phi,y] =0 \}.
\]
Consequently, we have a categorical isomorphism
\[
D^b_{H_\lambda}(\Lambda_C) \iso D^b_{Z_{H_\lambda}(x)}(\Lambda_x).
\]

We now define some subvarieties of $\Lambda_{C_\phi}$ that will play an important role in this paper.
To begin, for every $H_\lambda$-orbit $C\subseteq V_\lambda$, consider the intersection 
$
\Lambda_{C_\phi} \cap \overline{\Lambda_{C}}$.
This intersection is empty unless $C\geq C_\phi$, in which case 
$
\Lambda_{C_\phi} \cap \overline{\Lambda_{C}}
=
(C_\phi \times \overline{C^*} ) \cap \Lambda_\lambda.
$ 
The {\it singular} parts of $\Lambda_{C_\phi}$ are defined as follows:
\IN{$\Lambda_{C_\phi}^\mathrm{sing}$, singular part of conormal bundle} 
\begin{equation}\label{definition: reg}
\Lambda_{C_\phi}^\mathrm{sing}
\ceq
\bigcup_{C \ne C_\phi} \left( \Lambda_{C_\phi} \cap \overline{\Lambda_{C}} \right)
\qquad\textrm{and}\qquad
\Lambda_{C_\phi}^\reg
\ceq
\Lambda_{C_\phi} \setminus \Lambda_{C_\phi}^\mathrm{sing}.
\end{equation}
Note that $\Lambda_{C_\phi}^\reg$ is open in $\Lambda_{C_\phi}$ and carries the $H_\lambda$-action inherited from $\Lambda_{C_\phi}$.
Using \eqref{equation: induction space}, one can see that \IN{$\Lambda^\reg_{x_\phi}$, regular covectors at $x_\phi \in C_\phi$}
\[
\Lambda^\reg_{x_\phi} \ceq \{ y\in V^*_\lambda \tq (x_\phi,y) \in \Lambda^\reg_{C_\phi} \}
\]
satisfies 
\[
\Lambda^\reg_{C_\phi} \iso H_\lambda\times_{Z_{H_\lambda}(x_\phi)}\Lambda_{x_\phi}^\reg.
\]

\subsection{Equivariant local systems on conormal bundles}\label{subsection: ELS}

Here we study the equivariant fundamental group of $\Lambda_{C_\phi}^\reg$.
Since the action of $H_{\lambda}$ on $\Lambda_{C_\phi}^\reg$ may not transitive, the equivariant fundamental group of $\Lambda_{C_\phi}^\reg$ may not be the component group of an $H_\lambda$-orbit. 
For any $y \in \Lambda_{x_\phi}^\reg$, set\IN{$A^\reg_{\phi,y}$, equivariant fundamental group of conormal bundle at $(x_\phi,y)$}
\[
A^{\reg}_{\phi,y}
:= \pi_1^{H_{\lambda}}(\Lambda_{C_\phi}^{\reg},(x_\phi,y)).
\]
This is the \'etale fundamental group of the Galois category of $H_\lambda$-equivariant finite \'etale covers of $\Lambda^\reg_{C_\phi}$ with respect to the fibre functor determined by $(x_\phi,y)\in \Lambda_{C_\phi}^\reg$.
The choice of base point $y\in \Lambda_{x_\phi}^\reg$ determines the equivalence of categories
\begin{equation}\label{equation:gen}
\Loc_{H_{\lambda}}(\Lambda_{C_\phi}^{\reg}) 
\cong
\Rep(A_{\phi,y}^{\reg}).
\end{equation}
By contrast, for every $z\in \Lambda_\lambda$, the equivariant fundamental group of its $H_\lambda$-orbit is
\[
\pi_1^{H_\lambda}(H_\lambda\cdot z,z) 
= \pi_0( Z_{H_\lambda}(z) )
= Z_{H_\lambda}(z)/Z_{H_\lambda}(z)^0.
\]
For any $y\in {\Lambda^\reg_{x_\phi}}$, set \IN{$A^\mic_{\phi,y}$, microlocal fundamental group at $y$}
\[
A^\mic_{\phi,y} \ceq \pi_1^{H_\lambda}(H_\lambda\cdot (x_\phi,y),(x_\phi,y)) 
= \pi_0(Z_{H_\lambda}(x_\phi,y)).
\]

The groups $A_{\phi,y}^{\reg}$ and $A_{\phi,y}^{\mathrm{mic}}$ are related by restriction of local systems, specifically, by restriction from $\Lambda_{C_\phi}^{\reg}$ to the $H_{\lambda_\phi}$-orbit $H_{\lambda_\phi}\cdot(x_{\phi},y)$.
Moreover, the choice of $y\in \Lambda_{x_\phi}^{\reg}$
determines the following commuting diagram.
\[
\begin{tikzcd}
\arrow{d}[swap]{\text{restriction}}
\Loc_{H_{\lambda_\phi}}(\Lambda_{C_\phi}^{\reg}) 
\arrow{r}{\cong} & \Rep(A_{\phi,y}^{\reg})
\arrow{d}{\text{restriction}} \\
\Loc_{H_{\lambda_\phi}}({H_{\lambda_\phi}}\cdot(x_\phi,y)) 
\arrow{r}{\cong} &  \Rep(A_{\phi,y}^{\mathrm{mic}})
\end{tikzcd}
\]
    We say $y\in \Lambda_x^\reg$ is {\it microlocally generic} if it appears in a $Z_{H_\lambda}(x)$-invariant subvariety appearing in \cite{ABV92}*{Lemma 24.3}.
    We write $\Lambda^\mic_x$\IN{$\Lambda^\mic_x$, microlocally generic covectors at $x$} for set of microlocally generic $y\in \Lambda_x^\reg$ and set $\Lambda^\mic_C \ceq H_\lambda\times_{Z_{H\lambda}(x)} \Lambda_x^\mic$ for $x\in C$.\IN{$\Lambda^\mic_C$, microlocally generic part of $\Lambda_C$}
If the choice of $y\in \Lambda^\reg_{x_\phi}$ is implicit or irrelevant, we may suppress it from the notation.
In fact, by \cite{ABV92}*{Lemma 24.3}, $A^\mic_{\phi,y}$ is independent of  $y\in \Lambda^\mic_{x_\phi}$.

\begin{remark}\label{remark: generic KS}
In \cite{CFK22} one finds an example of Langlands parameter $\phi=\phi_\KS$ for $\GL_{16}$ such that the group action on $\Lambda_{C_{\phi}}^\reg$ is not transitive and, moreover, $A_{\phi,y}^{\mathrm{mic}}$ is different from $A_{\phi,y}^{\reg}$; indeed, the former is trivial while the latter is non-trivial for $y\in \Lambda^\mic_{x_\phi}$.
\end{remark}

\subsection{Normalized vanishing cycles}\label{subsection: NEvs}

Suppose that $G$ satisfies Definition~\ref{assu: groups we consider}.
    Let $\lambda$ be an infinitesimal parameter for $G$, and let $V_{\lambda}$ and $H_{\lambda}$ be as in Definition \ref{definition: Vogan variety}(1).
    Let $C_\phi$ be any $H_\lambda$-orbit in $V_\lambda$.
    In Appendix~\ref{section: global vanishing cycles} we make some improvements on the functor $\Evs_C$, originally introduced in \cite{CFMMX22} and summarized here:
    \IN{$\Evs_{C}$, vanishing cycles functor}
    \begin{align*}
    \Evs_{C}: 
    D^b_{c,H_\lambda}(V_\lambda) &\to D^b_{c,H_\lambda}(\Lambda_C^\reg)\\
    \mathcal{F} &\mapsto j^*_{\Lambda_{C^*}^\reg}\, \Rphi{f_{C^*}}[-1]\ p_{C^*}^* \mathcal{F} [-\codim C^*],
    \end{align*}
    where:
    \begin{enumerate}
        \item 
        $C^* \subseteq V_\lambda^*$ is the Pyasetskii dual of the orbit $C$, $\codim C^*=\dim  V_\lambda^*-\dim C^*$, and $p_{C^*} : V_\lambda \times C^* \to V_\lambda$ is the restriction of the projection $p : V_\lambda \times V_{\lambda}^* \to V_\lambda$;
        \item 
        $f_{C^*}:V_\lambda \times C^* \to \mathbb{A}_\CC^1$ is the restriction of the Killing form on $V_\lambda \times V_\lambda^* \subset \widehat{\mathfrak{g}}\times \widehat{\mathfrak{g}}^{\ast}$;
        \item  
        $\Lambda_{C^*}^\reg\subset V_\lambda \times V_\lambda^*$ is the regular part of the conormal bundle above $C^*$, recalled in Section~\ref{subsection: ELS}, and $j_{\Lambda_{C^*}^\reg} : \Lambda_{C^*}^\reg \hookrightarrow f_{C^*}^{-1}(0)$ is inclusion;
        \item the functor $\Rphi{f_{C^*}} :D^b_{c,H_{\lambda}}(V_\lambda\times C^*) \to D^b_{c,H_\lambda}(f_{C^*}^{-1}(0))$ is the equivariant vanishing cycles functor, as reviewed in Appendix~\ref{section: equivariant vanishing cycles}, see also  \cite{SGA7II}*{Exposés XIII et XIV}.
    \end{enumerate}
This is explained in greater detail in Theorem \ref{theorem: refinement} from the Appendix. The corresponding stratification $\mathscr{S}$ is constructed in Lemma \ref{lemma: Whitney}. 
Especially, the appendices to this paper provide more detail about the categories appearing as the target and source in the definition of $\Evs_C$, above.
We note that the cohomology sheaves of objects in $D^b_{c,H_\lambda}(V_\lambda)$ are constructible with respect to the stratification of $V_{\lambda}$ into $H_{\lambda}$-orbits, of which there are finitely many.
By contrast, the cohomology sheaves of objects in $D^b_{c,H_\lambda}(\Lambda_\lambda)$ are constructible with respect to a very different stratification, described in Appendix~\ref{section: global vanishing cycles}.

In Appendices~\ref{section: equivariant vanishing cycles} and \ref{section: global vanishing cycles} we also show that the restriction of $\Evs_C$ to equivariant perverse sheaves on $V_\lambda$ produces equivariant local systems on $\Lambda_{C}^\reg$.

\begin{definition}\label{definition: nevs}
Let $\lambda$ be an infinitesimal parameter for $G$ and let $V_\lambda$ be the Vogan variety for $\lambda$, equipped with the action of the centralizer group $H_\lambda$.
For an $H_\lambda$-orbit $C$ in $V_\lambda$, the functor $\NEvs_C$\IN{$\NEvs_C$, normalized $\Evs_C$}
is defined by
\begin{align*}
    \NEvs_C: D^b_{c,H_\lambda}(V_\lambda) &\to D^b_{c,H_\lambda}(\Lambda_C^\reg)\\
    \mathcal{F} &\mapsto \mathcal{T}_C^\vee \otimes \Evs_C(\mathcal{F}),
\end{align*}
where 
$
\mathcal{T}_C := \Evs_C\, \IC(\1_C)
$ is the {\it normalizing sheaf for $C$}\IN{$\mathcal{T_C}$, normalizing sheaf for $C$} and where $\mathcal{T}^\vee_C$ is the dual local system to $\mathcal{T}_C$, taken in $\Loc_{H_\lambda}(\Lambda_C^\reg)$.
\end{definition}

In Theorem~\ref{theorem: refinement} we prove that $\NEvs_{C} : D^b_{c,H_\lambda}(V_\lambda) \to D^b_{c,H_\lambda}(\Lambda_C^\reg)$ restricts to
\[
\NEvs_{C} : \Per_{H_\lambda}(V_\lambda) \to \Loc_{H_\lambda}(\Lambda_{C}^\reg).
\]
That result is used in Definition~\ref{definition: nevs} in the definition of the normalizing sheaf. 
It follows from \cite{CFMMX22}*{Section 7.10} that the rank of that sheaf is $1$.
This equivariant rank one local system doesn't depend  the on the choice of the nondegenerate and invariant symmetric bilinear form $\langle \cdot , \cdot\rangle$
by Remark \ref{rem-rank-one}.
Indeed, the normalization of the functor $\Evs_C$ appearing in Definition~\ref{definition: nevs} ensures that the analogue of condition \cite{CFMMX22}*{Theorem 7.22(d)} is satisfied: 
\[
\NEvs_C\IC(\locsys{\tau}{C}) 
= 
p_{\Lambda_{C}^\reg}^*
\
\stsh{\tau}{C} ,
\qquad \forall \tau\in \Rep(A_C).
\]
where $p_{\Lambda_{C^*}^\reg}$ is the restriction of the projection $p: V_\lambda\times V_\lambda^* \to V_\lambda$ to $\Lambda_{C}^\reg$.
In particular, this ensures the equality
\[
\NEvs_C \IC(\1_C) =  \1_{\Lambda_C^\reg}.
\]
This also follows from Remark \ref{local-sys}.

\begin{remark}\label{remark: Evs refinement}
    The definitions of $\Evs_C$ and $\NEvs_C$, above, have a subtle but important difference from the definitions that appeared in \cite{CFMMX22}*{Chapter 7}.
    Specifically, in \cite{CFMMX22} we introduced a subvariety $\Lambda_{C}^\gen \subseteq \Lambda_{C}^\reg$ defined tautologically so that the restriction of $\Evs_{C}$ to perverse sheaves gives
    \[
\NEvs_{C} : \Per_{H_\lambda}(V_\lambda) \to \Loc_{H_\lambda}(\Lambda_{C}^\gen).
\]
 At that time we did not know that $\Lambda_{C}^\gen = \Lambda_{C}^\reg$, but that is exactly the content of Theorem~\ref{theorem: refinement}. 
 This is an improvement since the definition of $\Lambda_{C}^\gen$, given in \cite{CFMMX22}*{Section 7.9}, was unwieldy.
\end{remark}

\subsection{ABV-packets}\label{subsection: ABV-packets}

Using the functor $\NEvs_C$, we may now define ABV-packets.

\begin{definition}\label{definition: ABV-packet}
    Let $\phi$ be a Langlands parameter for $G$ with infinitesimal parameter $\lambda$. Let $C_\phi$ be the $H_\lambda$-orbit in $V_\lambda$ corresponding to the parameter $\phi$.
    With reference to Section~\ref{subsection: Vogan's LLC}, for every $\pi \in \Pi_{\lambda}(G)$ we define $\IC(\pi)\in \Per_{H_\lambda}(V_\lambda)$ by
    \[
    \IC(\pi) = \IC(\locsys{\mathcal L}{C_\phi}) 
    \qquad\iff\qquad
    \LLC(\pi)=(\phi, \mathcal{L}).
    \]
For every pure inner form $\delta \in Z^1(F,G^0)$, the {\it ABV-packet}\IN{$\Pi^\ABV_{\phi}(G_\delta)$, ABV-packet} is     
\[
        \Pi_\phi^\ABV(G_\delta)
            =
            \{\pi \in \Pi_\lambda(G_\delta)
            \tq
            \NEvs_{C_\phi}\,
            \IC(\pi)\neq 0\},
    \]
while the {\it pure ABV-packet for $\phi$} is
    \IN{$\Pi^\ABV_{\phi}(G)^{\pure}$, pure ABV-packet}
    \[
    \Pi^\ABV_{\phi}(G)^{\pure} 
    \ceq
    \{\pi \in \Pi_\lambda(G)^\pure~\vert~ \NEvs_{C_\phi}\,
    \IC(\pi)\neq 0\}.\]
\end{definition}

\begin{remark}\label{remark: gen=reg}
In \cite{CFMMX22}*{Definition 1}, the ABV-packet is defined using the $\mathrm{Ev}_C$-functor. These definitions are equivalent since   $\NEvs_C(\mathcal{P})\neq 0$ if and only if $\mathrm{Ev}_C(\mathcal{P})\neq 0$; in this paper we use only $\NEvs_C$ for notational hygiene. 
Moreover, the Appendices to this paper make some important improvements regarding $\Ev_{C}$ and $\Evs_{C}$ as they appeared in \cite{CFMMX22}, the most notable is that now we know that the set characterized in \cite{CFMMX22} as the ``generic part of $\Lambda_C^\reg$" in fact coincides with the full regular part $\Lambda_{C}^\reg$, as explained in Remark~\ref{remark: Evs refinement}. 
This is equivalent to the assertion, above, that $\Evs_{C}$, as it is defined on this section, takes equivariant perverse sheaves to equivariant local systems.
\end{remark}

\begin{theorem}\label{theorem: ABV-LLC compatability}
  Let $G$ be a quasisplit group of the sort considered in Definition~\ref{assu: groups we consider} and let $\delta$ be a pure inner form for $G$.
  Let $\phi$ be a Langlands parameter for $G$ with infinitesimal parameter $\lambda$.
Then
\[
\Pi_{\phi}(G_\delta) \subseteq 
\Pi^\ABV_{\phi}(G_\delta) 
=
\Pi^\mic_{\phi}(G_\delta) 
\subseteq
\Pi_{\lambda}(G_\delta)
\]
and
\[
\pi \in \Pi^\ABV_{\phi}(G_\delta)
\qquad\implies\qquad
\phi_{\pi} \geq \phi
\]
with equality if and only if $\pi \in \Pi_\phi(G_\delta)$.
Moreover, the following diagram commutes
\[
\begin{tikzcd}
 \Pi^\ABV_{\phi}(G_\delta) 
 \arrow{rrr}{\NEvs_{C_\phi}}  
 &&& \Rep(A_{\phi,y}^\reg) \arrow{d} & A_{\phi,y}^\reg \\
 \Pi^\mic_\phi(G_\delta) \arrow[equal]{u}
 \arrow{rrr}{\NEvs_{(x_\phi,y)}}
  &&& \Rep(A_{\phi,y}^\mic) & A_{\phi,y}^\mic \arrow{d} \arrow{u} \\
 \Pi_\phi(G_{\delta}) \arrow{u} \arrow{rrr}{\LLC^{\mathfrak{e}}} &&& \Rep(A_\phi) \arrow{u} & A_\phi 
\end{tikzcd}
\]
where:
    \begin{enumerate}
    \item 
    $y\in \Lambda^\reg_{x_\phi}$;
    \item
    $\Rep(A_\phi) \to \Rep(A_{\phi,y}^\mic)$ is obtained by the pullback along $A_{\phi,y}^\mic \to A_\phi$ which is induced from the inclusion $Z_{H_\lambda}(x_\phi,y) \to Z_{H_\lambda}(x_\phi)$ ;
    and 
    \item $\Rep(A_{\phi,y}^{\reg}) \to \Rep(A_{\phi,y}^\mic)$ is obtained by the pullback along $A_{\phi,y}^\mic \to A_{\phi,y}^\reg$, which corresponds to the functor defined by restriction $\Loc_{H_\lambda}(\Lambda^\reg_{C_\phi}) \to \Loc_{H_\lambda}(H_\lambda\cdot (x_\phi,y))$.
\end{enumerate}
\end{theorem}

\begin{proof}
    For connected $G$, the first two statements were proved in \cite{CFMMX22}*{Theorem 7.22}.
    The proof there applies here because the results all follow from properties of vanishing cycles of simple objects in $\Per_{H_\lambda}(V_\lambda)$.
    The crucial observation here is that while $V_{\lambda} = V_{\lambda^0}$, even when $G$ is connected, $H_\lambda$ may not be connected so as far as vanishing cycles are concerned, nothing new happens for disconnected $G$ that did not already happen for $G^0$. 
    
    The maps between representations are all induced by pullback along the group homomorphisms on the right. 
    In Section~\ref{subsection: ELS} we saw the arrow $A_{\phi,y}^\mic \to A_{\phi,y}^\reg$.
    Consequently, the composition 
\[
\Pi^\ABV_\phi(G_\delta) \to \Rep(A_{\phi,y}^\reg) \to \Rep(A_{\phi,y}^\mic)
\]
is determined by 
\[
\NEvs_{(x_\phi,y)} : \Per_{H_\lambda}(V_\lambda) \to \Loc_{H_\lambda}(H_\lambda\cdot (x_\phi,y)) 
\iso
\Rep(A_{\phi,y}^\mic)
\]
for $(x_\phi,y)\in \Lambda_{C_\phi}^\reg$,
where $y$ is also chosen when determining the equivalence $\Loc_{H_\lambda}(\Lambda_{C_\phi}^\reg) \iso \Rep(A_{\phi,y}^\reg)$ in the statement of the theorem.
Now the commutativity of the diagram follows from 
\cite{CFMMX22}*{Theorem 7.22 (d)}.

See \cite{CDFZ25} for the role of Whittaker normalization in the horizontal arrows in the statement of the proposition.
\end{proof}

\section{ABV-distributions}\label{section: etaABV}

In this section we attach distributions to the ABV-packets that were defined in Section~\ref{section: ABV-packets} and we use these to find a basis for stable distributions on quasisplit connected groups (under the assumptions for Theorem \ref{thm: stable standards form a basis - quasisplit}).
The main result of this section is Theorem~\ref{theorem: Conjecture 2 qs}.

\subsection{Stable ABV-distributions}\label{subsection: stable ABV-distributions}

It is not $\NEvs_{C_\phi}$ but its restriction 
\IN{$\NEvs_{(x_\phi,y)}$, normalized vanishing cycles at $(x_\phi,y)\in \Lambda^\reg_{C_\phi}$}
\IN{$\NEvs_{z}$, normalized vanishing cycles at $z\in \Lambda_\lambda^\reg$}
\begin{equation}
\NEvs_{(x_\phi,y)} \ceq j^*_{{H_\lambda}\cdot(x_\phi,y)} \NEvs_{C_\phi} : D^b_{c,H_\lambda}(V_\lambda) \to D^b_{c,H_\lambda}({H_\lambda}\cdot(x_\phi,y))
\end{equation}
that appears in the main results here, because $D^b_{c,H_\lambda}({H_\lambda}\cdot(x_\phi,y))$ gives us access to representations of $A^\mic_\phi$, through the choice of the base point $y\in \Lambda^\reg_{x_\phi}$.
For any Langlands parameter, the restriction of $\NEvs_{(x_\phi,y)}$ to equivariant perverse sheaves gives
\[
\NEvs_{(x_\phi,y)} : \Per_{H_\lambda}(V_\lambda) \to \Rep(A^\mic_{\phi,y}).
\]
obtained by composing $\NEvs_{C_{\phi}} : \Per_{H_\lambda}(V_\lambda) \to \Loc_{H_\lambda}(\Lambda^\reg_{C_\phi})$ with the restriction 
$\Loc_{H_\lambda}(\Lambda_{C_\phi}) \to \Loc_{H_\lambda}({H_\lambda}\cdot(x_\phi, y))$
and the equivalence $\Loc_{H_\lambda}({H_\lambda}\cdot(x_\phi, y)) \simeq \Rep(A^\mic_{\phi})$ determined by $y \in \Lambda_{x_\phi}^\reg$.

\begin{lemma}\label{lemma: rank and mic}
    Let $G$ be a quasisplit group as considered in Section~\ref{subsection: assumptions on groups}. 
    Let $\phi$ be a Langlands parameter for $G$.
    Then, for any $\mathcal{F}\in D^b_{c,H_\lambda}(V_\lambda)$, the Euler characteristic of $\NEvs_{(x_\phi,y)} \mathcal{F}$ is independent  $y\in \Lambda_{x_\phi}^\reg$ and
\[
 \chi\left(\NEvs_{(x_\phi,y)} \mathcal{F} \right)
=
\chi^\mic_{C_\phi}(\mathcal{F}).
\]
The corresponding Euler characteristic $\chi$ is discussed in Equation
\eqref{equation: stalks 8}.
In particular, if $\mathcal{F}$ is a perverse sheaf then $\NEvs_{(x_\phi,y)} \mathcal{F}$ is concentrated in  degree zero where it is a local system and in this case, it's microlocal characteristic $\chi^\mic_{C_\phi}\mathcal{F}$ records the rank of the local system $\NEvs_{(x_\phi,y)} \mathcal{F}$.
For facts about $\chi^\mic_{C_\phi}$, see Appendix~\ref{subsection: KDF}.
\end{lemma}

\begin{proof}
Lemma~\ref{lemma: rank and mic} is proved in Appendix~\ref{subsection: normalization}, \eqref{equation: NMD}.
\end{proof}

\begin{definition}\label{definition: stable etaABV}
    Let $G$ be a quasisplit group of the sort considered in Section~\ref{subsection: assumptions on groups}.
    Let $\phi$ be a Langlands parameter for $G$.
    Set
    \IN{$\etaABV{\phi}{}, \etaABV{\phi}{,\pure}, \etaABV{\phi}{,\delta},$ stable ABV-distributions}
    \[
\etaABV{\phi}{}
\ceq
\sum_{\pi \in \Pi_{\phi}(G)} 
(-1)^{d(\phi) - d(\pi)}\ \chi^\mic_{C_\phi}(\IC(\pi)) \ [\pi].
\]
Also define
   \[
\etaABV{\phi}{,\pure}
\ceq
\sum_{\pi \in \Pi_{\phi}(G)^\pure} e(\pi) 
(-1)^{d(\pi)}\ 
\chi^\mic_{C_\phi}(\IC(\pi))  \ [\pi].
\]
Finally, for any pure inner form $\delta$ for $G$, set
  \[
\etaABV{\phi}{,\delta}
\ceq
\sum_{\pi \in \Pi_{\phi}(G_\delta)^\pure} (-1)^{d(\phi) - d(\pi)} 
\chi^\mic_{C_\phi}(\IC(\pi)) \ [\pi] .
\]
We refer to these as {\it stable ABV-distributions}.
\end{definition}

\begin{lemma}\label{lemma: NEvs and chimic}
Let $G$ and $\phi$ be as above.
    \[
\langle
\etaABV{\phi}{,\pure}, \mathcal{F} \rangle
=
\chi^\mic_{C_\phi}(\mathcal{F}),
\qquad\forall \mathcal{F}\in D^b_{c,H_\lambda}(V_\lambda).
\]
\end{lemma}

\begin{proof}
To prove the lemma we may take $\mathcal{F}$ from any basis for $K_0 D^b_{c,H_\lambda}(V_\lambda)$.
We use the basis of simple $H_\lambda$-equivariant perverse sheaves on $V_\lambda$. So, without loss of generality, we now take $\mathcal{F} = \IC(\sigma)$ for $\sigma \in \Pi_\lambda(G)^\pure$.
Then
\[
\begin{array}{rcl}
\langle
\etaABV{\phi}{,\pure}, \mathcal{F}
\rangle
&=& 
\langle
\etaABV{\phi}{,\pure}, \IC(\sigma)
\rangle \\
&=&
\mathop{\sum}\limits_{\pi\in \Pi_\phi^{\ABV}(G)^\pure} (-1)^{d(\pi)} e(\pi)\ \chi^\mic_{C_\phi}(\IC(\pi)) \ 
\langle \pi, \IC(\sigma)\rangle \\
&=& 
(-1)^{d(\sigma)}\, e(\sigma) \, (-1)^{d(\sigma)} e(\sigma)\ 
\chi^\mic_{C_\phi}(\IC(\sigma))  \\
&=&
\chi^\mic_{C_\phi}(\IC(\sigma)) 
\\
&=&
\chi^\mic_{C_\phi}(\mathcal{F}).
\end{array}
\]
\end{proof}

\begin{theorem}\label{theorem: Conjecture 2 qs}
Let $G$ be a {\bf connected} reductive algebraic group, quasisplit over $F$. 
Let $\lambda$ be an infinitesimal parameter for $G$.
Assume Hypothesis~\ref{hypothesis: KLH} for $G$ and $\lambda$.
Then the ABV-distributions $\etaABV{\phi}{}$ defined in Section~\ref{subsection: ABV-distributions} satisfies that
\[
\operatorname{span}
\{\mu_{\phi} \tq [\phi] \in \Phi_\lambda(G) \}
= 
\operatorname{span}
\{\eta_{\phi}^{\ABV} \tq \phi \in \Phi_\lambda(G) \}.
\]
Suppose further that the local Langlands correspondence for $G$ satisfies the assumptions of Theorem \ref{thm: stable standards form a basis - quasisplit}. Then $\{\etaABV{\phi}{} \tq \phi \in \Phi_\lambda(G) \}$ is a basis of $K_{\BC}\Pi_{\lambda}(G)^{\st}$.

\end{theorem}

\begin{proof}
We use the Kashiwara-Dubson Index Formula, which is recalled in \eqref{equation: stalks 7}.
We begin by explaining how that result applies to the matter at hand. 
From Section~\ref{subsection: KDF}, recall the definition of the linear functionals $\chi^\textrm{loc}_C$ and $\chi^\textrm{mic}_C$ on the Grothendieck group $K_0 D^b_{c,H_\lambda}(V_\lambda)$.
Now recall Lemmas~\ref{lemma: NEvs and chimic} and \ref{lemma: mu and chiloc}:
for all $\mathcal{F} \in K_0 D^b_{c,H}(V)$,
\begin{equation}\label{equation: mu loc, miceta}
\chi^\textrm{loc}_{C_\phi}(\mathcal{F})
= \langle \STD{\phi}{\pure}, \mathcal{F} \rangle
\qquad\textrm{and}\qquad
\chi^\textrm{mic}_{C_\phi}(\mathcal{F})
= \langle \etaABV{\phi}{,\pure}, \mathcal{F} \rangle,
\end{equation}
 where we have used the pairing of Definition~\ref{definition: pairing}. 
 Note that Lemma~\ref{lemma: mu and chiloc} assumes Hypothesis~\ref{hypothesis: KLH}.
Now, from Section~\ref{subsection: KDF} we have 
$
\chi^\textrm{mic}_{C_\phi}
= 
\sum_{C_{\phi_i}} c(\phi_i,\phi) \chi^\textrm{loc}_{C_{\phi_i}} 
$
as linear functionals on $K_0 D^b_{c,H_\lambda}(V_\lambda)$, 
where $c(\phi,\phi_i) =0$ unless $\phi_i \geq \phi$ and $c(\phi,\phi) =(-1)^{d(\phi)}$.
(Recall that $\phi \leq \phi_i$ means $C_\phi \subseteq \overline{C_{\phi_i}}$.)
Accordingly, these coefficients $c(\phi_i,\phi)$ define a change of basis matrix.
Thus,
\[
\chi^\textrm{mic}_{C_\phi}(\mathcal{F})
= 
\sum_{i} c(\phi_i,\phi) \chi^\textrm{loc}_{C_{\phi_i}}(\mathcal{F}),
\qquad
\forall \mathcal{F} \in K_0 D^b_{c,H_\lambda}(V_\lambda).
\]
By the work above, this means
\[
{\langle \etaABV{\phi}{,\pure}, \mathcal{F} \rangle}
=
\sum_{i} c(\phi_i,\phi) 
{\langle \STD{\phi_i}{\pure}, \mathcal{F} \rangle}
=
{\langle \sum_{i} c(\phi_i,\phi) \STD{\phi_i}{\pure}, \mathcal{F} \rangle}
\]
Since the pairing in Definition~\ref{definition: pairing} is non-degenerate, it follows now that
    \begin{equation}
        \etaABV{\phi}{,\pure}
        =
    \sum_{[\phi_i]\in \Phi_{\lambda}(G)} c(\phi_i,\phi) \STD{\phi_i}{\pure}.
    \end{equation}
Note that the equation about could be written in the form
\begin{equation}
        \etaABV{\phi}{,\pure}
        =
        (-1)^{d(\phi)} \STD{\phi}{\pure} + \sum_{\substack{\phi_i\in \Phi_{\lambda}(G) \\ \phi_i>\phi}} c(\phi_i,\phi) \STD{\phi_i}{\pure},
    \end{equation}
In particular, it follows that, using the projection 
$
K_\CC \Pi_\lambda(G)^{\pure} 
\to
K_\CC \Pi_\lambda(G),
$
that this gives
\begin{equation}
        \etaABV{\phi}{}
        =
        (-1)^{d(\phi)}\STD{\phi}{} + \sum_{\substack{\phi_i\in \Phi_{\lambda}(G) \\ \phi_i>\phi}} c(\phi_i,\phi) \STD{\phi_i}{}.
    \end{equation}
In Theorem~\ref{thm: stable standards form a basis - quasisplit}
we saw that stable standards form a basis:
\[
K_\CC \Pi_\lambda(G)^\st 
=
\operatorname{span}_\CC
\lbrace
\STD{\phi_j}{}
\tq
[\phi_j]\in \Phi_\lambda(G) \rbrace,
\qquad\text{and}\qquad
\dim K_\CC \Pi_\lambda(G)^\st  = \vert \Phi_\lambda(G) \vert.
\]
This implies the result, since the change of basis matrix is invertible.
\end{proof}

    In the above proof, we only impose the $p$-adic Kazhdan-Lusztig hypothesis for the local Langlands correspondence. If we further assume the local Langlands correspondence satisfies \textbf{(LC)} and \textbf{(St)} of Desiderata \ref{desideratum: LLC}, then $\mu_{\phi}$ is stable. 
    Hence, $\etaABV{\phi}{}$ is stable, that is why we refer to these as {\it stable ABV-distributions}.
    If we further includes the assumptions in Theorem \ref{thm: stable standards form a basis - quasisplit}, then $\{ \etaABV{\phi}{} 
    \tq 
    \phi \in \Phi_\lambda(G) \}$ forms a basis of $K_{\BC}\Pi_{\lambda}(G)^{\st}$. 

\subsection{Trace of equivariant vanishing cycles}\label{subsection: trace}

First, we recall some basic facts. 
Let $Z$ be any $H$-variety, where $H$ is a linear algebraic group, both taken here over the complex numbers.
For any $\mathcal{E}\in D^b_{c,H}(Z)$ and semisimple $h\in H$ and $z\in Z^h$, $h$ defines an automorphism of $\mathcal{H}^\bullet_{z}\mathcal{E}$; the trace of this automorphism is commonly denoted by 
\[
\trace(h;\mathcal{H}^\bullet_{z}\mathcal{E})
\ceq
\sum_{i\in \mathbb{Z}}(-1)^i \trace(h;\mathcal{H}^i_{z}\mathcal{E}).
\]
Suppose that $\mathcal{E}$ is an equivariant local system on $Z$. Again suppose $z\in Z^s$. Then $\mathcal{H}^\bullet_{z}\mathcal{E} = \mathcal{E}_{z}$ and $z\in Z^h$ determines a character of $\pi_1^H(Z,z)$ and also a representation of $\pi_1^H(H\cdot z,z)$ where, as above, $H\cdot z$ denotes the $H$-orbit of $z\in Z$.
This latter fundamental group is $Z_H(z)/Z_H(z)^0$; note that $h\in Z_{H}(z)$.
In this case, 
\[
\trace(h,\mathcal{H}^\bullet_{z}\mathcal{E}) = \trace(h, E),
\qquad
E\in \Rep(\pi_1^H(H\cdot z,z)),
\]
where $E$ is the representation of $\pi_1^H(H\cdot z,z)$ determined by $\mathcal{E}$ at the base point $z\in Z^h$.

We will need the specific case of this construction for $Z = \Lambda_{C_\phi}^{\reg}$ and $\mathcal{E}=\NEvs_{C_\phi}\IC(\pi)$, where $\IC(\pi)$ is the simple perverse sheaf corresponding to some $\pi\in \Pi_{\lambda}(G)$ under the local Langlands correspondence, as explained in Section~\ref{section: LLC}.
Suppose further that $h$ is the image of semisimple element $s\in \dualgroup{G}$, and $z=(x_\phi,y)\in (\Lambda_{C_\phi}^{\reg})^s$.
Now $\NEvs_{C_\phi}\, \IC(\pi)$ is concentrated in degree $0$, which is to say, $\mathcal{H}^i_{(x_\phi,y)}\NEvs_{C_\phi}\, \IC(\pi)$ is trivial except for $i=0$.
And in the appendix we show that $\NEvs_{C_\phi}\, \IC(\pi)$ is an $H_\lambda$-equivariant local system on $\Lambda_{C_\phi}^{\reg}$.
In Section~\ref{subsection: ELS} we saw that the category of $H_\lambda$-equivariant local systems on $\Lambda_{C_\phi}^{\textrm{reg}}$ is canonically equivalent to $\Rep(A^\reg_{\phi})$ where, as above, $A^\reg_{\phi}$ is simply the $H_\lambda$-equivariant fundamental group of $\Lambda^\reg_{C_\phi}$.
The vector space $\mathcal{H}^0_{(x_\phi,y)}\NEvs_{C_\phi}\, \IC(\pi)$ comes equipped with an action of $Z_{H_\lambda}(x_\phi,y)$.  
In this way we see that, for $z=(x_\phi,y)\in (\Lambda^\reg_{C_\phi})^s$, we may view $\NEvs_z\IC(\pi)$ as a representation of $A^\mic_{\phi,y}$. 
More generally, for any $\mathcal{F}\in D^b_{c,H_\lambda}(V_\lambda)$ we now see and that
\[
\trace(s,\NEvs_{(x_\phi,y)}\mathcal{F})
\]
is the trace at $h$ of $\NEvs_{(x_\phi,y)}\mathcal{F} \in D \Rep(A^\mic_{\phi,y})$, with $h\in A^\mic_{\phi,y}=Z_{H_\lambda}(x_\phi,y)/Z_{H_\lambda}(x_\phi,y)^0$ the image of $s$.

\begin{proposition}
\label{proposition: independence of y}
Let $G$ be a quasisplit group of the sort considered in Section~\ref{subsection: assumptions on groups} and let $\phi$ be a Langlands parameter for $G$.
Let $s\in \dualgroup{G}$ be semisimple such that $\,^s\phi = \phi$ and also $(\Lambda^\reg_{x_\phi})^s\neq \emptyset$.
Then, for every $\mathcal{F}\in D^b_{c,H_\lambda}(V_\lambda)$,
\[
\trace\left({s}; \NEvs_{(x_\phi,y_0)} \mathcal{F} \right) 
= \trace\left({h}; \NEvs_{(x_\phi,y_0)} \mathcal{F} \right) 
\]
is independent of the choice of base point $y_0\in (\Lambda_{x_\phi}^\reg)^s$.
\end{proposition}

\begin{proof}
First, note that the hypothesis $\,^s\phi=\phi$ implies $\,^s\lambda = \lambda$ and therefore $s\in H_\lambda$. 
We have also assumed $(\Lambda_{x_\phi}^\reg)^s \ne \emptyset$. But $(\Lambda_{x_\phi}^\reg)^s$ is also Zariski open in the complex vector space $(\Lambda_{x_\phi})^s$, so that $(\Lambda_{x_\phi}^\reg)^s$ is moreover \emph{connected}.
By additivity of the traces, we can also assume 
$\mathcal{F}=\IC(\pi)$, where $\IC(\pi)$ is the simple perverse sheaf corresponding to some $\pi\in \Pi_{\lambda}(G)$ under the local Langlands correspondence.
Pick $y_0\in (\Lambda_{x_\phi}^\reg)^s$; thus, $y_0$ is fixed by $s$.
By Section~\ref{theorem: refinement}, $\NEvs_{C_\phi}\IC(\pi)$ is an $H_\lambda$-equivariant local system on $\Lambda^\reg_{C_\phi}$.
Let $S$ be the diagonalizable group generated by $s$, which is to say, $S$ is the Zariski-closure of $\langle s \rangle$ in $H_\lambda$. Then $\trace\left({s}; \NEvs_{(x_\phi,y_0)} \IC(\pi) \right)$ only depends on 
the restriction of $\NEvs_{C_\phi}\IC(\pi)$ viewed as an $S$-equivariant local system on $(\Lambda_{x_\phi}^\reg)^s$. But $S$ is acting trivially on $(\Lambda_{x_\phi}^\reg)^s$, so that this can be viewed as a local system of $S/S^0$-modules
(e.g., see \cite{Sch03}*{Example 3.1.2 on pp. 153-154}), with
\[
\trace\left({s}; \NEvs_{(x_\phi,y_0)} \IC(\pi) \right) 
= \trace\left({h}; \NEvs_{(x_\phi,y_0)} \IC(\pi) \right) 
\]
viewed as the trace of this $S/S^0$-module at the image $h$ of $s$ in $S/S^0$. But this trace is then locally constant and therefore also constant on $(\Lambda_{x_\phi}^\reg)^s$, which is connected.
\end{proof}

\begin{definition}\label{definition: trace}
Let $G$ be a quasisplit group of the sort considered in Section~\ref{subsection: assumptions on groups} and let $\phi$ be a Langlands parameter for $G$.
Assume a local Langlands correspondence (Definition~\ref{definition: LLC}) for $G$.
Let $s\in \dualgroup{G}$ be semisimple such that $\,^s\phi = \phi$ and also $(\Lambda^\reg_{x_\phi})^s\neq \emptyset$.
Using Lemma~\ref{proposition: independence of y}, we define
\IN{$\langle \cdot,\pi\rangle_\phi^\ABV$, character of representation of microlocal fundamental group}
\[
 \langle {s},\pi\rangle_\phi^\ABV
\ceq
\trace\left({s}; \NEvs_{(x_\phi,y_0)} \IC(\pi) \right) 
\]
for any $y_0\in (\Lambda_{x_\phi}^\reg)^s$.
\IN{$\langle {s},\pi\rangle_\phi^\ABV= \trace\left({s}; \NEvs_{(x_\phi,y_0)} \IC(\pi) \right)$, trace paring}
\end{definition}

\subsection{ABV-distributions}\label{subsection: ABV-distributions}

In this section we define distributions $\etaABV{\phi,s}{\,\delta}$. If $\phi$ is of Arthur type, these distributions satisfy the endoscopic character identities in \cite{Art13}, which we will verify later (see Theorem \ref{theorem: transfer of ABV-distributions}).
Throughout this section we assume the local Langlands correspondence for $G$ in Definition~\ref{definition: LLC}, but not the Desiderata~\ref{desideratum: LLC}.

\begin{definition}\label{definition: etaABV}
Let $\delta$ be a pure inner form of a quasisplit group $G$ from Section~\ref{subsection: assumptions on groups}.
Let $\phi$ be a Langlands parameter for $G$ fixed by semisimple $s\in \dualgroup{G}$ and further such that $(\Lambda^\reg_{x_\phi})^s\neq \emptyset$.
Then
     $\etaABV{\phi,s}{,\delta}$ is the virtual representation of $G_{\delta}(F)$ defined by
    \[        
    \etaABV{\phi,s}{,\delta}
    \coloneqq
    (-1)^{d(\phi)}
        \sum_{\pi \in \Pi^\ABV_{\phi}(G_\delta)} (-1)^{d(\pi)} 
        \langle {s}, \pi\rangle_{\phi}^{\ABV} \ 
        [\pi].
\]
Also define
\[  
    \etaABV{\phi,s}{,\pure}
    \coloneqq
       (-1)^{d(\phi)} \operatorname{\sum}\limits_{\pi \in \Pi^\ABV(G)^\pure} e(\pi) (-1)^{d(\pi)} 
     \langle {s}, \pi\rangle_{\phi}^{\ABV} \ 
       [\pi] 
    \]
When $[\delta]$ is the class for the quasisplit form $G$ we elide $\delta$ from the notation above, thus defining $\etaABV{\phi,s}{}$.
Likewise, if $s$ is trivial we will elide it from the notation, thus defining $\etaABV{\phi}{,\delta}$ and $\etaABV{\phi}{}$.
Thus, in particular,
\[        
    \etaABV{\phi}{}
    \coloneqq
        (-1)^{d(\phi)}\sum_{\pi \in \Pi^\ABV_{\phi}(G)} (-1)^{d(\pi)} 
        \langle1, \pi\rangle_{\phi}^{\ABV} \ 
        [\pi].
\]
We refer to all these as {\it ABV-distributions}\IN{$\etaABV{\phi,s}{,\delta}$, ABV-distributions}.
\end{definition}

\begin{remark}
Definitions~\ref{definition: etaABV} and \ref{definition: stable etaABV} both define $\etaABV{\phi}{,\delta}$. 
These definitions are compatible because
    \[
    \langle 1,\pi\rangle_\phi^\ABV
    =
    (-1)^{d(\phi)}\, \chi_{C_\phi}^\mic(\IC(\pi)),
    \]
as follows directly from the definitions and from Lemma~\ref{lemma: rank and mic}.
\end{remark}

\begin{remark}\label{remark: trilinguality}
The functor $Q^\mic_C$ is indeed related to $\NEvs_C$, as anticipated in 
\cite{CFMMX22}*{Section 7.11}.
Appendix~\ref{section: Qmic} shows, unconditionally, that 
\[
\Pi^\ABV_\phi(G)
=
\Pi^\mic_\phi(G)
\qquad\textrm{and}\qquad
\etaABV{\phi}{}
=
\eta^\mic_{\phi}
\]
for all Langlands parameters $\phi$ and for all $s\in Z_{\dualgroup{G}}(\phi,y)$ for some microlocally generic $y\in \Lambda_{x_\phi}^\reg$.
Future work will even show the equalities $\etaABV{\phi,s}{,\delta} = \eta^{\mic,\delta}_{\phi,s}$.
This is discussed in more detail in Remark 
\ref{rem-comparison} from the Appendix.
\end{remark}

\begin{proposition}\label{prop L-packet inside ABV packet}
Let $G$ be a quasisplit group as in Section~\ref{subsection: assumptions on groups}.
    If  ${}^{s}\phi=\phi$ and $\pi \in \Pi_{\phi}(G)^\pure$ then 
    \[
    \langle {s},\pi\rangle_{\phi}^{\ABV} 
    = 
    \langle {s},\pi\rangle_{\phi}
    \]
    where the pairing on the right-hand side is defined in \eqref{equation: pairing phi}.
\end{proposition}

\begin{proof}
If $\pi \in \Pi_\phi(G)^\pure$ then $\pi=\pi(\phi,\tau)$ for a unique $\tau \in \widehat{A_\phi}$, by Desiderata~\ref{desideratum: LLC}.
The definition of $\NEvs_{C_\phi}$ is designed so that $\NEvs_{C_\phi} \IC(\pi(\phi,\tau))) = (p^*\tau)_{\Lambda^\reg_{C_\phi}}$.
Thus,
\[
\begin{array}{rcl}
\langle {s},\pi\rangle_{\phi}^{\ABV}
&=&
\trace\left({s}; \NEvs_{(x_\phi,y_0)}\IC(\pi(\phi,\tau))\right) \\
&=&
\trace\left({s};(p^*\tau)_{\Lambda^\reg_{C_\phi}} \right) \\
&=&
\trace\left({s};\tau\right) \\
&=& 
\langle {s},\pi\rangle_{\phi}. 
\end{array}
\]
where $p^*$ is the pullback under the corresponding projection map.
\end{proof}

Our next result generalizes Lemma~\ref{lemma: NEvs and chimic}; the method of proof is similar.

\begin{proposition}\label{proposition: etaABV,F}
With $\phi$ and $s$ as above, and $y_0\in (\Lambda^\reg_{x_\phi})^s$, 
\[
\langle
\etaABV{\phi,s}{,\pure}, \mathcal{F}
\rangle
=
(-1)^{d(\phi)}
\trace({s},\NEvs_{(x_\phi,y_0)} \mathcal{F}),
\qquad \forall \mathcal{F}\in  D^b_{c,H_\lambda}(V_\lambda).
\]
If we fix a pure inner form $\delta$ of $G$, then
\[
\langle
\etaABV{\phi,s}{,\delta}, \mathcal{F}
\rangle
=e(G_{\delta})
(-1)^{d(\phi)}
\trace({s},\NEvs_{(x_\phi,y_0)} {}^t \operatorname{proj}_{\delta}\mathcal{F}),
\qquad \forall \mathcal{F}\in  D^b_{c,H_\lambda}(V_\lambda),
\]
where $\operatorname{proj}_{\delta}$ is described in Lemma \ref{lem: transpose of projection IC} or Lemma \ref{lem: transpose of projection standard}.
\end{proposition}

\begin{proof}
Let us again emphasize that $\trace({s},\NEvs_{(x_\phi,y_0)} \mathcal{F})$ is, by the definitions above, the trace of the representation of the complex $\NEvs_{C_\phi} \mathcal{F}$
 restricted to $H_\lambda\cdot(x_\phi,y_0)$.
 To prove the proposition, we may take $\mathcal{F}$ from any basis for $K_0 D^b_{c,H_\lambda}(V_\lambda)$. 
Using the pairing of Grothendieck groups from Section~\ref{subsection: KLH}, we take $\mathcal{F}$ from the basis adapted to ABV-distributions, namely equivariant perverse sheaves, so we take $\mathcal{F} = \IC(\sigma)$.
Then
\[
\begin{array}{rcl}
\langle
\etaABV{\phi,s}{,\pure}, \mathcal{F}
\rangle
&=& 
(-1)^{d(\phi)} \mathop{\sum}\limits_{\pi\in \Pi_\phi^\ABV(G)^\pure} (-1)^{d(\pi)} e(\pi) \ \langle {s}, \pi\rangle_{\phi}^\ABV\, 
\langle \pi, \IC(\sigma)\rangle \\
&=& 
(-1)^{d(\phi)-d(\sigma)}\, e(\sigma) \, (-1)^{d(\sigma)} e(\sigma)
\langle {s}, \sigma\rangle_{\phi}^\ABV \\
&=&
(-1)^{d(\phi)}
\trace({s},\NEvs_{(x_\phi,y)} \mathcal{F}).
\end{array}
\]
Finally, the form by form version follows from $\operatorname{proj}_{\delta}( \etaABV{\phi,s}{,\pure})= e(G_{\delta})\etaABV{\phi,s}{,\delta}$.
\end{proof}

\subsection{Quasisplit transfer of stable ABV-distributions}\label{ssec: proof of prop QST}

Having introduced the distributions $\etaABV{\phi,s}{,\delta}$, here we study the transfer of $\etaABV{\phi}{}$ from $G$ to the forms $G_\delta$.
As above, in this section we assume the local Langlands correspondence for $G$ in Definition~\ref{definition: LLC}, but not the Desiderata~\ref{desideratum: LLC}.

We want to describe the geometric analogue of the following application of a special case of endoscopic transfer, namely, the endoscopic triple $(G^0,1,\xi)$ for $G$, where $\xi : \Lgroup{G^0} \to \Lgroup{G}$ is the obvious inclusion. See Section~\ref{section: passage to the identity} for related material.
\begin{definition}\label{definition: QST and z-strongly stable}
Let $G$ be a quasisplit group satisfying Definition~\ref{assu: groups we consider}.
Let $\QST_{G}$ be the linear transformation defined by
\begin{equation}\label{eqn:QST defined}
\begin{array}{rcl}
     \QST_{G}: K_{\BC}\Pi(G^0)^{\st} &\to& K_{\BC}\Pi(G)^{\pure}, \\
     \eta &\mapsto& \sum\limits_{[\delta] \in H^1(F, G)_{\theta}} \Transfer{G^0}{G_{\delta}}(\eta).
\end{array}
\end{equation}
We refer to $\QST$ as {\it quasisplit transfer}\IN{$\QST$, quasisplit transfer}.
If $G=G^0$, following \cite{ABV92}*{\S 18}, we say that a distribution in $K_{\BC}\Pi(G)^{\pure}$ is {\it strongly stable}\IN{strongly stable distribution} if it lies in the image of $\QST_{G}$. 
\end{definition}

Inspired by Theorems \ref{theorem: transfer of STD} and \ref{theorem: Conjecture 2 qs}, , we give the following geometric definition. Note that we do not assume $G$ is connected in this definition.

\begin{definition}\label{definition: QSTgeo}
Suppose that $G$ satisfies Assumption \ref{assu: groups we consider} and there exists a local Langlands correspondence for $G$ as in Definition \ref{definition: LLC}. 
Let $\lambda$ be an infinitesimal parameter of $G.$    \begin{enumerate}
        \item We let $K_{\BC}\Pi_{\lambda}(G)^{\std}$ denote the subspace of $K_{\BC}\Pi_{\lambda}(G)$ spanned by $\{\mu_{\phi}\tq \phi \in \Phi_{\lambda}(G)\}$.
        \item We define $\QST_G^{\geo}: K_{\BC}\Pi_{\lambda}(G)^{\std} \to K_{\BC}\Pi_{\lambda}(G)^{\pure}$ by
        \[ 
        \QST_G^{\geo}(\mu_{\phi}):= \mu_{\phi}^{\pure},
        \]
    and likewise define $\QST_G^{\geo,\delta}: K_{\BC}\Pi_{\lambda}(G)^{\std} \to K_{\BC}\Pi_{\lambda}(G_\delta)$ by
        \[ 
        \QST_G^{\geo,\delta}(\mu_{\phi}):= e(G_{\delta}) \ \mu_{\phi}^{\delta}.
        \]
    \end{enumerate}
We refer to $\QST^\geo$ as \IN{$\QST^\geo$, geometric quasisplit transfer}. 
\end{definition}

\begin{lemma}\label{lemma: QST} Assume that $G$ is connected.
Suppose that the local Langlands correspondence for $G$ satisfies the assumptions of Theorem \ref{thm: stable standards form a basis - quasisplit}. 
Then $\mu_{\phi}$ is stable 
and 
\[ 
K_{\BC}\Pi_{\lambda}(G)^{\std} 
=
K_{\BC}\Pi_{\lambda}(G)^{\st}
\qquad\mathrm{and}\qquad
\QST^{\geo}_G = \QST_G.
\]
\end{lemma}

\begin{proof}
The distribution $\STD{\phi}{}$ is stable since it is a parabolic induction of a stable essentially tempered distribution. 
The equality above follows from Theorem \ref{thm: stable standards form a basis - quasisplit}. For the last claim,
Theorem \ref{theorem: transfer of STD} implies that the two linear transformations agree on the basis $\{\mu_{\phi}\tq\phi \in \Phi_{\lambda}(G)\}$. Hence, the two linear transformation agree.
\end{proof}

Now we evaluate $\QST^{\geo}$ at the stable ABV-distributions.

\begin{proposition}\label{prop: QST of ABV}
    Let $G$ be a connected reductive algebraic group, quasisplit over $F$. 
For $\delta\in Z^1(F,G)$,
\[
\QST_{G}^{\geo,\delta}
    \left( \etaABV{\phi}{} \right) 
    = 
    e(G_\delta)\,
    \etaABV{\phi}{,\delta}.
\]
If the local Langlands correspondence satisfies the assumptions of Theorem \ref{thm: stable standards form a basis - quasisplit}, then
    \[
    \Transfer{G}{G_\delta}
    \left( \etaABV{\phi}{} \right) 
    = 
    e(G_\delta)\,
    \etaABV{\phi}{,\delta}.
    \]
\end{proposition}

\begin{proof}
We give the argument for assuming the desiderata listed above. The proof without them follows the same structure after replacing $\Transfer{G}{G_\delta}$ with $\QST_{G}^{\geo,\delta}$, with reference to Lemma~\ref{lemma: QST}.
In the proof of Theorem~\ref{theorem: Conjecture 2 qs} we saw that
\begin{equation}\label{equation: KD}
        \etaABV{\phi}{,\pure}
        =
        \sum_{[\phi_i]\in \Phi_{\lambda}(G)}
        c(\phi_i,\phi) \STD{\phi_i}{\pure}.
    \end{equation}
    By Theorem \ref{theorem: transfer of STD} (or the definition of $\QST_{G}^{\geo, \delta}$), we have
\[
    \left( \Transfer{G}{G_\delta} \circ \operatorname{proj}_{0} \right) (\STD{\phi_i}{\pure}) 
    = 
\operatorname{proj}_{\delta} (\STD{\phi_i}{\pure}),
    \]
for all $[\phi_i]\in \Phi_\lambda(G)$,
where 
\[
\operatorname{proj}_\delta  :
K_\CC \Pi_\lambda(G)^{\pure}
\to
K_\CC \Pi_\lambda(G_\delta)
\]
is projection to virtual representations of $G_\delta(F)$, and $ \operatorname{proj}_0$ is the projection to the quasisplit form.
Applying this to \eqref{equation: KD} gives 
\[
\begin{array}{rcl}
\Transfer{G}{G_\delta} \left( \etaABV{\phi}{} \right) 
&=&
\Transfer{G}{G_\delta}
\circ
\operatorname{proj}_0 \left( \etaABV{\phi}{,\pure} \right) \\
&=&
\Transfer{G}{G_\delta}
\circ
\operatorname{proj}_0 \left(
\sum_{i} c(\phi_i,\phi)\ \STD{\phi_i}{\pure}
\right)\\
&=&
\sum_{i} c(\phi_i,\phi)\
\Transfer{G}{G_\delta} \circ
\operatorname{proj}_0 
\left(
\STD{\phi_i}{\pure}
\right) \\
&=&
\sum_{i} c(\phi_i,\phi)\
\operatorname{proj}_\delta
\left(
\STD{\phi_i}{\pure}
\right)\\
&=&

\operatorname{proj}_\delta(\sum_{i} c(\phi_i,\phi)\
\STD{\phi_i}{\pure})
\\
&=& 
\operatorname{proj}_\delta(\etaABV{\phi}{,\pure})
\\
&=&
e(G_\delta)\ 
\etaABV{\phi}{, \delta}.
\end{array}
\] 
This gives the desired result. 
\end{proof}

\begin{corollary}\label{corollary: QST of ABV}
Let $G$ be a quasisplit connected reductive algebraic group over $F$. 
Assume the local Langlands correspondence  (Definition~\ref{definition: LLC}) for $G$.
Recall quasisplit transfer $\QST_{G}$ and $\QST_{G}^\geo$ from Definition~\ref{definition: QST and z-strongly stable} and \ref{definition: QSTgeo}, respectively. 
Then
    \[
    \QST_{G}^\geo (\etaABV{\phi}{}) 
    = 
    \etaABV{\phi}{,\pure}
    \qquad \forall \phi\in \Phi(G).
    \]
If the local Langlands correspondence satisfies  the assumptions of Theorem \ref{thm: stable standards form a basis - quasisplit}, then
    \[
    \QST_{G} (\etaABV{\phi}{}) 
    = 
    \etaABV{\phi}{,\pure}
    \qquad \forall \phi\in \Phi(G).
    \]
\end{corollary}

\section{Arthur packets are ABV-packets, and atomic stability}\label{subsection: Arthur packets revisited, 2}

We now prove the weak form of Vogan's conjecture for orthogonal and symplectic groups, under the assumption that all ABV-packets are preserved under Aubert-Zelevinsky involution (Conjecture \ref{conj AZ}). 
This proof can be regarded as a generalization of the proof of Vogan's conjecture for $\GL_n$ in \cite{Lo24}. Based on the proof, we propose a conjecture on local Langlands correspondence for general group $G$ (Conjecture \ref{conj alt A-packet}), which suggests a definition for Arthur packets for $G$.)

\subsection{ABV-packets and Aubert-Zelevinsky duality}\label{sec AZ}

In this section, we state our main assumption on ABV-packets. First, we recall the definition of Aubert-Zelevinsky involution from \cite{Aub95}. 
Let $G$ be a connected reductive algebraic group. For a representation  $\pi$ in the Grothendieck group $K_{\BC}\Pi(G)$, we define
\begin{align}\label{eq def AZ}
    \AZ(\pi):=\sum_P (-1)^{\mathrm{dim}(A_P)}\mathrm{Ind}_{P}^{G}(\textrm{Jac}_P(\pi))
\end{align}
Here, the sum is taken over all standard parabolic subgroups $P$ of $G$ and $A_P$ is the maximal split torus of the center of the Levi subgroup of $P$.
It is proved in \cite{Aub95} that if $\pi$ is irreducible, then there exists a $\varepsilon\in\{\pm 1\}$ and irreducible representation $\widehat{\pi}$ 
\IN{$\widehat{\pi}$, Aubert-Zelevinsky dual}
such that $\AZ(\pi)=\varepsilon \widehat{\pi}$.
Moreover, $\widehat{\widehat{\pi}}=\pi$.
We call $\widehat{\pi}$ the Aubert-Zelevinsky dual or Aubert-Zelevinsky involution of $\pi$.

From the discussion in \cite{CFMMX22}*{Section 10.3.4}, it is expected that the Aubert-Zelevinsky involution on representation side is closely related to the Fourier transform of perverse sheaves on the geometric side. This is recently established for unipotent representations for a large class of groups in \cite{AO26}. Since Fourier transform interacts nicely with vanishing cycle functors, the following conjecture is expected. 

\begin{conjecture}[Aubert duality preserves ABV-packets and distributions]
\label{conj AZ}
    For any 
    $\phi \in \Phi(G)$, 
    we have 
    \[\Pi_{\widehat{\phi}}^{\ABV}(G)=\{\widehat{\pi}\ | \ \pi \in \Pi_{\phi}^{\ABV}(G)\}.\]
     Moreover, there exists a non-zero complex number $\epsilon_{\phi}$ such that 
     \[
\AZ(\eta_{\phi}^{\ABV}(G))= \epsilon_{\phi}\ \eta_{\widehat{\phi}}^{\ABV}(G).
     \]
\end{conjecture}

\subsection{Desiderata for Arthur packets and ABV packets}
In this section, we state the rest of the desiderata for Arthur packets and ABV packets we need in our argument.
First, we give a technical definition for $L$-parameters. 

\begin{definition}\label{def weak equiv}
    Two admissible homomorphisms $\phi_1, \phi_2: W_F \times \SL_2(\BC) \to \Lgroup{G}$ are said to be {\it weakly equivalent} if, for every finite-dimensional representation $\xi: \Lgroup{G} \to \GL_n(\BC)$, the homomorphisms $\xi \circ \phi_1$ and $\xi \circ \phi_2$ are conjugate in $\GL_n(\BC)$. We write $[\phi]_w$
    to denote the weak equivalence class of a Langlands parameter $\phi$.\IN{$[\phi]_w$, weak equivalence class of Langlands parameters}
\end{definition}

\begin{remark}\label{rmk weakly equivalent}\ 
    \begin{enumerate}
        \item For $G= \SO_{2n+1}(F)$, $\Sp_{2n}(F)$, and $\mathrm{O}_{2n}(F)$, let $\std: \Lgroup{G} \to \GL_{N}(\BC)$ denote the standard embedding, where $N=2n+1$ if $G=\Sp_{2n}(F)$ and $N=2n$ otherwise. The induced map $\std_*\colon \Phi(G) \to \Phi(\GL_N)$ is injective \textnormal{(see \cite{GGP12}*{Theorem 8.1})}. Thus, for these groups, the weak equivalence class of any $L$-parameter $\phi$ consists of a single $L$-parameter. For the remainder of the paper, we will simply write $r(\phi) := r_{\std}(\phi)$ when $G$ is one of these classical groups.
        \item For $G=\RG_2(F)$, consider the sequence of embeddings
        \[
            \RG_2(\BC) \xrightarrow{\iota} \mathrm{Spin}_7(\BC) \xrightarrow{\mathrm{std}} \SO_{7}(\BC).
        \]
        The composition $\mathrm{std}_* \circ \iota_* : \Phi(\RG_2) \to \Phi(\Sp_6)$ is also injective \textnormal{(\cite{GS23}*{Lemma 2.3(iii)})}. Consequently, the weak equivalence class of any $L$-parameter $\phi$ of $\RG_2$ also consists of a single $L$-parameter.
        \item For $G=\SO_{2n}$ with $n > 2$, again let $\std: \SO_{2n}(\BC) \to \GL_{2n}(\BC)$ be the standard embedding. In this case, $\std(\phi_1)$ and $\std(\phi_2)$ are conjugate in $\GL_{2n}(\BC)$ if and only if $\phi_2 \in \{\phi_1, \phi_1^{c}\}$, where $\phi_1^c$ denotes the image of $\phi_1$ under the outer automorphism. Therefore, $[\phi]_w \subseteq \{\phi, \phi^c\}$, and there exist $L$-parameters $\phi$ for which $[\phi]_w = \{\phi, \phi^c\}$ \textnormal{(see \cite{Mat24}*{Proposition 7.3})}. It is worth noting that in this case, $\dim(C_{\phi}) = \dim(C_{\phi^c})$. 
    \end{enumerate}
\end{remark}

Now we state our working hypothesis.

\begin{assumption}\label{assu final}
We assume that $G$ is a connected reductive algebraic group defined and quasisplit over $F$ with the following properties.
\begin{enumerate}
    \item [(a)]For the group $G$, we assume the following holds.
    \begin{enumerate}
    \item [(a1)]There exists an enhanced local Langlands correspondence for pure inner forms of $G$.

    \item [(a2)]Suppose that $\phi_1, \phi_2 \in \Phi(G)$ are weakly equivalent (see Definition \ref{def weak equiv}) and share the same infinitesimal parameter $\lambda$. Then their corresponding $H_{\lambda}$-orbits of $V_{\lambda}$ have the same dimension.
    \end{enumerate}
    \item [(b)] 
    For the theory of ABV-packets, we assume the following holds.
    \begin{enumerate}
        \item [(b1)] 
        The associated distribution of $\eta_{\phi}^{\ABV}$ is stable for all $\phi$. Moreover, for any infinitesimal parameter $\lambda$, the collection $\{\eta_{\phi}^{\ABV}\tq \phi \in \Phi_{\lambda}(G)\}$ forms a basis of $K_{\BC}\Pi_{\lambda}(G)^{\st}$.
        \item [(b2)] 
        Conjecture \ref{conj AZ} holds.
    \end{enumerate}
    \item [(c)]
    We assume that there is a theory of Arthur packets for $G$. More precisely, we assume that
    for any local Arthur parameter $\psi$ of $G$, there exists a finite set of representations $\Pi_{\psi}(G) \subset \Pi(G)$, called the Arthur packets, with the following properties.
    \begin{enumerate}
        \item [(c1)]
 $\Pi_{\phi_{\psi}}(G) \subseteq \Pi_{\psi}(G) \subseteq \Pi_{\lambda}(G)$.
        \item [(c2)]
        There is some $\eta_{\psi}^{\Art}\in K_{\BC}\Pi_{\lambda}(G)^\st$ with support $\Pi_{\psi}(G)$. 
        \item [(c3)]
 $\Pi_{\widehat{\psi}}(G)= \{\widehat{\pi}\ | \ \pi \in \Pi_{\psi}(G)\}$.
        \item [(c4)]
        If $\pi \in \Pi_{\psi}(G)$ then $\phi_{\pi} \geq \phi_{\psi}$.
    \end{enumerate}
\end{enumerate}
\end{assumption}

\begin{remark} Here are the motivation and status of these Working hypotheses.
 \begin{enumerate}
     \item [(a)] The local Langlands correspondence for $p$-adic groups is currently known for $\GL_n$ by \cites{Hen00, HT01, Sch13}, for classical groups by \cites{Art13, Mok15, KMSW14, MR18, Ish24}, and for $G_2$ by \cites{AX22,GS23}. For all of these groups, Condition (a2) holds. See Remark \ref{rmk weakly equivalent}. It will be an interesting problem to classify the groups such that Condition (a2) holds.
     \item [(b)] Working hypothesis (b1) holds if the endoscopic character relations for tempered $L$-packets hold. See Section \ref{subsection: basis of stable standards Bin}. 
     For (b2), Conjecture \ref{conj AZ} can be deduced from \cite{AO26} for certain groups $G$, when combined with results from Part~\ref{part: VI Schurmann}. 
     \item [(c)] In the cases that Arthur packets are defined, Working Hypotheses (c1)-(c3) are all satisfied. For (c4), it is verified for $\Sp_{2n}(F)$ and $\SO_{2n+1}(F)$ in \cite{HLLZ25}, and for pure inner forms of $\O_{2n}(F)$ in \cite{HLL26}. As a corollary of Theorem \ref{theorem intro: Vogan's conjecture} that Arthur packets are ABV-packets, (c4) also holds for unitary groups.
 \end{enumerate}
\end{remark}

\subsection{Weak form of Vogan's conjecture}

We now show that the weak form of Vogan's conjecture follows from Working Hypothesis \ref{assu final}. For $\eta \in K_{\BC}\Pi_{\lambda}(G)$ we may write
\[
\eta = \sum_{\pi \in \Pi_{\lambda}(G)} a_{\pi} \pi
\]
with $a_{\pi} \in \BC$. Define the support $\supp(\eta) := \{\pi \in \Pi_{\lambda}(G) \mid a_{\pi} \neq 0\}$.

\begin{theorem}\label{thm Vogan}
    Suppose that $G$ satisfies Working Hypothesis \ref{assu final}. For any nonzero $\eta \in K_{\BC}^\st\Pi_{\lambda}(G)$,
    the following are equivalent.
    \begin{enumerate}
    \item [(i)] 
    $\supp(\eta) \subseteq \Pi_{\psi}(G)$.
    \item [(ii)] $\supp(\eta) \subseteq \{\pi\in \Pi_{\lambda}(G) \ | \ \phi_{\pi}\geq \phi_{\psi},\ \ \phi_{\widehat{\pi}} \geq \phi_{\widehat{\psi}} \}$.
    \item [(iii)] $\eta$ is a nonzero scalar multiple of $\eta_{\phi_{\psi}}^{\ABV}$.
    \end{enumerate}
    In particular, the {\it weak form of Vogan's Conjecture}\IN{Vogan's conjecture, weak form} holds:
    \[
    \Pi_{\psi}(G)
    =
    \Pi_{\phi_{\psi}}^{\ABV}(G).
    \]
\end{theorem}

\begin{proof}
We first show that Condition (i) implies Condition (ii). Indeed, it suffices to show that for any  $\pi\in \Pi_{\psi}$, we must have $\phi_{\pi}\geq \phi_{\psi}$ and $\phi_{\widehat{\pi}} \geq \phi_{\widehat{\psi}}$. The first inequality is exactly the Working Hypothesis \ref{assu final}(c4). By Working Hypothesis \ref{assu final}(c3), we have $\widehat{\pi}\in \Pi_{\widehat{\psi}}$, which gives us the second inequality again by Working Hypothesis \ref{assu final}(c4).

Next, we show that Condition (ii) implies Condition (iii). By Working Hypothesis \ref{assu final}(b1), we write
\begin{align}\label{eq Vogan 1}
     \eta= \sum_{\phi \in \Phi_{\lambda}} a_{\phi}\, \eta_{\phi}^{\ABV}
\end{align}
for some collection of complex numbers $a_{\phi}$. Let $\Sigma:= \{ \phi\in \Phi_{\lambda}\ | \ a_{\phi} \neq 0\}$ and let $\phi'$ be any minimal member in $\Sigma$ under the closure ordering. Then $\Pi_\phi(G)^\pure \subseteq \Pi^\ABV_\phi(G)^\pure$ implies that the $\supp(\eta_{\phi'}^{\ABV}) \supseteq \Pi_{\phi'}(G)$. Moreover, 
the minimality of $\phi'$ implies that $\Pi_{\phi'}(G) \cap \Pi_{\phi''}^{\ABV}(G)$ is empty for any $\phi'' \in \Sigma \setminus \{\phi'\}$. As a consequence, $\supp(\eta) \supseteq \Pi_{\phi'}(G)$ as well. Then Condition (ii) implies that $\phi' \geq \phi_{\psi}$. Since $\phi'$ is an arbitrary minimal element in $\Sigma$, we obtain
\begin{align}\label{eq Vogan 2}
    \Sigma \subseteq \{\phi \ | \ \phi \geq \phi_{\psi} \}.
\end{align}

On the other hand, we may apply the Aubert-Zelevinsky involution on $\eta$. Working Hypothesis \ref{assu final}(b2) implies that there are nonzero complex numbers $\epsilon_{\phi}$ such that
\begin{align*}
     \AZ(\eta)= \sum_{\phi \in \Phi_{\lambda}(G)} \epsilon_{\phi}\, a_{\phi}\  \eta_{\widehat{\phi}}^{\ABV}.
\end{align*}
Condition (ii) implies 
\[
\supp(\AZ(\eta))
\subseteq \{\pi\in \Pi_{\lambda}(G) \ | \ \phi_{\pi}\geq \phi_{\widehat{\psi}},\ \ \phi_{\widehat{\pi}} \geq \phi_{{\psi}} \}.
\]
Therefore, the same argument in the previous paragraph implies 
    \begin{align}\label{eq Vogan 4}
    \Sigma \subseteq \{\phi \ | \ \widehat{\phi} \geq \phi_{\widehat{\psi}} \}.
\end{align}
Combining \eqref{eq Vogan 2}, \eqref{eq Vogan 4}, we see that 
\[ \Sigma \subseteq \{ \phi\in \Phi_{\lambda}(G)\ | \ \phi\geq \phi_{\psi},\ \ \widehat{\phi} \geq {\phi_{\widehat{\psi}}}   \}=: \Sigma'.\]

We claim that $\Sigma'=\{\phi_{\psi}\}$ by our assumption on $G$. Indeed, if $\phi \in \Sigma'$, then for any algebraic representation $\xi: \Lgroup{G} \to \GL_N(\BC)$, we have the dual inequalities $\xi \circ \phi \geq \xi \circ \phi_{\psi}$ and $\xi \circ \widehat{\phi} \geq \xi \circ \phi_{\widehat{\psi}}$. Then \cite{Lo24}*{Lemma 1.3} implies that $\xi\circ \phi= \xi \circ \phi_{\psi}$ as an $L$-parameter of $\GL_N(F)$. Therefore, we conclude that $\phi$ is weakly equivalent to $\phi_{\psi}$. Finally, since $\phi \geq_{C} \phi_{\psi}$ by the definition of $\Sigma'$, we see that if $ \phi \neq \phi_{\psi}$, then $\dim(C_{\phi})> \dim(C_{\phi_{\psi}})$, which contradicts the Working Hypothesis \ref{assu final}(a2). We conclude that $\phi=\phi_{\psi}$ must hold. 

By the above claim, \eqref{eq Vogan 1} becomes 
\[ \eta= a_{\phi_{\psi}}\,  \eta_{\phi_{\psi}}^{\ABV}\]
for some nonzero complex number $a_{\phi_{\psi}}$. This completes the proof that Condition (ii) implies Condition (iii).

Finally, we show that Condition (iii) implies Condition (i). Take the distribution $\eta_{\psi}^{\Art}$ in Working Hypothesis \ref{assu final}(c2). Since we have verified that Condition (i) implies Condition (iii), we have $\eta_{\psi}^{\Art}= a_{\phi_{\psi}} \eta_{\phi_{\psi}}^{\ABV}$ for some  nonzero complex number $a_{\phi_{\psi}}$. 
 Therefore, we have
\[ \Pi_{\phi}^{\ABV}(G)= \supp( \eta_{\phi_{\psi}}^{\ABV}) =\supp( \eta_{\psi}^{\Art})= \Pi_{\psi}(G),\]
where the first equality follows from the definition of ABV-packets and $\eta_{\phi_{\psi}}^{\ABV}$ (Definitions \ref{definition: ABV-packet}, \ref{definition: etaABV}).
This completes the proof that Condition (iii) implies Condition (i), and Vogan's Conjecture.
\end{proof}

\begin{remark}\label{rmk CO26}
The above theorem is proved independently by Ciubotaru-Okada in \cite{CO26} in the unipotent setting. The key ingredient of the proof is the equality
\[ \{ \phi\in \Phi(G)\ | \ \phi\geq \phi_{\psi}, \widehat{\phi} \geq \phi_{\widehat{\psi}} \}= \{\phi_{\psi}\} \]
for any local Arthur parameter $\psi$ of $G$. In our proof, we use Working hypothesis \ref{assu final}(a2) to relate this problem to general linear groups, which is solved combinatorially in \cite{Lo24}. On the other hand, in \cite{CO26}*{Theorem 4.2}, they directly proved the above equality for general unipotent $\psi$ without relating to general linear groups. Hence, they do not need Working hypothesis \ref{assu final}(a2).
\end{remark}

\section{A conjectural characterization of Arthur packets}

We expect that the equivalence between Conditions (i) and (ii) in Theorem \ref{thm Vogan} hold for general $G$, which we state as the following conjecture.

\begin{conjecture}\label{conj alt A-packet}
Let $G$ be a connected reductive algebraic group defined and quasisplit over $F$. Suppose that the local Langlands correspondence $\LLC: \Pi(G) \to \Phi(G)$ is constructed. For any infinitesimal parameter $\lambda$ of $G$, recall the definition $\Pi_{\lambda}(G):= \{ \pi \in \Pi(G)\ | \ \lambda_{\LLC(\pi)}= \lambda\}$.
For any local Arthur parameter $\psi$ of $G$, the vector space of stable distributions supported on the set 
\begin{align}\label{eq conj alt A-packet}
    \{\pi \in \Pi_{\lambda_{\phi_{\psi}}}(G)\ | \ \LLC(\pi) \geq \phi_{\psi},\ \LLC(\widehat{\pi}) \geq \phi_{\widehat{\psi}}\}
\end{align}
is one dimensional.
\end{conjecture}
The above expectation suggests a conjectural definition of Arthur packets as a set for any connected reductive group $G$ with a local Langlands correspondence $\LLC: \Pi(G) \to \Phi(G)$. Namely, we may define $\Pi_{\psi}(G)$ as the support of any nonzero stable distribution supported on \eqref{eq conj alt A-packet}. 

\part{The geometrization of Kottwitz-Shelstad transfer}
\label{part: IV Lift}
To get the stronger form of Vogan's conjecture, we need to show that the ABV distributions for $L$-parameters of Arthur type defined in Section \ref{subsection: ABV-distributions} satisfy endoscopic character identities. To this end, we geometrize Kottwitz-Shelstad transfer in out setting. More precisely, we define the linear transformation 
\[\Lift{G',s}{G}: K_{\BC}\Pi_{\lambda'}^{\std}(G') \to K_{\BC}\Pi(G)\]
and evaluate it on the stable ABV-distributions; see Theorem \ref{theorem: Lift of ABV}. Assuming certain desiderata of the local Langlands Correspondence, we show that $\Lift{G',s}{G}$ coincides with $\Transfer{G',s}{G}$; see Theorem \ref{theorem: lift = transfer}.

\section{Endoscopic lift}\label{section: endoscopic lift}

In this section we use the geometric framework of Section~\ref{subsection: Vogan's LLC} to define a linear transformation called {\it endoscopic lift}, from stable invariant distributions on $G'(F)$ (as studied in Sections~\ref{section: stable standards} and \ref{section: etaABV}) to invariant distributions on $G_\delta(F)$ (as studied in Section~\ref{section: endoscopic transfer}). 
We show that endoscopic lift provides a geometrization of endoscopic transfer (both twisted and ordinary), by  calculating this linear transformation on the bases for stable distributions appearing in Section~\ref{subsection: basis of stable standards Bin}.

\subsection{Endoscopic lift}\label{subsection: Lift}

Let $G$ be a group of the sort considered in Definition~\ref{assu: groups we consider} and let
$(G',s, \xi)$ be an endoscopic triple for $G$ as discussed in Section~\ref{subsection: endoscopic triple}.
Let $S$ be the diagonalizable group generated by $s$, which is to say, $S$ is the Zariski-closure of $\langle s \rangle$ in $\dualgroup{G}$.
Then $S$ has a cyclic component group.
Let $\lambda'$ be an infinitesimal parameter for $G'$ and set $\lambda = \xi \circ \lambda'$. Then, since 
\[\textrm{Lie}(\widehat{G'})= \{ x \in \textrm{Lie}(\dualgroup{G})\ | \ \Ad(s)x=x\},\]
it follows that 
\[
\begin{array}{rcl}
V_{\lambda'}
&=& \{ x \in \textrm{Lie}(\dualgroup{G'})\ | \ \Ad(\lambda'(w))x=|w|x\} \\
&=& \{ x \in \textrm{Lie}(\dualgroup{G})\ | \ \Ad(\lambda'(w))x=|w|x,\ \Ad(s)x=x \} \\
&=& V_{\lambda}^{s}. 
\end{array}
\]
Although $H_{\lambda'}$ is a subgroup of $H_{\lambda}^{s}$, equality fails in general. Indeed, the diagram in Table~\ref{eq diagram centralizer endoscopic lift} shows the difference between these groups (see Lemma \ref{lem centralizer}).
Note that $\dualgroup{G'}= Z_{\dualgroup{G^0}}(s)^0$, and the component groups in the bottom row are all finite. Now taking the fixed points with respect to conjugation by $\lambda$, we see that 
\[ 
H_{\lambda'}= 
Z_{\dualgroup{G'}}(\lambda')
= Z_{Z_{{\dualgroup{G^0}}} (s)^0}(\lambda) \leq Z_{Z_{{\dualgroup{G^0}}} (s)}(\lambda)  \leq Z_{{\dualgroup{G}}} (s,\lambda) = H_{\lambda}^{s} \]
 is a sequence of finite index subgroups. Note that $s$ lies in the center of $ H_{\lambda}^{s}$.

\begin{table}[ht]
\caption{Relating the groups $H_{\lambda'}$ and $H_{\lambda}^{s}$.}
    \centering
\begin{equation}
    \begin{tikzcd}
& 0 \ar[d]&  0\ar[d]& & \\
& Z_{\dualgroup{G^0}}(s)^0 \ar[r, equal] \ar[d] &Z_{\dualgroup{G}}(s)^0  \ar[d]&&\\
    0 \ar[r]& Z_{{\dualgroup{G^0}}} (s) \ar[r] \ar[d] &  Z_{{\dualgroup{G}}}(s) \ar[r] \ar[d] & \textrm{im}( Z_{{\dualgroup{G}}}(s) \to \pi_0(\dualgroup{G})) \ar[r] \ar[d,equal] & 0\\
   0\ar[r]& \pi_0(Z_{{\dualgroup{G^0}}} (s)) \ar[r] \ar[d]&  \pi_0(Z_{{\dualgroup{G}}} (s))\ar[r] \ar[d] & \textrm{im}( Z_{{\dualgroup{G}}}(s) \to \pi_0(\dualgroup{G})) \ar[r] & 0 \\
   &0 &0&&
\end{tikzcd}
\end{equation}
\label{eq diagram centralizer endoscopic lift}
\end{table}

 To separate the ``central character" from the action of $s$, we introduce the homomorphism \IN{$H^+_{\lambda'}$} 
 \begin{equation}\label{eq m}
 \calm : H^+_{\lambda'}
 \ceq H_{\lambda'} \times S \to H_\lambda^s
  \end{equation}
  given by multiplication: $\calm(h,s):= hs$.
  \IN{$\calm$}
 This defines an equivariant functor 
\begin{align}\label{eq m star}
    \calm^* : \Per_{H_{\lambda}^s}(V_{\lambda}^s) \to \Per_{H^+_{\lambda'}}(V_{\lambda'}).
\end{align}
\IN{$\calm^*$}
 Note that 
 \[|H_{\lambda'} \backslash V_{\lambda'}| =|H_{\lambda'}^+ \backslash V_{\lambda'}| \geq |H_{\lambda}^{s} \backslash V_{\lambda}^{s}|,   
 \]
which are all finite.

\begin{lemma}\label{lemma: m}  
Let $G$, $(G',s,\xi)$ and $S$ be as above. 
Let $\calm^* : D^b_{c,H^\sigma_{\lambda}}(V_{\lambda}^{s}) \to D^b_{c,H^+_{\lambda}}(V_{\lambda'})$
be the equivariant derived functor induced by $\calm: H_{\lambda'}^+ \to H_{\lambda}^{s}$ from \eqref{eq m}.
For every $\mathcal{F'} \in \Per_{H_{\lambda}^{s}}(V_{\lambda}^{s})$ there is a unique decomposition
\[
\calm^*\mathcal{F'} 
= 
\mathop{\oplus}\limits_{\chi\in \Irrep(\pi_0(S))}
\mathcal{F}''_\chi\otimes \locsys{\chi}{V_{\lambda'}}.
\]
where $\mathcal{F}''_\chi\in \Per_{H_{\lambda'}}(V_{\lambda'})$, for each $\chi\in \Irrep(\pi_0(S))$.
\end{lemma}

\begin{proof}
$\Per_{H^+_{\lambda'}}(V_{\lambda'}) 
= 
\Per_{H_{\lambda'}\times S}(V_{\lambda'})  
\iso 
\Per_{H_{\lambda'}}(V_{\lambda'})\boxtimes \Per_{S}(0)
\iso
\Per_{H_{\lambda'}}(V_{\lambda'})\boxtimes \Rep(\pi_0(S))$
\end{proof}

\begin{definition}\label{definition: Lift}
Let
$(G', s,\xi)$ be an endoscopic triple for $G$, where $G$ satisfies Assumption~\ref{assu: groups we consider}.
Let $\delta$ be a pure inner form for $G$. 
Let $\lambda'$ be an infinitesimal parameter for $G'$
and set $\lambda = \xi \circ \lambda'$.
Define the linear transformation $\Lift{G',s}{G_\delta}$ \IN{$\Lift{G',s}{G_\delta}$, endoscopic lift} by the following diagram, 
\begin{equation*}   
    \begin{tikzcd}
        K_\CC \Pi_{\lambda'}(G')^\mathrm{std} \arrow{d}[swap]{\QST_{G'}^{\geo}} \arrow{rrrr}{\Lift{G',s}{G_\delta}} 
                &&&& K_\CC \Pi_{\lambda}(G_\delta) \\
        K_\CC \Pi^\mathrm{pure}_{\lambda'}(G') \arrow
        {rrrr}[swap]{\,^t(T\circ\, \calm^{\ast}\circ\, d\xi^*)} 
        &&&& K_\CC \Pi^\mathrm{pure}_{\lambda}(G) \arrow{u}[swap]{\mathrm{proj}_\delta} 
    \end{tikzcd}
\end{equation*}
where:
\begin{enumerate}
    \item $\QST_{G'}^{\geo} : K_\CC\Pi_{\lambda'}(G')^{\std} \to K_\CC \Pi_{\lambda'}(G')^\pure$ is defined in Definition \ref{definition: QSTgeo};
    \item 
    $
d\xi^* : 
D^b_{c,H_{\lambda}}(V_{\lambda}) 
\to
D^b_{c,H^{s}_{\lambda}}(V^s_{\lambda})
$
is the equivariant functor defined by restriction to $s$-fixed points;
\item 
$
\calm^* : \Per_{H_{\lambda}^s}(V_{\lambda}^s) \to \Per_{H^+_{\lambda'}}(V_{\lambda'})
$
is the equivariant functor defined by $\calm$ appearing in \eqref{eq m};
\item 
$
T : K_\CC \Per_{H^+_{\lambda'}}(V_{\lambda'}) \to 
    K_\CC \Per_{H_{\lambda'}}(V_{\lambda'})
$
\IN{$T: K_\CC \Per_{H^+_{\lambda'}}(V_{\lambda'}) \to 
    K_\CC \Per_{H_{\lambda'}}(V_{\lambda'})$}
is the linear transformation defined by 
\[
\mathcal{F''}\otimes \locsys{\chi}{V_\lambda} \mapsto \trace(s,\chi)\, \mathcal{F''};
\]
    \item
    $\,^t(T\circ\, \calm^*\, \circ\, d\xi^*)$ is the transpose of $T\circ\, f^*\, \circ\, d\xi^*$ using the pairing of \eqref{equation: pairing};
    \item $\operatorname{proj}_\delta$ is projection arising from the disjoint union \[
    \Pi_\lambda(G)^\pure = \mathop{\coprod}\limits_{[\delta]\in H^1(F,G)_{\theta}} \Pi_\lambda(G_\delta).
    \]
\end{enumerate}
\end{definition}

\begin{table}[h!]
\label{equation: Todxi}
\caption{Geometric operations defining endoscopic Lift}
\[
\begin{tikzcd}
K_\CC\Per_{H_{\lambda'}}(V_{\lambda'})
\arrow[bend right=20]{dd}[swap]{\mathrm{induction}}
&&&& 
K_\CC\Per_{H_{\lambda}}(V_{\lambda})
\arrow{ddll}{d\xi^*}
\arrow{llll}[swap]{T\circ\, \calm^*\,\circ\, d\xi^*}
\\
&&&&
\\
K_\CC\Per_{H^+_{\lambda'}}(V_{\lambda'})
\arrow{uu}[swap]{T}
&&
K_\CC\Per_{H^s_{\lambda}}(V^s_{\lambda})
\arrow{ll}{\calm^*}
\arrow{lluu}
&&
\end{tikzcd}
\]
\end{table}

\begin{theorem}\label{theorem: lift = transfer}
Suppose $G$ satisfies Assumption~\ref{assu: groups we consider}. 
Let $(G',s,\xi)$ be an endoscopic triple for $G$. Assume Hypothesis~\ref{hypothesis: KLH}. Then
\[
\Lift{G',s}{G_\delta}(\STD{\phi'}{}) = e(G_\delta)\, \STD{\xi\circ\phi',s}{G_\delta}
\qquad \forall\ 
\phi'\in \Phi_{\lambda'}(G').
\]
In particular, if the local Langlands correspondence satisfies the assumptions of Theorem \ref{thm: stable standards form a basis - quasisplit}, then for all pure inner forms $\delta$, we have
    \[
    \Transfer{G',s}{G_\delta} = \Lift{G',s}{G_\delta}.
    \]
on the space $\mathcal{D}(G')^\st$ of stable distributions on $G'(F)$ supported by characters of representations.
\end{theorem}

\begin{proof}
Let $\lambda'$ be any infinitesimal parameter of $G'$.
In Theorem~\ref{thm: stable standards form a basis - quasisplit} we saw that 
\[
K_\CC \Pi_{\lambda'}(G')^{\st}
=
\operatorname{span}_\CC\{ 
\STD{\phi'}{} \tq \phi'\in \Phi_{\lambda'}(G') \} 
\qquad\mathrm{and}\qquad
\dim K_\CC \Pi_{\lambda'}(G')^{\st}
= 
| \Phi_{\lambda'}(G') |.
\]
In Proposition~\ref{theorem: transfer of STD} we saw 
\[
\Transfer{G',s}{G_\delta}(\STD{\phi'}{}) 
= 
e(G_\delta)\, 
 \STD{\xi\circ\phi',s}{\delta}
\qquad \forall\ 
\phi'\in \Phi_{\lambda'}(G').
\]
Thus, the second part of the theorem is a consequence of the first part for the computation of $\Lift{G',s}{G_\delta}(\STD{\phi'}{})$, which we verify now.
It suffices to show that 
\[ \langle 
    \Lift{G',s}{G_\delta}(\STD{\phi'}{}), \stsh{\tau}{C_{\phi}} \rangle= \langle 
    e(G_\delta)\,
    \STD{\xi\circ\phi',{s}}{\delta}, \stsh{\tau}{C_{\phi}} \rangle\]
    for all $(\phi,\tau) \in \Phi_{\lambda}^{\mfr{e}}(G)$.

First, we compute the right hand side. It follows from the definition that
\begin{equation}\label{equation: Lift=Trans 0}
    \begin{array}{rcl}
    
    \langle 
    e(G_\delta)\,
    \STD{\xi\circ\phi',{s}}{\delta}, \stsh{\tau}{C_{\phi}} \rangle \vspace{0.5em}  
    &=& 
    e(G_\delta)\,
    \sum\limits_{\rho\in \widehat{A}_{\xi\circ\phi'}[\delta]} \trace(s, \rho)\, 
    \langle \SR(\xi\circ\phi',\rho), \stsh{\tau}{C_{\phi}}  \rangle.
    \end{array}
\end{equation}
By the Kazhdan-Lusztig Hypothesis~\ref{hypothesis: KLH}, in this sum (note that $\rho \in \widehat{A}_{\xi\circ\phi'}[\delta]$),
\[ 
\langle \SR(\xi\circ\phi',\rho), \stsh{\tau}{C_{\phi}} \rangle
= 
\begin{cases}
    e(G_\delta) &\text{ if } C_{\xi\circ\phi'}=C_{\phi} \text{ and } \tau=\rho,\\
    0 &\text{ otherwise.}
\end{cases}
\]
Thus,
\begin{equation}\label{eq Lift=trans RHS}
    \begin{array}{rcl}
   \langle 
    e(G_\delta)\,
    \STD{\xi\circ\phi',{s}}{\delta}, \stsh{\tau}{C_{\phi}} \rangle
    &=&
    \begin{cases}
     \trace(\sigma, \tau) & \text{if }  [\xi\circ\phi'] = [\phi] \text{ and }\tau \in \widehat{A}_{\xi\circ\phi'}[\delta]\\
     0 & \text{otherwise}.
     \end{cases}
    \end{array}
\end{equation}

We now calculate $\langle \Lift{G',\sigma}{G_\delta}(\STD{\phi'}{}), \stsh{\tau}{C_{\phi}} \rangle$. By Lemma \ref{lem: transpose of projection standard}, this pairing is nonzero only when $\tau \in \widehat{A}_{C_{\phi}}[\delta]$, which we assume now. It follows from Definition~\ref{definition: Lift} that
\begin{equation}\label{equation: Lift=Trans 1}
    \begin{array}{rcl}
    \langle \Lift{G',s}{G_\delta}(\STD{\phi'}{}), \stsh{\tau}{C_{\phi}} \rangle
    &=&
    \langle \QST_{G'}^{\geo}(\STD{\phi'}{}), (T \circ \calm^{\ast}\circ d\xi^*)\stsh{\tau}{C_{\phi}} \rangle,
    \qquad (\phi,\tau) \in \Phi^\mathfrak{e}_{\lambda}(G).
    \end{array}
\end{equation}
By definition,
$\QST_{G'}^{\geo}(\STD{\phi'}{}) = \STD{\phi}{\pure}$.
Thus, the right-hand side of \eqref{equation: Lift=Trans 1} becomes
\begin{equation}\label{equation: Lift=Trans 2}
    \begin{array}{rcl}
    &&
    \hskip-30pt \langle \QST_{G'}^{\geo}(\STD{\phi'}{}), (T\circ \calm^{\ast} \circ  d\xi^*)\stsh{\tau}{C_{\phi}} \rangle \\
    &=& 
    \langle \STD{\phi'}{\pure}, (T\circ \calm^{\ast} \circ  d\xi^*)\stsh{\tau}{C_{\phi}} \rangle \vspace{0.5em} \\
    &=&
    \sum\limits_{\tau'\in \widehat{A_{\phi'}}} e(\tau') \trace(1,\tau')\,
    \langle \SR(\phi',\tau'), (T\circ \calm^{\ast} \circ  d\xi^*)\stsh{\tau}{C_{\phi}}\rangle,
    \end{array}
\end{equation}
so we now turn our attention to the term $\langle \SR(\phi',\tau'), (T\circ \calm^{\ast} \circ  d\xi^*)\stsh{\tau}{C_{\phi}} \rangle$.
In the following, we introduce certain notations and then compute $(T \circ \calm^{\ast} \circ  d\xi^*)\stsh{\tau}{C_{\phi}}$ explicitly using \cite{ABV92}*{Proposition 23.7(b)}. 

Write $C=C_{\phi}$ for short. 
We decompose 
\begin{align}\label{eq lift=trans decomp of orbit}
    C\cap V_{\lambda}^{s }=  \bigsqcup_{i=1}^r C_i,\qquad C_i= \bigsqcup_{j=1}^{r_i} C_{ij},
\end{align}
where each $C_i$ is an 
$H_{\lambda'}^{s}$-orbits inside $C\cap V_{\lambda}^{s}$, and each $C_{ij}$ is an $H_{\lambda'}^+=H_{\lambda'} \times S$-orbit inside $C_i$. 
Note that $C_{ij}$ is also a single $H_{\lambda'}$-orbit since $s$ acts trivially on $V_{\lambda}^{s}$. 
We will write $C_{ij}'$ for $C_{ij}$ when we regard it as an $H_{\lambda'}$-orbit.
Fixing any $x \in C_{ij}$, we obtain homomorphisms between component groups
\[ 
A_{C_{ij}'}=  \pi_0(Z_{H_{\lambda'}} (x)) \to A_{C_{ij}}= \pi_0(Z_{H_{\lambda'} \times S}(x)) \to A_{C_i}=\pi_0(Z_{H_{\lambda}^{s} }(x)) \to A_{C}=\pi_0(Z_{H_{\lambda} }(x)). \]
Again, since $s$ acts trivially on $V_{\lambda'}= V_{\lambda}^{s}$, we see that $Z_{H_{\lambda'}^+}(x)= Z_{H_{\lambda'}}(x) \times S$; hence,
\[ A_{C_{ij}}= \pi_0(Z_{H_{\lambda'}}(x)) \times \pi_0(S)= A_{C_{ij}'} \times B, \]
where we write $B:= \pi_0(S)$, a finite cyclic group generated by the image of $s$. 
As a consequence, we have a bijection $\widehat{A}_{C_{ij}} \simeq \widehat{A}_{C_{ij}'} \times \widehat{B}$, and any irreducible representation $\rho_{ij}$ of $A_{C_{ij}}$ can be written as $\rho_{ij}= \rho_{ij}' \otimes \rho_{ij,B}$, where $\rho_{ij}'= \rho_{ij}|_{A_{C_{ij}'}}$ is an irreducible representation of $A_{C_{ij}'},$ and $\rho_{ij,B} \in \widehat{B}$ is the central character of $\rho_{ij}$ restricted to $B$. Also, if $\rho_{i} \in \widehat{A}_{C_i}$, then since the image of $B$ lies in the center of $A_{C_i}$, we see that $\rho_i|_{B}$ is an isotypic direct sum of the restriction of the central character of $\rho_i$ to $B$, which we denote by $\rho_{i,B}$.

Now, we decompose the restriction of representations of these component groups. Let $\tau \in \widehat{A}_{C}$. We write
\[ 
\tau|_{A_{C_i}}= \bigoplus_{\rho_i \in \widehat{A}_{C_{i}}} m(\tau, \rho_i) \rho_i, 
\]
and
\[
\begin{array}{rcl}
\rho_i|_{A_{C_{ij}}}
&=& \bigoplus_{ \rho_{ij}\in \widehat{A}_{C_{ij}}} m(\rho_i, \rho_{ij}) \rho_{ij} \\
&=& \bigoplus_{\rho_{ij}' \in \widehat{A}_{C_{ij}'}}  m(\rho_i, \rho_{ij}' \otimes \rho_{i,B}) \rho_{ij}' \otimes \rho_{i,B}, 
\end{array}
\]
where $m(\tau, \rho_i)$ and $m(\rho_i, \rho_{ij}) $ are non-negative integer represents the multiplicity. The trace of $\tau$ at $\sigma$ can be expressed in terms of this decomposition: fixing any pair $(i,j)$, we have
\begin{equation}\label{equation: Lift of ABV, 4}
\begin{array}{rcl}
 \trace(s, \tau)
 &=& \trace(s, \tau|_{A_{C_{i}}}) \\
 &=& \sum_{\rho_i\in \widehat{A}_{C_{i}}} m(\tau, \rho_{i}) \trace(s, \rho_i) \\
&=&\sum_{\rho_i\in \widehat{A}_{C_{i}}} \rho_{i,B}(s)   m(\tau, \rho_{i}) \dim(\rho_i)\ \\
&=&\sum_{\rho_i\in \widehat{A}_{C_{i}}}  \rho_{i,B}(s)  m(\tau, \rho_{i})   \dim(\rho_i |_{A_{C_{ij}'}})\\
&= &\sum_{\rho_i\in \widehat{A}_{C_{i}}}  \sum_{\rho_{ij}'\in \widehat{A}_{C_{ij}'}}
\rho_{i,B}(s)  m(\tau, \rho_{i}) m(\rho_i, \rho_{ij}' \otimes \rho_{i,B})  \dim( \rho_{ij}')\\
&=& \sum_{\rho_i\in \widehat{A}_{C_{i}}}  \sum_{\rho_{ij}'\in \widehat{A}_{C_{ij}'}} \rho_{i,B}(s) \trace(1, \rho_{ij}')  m(\tau, \rho_{i}) m(\rho_i, \rho_{ij}' \otimes \rho_{i,B}). 
\end{array}
\end{equation}

With the above preparation, we are ready to apply \cite{ABV92}*{Proposition 23.7(b)}. First, we pull back along $d{\xi^{\ast}}$:
\begin{align*}
    d\xi^{\ast} ( \stsh{\tau}{C_{\phi}} )
    &=  
    \sum_{ i=1 }^r 
    \hskip5pt 
    \sum_{\rho_i \in \widehat{A}_{C_i}} 
    \hskip5pt
    m(\tau, \rho_i) \stsh{\rho_i}{C_i}.
\end{align*}
Next, we apply $\calm^{\ast}$:
\begin{align*}
    ( \calm^{\ast} \circ d\xi^{\ast}) ( \stsh{\tau}{C_{\phi}} )
    &=  
    \sum_{ i=1 }^r 
    \hskip5pt \sum_{ j=1 }^{r_i} 
    \hskip5pt 
    \sum_{\rho_i \in \widehat{A}_{C_i}}  
    \hskip5pt
     \sum_{\rho_{ij}' \in \widehat{A}_{C_{ij}'}} 
    \hskip5pt
    m(\tau, \rho_i) m(\rho_i, \rho_{ij}' \otimes \rho_{i,B})\stsh{\rho_{ij}'\otimes \rho_{i,B}}{C_{ij}}.
\end{align*}
Finally, we apply $T$. Note that $T(\stsh{\rho_{ij}'\otimes \rho_{i,B}}{C_{ij}})= \rho_{i,B}(\sigma) \stsh{\rho_{ij}'}{C_{ij}'}. $ 
Therefore,
\begin{align}\label{equation: Lift=Trans 3}
    (T \circ  \calm^{\ast} \circ d\xi^{\ast}) ( \stsh{\tau}{C_{\phi}} )
    &=  
    \sum_{ i=1 }^r 
    \hskip5pt \sum_{ j=1 }^{r_i} 
    \hskip5pt 
    \sum_{\rho_i \in \widehat{A}_{C_i}} 
    \hskip5pt
     \sum_{\rho_{ij}' \in \widehat{A}_{C_{ij}'}}  
    \hskip5pt
   \rho_{i,B}({\sigma}) m(\tau, \rho_i) m(\rho_i, \rho_{ij}' \otimes \rho_{i,B}) \stsh{\rho_{ij}'}{C_{ij}'}.
\end{align}

Now we return to the computation of $\langle \Lift{G',\sigma}{G_\delta}(\STD{\phi'}{}), \stsh{\tau}{C_{\phi}} \rangle$. Recall the Kazhdan-Lusztig Hypothesis~\ref{hypothesis: KLH} for $G'$ gives
\[ 
\langle \SR(\phi', \tau'), \stsh{\rho'}{C'} \rangle
= 
\begin{cases}
    e(\tau') &\text{ if } C'=C_{\phi'} \text{ and } \tau'=\rho',\\
    0 &\text{ otherwise.}
\end{cases}
\]
It now follows that
$\langle \Lift{G',s}{G_\delta}(\STD{\phi'}{}), \stsh{\tau}{C_{\phi}} \rangle$ is zero unless $C_{\phi'}\subseteq C_{\phi}$, so $C_{\phi'}=C_{ij}'$ for some pair $(i,j)$, which we now assume. Then, with \eqref{equation: Lift of ABV, 4}, we have
\begin{equation*}
    \begin{array}{rcl}
     \langle \Lift{G',s}{G_\delta}(\STD{\phi'}{}), \stsh{\tau}{C_{\phi}} \rangle
     &=& 
     \sum\limits_{\tau'\in \widehat{A}_{\phi'}} e(\tau') \trace(1,\tau') \langle \SR(\phi',\tau'), (T\circ \calm^{\ast} \circ  d\xi^*)(\stsh{\tau}{C_{\phi}}) \rangle\\
    &=& 
    \sum_{\rho_i \in \widehat{A}_{C_i}} 
    \hskip5pt
     \sum_{\rho_{ij}' \in \widehat{A}_{C_{ij}'}}  
    \hskip5pt
  \trace(1,\rho_{ij}') \rho_{i,B}(\overline{s}) m(\tau, \rho_i) m(\rho_i, \rho_{ij}' \otimes \rho_{i,B}) \\
    &=& \trace(s, \tau).
    \end{array}
\end{equation*}
We conclude that
\begin{equation}\label{eq Lift=trans LHS}
    \langle \Lift{G',s}{G_\delta}(\STD{\phi'}{}), \stsh{\tau}{C_{\phi}} \rangle
= 
\begin{cases}
    \trace (s, \tau ) & \text{ if } C_{\phi'} \subseteq C_\phi \cap V_{\lambda}^{s} \text{ and }\tau \in \widehat{A}_{C_{\phi}}[\delta] \\
    0 & \text{otherwise}.
\end{cases}
\end{equation}
Note that $C_{\phi'} \subseteq C_\phi \cap V_{\lambda}^{s}$ if and only if $[\xi\circ\phi'] = [\phi] $. Comparing \eqref{eq Lift=trans LHS} and \eqref{eq Lift=trans RHS}, we conclude that 
\[\Lift{G',s}{G_\delta}(\STD{\phi'}{}) = e(G_\delta)\, \STD{\xi\circ\phi',{s{\hat\theta}}}{G_\delta}
\qquad \forall\ 
\phi'\in \Phi_{\lambda'}(G').\]
This completes the proof of the theorem.
\end{proof}

\subsection{An application of the fixed-point formula}\label{subsection: FPF}

Our next goal is to evaluate the geometric lift on the stable ABV-distribution (Theorem~\ref{theorem: Lift of ABV} below). In this section, we explain a key step in the proof of the proof and explain how it follows from a fixed-point formula for vanishing cycles, Theorem~\ref{theorem: Schurmann}, established by Jörg Schürmann in the appendix to this paper.

\begin{proposition}\label{theorem: FPF}
Suppose that $G$ satisfies Definition~\ref{assu: groups we consider} and let $(G',s,\xi)$ be and endoscopic triple for $G$.
Let $\phi'$ be a Langlands parameter for $G'$.
Then
\[
(-1)^{d(\phi')} \trace({s},\NEvs_{z'}\ d\xi^*\, \mathcal{F})
=
(-1)^{d(\xi\circ\phi')}
\trace({s},\NEvs_{d\xi(z')} \mathcal{F})
\] 
for all $z'\in \Lambda^\mic_{\lambda'}$ such that $d\xi(z')\in \Lambda^\mic_{\lambda}$.
\end{proposition}

\begin{proof}
Recall $\NEvs_{C}$ from Definition~\ref{definition: nevs} and \eqref{equation: appendix NEvs}:
\[
\NEvs_{C} = \Evs_{C}\IC(\locsys{\1}{C})^\vee \otimes \Evs_{C} = \mathcal{T}_{C}^\vee\otimes \Evs_{C}
\]
so
\[
\NEvs_{z} = j_{z}^* \NEvs_{C_x} = \mathcal{T}_{z}^\vee \otimes \Evs_{z}, 
\qquad
\forall z\in \Lambda^\reg_{\lambda}.
\]
Thus, Lemmas~\ref{lemma: eccentricity} and \ref{lemma: Evs FPF} now give
\[
\begin{array}{rcl}
&& \hskip-20pt (-1)^{d(\phi')}
\trace(s,\NEvs_{z'}d\xi^*\mathcal{F}) \\
&=&
(-1)^{d(\phi')}
\trace(s,\mathcal{T'}^\vee_{z'})\times \trace(s, \Evs_{z'} d\xi^*\mathcal{F})\\
&=&
(-1)^{e(\phi')}
\trace(s,\mathcal{T'}^\vee_{z'})\times 
(-1)^{d(\phi')-e(\phi')} \trace(s, \Evs_{z'} d\xi^*\mathcal{F})\\
&=&
(-1)^{e(\phi)}
\trace(s,\mathcal{T}^\vee_{z})\times 
(-1)^{d(\phi)-e(\phi)} \trace(s, \Evs_{z} \mathcal{F})\\
&=&
(-1)^{d(\phi)}
\trace(s,\mathcal{T}^\vee_{z}\otimes \Evs_{z} \mathcal{F})\\
&=&
(-1)^{d(\phi)}
\trace(s,\NEvs_{z} \mathcal{F}),
\end{array}
\]
where, as above, $z = d\xi(z')$.
This completes the proof of Proposition~\ref{theorem: FPF}. 
\end{proof}
As an application, this describes $\Lift{G',s}{G}(\etaABV{\phi'}{})$ when $(G',s,\xi)= (G^0, 1,\xi)$; see Section \ref{section: passage to the identity}.
We now prove the lemmas used in the proof of Proposition~\ref{theorem: FPF}.

\begin{lemma}
\label{lemma: EFPF}
Suppose $G$ satisfies Assumption~\ref{assu: groups we consider} and let $(G',s,\xi)$ be an endoscopic triple for $G$.
Let $\phi'$ be a Langlands parameter for $G'$.
Then
\[
\trace({s},\Ev_{z'}\ d\xi^*\, \mathcal{F})
=
\trace({s},\Ev_{d\xi(z')} \mathcal{F})
\qquad\forall
z'\in \Lambda_{\lambda'}= \Lambda^s_{\lambda},
\]
where the functor $\Ev$ is defined and studied in Appendix~\ref{section: global vanishing cycles}, see especially Definition~\ref{definition: Appendix ev}.
\end{lemma}

\begin{proof}   
Consider Theorem~\ref{theorem: Schurmann} with $M = V_{\lambda}$ and $f=d\xi(y')$.
Since $H_\lambda$ is a linear algebraic group, the Zariski closure of the group generated by $s$ is a diagonalizable group - this group is played by $H$ in the statement of Theorem~\ref{theorem: Schurmann} and was denoted by $S$ in Lemma~\ref{lemma: m}.
Now, this gives
\[ \trace({s},(\Rphi{y'}[-1](d\xi^*\, \mathcal{F})_x
=
\trace({s},(\Rphi{d\xi(y')}[-1] \mathcal{F})_{d\xi(x')},
\qquad
\forall (x',y')\in \Lambda_{\lambda'}.
\]
By Proposition~\ref{theorem: family ev},
\[
\Ev_{(x,y)}\mathcal{F} 
\ceq 
(\ev{C^*_x}\mathcal{F})_{x,y}
=
(\Rphi{y}[-1]\mathcal{F})_x \in D^b_{Z_{H_\lambda}(x,y)}(\{ (x,y)\}),
\qquad
\forall y\in \Lambda_{x},
\]
and
\[
\Ev_{(x',y')}d\xi^* \mathcal{F} 
=
(\Rphi{y'}[-1]d\xi^*\mathcal{F})_{x'} \in D^b_{Z_{H_{\lambda'}}(x',y')}(\{ (x',y')\}),
\qquad
\forall y'\in \Lambda_{x'}.
\]
Lemma~\ref{lemma: EFPF} now follows immediately.
\end{proof}

\begin{remark}
  Note that the target of the functor $\Ev_{C^*}: D^b_{c,H_\lambda}(V_\lambda) \to D^b_{\mathscr{S}-c,H_\lambda}(\Lambda_{C^*})$ is $H_\lambda$-equivariant perverse sheaves on $\Lambda_{C^*} = \Lambda_\lambda \cap (V_\lambda\times C^*)$ constructible  with respect to a stratification $ \mathscr{S}$ defined in Appendix~\ref{section: global vanishing cycles}, which is not a stratification by $H_\lambda$-orbits. 
  This is explained in Theorem \ref{theorem: refinement} from the Appendix. The corresponding stratification $\mathscr{S}$ is constructed in Lemma~\ref{lemma: Whitney}.  
  Moreover, Theorem~\ref{theorem: Schurmann} applies to all $z'\in \Lambda_{C^*}$, not just $\Lambda^\reg_{C^*}$. 
  The next two lemmas, and Proposition~\ref{theorem: FPF}, require $z'\in \Lambda_{\lambda'}$ to be regular and also $z = d\xi(z')\in \Lambda_{\lambda}$ to be regular only because of the definition of $\Evs$ and also because of how we will use it, and further that these elements lie in $\Lambda^\mic$, as used in the definition of $\etaABV{\phi'}{}$ and $\etaABV{\phi}{}$. The larger target of $\Ev$ and larger scope of Theorem~\ref{theorem: Schurmann} will play a role in future work; see also Section~\ref{section: Qmic}.
\end{remark}

Define $e(\phi) \ceq \dim(C_\phi) + \dim(C_\phi^*) - \dim \Lambda_\lambda$; this is called the {\it eccentricity}\IN{$e(\phi)$, eccentricity} of the parameter $\phi$ in \cite{CFMMX22}*{(7.17)}.
Likewise define $e(\phi')$.
    Let $\mathcal{T}\in \Loc_{H_\lambda}(\Lambda_\lambda^\reg)$ be the $H_\lambda$-equivariant local system defined component-by-component by $\mathcal{T}\vert_{\Lambda^\reg_{C}} \ceq \Evs_{C}\IC(\locsys{\1}{C})$; we refer to $\mathcal{T}$ as the {\it twisting sheaf}\IN{$\mathcal{T}$, twisting sheaf} on $\Lambda^\reg_\lambda$.
    Likewise define $\mathcal{T}'\in \Loc_{H^s_{\lambda}}(\Lambda_{\lambda'}^\reg)$ or $\mathcal{T}'\in \Loc_{H_{\lambda'}}(\Lambda_{\lambda'}^\reg)$,
    depending on working with the action of $H^s_{\lambda} $ or $H_{\lambda'} $ on $V_{\lambda}^s=V_{\lambda'}$ (fitting with the corresponding restriction of group functor).

\begin{lemma}\label{lemma: eccentricity}
Let $(G',s,\xi)$ be an endoscopic triple for $G$.
Let $\phi'$ be a Langlands parameter for $G'$ and set $\phi = \xi\circ\phi'$.
Then
\[(-1)^{e(\phi')}=
(-1)^{e(\phi')}\trace({1},\mathcal{T}'_{z'}) 
=
(-1)^{e(\phi')}\trace({s},\mathcal{T}'_{z'}) 
=
(-1)^{e(\phi)}\trace({s},\mathcal{T}_{d\xi(z')}),
\]   
for all $z'\in \Lambda^\reg_{\lambda}$ such that $d\xi(z')\in \Lambda^\reg_{\lambda'}$.
\end{lemma}

\begin{proof}
Write $z'=(x',y')$ and $z=d\xi(z') = (x,y)$. 
Observe that $e(\phi) = \dim C_\phi - \codim C^*_{\phi}$.
Recall $\Evs_{C}$ from Section~\ref{section: ABV-packets} and Definition~\ref{theorem: refinement}: $\Evs_{C} \ceq j^*_{\Lambda^\reg_{C}}\, \Ev_{C}[-\codim C^*]$, so
$\Evs_{z} = \Ev_{z}[-\codim C^*]$
for $z=(x,y)\in \Lambda_\lambda^\reg$.
Then, by Lemma~\ref{lemma: EFPF},
\[
\begin{array}{rcl}
(-1)^{e(\phi)}
\trace(s, \mathcal{T}_{z})
&=&
\trace(s, \mathcal{T}_{z}[-e(\phi)])\\
&=&
\trace(s, \Evs_{z}\IC(\locsys{\1}{C_{\phi}})[-\dim C_\phi+\codim C^*_\phi])\\
&=& \trace(s, \Ev_z\locsys{\1}{C_{\phi}})\\
&=&
 \trace(s, \Ev_{z'}\locsys{\1}{C_{\phi'}})\\
&=&
 \trace(s, \Ev_{z'}\IC(\locsys{\1}{C_{\phi'}})[-\dim C_{\phi'}]) \\
&=&
 \trace(s, \Evs_{z'}\IC(\locsys{\1}{C_{\phi'}})[-\dim C_{\phi'}+\codim C^*_{\phi'}])\\
&=&
 \trace(s, \Evs_{z'}\IC(\locsys{\1}{C_{\phi'}})[-e(\phi')])\\
&=& (-1)^{e(\phi')}
 \trace(s, \Evs_{z'}\IC(\locsys{\1}{C_{\phi'}}))\\
&=& (-1)^{e(\phi')}
 \trace(s, \mathcal{T}'_{z'}).
\end{array} 
\]
Similarly $$\trace(s, \Ev_{z'}\locsys{\1}{C_{\phi'}})=
\trace(1, \Ev_{z'}\locsys{\1}{C_{\phi'}}) =(-1)^{e(\phi')}
 \trace(1, \mathcal{T}'_{z'}) =(-1)^{e(\phi')} \:,
$$
since $s$ is acting trivially on $V_{\lambda}^s=V_{\lambda'}$.
\end{proof}

\begin{lemma}\label{lemma: Evs FPF}
    Let $G$, $G'$, $s$, $\xi$, $\phi'$ and $\phi$ be as above.
    Then, for all $\mathcal{F}\in D^b_{c,H_\lambda}(V_\lambda)$
\[
(-1)^{d(\phi')-e(\phi')} \trace(s,\Evs_{z'}\ d\xi^*\, \mathcal{F})
=
(-1)^{d(\phi)-e(\phi)}
\trace({s},\Evs_{d\xi(z')} \mathcal{F})
\] 
for all $z'\in \Lambda^\reg_{\lambda'}$ such that $d\xi(z')\in \Lambda^\reg_{\lambda}$.
\end{lemma}

\begin{proof}
As above, set $\phi = \xi\circ\phi'$.
Observe that $d(\phi') - e(\phi') = \codim C^*_{\phi'}$.
By Lemma~\ref{lemma: EFPF},
\[
\begin{array}{rcl}
(-1)^{d(\phi')-e(\phi')} \trace(s,\Evs_{z'}\ d\xi^*\, \mathcal{F})
&=&
\trace(s,\Evs_{z'}\ d\xi^*\, \mathcal{F}[\codim C^*_{\phi'}])\\
&=&
\trace(s,\Ev_{z'} d\xi^*\, \mathcal{F})\\
&=&
\trace(s,\Ev_{d\xi(z')} \mathcal{F})\\
&=&
\trace(s,\Evs_{d\xi(z')}\mathcal{F}[\codim C^*_{\phi}])\\
&=&
(-1)^{d(\phi)-e(\phi)}
\trace(s,\Evs_{d\xi(z')}\mathcal{F}).
\end{array}
\]   
\end{proof}

\subsection{Trace and T and twisting}\label{subsection: trace of sigma}

In this section we calculate $\trace(s,\NEvs_{z'}(T\circ\calm^*)\mathcal{G})$ for all $z'\in \Lambda^\reg_{\lambda'}$ and for all $\mathcal{G}\in D^b_{c,H_\lambda^s}(V_\lambda^s)$, for use in Section~\ref{subsection: lift of ABV}, in particular.

\begin{proposition}\label{proposition: T}
For all $\mathcal{F'}\in D^b_{c,H_\lambda^s}(V_\lambda^s)$ and $z'\in \Lambda_{\lambda'}^\reg$,
\[
\begin{array}{rcl}
\trace(s,\NEvs_{z'}  \mathcal{F'})
&=& 
\trace(1,\NEvs_{z'}(T\circ\calm^*)\mathcal{F'}).
\end{array}
\]
\end{proposition}

\begin{proof}
Recall that  $\NEvs_{z'}\mathcal{F'} = \mathcal{T'}_{z'}^{\vee} \otimes \Ev_{z'} \mathcal{F'}[-\codim C_{x'}^*]$
where 
$\mathcal{T'}_{z'} = \Ev_{z'}\IC(\1_{C_x'})$.
Note that on the left (resp. right)  hand side 
we are considering $\NEvs_{z'}$ for $H_\lambda^s $-
(resp $H_{\lambda'}$)-equivariant coefficients.
Then, using Lemma~\ref{lemma: eccentricity} and
Lemma~\ref{lemma: T}, we have
\[
\begin{array}{rcl} (-1)^{\codim C_{x'}^*}\cdot 
\trace(1,\NEvs_{z'}(T\circ\calm^*)\mathcal{F'})
&=&
\trace(1,\mathcal{T'}_{z'}^{\vee} \otimes \Ev_{z'}(T\circ\calm^*)\mathcal{F'})\\
&=&
\trace(1,\mathcal{T'}_{z'}^{\vee}) \trace(1,\Ev_{z'}(T\circ\calm^*)\mathcal{F'})\\
&=&
\trace(1,\mathcal{T'}_{z'}^{\vee}) \trace(s,\Ev_{z'}\mathcal{F'})\\
&=&
\trace(s,\mathcal{T'}_{z'}^{\vee}) \trace(s,\Ev_{z'}\mathcal{F'})\\
&=&
\trace(s,\mathcal{T'}_{z'}^{\vee}\otimes \Ev_{z'}\mathcal{F'})\\
&=&
(-1)^{\codim C_{x'}^*}\cdot \trace(s,\NEvs_{z'}\mathcal{F'}).
\end{array}
\]
\end{proof}

\begin{lemma}\label{lemma: T}
For all $\mathcal{F'}\in D^b_{c,H_\lambda^s}(V_\lambda^s)$ and $z'\in \Lambda_\lambda^s$,
\[
\begin{array}{rcl}
\trace(s,\Ev_{z'}  \mathcal{F'})
&=& 
\trace(1,\Ev_{z'}(T\circ\calm^*)\mathcal{F'}).
\end{array}
\]
\end{lemma}

\begin{proof}
Write $z'=(x',y')\in \Lambda_{C_{\phi'}}$.
Lemma~\ref{lemma: m} gives the decomposition
\[
\calm^*\mathcal{F'} = \oplus_{\chi \in \Irrep(\pi_0(S))} \mathcal{F}_\chi''\otimes\locsys{\chi}{V_{\lambda'}}
\]
for $\mathcal{F'}\in D^b_{c,H_{\lambda'}}(V_{\lambda'})$ and $\chi \in \Irrep(\pi_0(S))$.
Then
\[
\begin{array}{rcl}
\trace(1,\Ev_{z'}(T\circ \calm^*)\mathcal{F'})
&=&\oplus_{\chi \in \Irrep(\pi_0(S))} 
\trace(1,\Ev_{z'}T(\mathcal{F''}\otimes\locsys{\chi}{V_{\lambda'}})) \\
&=&\oplus_{\chi \in \Irrep(\pi_0(S))} 
\trace(1,\Ev_{z'}T(\mathcal{F''}\otimes\locsys{\chi}{V_{\lambda'}})) \\
&=& \oplus_{\chi \in \Irrep(\pi_0(S))} 
\trace(1,(\Rphi{y'}[-1]T(\mathcal{F''}\otimes\locsys{\chi}{V_{\lambda'}}))_{x'}) \\
&=&\oplus_{\chi \in \Irrep(\pi_0(S))} 
\trace(1,(\Rphi{y'}[-1]\mathcal{F''})_{x'}) \trace(s,\chi)\\
&=&\oplus_{\chi \in \Irrep(\pi_0(S))} 
\trace((1,s), (\Rphi{y'}[-1]\mathcal{F''})_{x'} \otimes \locsys{\chi}{{x'}}) \\
&=&
\trace((1,s),\calm^*(\Rphi{y'}[-1]\mathcal{F'})_{x'}) \\
&=&
\trace(\calm(1,s),(\Rphi{y'}[-1]\mathcal{F'})_{x'}) \\
&=&
\trace(s,(\Rphi{y'}[-1]\mathcal{F'})_{x'})\\
&=&
\trace(s,\Ev_{z'}\mathcal{F'}).
\end{array}
\]
\end{proof}

\subsection{Lift of ABV-distributions under certain conditions}\label{subsection: lift of ABV}

In this section we calculate lift of ABV-distributions under the following condition.

\begin{definition}\label{definition: regular parameter}
Let $(G',s,\xi)$ be an endoscopic triple for $G$. 
Let $\phi'$ be a Langlands parameter for $G'$. 
We say that $\phi'$ is  {\it a regular parameter for $G$}\IN{regular parameter for connected $G$} if  $d\xi(\Lambda^\reg_{C_{\phi'}}) \cap \Lambda^\reg_{C_{\xi\circ\phi'}}\ne \emptyset$ or, equivalently, if there is some {\it regular} covector $y'\in \Lambda_{x_{\phi'}}$ such that its image $d\xi(y')\in \Lambda_{x_{\xi\circ\phi'}}$ is also regular.
    \end{definition}

Note that this definition can only be made after an endoscopic triple $(G',s,\xi)$ for $G$ has been specified.

\begin{remark}\label{remark: irregular}
In Lemma~\ref{lemma: Arthur is regular} we show that if $\phi'$ is of Arthur type then $\phi'$ is regular for $G$.
However, there are many examples of parameters $\phi'$ that are regular for $G$ but not of Arthur type; one such example is given in \cite{CFK22}. 
Likewise, parameters $\phi'$ that are not regular for $G$ are not exotic.
For example, take $G= \GL_2$ and let $G'=T$ be a maximal split torus and let $s$ be any regular element of $\dualgroup{T}$; then $\xi : \Lgroup{G'} \to \Lgroup{G}$ is induced from $\dualgroup{T} \hookrightarrow \dualgroup{G}$.
For the standard Borel $B$, let $\phi' = \lambda'$ be the Langlands parameter for the character $\delta_B^{1/2}$ of $T(F)$.
Then $\phi = \xi\circ \phi'$ is the Langlands parameter for the trivial representation.
In this case, $\phi'$ is not regular for $G$: although $x_{\phi'} = 0$ is regular in the conormal for $\phi'$, $T$, its image $x_\phi = d\xi(x_{\phi'}) = 0$ is not regular in the conormal for $G$.
\end{remark}

\begin{remark}
    For any Langlands parameter $\phi$ for $G$, there is always some endoscopic $(G',s,\xi)$ and Langlands parameter $\phi'$ for $G'$ that is regular for $G$ and $\phi = \xi \circ \phi'$.
    The minimal such $G'$ is relevant to the endoscopic classification of representations in $\Pi^\ABV_\phi(G)$.
    In the example of Remark~\ref{remark: irregular}, if $\phi$ is for the trivial representation, then $G$ itself is the minimal such endoscopic group.
\end{remark}

\begin{theorem}\label{theorem: Lift of ABV}
Let $G$ be be a quasisplit group satisfying Assumption~\ref{assu: groups we consider} and let $\delta$ be a pure inner form of $G$.
Let $(G',s,\xi)$ be an endoscopic triple for $G$.
Let $\phi'$ be a Langlands parameter for $G'$ and set
$\phi = \xi\circ\phi'$.
Suppose the $\phi'$ is regular for $G$. 
Then
    \[
    \Lift{G',s}{G_\delta}
    (\etaABV{\phi'}{})
    =
    e(G_\delta)\,
    \etaABV{\phi,s}{,\delta}.   
    \]
\end{theorem}

\begin{proof}
Let $\lambda'$ be the infinitesimal parameter for $\phi'$ and let $\lambda$ be the infinitesimal parameter for $\phi$.
    We will prove
\begin{equation}\label{equation: Lift of ABV, 0}
    \langle \Lift{G',s}{G_\delta}( \etaABV{\phi'}{}), \mathcal{F} \rangle
    =
    e(G_\delta)\,
    \langle
    \etaABV{\phi,s}{,\delta} , \mathcal{F} \rangle,
    \qquad\forall \mathcal{F}\in K_\CC\Per_{H_{\lambda}}(V_{\lambda}).  
\end{equation}
By the non-degenerateness of the pairing \eqref{equation: pairing}, this will imply the claimed result.

First, consider the left-hand side of \eqref{equation: Lift of ABV, 0}.
By Definition~\ref{definition: Lift} and Lemma~\ref{proposition: etaABV,F}
\begin{equation*}
    \begin{array}{rcl}
\langle
    \Lift{G',s}{G_\delta}( \etaABV{\phi'}{}) , \mathcal{F} 
\rangle
    &=& 
    \langle
    \QST^\geo_{G'}(\etaABV{\phi'}{}) , (T\circ \calm^*\circ d\xi^*) 
    {}^t\operatorname{proj}_{\delta} \mathcal{F} 
\rangle
\vspace{0.5em} \\
    &=&
    \langle
    \etaABV{\phi'}{,\pure} , (T\circ \calm^*\circ d\xi^*)
{}^t\operatorname{proj}_{\delta} \mathcal{F} 
\rangle
\vspace{0.5em} \\
    &=&
 (-1)^{d(\phi')}
\trace(1,\NEvs_{z'} (T\circ\calm^*\circ d\xi^* \circ
{}^t\operatorname{proj}_{\delta} \mathcal{F}))
\vspace{0.5em} 
    \end{array}
\end{equation*}
for $z'\in \Lambda^\mic_{C_{\phi'}}$.
Now, consider the right-hand side of \eqref{equation: Lift of ABV, 0}.
By Lemma~\ref{proposition: etaABV,F},  
\begin{equation}\label{equation: Lift of ABV, 1}
    \begin{array}{rcl}
\langle
    \etaABV{\phi,s}{,\delta} , \mathcal{F} 
\rangle
    &=&
 e(G_\delta)\,
(-1)^{d(\phi)}
\trace(s,\NEvs_{z}  
{}^t\operatorname{proj}_{\delta} \mathcal{F})
    \end{array}
\end{equation}
for $z\in \Lambda^\mic_{C_\phi}$,
so
\begin{equation}\label{equation: Lift of ABV, 2}
    \begin{array}{rcl}
e(G_\delta)\,
\langle
    \etaABV{\phi,s}{,\delta} , \mathcal{F} 
\rangle
    &=&
(-1)^{d(\phi)}
\trace(s,\NEvs_{z}  
{}^t\operatorname{proj}_{\delta} \mathcal{F}).
    \end{array}
\end{equation}
By Proposition~\ref{theorem: FPF}, we have
\[
(-1)^{d(\phi)}
\trace(s,\NEvs_{z} {}^t\operatorname{proj}_{\delta}\mathcal{F})
=
(-1)^{d(\phi')} \trace(s,\NEvs_{z} d\xi^* \circ
{}^t\operatorname{proj}_{\delta} \mathcal{F}).
\]
By Proposition~\ref{proposition: T},
\[
\trace(s,\NEvs_{z'} d\xi^*\circ {}^t\operatorname{proj}_{\delta}\mathcal{F})
=
\trace(1,\NEvs_{z'} (T\circ \calm^*\circ d\xi^*){}^t\operatorname{proj}_{\delta}
\mathcal{F}).
\]
Thus
\[
\begin{array}{rcl}
e(G_\delta)\,
\langle
    \etaABV{\phi,s}{,\delta} ,  \mathcal{F} 
\rangle
    &=&
    (-1)^{d(\phi')}
   \trace(1,\NEvs_{z'} (T\circ \calm^*\circ d\xi^*)
  {}^t\operatorname{proj}_{\delta} \mathcal{F}).
\end{array}
\]
We conclude that 
\[ \Lift{G',s}{G_\delta}( \etaABV{\phi'}{}) , \mathcal{F} 
\rangle= 
\langle
   e(G_\delta) \etaABV{\phi,s}{,\delta} ,  \mathcal{F} 
\rangle,\]
and therefore completes the proof.
\end{proof}

\begin{theorem}\label{theorem: transfer of ABV-distributions}
Let $(G',s,\xi)$ be an endoscopic triple for $G$ as in defined in Section~\ref{subsection: endoscopic triple}.
Suppose the local Langlands correspondence for $G$ satisfies the assumptions of Theorem \ref{thm: stable standards form a basis - quasisplit}. 
Let $\phi'$ be a Langlands parameter for $G'$ and set $\phi = \xi\circ\phi'$.
Suppose $\phi'$ is regular for $G$ as in Definition~\ref{definition: regular parameter}.
Then, for every pure inner inner form $\delta$ for $G$,
\[
\Transfer{G',s}{G_\delta}
(\etaABV{\phi'}{})
=
e(G_\delta)\ 
\etaABV{\phi,s}{,\delta}.
\]
\end{theorem}

\begin{proof}
    This is a consequence of Theorems \ref{theorem: lift = transfer} and \ref{theorem: Lift of ABV}. 
\end{proof}

\subsection{Lift from a disconnected group}
In this subsection, we consider the triple of groups $(G,G',G'')$ satisfies the following assumptions.
\begin{assu}\label{assu G G' G''}
Let $G$ and $G''$ be two groups satisfying Definition~\ref{assu: groups we consider}; we do not require $|\pi_0(G)|=|\pi_0(G'')|$. 
Further, suppose there exists an endoscopic triple $(G',s, \xi)$ of $G$ such that $G'= (G'')^{0}$, and the following conditions hold.
\begin{enumerate}
    \item [(i)] There exists a local Langlands correspondence (Definition \ref{definition: LLC}) for the groups $G$ and $G''$ such that Kazhdan-Lusztig hypothesis (Hypothesis \ref{hypothesis: KLH}) holds.
    \item [(ii)] If two Langlands parameters $\phi_1', \phi_2'$ of $G'$ are conjugate under $\dualgroup{G''}$, then $\xi \circ \phi_1'$ and $\xi\circ \phi_2'$ are conjugate under $\dualgroup{G}$.
\end{enumerate}
\end{assu}
\begin{remark}
    It is crucial that we do not assume the local Langlands correspondence for $G'$. Also, we do not assume the local Langlands correspondences for $G$ and $G''$ satisfy Desiderata \ref{desideratum: LLC}.
\end{remark}

Here are two important examples of $(G,G',G'')$ satisfying Assumption~\ref{assu G G' G''}.
\begin{itemize}
 \item $(\Sp_{2n}, \Sp_{2n_1}\times \SO_{2n_2}^{\alpha}, \Sp_{2n_1}\times \O_{2n_2}^{\alpha})$; see Section \ref{subsubsection: C} for details.
    \item $(\GL_{2n}\rtimes \langle \theta \rangle, \SO_{2n}^{\alpha}, \O_{2n}^{\alpha})$; see Section \ref{subsubsection: D} for details.
\end{itemize}
These triples come from Arthur's theory. The endoscopic classification only gives us a local Langlands correspondence for the even orthogonal group $ \O_{2n}^{\alpha}$. For the identity component $\SO_{2n}^{\alpha}$, we only have a ``weak" local Langlands correspondence; see \cite{AG17}. Therefore, we do not have an isomorphism $ K_{\BC}\Pi_{\lambda'}(G')^\pure\simeq K_{\BC} \Per_{H_{\lambda'}}(V_{\lambda'})^{\ast}$. In fact, $K_{\BC}\Pi_{\lambda'}(G')^\pure$ is not defined. The goal of this subsection is to define a geometric analogue of the composition
\[ \Transfer{G',s}{G} \circ \Res^{G''}_{G'},\]
and study its properties.
We will show in Sections \ref{sec step 2}, \ref{sec step 3} that the analogue of Theorem \ref{theorem: lift = transfer} holds for the two triples above.
\begin{remark}
    The triple $(\O_{2n}, \SO_{2n_1}^{\alpha_1}\times \SO_{2n_2}^{\alpha_2}, \O_{2n_1}^{\alpha_1}\times \O_{2n_2}^{\alpha_2} )$ does not satisfy Definition~\ref{assu G G' G''} since $\pi_0( \O_{2n_1}^{\alpha_1}\times \O_{2n_2}^{\alpha_2})\simeq (\Z/2\Z)^2$. However, the results in this section also hold for this triple with minor adaptations to the same argument.
\end{remark}

In the rest of this subsection we write $G''= G' \rtimes \langle \theta \rangle$ as in Definition~\ref{assu: groups we consider} and let $\ell:= |\pi_0(G'')|$.
Before we give the explicit definition, let us explain the motivation. Let $\phi''$ be a Langlands parameter of $G''$ and set $\ell_{\phi''}:=1$ if ${}^{\hat{\theta}}\phi''$ is conjugate to $\phi''$ under $\dualgroup{G'}$, and set $\ell_{\phi''}:=\ell$ otherwise. Regard $\phi''$ as a Langlands parameter of $G'$ and define a set of Langlands parameters of $G^0$
\begin{align}\label{eq def Omega phi''}
    \Omega_{\phi''}=\{{}^{\hat{\theta}^k}\phi'' \tq 1\leq k \leq \ell_{\phi''} \},
\end{align}
which is a set of representatives of the set of $L$-parameters $\phi'$ of $G'$ that are conjugate to $\phi''$ under $\dualgroup{G''}.$ For each Langlands parameter $\phi' \in \Omega_{\phi''}$,
$A_{\phi'}$ is a subgroup of $A_{\phi''}$ of index $1$ or $\ell$, by Corollary \ref{cor: component group  l pardisc}. For $\tau'' \in \widehat{A}_{\phi''}$, we  decompose
\[ \tau''|_{A_{\phi'}}= \bigoplus_{\tau' \in A_{\phi'}} m(\tau',\tau'') \tau', \]
which is multiplicity free.

Assume the local Langlands correspondence for $G'$ for the moment. Then it is expected that (see Desideratum \ref{desideratum: LLC} \RES)
\[ \Res_{G'}^{G''} (\SR(\pi(\phi'', \tau'')))= \bigoplus_{\phi' \in \Omega_{\phi''} } \bigoplus_{\tau' \in  \widehat{A}_{\phi'}} m(\tau',\tau'') \SR(\pi(\phi', \tau')). \]
We would like to define a linear transformation on the Grothendieck groups of perverse sheaves directly without the local Langlands correspondence for $G'$ that is transpose to the above expected description $\Res_{G'}^{G''}$.

Fix an infinitesimal parameter $\lambda''$ of $G''$, which is a homomorphism $W_F \to \Lgroup{G''}$ with image in $\Lgroup{G'}$, and let $\lambda= \xi \circ \lambda''$ be the corresponding infinitesimal parameter of $G$.
Define $K_{\BC}\Pi_{\lambda''}(G'')^{\std}$ (resp. $K_{\BC}\Pi_{\lambda''}(G'')^{\pure,\std}$) to be the vector space spanned by $\{\mu_{\phi''} \tq \phi'' \in \Phi_{\lambda''}(G'')\}$ (resp. $\{\mu_{\phi''}^{\pure}\tq \phi'' \in \Phi_{\lambda''}(G'')\}$). Let $\Omega_{\lambda''}$ be defined as in \eqref{eq def Omega phi''}. 
Define the linear transformation 
\[ R: K_{\BC} \Per_{H_{\lambda''}}(V_{\lambda''})^{\ast} \to     \bigoplus_{\lambda' \in \Omega_{\lambda''}} K_{\BC} \Per_{H_{\lambda'}}(V_{\lambda'})^{\ast} \]
as follows. Recall that $\{\stsh{\tau''}{C_{\phi''}}\tq (\phi'',\tau'')\in \Phi_{\lambda''}^{\mfr{e}}(G'')\}$ forms a basis of $K_{\BC} \Per_{H_{\lambda''}}(V_{\lambda''})$. We let $\{ (\stsh{\tau''}{C_{\phi''}})^{\ast} \tq (\phi'',\tau'')\in \Phi_{\lambda''}^{\mfr{e}}(G'')\}$ be the dual basis for $K_{\BC} \Per_{H_{\lambda''}}(V_{\lambda''})^{\ast}$ and similarly for $K_{\BC} \Per_{H_{\lambda'}}(V_{\lambda'})^{\ast}$.
Note that Kazhdan-Lusztig hypothesis for $G''$ implies that
    \[  (\stsh{\tau''}{C_{\phi''}})^{\ast}= e(\tau'') \langle \SR(\pi(\phi'',\tau'')) ,-  \rangle_{\lambda''}.   \]
    Now we define $R$ on this basis by
\begin{align}\label{eq def of geometrization of finite index res}
    R( (\stsh{\tau''}{C_{\phi''}})^{\ast}):=  \sum_{\phi' \in \Omega_{\phi''}} \sum_{\tau' \in A_{\phi'}}  m(\tau',\tau'') (\stsh{\tau'}{C_{\phi'}})^{\ast}.
\end{align}

With the definition of $R$, we give the definition of the linear transformation $\Lift{G'',s}{G_{\delta}}$. 

\begin{definition}
Suppose that $(G, G', G'')$ satisfies Definition~\ref{assu G G' G''} and let $\delta$ be a pure inner form of $G$.  We define the linear transformation $\Lift{G'',s}{G_{\delta}}$ via the following diagram.\IN{$\Lift{G'',s}{G_{\delta}}$, lift from disconnected group}
  
\begin{equation*}
    \begin{tikzcd}
        K_{\BC}\Pi_{\lambda''}(G'')^{\std} \ar[rrd, "\Lift{G'', s}{G_{\delta}}"'] \ar[rr,"\frac{1}{\ell} \QST^\geo"] && K_{\BC}\Pi_{\lambda''}(G'')^\pure \ar[r,"R"] & 
      \bigoplus_{\lambda' \in \Omega_{\lambda''}} K_{\BC} \Per_{H_{\lambda'}}(V_{\lambda'})^{\ast}
        \ar[d," {}^t (T\, \circ\, m^{\ast} \, \circ\, d\xi^{\ast})"]\\
        &&K_{\BC}\Pi_{\lambda}(G_{\delta}) &\ar[l, "\mathrm{proj}_{\delta}"']  K_{\BC}\Pi_{\lambda}(G)^\pure
    \end{tikzcd}
\end{equation*}
where:
\begin{enumerate}
    \item the linear transformation $\QST^\geo$ is defined by sending $\mu_{\phi}$ to $\mu_{\phi}^{\pure}$ as in Definition \ref{definition: QSTgeo};
    \item the linear transformation $R$ is defined as in \eqref{eq def of geometrization of finite index res}, where we identify $K_{\BC}\Pi_{\lambda''}(G'')^\pure$ as $K_{\BC} \Per_{H_{\lambda''}}(V_{\lambda''})^{\ast} $ via the perfect pairing in Definition \ref{definition: pairing} from the local Langlands correspondence of $G''$;
    \item the linear transformation  ${}^t (T\circ m^{\ast} \circ  d\xi^{\ast})$ is given in Definition \ref{definition: Lift}, applied to each summand $K_{\BC} \Per_{H_{\lambda'}}(V_{\lambda'})^{\ast}$ of $    \bigoplus_{\lambda' \in \Omega_{\lambda''}} K_{\BC} \Per_{H_{\lambda'}}(V_{\lambda'})^{\ast}$.
\end{enumerate}
\end{definition}

The first goal of this section is to describe the linear transformation $\Lift{G'',s}{G}$ on the basis $\{\mu_{\phi''} \tq \phi'' \in \Phi_{\lambda''}(G'')\}$ of $ K_{\BC}\Pi_{\lambda''}(G'')^{\std}$.
\begin{proposition}\label{prop lift=trans disconnected}
    Suppose that $(G, G', G'')$ satisfies Definition~\ref{assu G G' G''}. Let $\phi'' \in \Phi_{\lambda''}(G'')$ and suppose $\phi''$ is microlocally generic for $G$ as in Definition~\ref{definition: regular parameter}.
    Then 
    \[ 
    \Lift{G'',s}{G_{\delta}}( \mu_{\phi''} ) = 
    e(G_\delta) \mu_{\xi \circ \phi'',s}^{\delta}
    \]
    for all pure inner forms $\delta$ for $G$.
\end{proposition}

\begin{proof}
    We evaluate the pairing
    \begin{align}\label{eq trans=lift disconnected 1}
        \langle \Lift{G'',s}{G_{\delta}} \mu_{\phi''}, \stsh{\tau}{C_{\phi}}\rangle= \frac{1}{\ell}\langle R \circ   \QST^\geo (\mu_{\phi''}),   (T \circ m^{\ast} \circ d\xi^{\ast})  \circ {}^t\textrm{proj}_{\delta} (\stsh{\tau}{C_{\phi}})\rangle.
    \end{align}
    By definition, we have
    \begin{align*}
        R \circ   \QST^\geo (\mu_{\phi''})&= R \left( \sum_{\tau'' \in \widehat{A}_{\phi''}} e(\tau'') \trace(\tau'',1) \SR(\pi(\phi'', \tau''))  \right)\\
        &= R \left( \sum_{\tau'' \in \widehat{A}_{\phi''}} \dim(\tau'') (\stsh{\tau''}{C_{\phi''}})^{\ast}  \right)\\
        &= \sum_{\phi'\in \Omega_{\phi''}}  \sum_{\tau'' \in \widehat{A}_{\phi''}} \sum_{\tau' \in \widehat{A}_{\phi'}} m(\tau',\tau'')  \dim(\tau'') (\stsh{\tau''}{C_{\phi''}})^{\ast},
    \end{align*}
    where we use Kazhdan-Lusztig hypothesis for $G''$ for the second equality.

For each $\lambda' \in \Omega_{\lambda''}$, the $K_{\BC} \Per_{H_{\lambda'}}(V_{\lambda'})^{\ast}$ component in $    \bigoplus_{\lambda' \in \Omega_{\lambda''}} K_{\BC} \Per_{H_{\lambda'}}(V_{\lambda'})^{\ast}$ of  $(T \circ m^{\ast} \circ d\xi^{\ast})  \circ {}^t\textrm{proj}_{\delta} (\stsh{\tau}{C_{\phi}})$ has already been computed in the proof of Theorem \ref{theorem: lift = transfer}. It is zero unless $e(\tau)=\delta$, which we assume from now on. Then \eqref{equation: Lift=Trans 3} gives that 
    \[
    \begin{array}{rcl}
    && (T \circ  \calm^{\ast} \circ d\xi^{\ast}) ( \stsh{\tau}{C_{\phi}} )  \\
    & = &  
\mathop{\sum}\limits_{\lambda' \in \Omega_{\lambda''}} \left(\mathop{\sum}\limits_{ i=1 }^r 
    \hskip5pt \mathop{\sum}\limits_{ j=1 }^{r_i} 
    \hskip5pt 
\mathop{\sum}\limits_{\rho_i \in \widehat{A}_{C_i}} 
    \hskip5pt
\mathop{\sum}\limits_{\rho_{ij}' \in \widehat{A}_{C_{ij}'}}  
    \hskip5pt
   \rho_{i,B}({\sigma}) m(\tau, \rho_i) m(\rho_i, \rho_{ij}' \otimes \rho_{i,B}) \stsh{\rho_{ij}'}{C_{ij}'}\right).
   \end{array}
   \]
Here, the orbits $C_{ij}'$ in the summation are the $H_{\lambda'}$-orbit of $V_{\lambda'} \cap C_{\phi}$, which are in bijection with $\Omega_{\phi} \cap \Phi_{\lambda'}(G')$, where
\[\Omega_{\phi}:= \{ \phi' \in \Phi(G')\ | \ \xi \circ \phi'= \phi  \}. \]
Thus, we may rewrite 
\[(T \circ  \calm^{\ast} \circ d\xi^{\ast}) ( \stsh{\tau}{C_{\phi}} )
    =  
    \sum_{\phi' \in \Omega_{\phi} } 
    \hskip5pt 
     \sum_{\rho' \in \widehat{A}_{C_{\phi'}}}  
    \hskip5pt
   c_{\rho' }  \cdot \stsh{\rho'}{C_{\phi'}},\]
   where if $C_{\phi'}= C_{ij}'$ and $\rho'=\rho_{ij}'$, then (note that $C_{ij}'$ determines $C_i$ by \eqref{eq lift=trans decomp of orbit}) 
   \[ c_{\rho' }=\sum_{\rho_i \in \widehat{A}_{C_i}}  
   \rho_{i,B}({\sigma}) m(\tau, \rho_i) m(\rho_i, \rho_{ij}' \otimes \rho_{i,B}). \]
   Recall from \eqref{equation: Lift of ABV, 4} that for each $\phi' \in \Omega_{\phi},$
   \begin{align}\label{eq trans=lift disconnected 2}
       \trace(s,\tau)= \sum_{\rho' \in \widehat{A}_{\phi'}} \dim(\rho') c_{\rho'}.
   \end{align}

Now we see that the pairing \eqref{eq trans=lift disconnected 1} is zero unless $\Omega_{\phi''} \cap \Omega_{\phi}$ is non-empty. This is equivalent to $\Omega_{\phi''} \subseteq  \Omega_{\phi}$ by Definition~\ref{assu G G' G''}(ii), and is equivalent to $\xi \circ \phi''=\phi$, which we also assume now. In this situation, \eqref{eq trans=lift disconnected 1} becomes
   \begin{align}\label{eq trans=lift disconnected 3}
        \frac{1}{\ell} \sum_{\phi'\in \Omega_{\phi''}}  \sum_{\tau'' \in \widehat{A}_{\phi''}} \sum_{\tau' \in \widehat{A}_{\phi'}} m(\tau',\tau'') \dim(\tau'')  c_{\tau'}.
   \end{align}
   We claim that \eqref{eq trans=lift disconnected 3} is equal to $\trace(s,\tau)$. If the claim holds, then we see that 
   \[ \langle \Lift{G'',s}{G_{\delta}} (\mu_{\phi''}), \stsh{\tau}{C_{\phi}}\rangle= \begin{cases}
       \trace(s,\tau)  & \text{ if }\tau \in \widehat{A}_{\phi}[\delta] \text{ and }\ \xi \circ \phi''=\phi,\\
        0 &\text{ otherwise.}
   \end{cases}
\]
Thus, by Kazhdan-Lusztig hypothesis for $G$, we conclude that
   \[
   \Lift{G'',s}{G_{\delta}}
   (\mu_{\phi''}) = \sum_{\tau \in \widehat{A}_{\phi} [\delta]} e(\tau) \trace(s,\tau) \SR(\pi(\xi\circ \phi'', \tau))= e(G_{\delta}) \mu_{\xi \circ \phi''}^{\delta}. \]
This gives the desired conclusion.

  Finally, we verify the claim for $|\Omega_{\phi''}|=\ell$ or $1$ separately. First, suppose that $|\Omega_{\phi''}|=\ell$. Then we have $A_{\phi'}\simeq A_{\phi''}$ for all $\phi' \in \Omega_{\phi''}$ by Corollary \ref{cor: component group  l pardisc}. Thus, by \eqref{eq trans=lift disconnected 2},  \eqref{eq trans=lift disconnected 3} becomes
   \[ \frac{1}{\ell}\sum_{\phi'\in \Omega_{\phi''}} \sum_{\tau' \in \widehat{A}_{\phi'}} \dim(\tau')  c_{\tau'}= \frac{1}{\ell}\sum_{\phi'\in \Omega_{\phi''}} \trace(s,\tau)= \trace(s,\tau). \]
   Next, suppose that $|\Omega_{\phi''}|=1$. By \eqref{eq trans=lift disconnected 2}, it suffices to show that 
   \begin{align*}
       \sum_{\tau'' \in \widehat{A}_{\phi''}}  m(\tau',\tau'') \dim(\tau'') = \ell \dim(\tau')
   \end{align*}
    for any $\tau' \in \widehat{A}_{\phi'}$. Let $\langle \cdot, \cdot \rangle_{A_{\phi''}}$ and $\langle \cdot, \cdot \rangle_{A_{\phi'}}$ be the inner product for characters on these finite groups. Then by Frobenius reciprocity, we have
    \begin{align*}
        \sum_{\tau'' \in \widehat{A}_{\phi''}}  m(\tau',\tau'') \dim(\tau'') 
        &= \sum_{\tau'' \in \widehat{A}_{\phi''}}  \langle \Res_{A_{\phi'}}^{A_{\phi''}}\tau'', \tau' \rangle_{A_{\phi''}} \dim(\tau'') \\
        & = \sum_{\tau'' \in \widehat{A}_{\phi''}}  \langle \tau'', \Ind_{A_{\phi'}}^{A_{\phi''}}\tau' \rangle_{A_{\phi'}} \dim(\tau'').
    \end{align*}
    The last term is exactly $\dim(\Ind_{A_{\phi'}}^{A_{\phi''}}\tau')= \ell \dim(\tau')$.
    This completes the proof of the claim and the proposition.
\end{proof}

The second goal of this section is to show the analogue of Theorem \ref{theorem: Lift of ABV} for $\Lift{G'',s}{G_{\delta}}$. 
First, we compute $R(\etaABV{\phi''}{, \pure})$ explicitly. For an $L$-parameter $\phi'$ of $G'$ with infinitesimal parameter $\lambda'$, recall from Section \ref{subsection: stable ABV-distributions} that $\chi_{C_{\phi'}}^{\mic}$ is the linear functional in $K_\BC \Per_{H_{\lambda'}}(V_{\lambda'})^{\ast}$ such that for all $\mathcal{F} \in K_\BC \Per_{H_{\lambda'}}(V_{\lambda'})$, 
\[ \chi_{C_{\phi'}}^{\mic}(\mathcal{F})= (-1)^{d(\phi')} \trace(1, \NEvs_{z'} \mathcal{F}), \]
for $z'\in \Lambda^\mic_{C_{\phi'}}$. Also, $\chi_{C_{\phi'}}^{\textrm{loc}}$ is the unique functional in $K_\BC \Per_{H_{\lambda'}}(V_{\lambda'})^{\ast}$ such that 
\[ \chi_{C_{\phi'}}^{\textrm{loc}}( \stsh{\tau'}{C_{\widetilde{\phi'}}} ) =\begin{cases}
   \rank(\tau') & \text{ if }\widetilde{\phi'}=\phi',\\
    0 &\text{ otherwise.}
\end{cases}\]
By Kashiwara-Dubson Index Formula, which is recalled in \eqref{equation: stalks 7}, we may write
\begin{align}\label{eq: Kashiwara-Dubson G'}
    \chi^\textrm{mic}_{C_{\phi'}} = \sum_{i} c(\phi_i',\phi')\, \chi^\textrm{loc}_{C_{\phi_i'}},
\end{align}
where the sum is taken over $\phi_i' \in \Phi_{\lambda'}(G')$ and $ c(\phi_i',\phi')$ are some constants.

\begin{lemma}\label{lem disconnected R ABV}
Let $\phi'' \in \Phi_{\lambda''}(G'')$.    With the above notation, we have 
    \[ R(\etaABV{\phi''}{, \pure})= \frac{\ell}{|\Omega_{\phi''}|} \sum_{\phi' \in \Omega_{\phi''}} \chi_{C_{\phi'}}^{\mic}.\]
\end{lemma}
\begin{proof}
For $\lambda' \in \Omega_{\lambda''}$, write $d\xi_{\lambda'}$ for the embedding from $V_{\lambda'}$ to $V_{\lambda''}$.    First, we compute ${}^td\xi_{\lambda'}^{\ast}(\chi_{C_{\phi'}}^{\textrm{loc}}) $ for any $\phi' \in \Phi_{\lambda'}(G')$. By definition and \cite{ABV92}*{Proposition 23.7(b)}, for arbitrary $\widetilde{\phi}''$ and $\widetilde{\tau}'' \in \widehat{A}_{\widetilde{\phi}''}$,
we have
\begin{align*}
{}^td\xi_{\lambda'}^{\ast}\chi_{C_{\phi'}}^{\textrm{loc}}(
  \stsh{\widetilde{\tau}''}{C_{\widetilde{\phi}''}})&= \chi_{C_{\phi'}}^{\textrm{loc}}(d\xi_{\lambda'}^{\ast}\stsh{\widetilde{\tau}''}{C_{\widetilde{\phi}''}})\\
    &= \chi_{C_{\phi'}}^{\textrm{loc}}\left( \sum_{\widetilde{\phi'}\in \Omega_{\widetilde{\phi}''} \cap \Phi_{\lambda'}(G')} \sum_{\tau' \in \widehat{A}_{\widetilde{\phi'}}}m(\tau', \widetilde{\tau}'')\stsh{\tau'}{C_{\widetilde{\phi'}}} \right)\\
    &= \begin{cases}
          \rank(\widetilde{\tau}'') & \text{ if }\phi' \in \Omega_{\widetilde{\phi}''}\\
        0& \text{ otherwise.}
    \end{cases}
\end{align*}
Therefore, we obtain
\[{}^td\xi_{\lambda'}^{\ast}( \chi_{C_{\phi'}}^{\textrm{loc}}) = \chi_{C_{\xi_{\lambda'} \circ \phi'}}^{\textrm{loc}}= \mu_{\xi_{\lambda'} \circ \phi'}^{\pure},\]
where we use Kazhdan-Lusztig hypothesis for $G''$ for the second equality; see Lemma \ref{lemma: mu and chiloc}.

Now we fix a $\phi' \in \Omega_{\phi''}$. By Proposition \ref{theorem: FPF}, for any $\mathcal{F} \in K_{\BC}\Pi_{\lambda''}(G'')$, we have
\begin{align*}
    {}^td\xi_{\lambda'}^{\ast}\chi_{C_{\phi'}}^{\textrm{mic}}(\mathcal{F})&= \chi_{C_{\phi'}}^{\textrm{mic}}(d\xi_{\lambda'}^{\ast}\mathcal{F})\\
    &=(-1)^{d(\phi')} \trace(1, \NEvs_{z'} d\xi_{\lambda'}^{\ast}\mathcal{F})\\
    &= (-1)^{d(\phi'')} \trace(1, \NEvs_{d\xi_{\lambda'}(z')} \mathcal{F})\\
    &= \langle \etaABV{\phi''}{,\pure}, \mathcal{F} \rangle.
\end{align*}
Namely, $ {}^td\xi_{\lambda'}^{\ast}( \chi_{C_{\phi'}}^{\textrm{mic}})= \etaABV{\phi''}{,\pure}$. Also, from \eqref{eq: Kashiwara-Dubson G'}, we write
\[\chi^\textrm{mic}_{C_{\phi'}} = \sum_{i} c(\phi_i',\phi')\, \chi^\textrm{loc}_{C_{\phi_i'}}\]
for some coefficient $c(\phi_i', \phi').$ For each $i$, write $\phi_i'':= \xi_{\lambda'} \circ \phi_i'$ for the corresponding $L$-parameter of $G''$.
Then by Proposition \ref{theorem: FPF}, we have
\begin{align*}
    R(\etaABV{\phi''}{,\pure})&= R( {}^t d \xi_{\lambda'}^{\ast} \chi_{C_{\phi'}}^{\textrm{mic}})\\
    &=R( {}^t d \xi_{\lambda'}^{\ast} \sum_{i}c(\phi_i', \phi') \chi^\textrm{loc}_{C_{\phi_i'}} )\\
    &= R( \sum_{i}c(\phi_i', \phi') \mu_{\phi_i''}^{\pure} )\\
    &=\sum_{i}c(\phi_i', \phi') T_{\phi_{i}'},
\end{align*}
where 
\[T_{\phi_{i}'}:= \sum_{\widetilde{\phi}_i'\in \Omega_{\phi_i''}} \sum_{\widetilde{\tau}_i' \in \widehat{A}_{\widetilde{\phi}_i'}} \sum_{{\tau}_i'' \in \widehat{A}_{\phi_i''}} m(\widetilde{\tau}_i', \tau_i'') \rank(\tau_i'') (\stsh{\widetilde{\tau}_i'}{C_{\widetilde{\phi}_i'}})^{\ast}.\]

Next, we claim that 
\[ T_{\phi_i'}= \frac{\ell}{ |\Omega_{\phi_i''}|} \sum_{\widetilde{\phi}_i' \in \Omega_{\phi_i''}} \chi^\textrm{loc}_{C_{\widetilde{\phi}_i'}} .\]
Indeed, if $|\Omega_{\phi_i''}|= \ell$, then $A_{\phi_i''} \simeq  A_{\widetilde{\phi}_i'}$ for all $\widetilde{\phi}_i' \in \Omega_{\phi_i''}$. Thus, 
\[ T_{\phi_{i}'}= \sum_{\widetilde{\phi}_i' \in \Omega_{\phi_i''}} \sum_{\widetilde{\tau}_i' \in \widehat{A}_{\widetilde{\phi}_i'}} \rank(\widetilde{\tau}_i') (\stsh{\widetilde{\tau}_i'}{C_{\widetilde{\phi}_i'}})^{\ast} = \sum_{\widetilde{\phi}_i' \in \Omega_{\phi_i''}} \chi^\textrm{loc}_{C_{\widetilde{\phi}_i'}}.\]
If $|\Omega_{\phi_i''}|= 1$, then $\Omega_{\phi_i''}=\{\phi_i'\}$. By the same argument at the end of the proof of Proposition \ref{prop lift=trans disconnected}, using Frobenius reciprocity of finite groups, we get   
\begin{align*}
   T_{\phi_{i}'}&= \sum_{\tau_i' \in \widehat{A}_{\phi_i'}} \sum_{{\tau}_i'' \in \widehat{A}_{\phi_i''}} m({\tau}_i', \tau_i'') \dim(\tau_i'') (\stsh{{\tau}_i'}{C_{{\phi}_i'}})^{\ast}\\
   &= \sum_{{\tau}_i' \in \widehat{A}_{{\phi}_i'}} \dim(\Ind_{A_{\phi_i'}}^{A_{\phi_i''}}\tau_i')(\stsh{{\tau}_i'}{C_{{\phi}_i'}})^{\ast}\\
   &= \ell \chi^\textrm{loc}_{C_{\phi_i'}}.
\end{align*}
This proves the claim.

In summary, we obtain
\begin{align*}
    R( \etaABV{\phi''}{,\pure})& = \sum_{i}c(\phi_i', \phi') \frac{\ell}{ |\Omega_{\phi_i''}|} \sum_{\widetilde{\phi}_i' \in \Omega_{\phi_i''}} \chi^\textrm{loc}_{C_{\widetilde{\phi}_i'}}\\
    &= \sum_{i}c(\phi_i', \phi') \sum_{1 \leq k \leq \ell}\chi^\textrm{loc}_{C_{{}^{\widehat{\theta}^k} \phi_i'}} \\
    &=  \sum_{1 \leq k \leq \ell} {}^{\widehat{\theta}^k}\left(\sum_{i}c(\phi_i', \phi')\chi^\textrm{loc}_{C_{ \phi_i'}}\right)\\
    &= \sum_{1 \leq k \leq \ell}  \chi^{\mic}_{C_{{}^{\widehat{\theta}^k}\phi_i'}}\\
    &= \frac{\ell}{|\Omega_{\phi''}|} \sum_{\phi' \in \Omega_{\phi''}} \chi^{\mic}_{C_{\phi_i'}}.
\end{align*}
This completes the proof of the lemma.
\end{proof}

Similarly to Definition~\ref{definition: regular parameter}, we say that {\it $\phi''$ is reguar for $G$}\IN{regular parameter of $G''$ for $G$} if $d\xi(\Lambda_{x_{\phi''}}^\mic) \cap \Lambda_{\xi\circ\phi''}^\mic$ is non-empty.
This is well-defined by Assumption~\ref{assu G G' G''}(ii).
\begin{proposition}\label{prop: Lift of ABV disconnected}
    Let $(G,G',G'')$ satisfy  Assumption~\ref{assu G G' G''}.
Suppose  $\phi''$ is microlocally generic for $G$, as defined above. 
Set $\phi= \xi \circ \phi''$. 
Then
    \[
    \Lift{G'',s}{G_\delta}
    (\etaABV{\phi''}{})
    =
    e(G_\delta)\,
    \etaABV{\phi,s}{,\delta}  ,
    \]
for all pure inner forms $\delta$ for $G$.
\end{proposition}

\begin{proof}
Lemma \ref{lem disconnected R ABV} gives us
\[ R \circ \frac{1}{\ell} \QST^{\geo}(\etaABV{\phi''}{})= \frac{1}{|\Omega_{\phi''}|}\sum_{\phi' \in \Omega_{\phi''}} \chi^{\mic}_{C_{\phi_i'}}.\]
For each $\phi'$ in the summation, the proof of Theorem~\ref{theorem: Lift of ABV} gives us ($\lambda':=\lambda_{\phi'}$)
\[ \mathrm{proj}_{\delta} \circ {}^t(T \circ m^{\ast} \circ d \xi_{\lambda'}^{\ast}) (\chi^{\mic}_{C_{\phi_i'}})= \mathrm{proj}_{\delta}(\chi_{C_{\phi}}^{\mic})= \mathrm{proj}_{\delta} (\etaABV{\phi}{,\pure})=e(G_{\delta})\etaABV{\phi}{,\delta}. \]
Here, the second equality comes from Lemma \ref{lemma: NEvs and chimic}. 
The rest of the argument follows the proof of Theorem~\ref{theorem: Lift of ABV}.
\end{proof}

\section{Passage to the identity component}\label{section: passage to the identity}

In this section we explain how ABV-packets and distributions on $G^0$ related to those on $G$.

\begin{proposition}\label{proposition: induction from identity}
    Let $G$ be a quasisplit group satisfying Assumption~\ref{assu: groups we consider} and let $\delta$ be a pure inner form for $G$ as defined in Section~\ref{subsection: pure inner forms}. 
    Let $\xi : \Lgroup{G^0} \to \Lgroup{G} = \Lgroup{G^0} \rtimes {\langle{\hat\theta}\rangle}$ be the obvious inclusion.
    Assume Hypothesis~\ref{hypothesis: KLH} for $G$ and $G^0$.
    Suppose $G$ and $G^0$ have compatible Langlands  correspondences. Namely, Desideratum \ref{desideratum: LLC} \RES \ holds. 
    Then, with respect to the pairings introduced in Section~\ref{subsection: pairing},
\[
    \langle \Ind_{G_\delta^0}^{G_\delta}(\sigma), \mathcal{F} \rangle
    =
    \langle \sigma , d\xi^* \mathcal{F}\rangle
\]
for all $\sigma \in K_0\Rep_{\lambda^0}(G^0_\delta)$ and for all $\mathcal{F}\in K_0 \Per_{H_{\lambda}}(V_{\lambda})$, for any infinitesimal parameter $\lambda$ for $G$.
\end{proposition}

\begin{proof} 
By definition, every Langlands parameter for $G$ factors through $G^0$, which is to say, $\phi = \xi\circ \phi^0$. 
Note that, passing to Lie algebras, $\xi$ induces an equivariant morphism $d\xi : V_{\lambda^0} \to V_{\lambda}$ of varieties with group action, for every  infinitesimal parameter $\lambda$ for $G$. 
Indeed, letting $\lambda$ (resp. $\lambda^0$) denote the infinitesimal parameter of $\phi$ (resp. $\phi^0$), it is clear that $d\xi : V_{\lambda^0}\to V_{\lambda}$ is trivial on the underlying space though not trivial on the group actions, since $H_{\lambda^0}$ is a finite-index subgroup of $H_{\lambda}$.
\[
\begin{tikzcd}
    H_{\lambda} &&  V_{\lambda} \ar[loop left]{r}{ H_{\lambda} } 
     \\
    \arrow[>->]{u}
    H_{\lambda^0} 
    && 
    V_{\lambda^0} \arrow[loop left]{r}{H_{\lambda^0}} \arrow[>->>]{u} {d\xi}
\end{tikzcd}
\]

Now, without loss of generality, we may assume that $\sigma$ is a standard representation of $G_\delta^0(F)$ and $\mathcal{F}$ is a standard sheaf on $V_{\lambda}$: 
take $\sigma = \SR(\pi^0)$ for an irreducible representation $\pi^0=\pi(\phi^0, \tau^0)$ of $G_\delta^0(F)$ with infinitesimal parameter $\lambda^0$ and take $\mathcal{F} = \stsh{\tau'}{C_{\phi'}}$ where $\phi'$ is a Langlands parameter for $G$ with infinitesimal parameter $\lambda$.
By Desideratum \ref{desideratum: LLC} \RES, we write
\[ \Ind_{G_\delta^0}^{G_\delta}(\pi^0) = \bigoplus_{\tau \in \widehat{A}_{\phi}}m(\tau, \Ind_{A_{\phi^0}}^{A_{\phi}}\tau^0) \pi(\phi,\tau).\]
From Definition \ref{def: induction of standards to disconnected}, we have 
\[ \Ind_{G_\delta^0}^{G_\delta}(\SR(\pi^0)) = \bigoplus_{\tau \in \widehat{A}_{\phi}}m(\tau, \Ind_{A_{\phi^0}}^{A_{\phi}}\tau^0) \SR(\pi(\phi,\tau)). \]
On the other hand, by \cite{ABV92}*{Proposition 23.7(b)}, write
\[ d\xi^* \stsh{\tau'}{C_\phi} = \sum_{(\phi^0)'}\sum_{(\tau^0)' \in \widehat{A}_{ (\phi^0)'}} m((\tau^0)', \tau'|_{A_{(\phi^0)'}} ) \stsh{(\tau^0)'}{C_{(\phi^0)'}}, \]
where the sum is taken over all $(\phi^0)'$ such that $\phi'= \xi^{\ast} \circ (\phi^0)'$.
Thus,
\[
\begin{array}{rcl}
\langle \Ind_{G_\delta^0}^{G_\delta}(\sigma), \mathcal{F} 
\rangle
&=& 
\langle \Ind_{G_\delta^0}^{G_\delta}\SR(\pi^0), \stsh{\tau'}{C_{\phi'}} 
\rangle \\
&=& 
\sum_{\tau \in \widehat{A}_{\phi}} 
m(\tau, \Ind_{A_{\phi^0}}^{A_{\phi}}\tau^0) \langle \SR(\pi(\phi,\tau)), \stsh{\tau'}{C_{\phi'}}  
\rangle\\
&=& \begin{cases}
   e(\tau') m(\tau', \Ind_{A_{\phi^0}}^{A_{\phi}}\tau^0) & \text{ if }\phi'=\phi,\\
   0 & \text{ otherwise.}
\end{cases}
\end{array}
\]    
On the other hand,
\[
\begin{array}{rcl}
\langle \sigma , d\xi^* \mathcal{F} \rangle
&=&
\langle \SR(\pi^0) , d\xi^* \stsh{\tau}{C_\phi} \rangle 
\\
&=& 
\sum_{(\phi^0)'}\sum_{(\tau^0)' \in \widehat{A}_{ (\phi^0)'}} m((\tau^0)', \tau'|_{A_{(\phi^0)'}} ) \langle \SR(\pi^0) , \stsh{(\tau^0)'}{C_{(\phi^0)'}} \rangle\\
&=& \begin{cases}
  e(\tau^0) m(\tau^0, \tau'|_{A_{\phi^0}}) & \text{ if }\phi'=\phi,\\
   0 & \text{ otherwise.}
\end{cases}
\end{array}
\]
By our definition of Kottwitz sign for our disconnected groups, the desired equality
\[
    \langle \Ind_{G_\delta^0}^{G_\delta}(\sigma), \mathcal{F} \rangle
    =
    \langle \sigma , d\xi^* \mathcal{F}\rangle
\]
follows from Frobenius reciprocity of finite groups. This completes the proof of the proposition.
\end{proof}

The previous proposition gives the following description of $\Lift{G^0}{G_\delta}$.

\begin{corollary}\label{corollary: induction from identity}
Let $G$ be a quasisplit group satisfying Assumption~\ref{assu: groups we consider}.
Let $\xi : \Lgroup{G^0} \to \Lgroup{G} = \Lgroup{G^0} \rtimes {\langle{\hat\theta}\rangle}$ be the obvious inclusion.
Let $\delta$ be a pure inner form for $G$ as defined in Section~\ref{subsection: pure inner forms}.
Suppose $G$ and $G^0$ have compatible local Langlands correspondences.
Then
\[
\Lift{G^0}{G_\delta}
    =
\Ind_{G_\delta^0}^{G_\delta} \ \Lift{G^0}{G^0_\delta}
    =
\mathrm{proj}_{\delta} \circ    \,^td\xi^* \circ \QST^{\geo}.
 \]
\end{corollary}

\begin{proof}
Note that $(G^0,1,\xi)$ is an endoscopic triple for $G$.
    Now the result is a direct consequence of the definitions and Proposition~\ref{proposition: induction from identity}.
\end{proof}

\begin{theorem}\label{theorem: tranfer from the identity}
Let $G$ be a quasisplit group satisfying Assumption~\ref{assu: groups we consider} and let $\delta$ be a pure inner form for $G$ as defined in Section~\ref{subsection: pure inner forms}. Let $\phi^0 \in \Phi(G^0)$ and let $\phi:= \xi \circ \phi^0$.
Then
\[
    \etaABV{\phi}{,\pure} = {}^t d \xi^{\ast}  (\etaABV{\phi_0}{,\,\pure}).
\]
Assume Hypothesis~\ref{hypothesis: KLH} for $G$ and $G^0$.
    Suppose $G$ and $G^0$ have compatible Langlands  correspondences. Namely, Desideratum \ref{desideratum: LLC} \RES\ holds. 
    Then as a set
\[
    \Pi^\ABV_{\phi}(G_\delta) = 
    \{ \pi\leq \Ind_{G_\delta^0}^{G_\delta} (\sigma) \tq \sigma\in \Pi^\ABV_{\phi^0}(G_\delta^0) \} .
\]
\end{theorem}

\begin{proof}
The first assertion follows from Proposition \ref{theorem: FPF}. The second assertion follows from Proposition \ref{proposition: induction from identity}.
\end{proof}

\begin{remark}
In \cite{CS2} it is shown that that $\Per_{H_\lambda}(V_\lambda)$ admits a categorical decomposition
\[
\Per_{H_\lambda}(V_\lambda)
=
\mathop{\bigoplus}\limits_{[\delta]\in H^1(F,G^0)}\Per^{[\delta]}_{H_\lambda}(V_\lambda).
\]
For every $[\delta]\in H^1(F,G^0)$, the Langlands correspondence for $G$ determines a perfect pairing
\[
K_0 \Rep_{\lambda}(G_\delta) \times K_0 \Per^{[\delta]}_{H_\lambda}(V_{\lambda}) \to \Z,
\]
and likewise, the Langlands correspondence for $G^0$ determines a perfect pairing
\[
K_0 \Rep_{\lambda^0}(G^0_\delta) \times K_0 \Per^{[\delta]}_{H_{\lambda^0}}(V_{\lambda^0}) \to \Z.
\]
Using these pairings, Proposition~\ref{proposition: induction from identity} may be expressed in the form
    \[
    \Ind_{G_\delta^0}^{G_\delta}
    =
     \,^td\xi^* .
    \]
\end{remark}

\part{Application to Arthur packets}\label{part: V Arthur}

In this section we prove Vogan's conjecture on Arthur packets for quasisplit classical groups and their pure inner forms.

\section{Arthur packets}\label{subsection: Arthur packets classical}

In this section we recall the theory of Arthur packets as given by \cite{Art13}*{Theorem 2.2.1}, \cite{Mok15}*{Theorem 3.2.1}, and \cite{KMSW14}*{Theorem* 1.6.1} (see Section \ref{subsection: Arthur conjectures}) in our framework, i.e., by considering twisted endoscopy as transfer to a potentially disconnected reductive groups. We take the view that $\O^{\alpha}_{2n}$ is a classical group.
Table~\ref{table: endo summary} summarizes the cases of endoscopy that appear in \cites{Art13,Mok15}; they are all special cases of twisted endoscopy (\cite{KS99}) which we recall explicitly in Section \ref{subsection: classical groups}. 

\subsection{Local Arthur parameters }\label{subsection: Arthur parameters}

Let $G$ be a reductive group over $F$ satisfying Definition \ref{assu: groups we consider}.
A local Arthur parameter for a connected reductive group $G$ is a homomorphism
$$\psi: W_F \times \SL_2^D(\BC) \times \SL_2^A(\BC) \longrightarrow {}^L{G},$$
for which the restriction $\psi|_{W_F \times \SL_2^D(\BC)}$ is a tempered $L$-parameter and the restriction $\psi|_{\SL_2^A(\BC)}$ is algebraic.
The first copy of $\SL_2(\BC)$ is called the Deligne-$\SL_2(\BC)$ and is denoted by $\SL_2^D(\BC)$; the second copy of $\SL_2(\BC)$ is called the Arthur-$\SL_2(\BC)$ and is denoted by $\SL_2^A(\BC)$. 
A local Arthur parameter $\psi$ is called {\it generic} or {\it tempered} if $\psi|_{\SL_2^A}$ is trivial. 

Two local Arthur parameters are equivalent if they are conjugated under $\widehat{G}$. We call such an equivalence class an $A$-parameter and we write $\Psi(G)$ for the set of $A$-parameter of $G$. From a local Arthur parameter $\psi$ of $G$, we associate an $L$-parameter $\phi_{\psi}$ of $G$ by
\begin{equation}\label{eqn: psi to phi psi}
\phi_\psi(w,x) = \psi(w,x,d_w), 
\qquad\text{for}\quad
d_w \ceq \begin{pmatrix}
    |w|^{1/2} & 0 \\
    0 & |w|^{-1/2}
\end{pmatrix}.
\end{equation}
From the condition that $\psi|_{W_F}$ is bounded, it follows that the mapping $\psi \mapsto \phi_{\psi}$ is injective.

For any local Arthur parameter $\psi$, define
\[ S_{\psi}:= Z_{\dualgroup{G}}(\psi),\ \ 
A_\psi := S_{\psi}/S_{\psi}^0.
\]
By an argument similar to Proposition~\ref{proposition: enhancement by purity}, there is a canonical function
\[
\widehat{A_\psi} \to H^1(F,G)
\]
defined by pullback along $\pi_0(Z_{\dualgroup{G}}^{\Gamma_F}) \to \pi_0(Z_{\dualgroup{G}}(\psi))$ together with Kottwitz isomorphism and this gives the decomposition
\[
\widehat{A_\psi} = \mathop{\bigsqcup}\limits_{[\delta]\in H^1(F,G)} \widehat{A_\psi}[\delta].
\]
Finally, we let $s_{\psi}$ be the following element in $S_{\psi}$, or its image in $A_{\psi}$ by abuse of notation.\IN{$s_{\psi}$}
\begin{align}\label{eq def of s psi}
    s_{\psi}:= \psi\left(1,1, -1\right).
\end{align}

\subsection{General linear groups}\label{subsection: general linear}

Consider the group $\GL_N(F)$. In this case, we have the local Langlands correspondence by \cites{Hen00, HT01, Sch13}, satisfying all parts of Desiderata \ref{desideratum: LLC}.
For each local Arthur parameter $\psi$ of $\GL_N$, 
\[ 
\Pi_{\psi}(\GL_N(F)) \ceq \Pi_{\phi_{\psi}}(\GL_N) = \{\pi_\psi\}
,
\]
the singleton $L$-packet containing the product of generalized Speh representations $\pi_\psi$, which is unitary, by \cite{Tad86}. 
Define $\etaArt{\psi}{} \in K_0 \Pi(G)$ by $\etaArt{\psi}{} \ceq \pi_\psi$.
Now set
\[
J_{N} \ceq \begin{pmatrix}
&&&1 \\
&&-1& \\
&\iddots&&\\
(-1)^{N-1}&&&
\end{pmatrix}
\]
and let
$\theta_N : \GL_N \to \GL_N$ be the automorphism of $F$-groups defined by
\begin{equation}\label{equation: thetaN}
\theta_N(g) 
= 
J_N^{-1}\ {}^{t}g^{-1}\ J_N.
\end{equation}
Note that $\theta_N$ is an outer automorphism that preserves the standard pinning.
Now define
\[
\GL_N^+\ceq \GL_N\rtimes {\langle\theta_N\rangle}.
\]
Then
\[
\Ind_{\GL_N}^{\GL_N^+}(\pi_{\psi}) = \pi_{\psi}^+ \oplus \pi_{\psi}^-,
\]
where $\pi_{\psi}^+$ is the Whittaker normalized extension of $\pi_{\psi}$ from $\GL_N(F)$ to $\GL_N^+(F)$ (see Section \ref{sec Whittaker normalization}).
We extend $\psi$ to a local Arthur parameter $\psi^+$ for $\GL_N^+$.
Motivated by our LLC, we set
\[
\Pi^\Art_{\psi^+}(\GL_N^+) 
= 
\{ \pi_{\psi}^+,\pi_{\psi}^- \}.
\]
When combined with our work in Part~\ref{part: II Transfer}, we get
\begin{align}\label{eq Arthur distribution disconnected GL}
    \Transfer{G^0,{\hat\theta_N}}{\GL_N^+}\ \frac{1}{2} \Res^{G}_{G^0} (\etaArt{\psi}{}) 
=
\pi_{\psi}^+ - \pi_{\psi}^- =:\etaArt{\psi^+,{\hat\theta_N}}{}.
\end{align}

Note that $\dualgroup{\GL}_N=\GL_N(\BC)$ and $\hat\theta_N(g)=J_N^{-1}\ {}^tg^{-1}\ J_N$ for $g\in\GL_N(\BC)$. 
Furthermore, we have that $H^1(F,\GL_n)=\{1\}$ via Kottwitz isomorphism and so $\GL_N$ has no nontrivial pure inner forms.

\subsection{Classical groups}\label{subsection: classical groups}

Having found Arthur packets for general linear groups, in the rest of this section we turn to quasisplit classical groups, by type. 

In this section we consider the case when $G$ is one of the quasisplit groups $\U^E_n$, $\SO_{2n+1}$, $\Sp_{2n}$ and $\O^{\alpha}_{2n}$. 
Here we recall the definitions of these groups, their pure inner forms, how they appear as twisted endoscopic groups, and their elliptic endoscopic data.

\subsubsection{Unitary groups (type A)}
\label{subsubsection: A}

In this section, $G = \U^E_n$, the quasisplit unitary $F$-group split over the quadratic extension $E/F$; thus,
\[
G(F) = \U^{E}_n(F) = \{ g\in \GL_n(E) \tq \,^tc(g) J g = g\}
\]
where $c : \GL_n(E) \to \GL_n(E)$ is the component-wise extension of the non-trivial element $c : E\to E$ in $\Gal(E/F)$.
The dual group to $G$ is 
\[
\dualgroup{G}= \GL_n(\BC),
\]
which comes equipped with a non-trivial action by $W_F$, defined as follows: 
if the image of $w\in W_F$ under $W_F \to W_F/W_E \iso \Gal(E/F)$ is the non-trivial element $c$, then
\begin{equation}\label{equation: theta A}
w(g) = {\hat\theta}_n(g) = J_n^{-1}\ \,^tg^{-1}\ J_n. 
\end{equation}
Note that $Z(\GL_N(\BC))^{c}=\{\pm I_{n}\}$, where $I_n$ denotes the identity. From Kottwitz isomorphism, we obtain $H^1(F,G)=\mathbb{Z}/2\mathbb{Z}$, {\it i.e.}, $G$ has a unique nontrivial pure inner form.

Now we have
$\Lgroup{G} = \GL_n(\BC) \rtimes W_F$.
It follows from the above discussion that the standard representation of $\Lgroup{G}$ is on $\BC^{2n}$ for the action of $\dualgroup{G}$ given by
$
(g\rtimes 1)(v_1,v_2) = (gv_2,{\hat\theta_n}(g)v_2)
$
and for the action of $W_F$ through the quotient to $\Gal(E/F)$ given by
$
(1\rtimes c)(v_1,v_2) = (v_2,v_1).
$

Now we define $G^+$ and explain how to see $G$ as a twisted endoscopic group.
In this case,
\[
G^+ \ceq \Res_{E/F}\GL_n \rtimes{\langle\theta\rangle},
\]
where  
\[
\theta(g):=J_n^{-1}\ {}^tc(g)^{-1}J_n,
\]
which is an outer automorphism that preserves the standard pinning. 
Then 
\[\dualgroup{G^+} =
\GL_n^2(\BC)\rtimes {\langle{\hat\theta}\rangle},
\]
where $\hat{\theta}(g_1,g_2)=(\hat{\theta}_n(g_2),\hat{\theta}_n(g_1)).$
Now let
$
\xi : \Lgroup{G} \to \Lgroup{G^+} 
$
be induced from 
\[
\dualgroup{G} 
= 
\GL_n(\BC)
\hookrightarrow
\dualgroup{G^+}
=
\GL_n^2(\BC) \rtimes{\langle{\hat\theta}\rangle}.
\]
Then $(G,{\hat\theta},\xi)$ is an endoscopic triple for $G^+$.

We now review elliptic endoscopic data for $G$.
For $n=2k+1$ and $k=k_1+k_2$ 
\begin{equation}\label{equation: s A odd}
G' = \U^{E}_{2k_1+1} \times \U^{E}_{2k_2}
\qquad\text{and}\qquad
s = \begin{pmatrix}
-I_{k_2} &&\\
&I_{2k_1+1}&\\
&&-I_{k_2}
\end{pmatrix}.
\end{equation}
Then the stabilizer of $s$ in $\dualgroup{G}$ is $\GL_{2k_1+1}(\BC) \times \GL_{2k_2}(\BC)$ and $(G',s,\xi')$ is an endoscopic triple for $G$, where $\xi' : \Lgroup{G'} \to \Lgroup{G}$ is the map induced from  
\[
\dualgroup{G'} = \GL_{2k_1+1}(\BC) \times \GL_{2k_2}(\BC) \hookrightarrow \GL_{2k} = \dualgroup{G}.
\]
For $n=2k$ and $k=k_1+k_2$, set 
\begin{equation}\label{equation: s A even}
G' = \U^{E}_{2k_1} \times \U^{E}_{2k_2}
\qquad\text{and}\qquad
s=\begin{pmatrix}
-I_{k_2}&&\\
&I_{2k_1}&\\
&&-I_{k_2}
\end{pmatrix};
\end{equation}
then the stabilizer of $s$ in $\dualgroup{G}$ is $\GL_{2k_1+1}(\BC) \times \GL_{2k_2}(\BC)$ and $(G',s,\xi')$ is an endoscopic triple for $G$, where $\xi' : \Lgroup{G'} \to \Lgroup{G}$ is the map induced from  
\[
\dualgroup{G'} = \GL_{2k_1}(\BC) \times \GL_{2k_2}(\BC) \hookrightarrow \GL_{2k} = \dualgroup{G}.
\]
Any elliptic endoscopic triple for $G$ is equivalent to some $(G',s,\xi')$ as above.

\subsubsection{Odd special orthogonal groups (type B)}
\label{subsubsection: B}

In this section $G = \SO_{2n+1}$, the split even orthogonal group over $F$, so
\[
G(F) 
=
\left\{
g\in \GL_{2n+1}(F) \tq 
\theta_{2n+1}(g) = g, \ \det(g) =1
\right\}.
\]
The dual group for $G = \SO_{2n+1}$ is
$\dualgroup{G}=\Sp_{2n}(\BC)$. Since $G_n$ is split, the $L$-group is $\Lgroup{G}=\Sp_{2n}(\BC) \times W_F$ and $Z(\dualgroup{\SO}_{2n+1})=\{\pm I_{2n}\}$. By Kottwitz isomorphism, $H^1(F,\SO_{2n+1})=\mathbb{Z}/2\mathbb{Z}$, {\it i.e.}, $G=\SO_{2n+1}$ has a unique nontrivial pure inner form.

Note that we have the natural embedding
$$\xi: \Lgroup{G}=\Sp_{2n}(\BC) \times W_F \hookrightarrow \GL_{2n}(\BC) \times W_F = \Lgroup{\GL_{2n}}.$$
In this case we set 
\[
G^+ = \GL_{2n}\rtimes \langle\theta \rangle
\]
for $\theta = \theta_{2n}$ as defined in \eqref{equation: thetaN}, so
\[
\dualgroup{G^+} = \GL_{2n}(\BC)\rtimes{\langle{\hat\theta}\rangle}
\]
for 
\begin{equation}\label{equation: theta B}
    {\hat\theta}(g) 
    =
    J_{2n}^{-1}\ \,^tg^{-1}\ J_{2n}.
\end{equation}
Then $G^+$ is a disconnected group of the sort considered in Definition~\ref{assu: groups we consider} and $(G, {\hat\theta} ,\xi)$ is an endoscopic triple for $G^+$ as defined in Section~\ref{subsection: endoscopic triple}.

The elliptic endoscopic data of $G=\SO_{2n+1}$ are as follows. 
Set
\begin{equation}\label{equation: s B}
G'= \SO_{2n_1+1} \times \SO_{2n_2+1}
\qquad\text{and}\qquad
s=\begin{pmatrix}
-I_{n_2}&&\\
&I_{2{n_1}}&\\
&&-I_{n_2}
\end{pmatrix},
\end{equation}
for $n=n_1+n_2$.
Then 
\[
Z_{\dualgroup{G}}(s) = \Sp_{2n_1}(\BC) \times \Sp_{2n_2}(\BC)
\]
and $(G',s,\xi')$ is an elliptic endoscopic triple for $G$ where $\xi' : \Lgroup{G'} \to \Lgroup{G}$ induced from 
\[
\dualgroup{G'} = \Sp_{2n_1}(\BC) \times \Sp_{2n_2}(\BC) \hookrightarrow \Sp_{2n}(\BC) = \dualgroup{G}.
\]

\subsubsection{Symplectic groups (type C)}\label{subsubsection: C}

In this section, $G = \Sp_{2n}$ over $F$, so 
\[
G(F) = \Sp_{2n}(F)
=
\{
g\in \GL_{2n}(F) \tq \theta_{2n}(g) = g, \det(g) =1
\}.
\]
The dual group of $G$ is
$\dualgroup{G} =\SO_{2n+1}(\BC)$ and the $L$-group is $\Lgroup{G}=\SO_{2n+1}(\BC) \times W_F$ since $\Sp_{2n}(F)$ is split. We obtain that $Z(\dualgroup{\Sp_{2n}})=\{I_{2n+1}\}.$ By Kottwitz isomorphism, $H^1(F,G)=\{1\}$ and so $\Sp_{2n}(F)$ has no nontrivial pure inner forms.

We have the natural embedding
\[
\Lgroup{G} =\SO_{2n+1}(\BC) \times W_F 
\hookrightarrow 
\Lgroup{\GL_{2n+1}} =
\GL_{2n+1}(\BC) \times W_F.
\]
Here we have 
\[
G^+= \GL_{2n+1}\rtimes \langle\theta_{2n+1}\rangle
\]
for $\theta_{2n+1}$ as in \eqref{equation: thetaN}.
Then $G^+$ is a disconnected group of the sort considered in Definition~\ref{assu: groups we consider}, and 
\[
\Lgroup{G^+} = \GL_{2n+1}(\BC)\rtimes{\langle{\hat{\theta}_{2n+1}}\rangle}
\]
for
\begin{equation}\label{equation: theta C}
{\hat\theta}_{2n+1}(g) 
=
J_{2n+1}^{-1}\ \,^tg^{-1}\ J_{2n+1}.
\end{equation}
Then $(G,{\hat\theta}_{2n+1},\xi)$ is an endoscopic triple for $G^+$ as defined in Section~\ref{subsection: endoscopic triple}.

The elliptic endoscopic data of $G=\Sp_{2n}$ are as follows. 
Set
\begin{equation}\label{equation: s C}
G' = \Sp_{2n_1} \times \SO_{2n_2}^{\alpha}
\qquad\text{and}\qquad
s=\begin{pmatrix}
-I_{n_2}&&\\
&I_{2{n_1}}&\\
&&-I_{n_2}
\end{pmatrix}
\end{equation}
for $n=n_1+n_2$, where $\alpha$ is a square class in $F^\times$ (see \S\ref{subsubsection: D} for the notation $\SO_{2n_2}^\alpha$). 
Then $(G',s,\xi')$ is an endoscopic triple for $G$ where $\xi': \Lgroup{G'} \hookrightarrow \Lgroup{G}$ is induced from
\[
\dualgroup{G'} 
=
Z_{\dualgroup{G}}(s)^0
=
\SO_{2n_1+1}(\BC) \times \SO_{2n_2}(\BC)
\hookrightarrow
\SO_{2n+1}(\BC)
=
\dualgroup{G}.
\]

\subsubsection{Even special orthogonal groups (type D)}\label{subsubsection: D}

Let $\alpha \in F^\times$ be a square class over $F^\times$, possibly trivial, and let $G = \O_{2n}^\alpha$ be the corresponding quasisplit even orthogonal group; then $G= \O_{2n}^\alpha$ is split over $F$ if and only if $\alpha$ is a square.
{More specifically, if $\alpha$ represents the trivial square class, then 
\[
\O_{2n}^\alpha=\{g\in\GL_{2n} \ | \ \tilde{J}_{2n}{}^tg=g\tilde{J}_{2n}\},
\]
where
\[
\tilde{J}_{2n}= \begin{pmatrix}
&&&1 \\
&&1& \\
&\iddots&&\\
1&&&
\end{pmatrix},
\]
while if $\alpha$ represents a nontrivial square class, then we may set
\[
\O_{2n}^\alpha=\{g\in\GL_{2n} \ | \ \tilde{J}_{2n,\alpha}{}^tg=g\tilde{J}_{2n,\alpha}\},
\]
where 
\[
\tilde{J}_{2n,\alpha}=\mathrm{diag}(I_{n-1},\begin{pmatrix}
    0 & 1 \\
    -\alpha & 0
\end{pmatrix}, I_{n-1}) \tilde{J}_{2n}.
\]
We note that these definitions are consistent when $\alpha$ represents the trivial square class since $F$ has characteristic 0, {\it i.e.}, both $\tilde{J}_{2n}$ and $\tilde{J}_{2n,\alpha}$ define equivalent quadratic forms over $F$. However, we set 
\[
\O_{2n}=\{g\in\GL_{2n} \ | \ \tilde{J}_{2n}{}^tg=g\tilde{J}_{2n}\}
\]
and so our convention is that $\O_{2n}=\O_{2n}^1$. 

Then $G^0 = \SO_{2n}^\alpha:=\{g\in\O_{2n}^\alpha \ | \ \det(g)=1\}$ (and similarly define $\SO_{2n}$). This is a connected reductive group and its complex dual group is $\dualgroup{\SO_{2n}^\alpha} = \SO_{2n}(\BC)$. When $\alpha$ represents the trivial square class, the Galois group acts entry-wise on $\SO_{2n}^\alpha$. Thus, in this case, we have $\SO_{2n}^\alpha(F)=\GL_{2n}(F) \cap \SO_{2n}^\alpha.$ In particular, $\SO_{2n}^\alpha$ splits over $F$ and thus $\Lgroup{\SO_{2n}^\alpha}=\SO_{2n}(\BC)\times W_F$ carries the trivial action of $W_F$ when $\alpha$ represents the trivial square class. Note that $\Lgroup{G} = \Lgroup{\O^\alpha_{2n}}=\O_{2n}(\BC)\times W_F$ in this case. 

In this paragraph, we assume that $\alpha$ represents a nontrivial square class. Then the Galois action is twisted (since it arises from an outer form of $\SO_{2n}$) which we make explicit below. Note that 
\[
{}^tP\tilde{J}_{2n,\alpha}P=\tilde{J}_{2n},
\]
where
\[
P=\mathrm{diag}(I_{n-2}\begin{pmatrix}
    1 & \frac{1}{2} \\
    \frac{1}{\sqrt{\alpha}} & \frac{-1}{2\sqrt{\alpha}}
\end{pmatrix}, I_{n-2}).
\]
This defines an isomorphism over $E=F[\sqrt{\alpha}]$, via conjugation by $P$, of $\SO_{2n}^\alpha(E)$ with $\SO_{2n}(E).$ In particular, we obtain that $\SO_{2n}^\alpha$ splits over $E$ and so the Galois action factors through $\Gal(E/F).$ Let $\sigma\in \Gal(E/F)$ denote the nontrivial element.
Via the above isomorphism, we obtain the action of $\Gal(E/F)$ on $\SO_{2n}^\alpha$. Specifically, for $g\in \SO_{2n}^\alpha$, the action of $\sigma$ on $\SO_{2n}^\alpha$ is given by
\[
\sigma_\alpha(g):=P^{-1}\sigma(PgP^{-1})P=w_\alpha\sigma(g)w_\alpha^{-1}
\]
where
\[
w_\alpha=P^{-1}\sigma(P)=\mathrm{diag}\left(
    I_{n-2},\begin{pmatrix}
     0 & \frac{1}{2} \\
     2 & 0
    \end{pmatrix}, I_{n-2}\right).
\]
Note that $w_\alpha\in \O_{2n} \setminus \SO_{2n}.$ We obtain that $\Lgroup{\SO_{2n}^\alpha}=\SO_{2n}(\BC)\rtimes W_F,$ where the action of $W_F$ factors through the quotient $W_F/W_E \simeq \Gal(E/F)$  since $\SO_{2n}^\alpha$ splits over $E$ that. In particular, we then identify $\sigma\in \Gal(E/F)$ with $w_\alpha$  which determines an isomorphism $\SO_{2n}(\BC)\rtimes\Gal(E/F)\cong\O_{2n}(\BC).$ 

From Kottwitz isomorphism, we determine the number of pure inner forms of $G^0 = \SO_{2n}^\alpha$ and hence of $G = \O_{2n}^\alpha$. First, we have $H^1(F,\SO_{2})=\{1\}$ since $Z(\dualgroup{\SO}_2)=\SO_{2}.$ For $n\geq 2$, we have that $H^1(F,\SO_{2n})=\mathbb{Z}/2\mathbb{Z}$ since $Z(\dualgroup{\SO}_{2n})=\{\pm I_{2n}\}$. Finally, if $\alpha$ represents a nontrivial square class, then $H^1(F,\SO_{2n}^\alpha)=\mathbb{Z}/2\mathbb{Z}$. Indeed, we have $\dualgroup{\SO_{2}}=\{\mathrm{diag}(t,t^{-1}) \ | \ t\in\BC^\times\}$ and so $Z(\dualgroup{\SO_2})^{w_\alpha}=\{\pm I_2\}.$ Furthermore,  $Z(\dualgroup{\SO}_{2n})=\{\pm I_{2n}\}$ for $n\geq 2$ from which the claim follows.

Next, we discuss how $G_n^0=\SO_{2n}^\alpha$ arises as a twisted endoscopic group of $\GL_{2n}$. First, we note that $\Lgroup{\SO_{2n}^\alpha}=\SO_{2n}(\BC)\rtimes \Gal(E/F)$ can be regarded as a subgroup of $\O_{2n}(\BC)$, where, if $\sigma\in\Gal(E/F)$ is nontrivial, then it is identified with $w_\alpha.$ This provides 
the $L$-embedding
\[
\xi^0:\Lgroup{\SO_{2n}^\alpha} \rightarrow \O_{2n}(\BC)\rtimes W_F\rightarrow \GL_{2n}(\BC)\times W_F=\Lgroup{\GL_{2n}}.
\]
We note that $\xi^0$ does depend on whether the square class $\alpha$ is trivial. Indeed, if $\alpha$ represents the trivial square class, then $\xi^0(1\rtimes w)=1\rtimes w$ for any $w\in W_F$, while if $\alpha$ represents a nontrivial square class, then $\xi^0(1\rtimes \sigma)=w_\alpha\rtimes \sigma$ via the identification of $\sigma\in\Gal(E/F)$ with $w_\alpha\in \O_{2n} \setminus \SO_{2n}$. 
We set \begin{equation}\label{eqn: s for SO even}
    s=J'_{2n}\rtimes \hat{\theta}_{2n},
\end{equation} where $J'_{2n}=\mathrm{diag}(1,-1,1,\dots).$
Then $\mathrm{Int}(s)(g)=\tilde{J}_{2n}{}^tg^{-1}\tilde{J}_{2n}$ for any $g\in \GL_{2n}.$ Let $\xi_\alpha$ denote the composition of $\xi^0$ with the embedding $\Lgroup{\GL_{2n}}\rightarrow\Lgroup{\GL_{2n}^+}$.
Then we obtain that $(\SO_{2n}^\alpha, s, \xi_\alpha)$ is an endoscopic triple for $\GL_{2n}\rtimes \langle\theta_{2n}\rangle$.
}

The elliptic endoscopic data of $G=\O^{\alpha}_{2n}$ are as follows. 
First, consider
\begin{equation}\label{equation: s D}
G' = \SO_{2n_1}^{\alpha_1} \times \SO_{2n_2}^{\alpha_2}
\qquad\text{and}\qquad
s=\begin{pmatrix}
-I_{n_2}&&\\
&I_{2{n_1}}&\\
&&-I_{n_2}
\end{pmatrix},
\end{equation}
for $n=n_1+n_2$, where $\alpha_1 \alpha_2$ is in the same square class as $\alpha$. We note that the case where $n_i=2$ and $\eta_i=1$ for some $i$ is excluded.
Then $(G',s,\xi')$ is an endoscopic triple for $G$ where, for $g=(g_1,g_2)\in \SO_{2n_1}(\BC)\times\SO_{2n_2}(\BC)$ and $w\in W_F$, we have that $\xi'$ is given by
\begin{align*}
\xi'(g\rtimes 1_{W_F})&= \mathrm{diag}(g_1,g_2) \rtimes 1_{W_F}, \\
\xi'(1\rtimes w)&= \mathrm{diag}(w_{\alpha_1},w_{\alpha_2}) \rtimes w.
\end{align*}
It should be noted that $G'$ is also an endoscopic group for $G^0 \SO^{\alpha}_{2n}$ for the same $s$.

There is one other kind of elliptic endoscopic data for $G = \O^{\alpha}_{2n}$, which is not also endoscopic for $G^0 = \SO^{\alpha}_{2n}$:
\begin{equation}\label{equation: s III}
G' = \Sp_{2n_1}\times \Sp_{2n_2}
\qquad\text{and}\qquad
s = t{\hat\theta}
\end{equation}
for 
\[
t =
\begin{pmatrix}
 -I_{n_2}&&\\
&I_{2{n_1}}&\\
&&-I_{n_2}    
\end{pmatrix},
\]
where $n= n_1+n_2$.
Here, $\theta$ agrees with $G^0 = \O^{\alpha}_{2n} = \SO^{\alpha}_{2n}\rtimes{\langle{\theta} \rangle}$.
Let $\xi'$ be induced from
\[
Z_{\dualgroup{G}}(s)^0 
=
\SO_{2n_1+1}(\BC)\times\SO_{2n_2+1}(\BC)
\hookrightarrow
\SO_{2n}(\BC)
\to
\dualgroup{G}.
\]
Then $(G',s,\xi')$ is an endoscopic triple for $G$. Arthur refers to this as the "third kind" of endoscopy in \cite{Art13}*{p.~14}.

\begin{table}[ht]
\label{table: endo summary}
\caption{Endoscopic triples $(G',s,\xi)$ for $G$}
\begin{tabular}{|c|c|c| } 
 \hline
Endoscopic group $G'$ & Ambient group $G$ & Semisimple $s\in \dualgroup{G}$ \\
\hline
\hline
$\U^E_n$ & $\Res_{E/F}\GL_n\rtimes{\langle\theta\rangle}$ &   
\eqref{equation: theta A}
\\
$\U^{E}_{n_1}\times\U^{E}_{n_2}$ & $\U^{E}_{n_1+n_2}$ & \eqref{equation: s A odd} \eqref{equation: s A even} 
\\
\hline
$\SO_{2n+1}$ & $\GL_{2n}\rtimes{\langle\theta\rangle}$ & 
\eqref{equation: theta B}
\\
$\SO_{2n_1+1} \times \SO_{2n_2+1}$ & $\SO_{2(n_1+n_2)+1}$ &
\eqref{equation: s B}
\\
\hline
$\Sp_{2n}$ & $\GL_{2n+1}\rtimes{\langle\theta\rangle}$ & 
\eqref{equation: theta C}
\\
$\Sp_{2n_1} \times \SO^{\alpha}_{2n_2}$ & $\Sp_{2(n_1+n_2)}$ & 
\eqref{equation: s C}
\\
\hline
$\SO^\alpha_{2n}$ & $\GL_{2n}\rtimes{\langle\theta\rangle}$ & 
\eqref{eqn: s for SO even}
\\
$\SO^{\alpha}_{2n_1} \times \SO^{\alpha}_{2n_2}$ & $\O^{\alpha}_{2(n_1+n_2)}$ &  
\eqref{equation: s D}
\\
$\Sp_{2n_1} \times \Sp_{2n_2}$ & 
$\O^\alpha_{2(n_1+n_2)}$
& 
\eqref{equation: s III} 
\\
\hline
\end{tabular}
\end{table}

\subsection{Non-elliptic endoscopic triples for classical groups}\label{subsection: G'}

In this subsection, we describe the endoscopic triple $(G', s, \xi')$ for the groups $G$, which are $\U^E_n$, $\SO_{2n+1}$, $\Sp_{2n}$ and $\O^{\alpha}_{2n}$. 

In the previous subsection, we have described an embedding $\std: \dualgroup{G} \hookrightarrow \GL_N(\BC)$ case by case. For any semisimple element $s \in \dualgroup{G}$, we let $\Lambda$ denote the multi-set of eigenvalues of $\std(s) \in \GL_N(\BC)$. We may write 
\[ \Lambda= \{\lambda_1^{m_1},\ldots, \lambda_{r}^{m_r}\} \sqcup \{(\lambda_1^{-1})^{m_1},\ldots, (\lambda_{r}^{-1})^{m_r}\} \sqcup \{1^{n_1}, (-1)^{n_2}\}, \]
where $\lambda_i \neq \pm 1$ and $m_1,\ldots, m_r, n_1, n_2$ indicate the multiplicity since 
these eigenvalues come in pairs (for unitary groups, notice that $s$ commutes with $\Gamma$). 
Then, the identity component of the centralizer of $s$ in $\dualgroup{G}$ is isomorphic to
\[ 
\GL_{m_1}(\BC) \times \cdots \times \GL_{m_r}(\BC) \times \dualgroup{G_1^0} \times \dualgroup{G_2^0}, 
\]
where $G_1$, $G_2$ are the groups of the same type as in the elliptic endoscopic triple in previous subsection. See \cite{LLS24}*{\S 6.1} for details of the computation. As a consequence, let $(G', s, \xi')$ be the endoscopic triple determined by $s$. Then 
\begin{align}\label{eq non-elliptic endoscopic group}
    G' \cong \GL_{m_1}(\BC) \times \cdots \times \GL_{m_r}(\BC) \times G_1^0 \times G_2^0,
\end{align}
where $G_1, G_2$ are classical groups from the list $\U^E_{n_i}$, $\SO_{2n_i+1}$, $\Sp_{2n_i}$ or $\O^{\alpha}_{2n_i}$. 

\subsection{Definition of Arthur packets for classical groups}\label{subsection: Arthur conjectures}

In this subsection, we state the definition of Arthur packets for (pure inner forms) of classical groups uniformly. Thus, we let $G$ is one of the quasisplit $\U_{n}^E$ , $\SO_{2n+1}$, $\Sp_{2n}$, or $\O_{2n}^{\alpha}$. 
Let $G^{0}$ be the identity component of $G$; then $G=G^{0}$ if $G$ is not $\O_{2n}^{\alpha}$, and $G^{0}= \SO_{2n}^{\alpha}$ otherwise. 
We let $G^+:=\Res_{E/F}\GL_N\rtimes \langle \theta \rangle$ be the disconnected group such that $(G, s, \xi)$ is an endoscopic triple for $G^+$ ($E=F$ unless $G = \U_n$).
For a local Arthur parameter $\psi$ of $G$, we let $\psi^0$ be the corresponding local Arthur parameter of $G^0$, and set 
\[ \psi^+:= \xi \circ \psi,\]
which is a local Arthur parameter of $G^+$.

\begin{conjecture}[{\cite{Art89}}]\label{conj: Arthur}
   In the above setting, the following holds for any local Arthur parameter $\psi$ of $G$. For Parts (2), (3) and (4), we fix a pure inner form $\delta \in H^1(F,G)_{\theta}$.
    \begin{enumerate}
        \item There exists a unique distribution ${\eta}_{\psi}$ of $G$ that lies in the image of $\Ind_{G^0}^G: K_{\BC} \Pi(G^0) \to K_{\BC}\Pi(G)$        
        such that $ 
\Res_{G^0}\eta_{\psi}$ is a stable distribution on $G^0$, and it satisfies the twisted endoscopic character relation
        \[\Transfer{G^0,s}{G^+}\left( {\frac{1}{[G:G^0]}} \Res_{G^0}\eta_{\psi} \right) =  \etaArt{\psi^+,s}{}, \]
        where $\etaArt{\psi_{\GL^+,s}}{}$ is defined in \eqref{eq Arthur distribution disconnected GL}.
        \item For $s \in S_{\psi}$, let $(G',s,\xi')$ be the corresponding endoscopic triple such that $ \psi$ factors through $\Lgroup{G}'$ as a local Arthur parameter $\psi'$. Let $G''$ be the product of general linear groups with at most two classical groups such that $G'= (G'')^{0}$ (see \eqref{eq non-elliptic endoscopic group}), and let $\eta_{\psi''}$ be the unique distribution on $G''$ defined by Part (1); see Remark \ref{rmk Arthur packet def}, below.
        Then, the distribution
          \begin{align*}
            \eta_{\psi,s}^{\Art,\delta}:= e(G_\delta) \Transfer{G',s}{G_{\delta}} \left( \frac{1}{[G'':G']} \Res_{G'} \eta_{\psi''}\right)
        \end{align*}
        only depends on the image of $s$ in $A_{\psi}$. Thus, for $x \in A_{\psi}$, we define $\eta_{\psi,x}^{\Art,\delta}:= \eta_{\psi,s}^{\Art,\delta}$ for any choice of $s\in S_{\psi}$ that maps to $x$.
    
        \item For any $\varepsilon \in \widehat{A_{\psi}}[\delta]$, the following element in $K_{\BC}\Pi(G_{\delta})$
        \[ \pi(\psi,\varepsilon):= \frac{1}{|A_{\psi}|} \sum_{x \in A_{\psi}} \overline{\mathrm{trace}(\varepsilon(x s_{\psi}))} \eta_{\psi,x}^{\Art, \delta}\]
        is a non-negative linear combination of irreducible representations of $G_{\delta}$. 
        \item When $\psi$ is tempered, $\pi(\psi,\varepsilon)$ is an irreducible tempered representation, and all tempered representations arises in this way. Moreover, $\pi(\psi, \varepsilon)$ is generic with respect to the Whittaker datum of $G_n$ we fixed if and only if $\varepsilon$ is trivial.
        \item Define $\Pi_{\psi}(G_{\delta})$ as a multi-set such that
        \[ \bigoplus_{\pi \in \Pi_{\psi}(G_{\delta})} \pi= \bigoplus_{\varepsilon \in \widehat{A_{\psi}}[\delta]} \pi(\psi,\varepsilon).\]
        Then $\Pi_{\psi}(G_{\delta})$ consists of unitary representations.
    \end{enumerate}
\end{conjecture}

The endoscopic relation in Part (1) of the conjecture is equivalent to the one in \cite{Art13}*{Theorem 2.2.1(a)}, as we now explain.

    First, recall that the transfer $\Transfer{G^0,\theta}{G^+}$ in Part (1) is scaled from Kottwitz-Shelstad transfer by $2=|\pi_0(G^+)|$ (see Definition \ref{definition: modified transfer}). Also,
\[ 
(\eta_{\psi^+})_{\theta}
= 
(\pi_{\psi}^+ - \pi_{\psi}^-)_{\theta}= 2 (\pi_{\psi}^+)_{\theta}.\]
Thus, our twisted endoscopic character relation agree with Arthur's.

Second, for even orthogonal groups, the twisted endoscopic transfer from $K_{\BC}\Pi(\SO_{2n})^{\st}$ to $K_{\BC}\Pi(G^+)$ is not an injection. More precisely, let $(\cdot)^c$ denote the outer conjugation of irreducible representations and extend it to a linear transformation of the Grothendieck group. Then the transfer factors through the quotient
\[K_{\BC}\Pi(\SO_{2n})^{\st}/\mathrm{span}(\{ \eta -\eta^c\ | \ \eta \in K_{\BC}\Pi(\SO_{2n})^{\st}\}).\]
In \cite{Art13}*{Theorem 2.2.1(a)}, Arthur considered a member $[\eta_{\psi^0}]$ in this quotient
that has the correct image under twisted endoscopic transfer. In our formulation, we consider $\eta_{\psi}:= \Ind_{G^0}^G(\eta_{\psi^0})$, which is a well-defined distribution in $K_{\BC}\Pi(\O_{2n}(F))$. Then it gives a unique representative $\frac{1}{[G:G^0]}\Res_{G^0}^{G}\ \eta_{\psi}$ of the class $[\eta_{\psi^0}]$ that lies in the image of $\Res_{G^0}$.

Since we require that $\eta
_{\psi}$ lies in the image of $\Ind_{G^0}^G$, if we consider the endoscopic triple $(G^0, 1 , \xi)$ of $G$ in Part (2), we get 
\[ \eta_{\psi,1}^{\Art}= \Trans_{G^0,1}^{G} (\frac{1}{[G:G^0]}\Res_{G^0}\, \eta_{\psi})=\frac{1}{[G:G^0]} \Ind_{G^0}^{G} \circ \Res_{G^0}\,  \eta_{\psi}= \eta_{\psi}.\]

\begin{remark}
\label{rmk Arthur packet def}\ Recall from \eqref{eq non-elliptic endoscopic group} that
\[G' \simeq \GL_{m_1} \times \cdots \times \GL_{m_r} \times G_1^{0} \times G_2^0,\]
where $G_1, G_2$ are certain classical groups or trivial groups. For a local Arthur parameter $\psi'= \psi_1 \boxtimes \cdots \boxtimes \psi_r \boxtimes \psi_{G_1} \boxtimes \psi_{G_2}$, the distribution $\eta_{\psi''}$ in the statement is 
\[ \pi_{\psi_1} \otimes \cdots \otimes \pi_{\psi_r} \otimes \eta_{\psi_{G_1}} \otimes \eta_{\psi_{G_2}} , \]
where $\pi_{\psi_i}$ is the unique element in $\Pi_{\psi_i}(\GL_{m_i})$, and $\eta_{\psi_{G_i}}$ is the distribution of $G_i$ defined in Part (1) of Conjecture \ref{conj: Arthur}.
\end{remark}

From Part (3) of Conjecture \ref{conj: Arthur}, we have
        \[ 
        \eta_{\psi,x}^{\Art, \delta}= \sum_{\varepsilon \in \widehat{A}_{\psi}[\delta]} \mathrm{trace}(\varepsilon(x s_\psi))\ \pi(\psi, \varepsilon). \]
This conjecture is proved by Arthur (\cite{Art13}) for quasisplit $G=\SO_{2n+1}, \Sp_{2n}, \mathrm{O}_{2n}$ and by \cite{Mok15} for quasisplit $G=\U_n$. For type $D$, Arthur considered outer conjugation classes of representations of $\SO^{\alpha}_{2n}$, and we reformulated in terms of representations of $\O_{2n}^{\alpha}$ following  \cite{AG17}, see Definition \ref{def A packet SO 2n}(2) below.
For pure inner form of quasisplit classical groups, the above conjecture is proved when $\psi$ is tempered for $\U_n$ by \cite{KMSW14}, for $\SO_{2n+1}$ by Ishimoto \cite{Ish24}, and for $\SO_{2n+1}$, $\Sp_{2n}$, $\mathrm{O}_{2n}$, $\U_n$ by M{\oe}glin-Renard (\cite{MR18}). 
See also \cite{Pen25} for more recent progress. We record the cases relevant to this paper below.

\begin{theorem}[\cites{Art13,KMSW14,Mok15,AG17,MR18,Ish24}]\label{thm Arthur's theory for quasisplit and tempered}
Conjecture \ref{conj: Arthur}(1) holds for the (quasisplit) classical groups $G$. Part (2), (3) and (4) of Conjecture \ref{conj: Arthur} holds when $\delta=1$ or $\psi$ is tempered.
\end{theorem}

Assume that Conjecture \ref{conj: Arthur} holds for $G_{\delta}$, we rephrase the definition of Arthur packets and emphasize the induced mapping $\Pi_{\psi}(G_{\delta}) \to \widehat{A_{\psi}}[\delta]$. We also record the definition of the Arthur packets for $\SO_{2n}^{\alpha}$.

\begin{definition}\label{def A packet SO 2n}
Assume that Conjecture \ref{conj: Arthur} holds for $G_{\delta}$.
\begin{enumerate}
    \item The multi-set $\Pi_{\psi}(G_{\delta})$ is called the Arthur packet of $G$ for the local Arthur parameter $\psi$. It is equipped with a mapping  from the multi-set $\Pi_{\psi}(G_{\delta})$ to $\widehat{A_{\psi}}[\delta]$ that sends an irreducible subquotient of $\pi(\psi, \varepsilon)$ to $\varepsilon$. \IN{$\Pi_{\psi}(G_{\delta})$, Arthur packet}

    For $\pi\in \Pi_{\psi}(G_{\delta})$, let 
    \[ \varepsilon_{\pi}:= \bigoplus_{\varepsilon \in \Irr{A_{\psi}}} m(\pi, \varepsilon) \ \varepsilon  \in \Rep(A_{\psi}),\]
where $m(\pi, \varepsilon)$ is the multiplicity of $\pi$ in $\pi(\psi, \varepsilon).$ Then for any $s \in A_{\psi}$, we define
\[ \langle s ,\pi \rangle_{\psi}^{\Art}:= \trace(s, \varepsilon_{\pi}).\]    
    \item Let  $G^0=\SO_{2n}^{\alpha}$. For $\pi \in \Pi(G^0)^\pure$, we write $\pi^c$ for the outer conjugation of $\pi$, and denote its outer conjugacy class by
    \[ [\pi]= \begin{cases}
        \{\pi, \pi^c\} & \text{ if }\pi\not\cong \pi^c,\\
        \{\pi\} & \text{ otherwise.}
    \end{cases} \]
    Let  $\psi^0$ be  a local Arthur parameter of $G^0$, and let $\psi$ denote the corresponding local Arthur parameter of $G= \O_{2n}^{\alpha}$. The Arthur packet $\Pi_{\psi^0}^{\Sigma}(G^0_{\delta})$ is a multi-set of outer conjugacy classes of irreducible representations together with a mapping $\Pi_{\psi^0}^{\Sigma}(G_{\delta}^0) \to \widehat{A_{\psi^0}}[\delta]$ characterized as follows: For any $\varepsilon \in \widehat{A_{\psi}}[\delta]$, let $\varepsilon^0:= \varepsilon|_{A_{\psi^0}}$ and write $\pi^{\Sigma}(\psi^0, \varepsilon^0)$ for the fibre of the mapping. Then
    \[ \bigoplus_{[\pi] \in \pi^{\Sigma}(\psi^0, \varepsilon^0)}  \bigoplus_{\sigma \in [\pi]}\sigma=  \Res_{G_{\delta}^0}^{G_{\delta}} \ \pi(\psi,\varepsilon).\]
\end{enumerate}
\end{definition}

\subsection{Local Langlands correspondence for classical groups}\label{sec verification for hypotehsis on LLC for classical groups}
Theorem \ref{thm Arthur's theory for quasisplit and tempered} gives a local Langlands correspondence of pure inner forms of classical groups that satisfies many parts of Desiderata \ref{desideratum: LLC}. Especially, \textbf{(LC)}, \textbf{(St)}, \textbf{(ECR)}, and \textbf{(Sha)} hold.  We describe the construction for $L$-packets of $\O_{2n}^{\alpha}$ in details, and the other groups follow similarly.

Here we write $G = G_n \O_{2n}^{\alpha}$.
According to \cite{GGP12}*{Theorem 8.1}, an $L$-parameter $\phi$ of $G$ is determined by its composition with the standard embedding $\std: \Lgroup{G} \hookrightarrow \GL_{2n}(\BC)$, which we denote by $\phi_{\GL}:= \std \circ \phi$. Thus, we may decompose $\phi_{\GL}$ into direct sum of irreducible representations of $W_F \times \SL_2(\BC):$
\begin{align*}
\phi_{\GL}= \bigoplus_{i \in I} \phi_{i} \lvert \cdot \rvert^{x_i} \otimes S_{a_i},
\end{align*}
where $\phi_i$ is an irreducible unitary representation of $W_F$, $x_i \in \R$, and $S_{a_i}$ is the unique $a_i$-dimensional irreducible representation of $\SL_2(\BC)$. Since $\phi_{\GL}$ factors through $\Lgroup{G}$, we may fix a choice of the decomposition 
\[ I= I_{>0} \sqcup I_{ngp,1} \sqcup I_{gp} \sqcup I_{ngp,2} \sqcup I_{<0},\]
where 
\begin{itemize}
    \item $I_{>0}= \{i \in I\ | \ x_i >0 \}$, and $I_{<0}= \{i \in I\ | \ x_i <0 \}$;
    \item $I_{gp}=\{i \in I \ | \ x_i=0 \text{ and } \phi_i \otimes S_{a_i} \text{ preserves an orthogonal form}\}$;
    \item there is a bijection $ (\cdot)^{\vee}: I_{>0} \sqcup I_{ngp,1} \to I_{<0} \sqcup I_{ngp,2}$ such that $\phi_{i^{\vee}}= \phi_{i}^{\vee}$ and $x_{i^{\vee}}= -x_i$.
\end{itemize}
Write $I_{>0}=\{i_1,\ldots, i_r\}$ with $x_{i_1}\geq x_{i_2} \geq \cdots \geq x_{i_r} >0$, and let $d_i:= \dim(\phi_i) a_i$. Now we describe the Langlands classification of $\phi$. Write $\phi \leftrightarrow (P=MN, \phi_t, \nu)$. Then
\[ M \simeq \GL_{d_1} \times \cdots \times \GL_{d_r} \times G_{m},\]
where $m= n- \sum_{i=1}^r d_i$. If $m=0$, then we set $G_m$ to be the trivial group.
With this isomorphism, we have 
\[ \phi_t= \phi_{1}\otimes S_{a_1} \boxtimes \cdots \boxtimes \phi_{r}\otimes S_{a_r} \boxtimes \phi_{=0},\]
where $\phi_{=0}$ is the $L$-parameter of $G_m$ such that
\[\std \circ \phi_{=0}= \bigoplus_{i \in I_{ngp,1} \sqcup I_{gp} \sqcup I_{ngp,2}} \phi_i \otimes S_{a_i}, \]
and the unramified character $\nu$ of $M$ is described by
\[ \nu(  (g_1,\ldots, g_r, h) ):= \prod_{i=1}^r |\det(g_i)|^{x_i}.\]

According to Lemma \ref{lemma: LC for enhanced parameters}, we have an isomorphism
\[ A_{\phi} \simeq A_{\phi_t} \simeq A_{\phi_{=0}}.\]
If $m=0,$ this is the trivial group.
\begin{remark}
    Note that this Levi subgroup is connected only when $m=0$. 
    Indeed, for any standard Levi $M=\GL_{d_1}\times\cdots\times\GL_{d_r} \times G_m^0$ with $m>0$, the unramified character $\nu$ of $M$ is automatically fixed by $\theta$.
\end{remark}

Finally, for $i \in I_{>0}$, we let $\St(\rho_i,a_i)$ denote the tempered representation of $\GL_{d_i}(F)$ whose $L$-parameter is $\phi_i \otimes S_{a_i}$. 
\begin{definition}\label{def LLC even orthogonal}
    With the above notation, the $L$-packet $\Pi_{\phi} (G)^\pure$ is given as follows:
    \begin{enumerate}
        \item If $m \neq 0$, then $\Pi_{\phi}(G)^\pure:=\{\pi(\phi,\varepsilon)\ | \ \varepsilon \in \widehat{A}_{\phi}= \widehat{A}_{\phi_{=0}}\}$, where $\pi(\phi, \varepsilon)$ is the unique irreducible quotient of the following induction
        \[ \SR(\pi(\phi,\varepsilon)):= \Ind_{P}^G  \ \St(\rho_1,a_1)\lvert\cdot\rvert^{x_1} \otimes \cdots \otimes \St(\rho_r,a_r)\lvert\cdot\rvert^{x_r} \otimes \pi(\phi_0, \varepsilon),  \]
        where we regard $\phi_0$ as a local Arthur parameter of $G_m$ trivial on $\SL_2$, and $\pi(\phi_0,\varepsilon)$ is the tempered representation as in Conjecture \ref{conj: Arthur}(4).
        \item If $m=0$, then $\alpha$ must be $1$ and
        $\Pi_{\phi}(G)^\pure=\{\pi(\phi,1)\}$ is a singleton consisting of the unique irreducible quotient of 
        \[\SR(\pi(\phi,1)):=\Ind_{P}^G \ \St(\rho_1,a_1)\lvert\cdot\rvert^{x_1} \otimes \cdots \otimes \St(\rho_r,a_r)\lvert\cdot\rvert^{x_r}.\]
    \end{enumerate}
\end{definition}
The local Langlands correspondence constructed above satisfies many parts of Desiderata \ref{desideratum: LLC}. Especially, \textbf{(LC)} follows from the construction, and \textbf{(St), (ECR), (Sha)} follows from Theorem \ref{thm Arthur's theory for quasisplit and tempered}. 

We end this section by describing the restriction of the standard modules in Definition \ref{def LLC even orthogonal} to $\SO_{2n}^{\alpha}$.

\begin{lemma}\label{lem restriction of std mod even orthogonal}
    In the notation of Definition \ref{def LLC even orthogonal}, we have 
    \begin{enumerate}
        \item If $m\neq 0$, then 
        \[ \Res_{G_n^0} ( \SR(\pi(\phi,\varepsilon)))= \Ind_{P^0}^{G_n^0}  \St(\rho_1,a_1)\lvert\cdot\rvert^{x_1} \otimes \cdots \otimes \St(\rho_r,a_r)\lvert\cdot\rvert^{x_r} \otimes \Res_{G_m^0} (\pi(\phi_0, \varepsilon)). \]
        \item If $m=0$, let $P_1$ and $P_2$ be two distinct standard parabolic subgroups of $\SO_{2n}^{\alpha}$ such that $P_i$ is conjugate under $G_n$ to $P$, the standard parabolic of $G_n$ in the Langlands classification of $\phi$. Then
        \begin{align*}
            \Res_{G_n^0} ( \SR(\pi(\phi,1)))&= \Ind_{P_1}^{G_n^0}  \St(\rho_1,a_1)\lvert\cdot\rvert^{x_1} \otimes \cdots \otimes \St(\rho_r,a_r)\lvert\cdot\rvert^{x_r} \\
      &  \oplus \Ind_{P_2}^{G_n^0}  \St(\rho_1,a_1)\lvert\cdot\rvert^{x_1} \otimes \cdots \otimes \St(\rho_r,a_r)\lvert\cdot\rvert^{x_r} . 
        \end{align*} 
    \end{enumerate}
\end{lemma}

\section{Proof of Vogan's conjecture based on geometric transfer}\label{section: endoclassical}

In this section we prove \cite{CFMMX22}*{Conjecture 1}: Arthur packets are ABV-packets and Arthur's coefficients are calculated by vanishing cycles.

\begin{proposition}
    If $\psi$ is a local Arthur parameter for $G$ then the conormal bundle above $C_{\phi_\psi}$ admits an open $H_\lambda$-orbit, contained in $\Lambda^\mic_{C_{\phi_\psi}}$, where $\lambda$ is the infinitesimal parameter of $\phi_\psi$.
    This orbit admits a canonical base-point $
    z_\psi \ceq (x_\psi,y_\psi) \in \Lambda^\mic_{C_{\phi_\psi}},
    $
    for $x_\psi \ceq x_{\phi_\psi}$, where $y_\psi$ is defined in the proof. 
    With respect to this base-point, 
    \[
    A_\psi =  A^\mic_{z_\psi}.
    \]
\end{proposition}

\begin{proof}
    Set
\begin{equation}
\label{equation: x psi y psi regular}
 x_{\psi}:=d\psi_1\begin{pmatrix}
        0&1\\0&0
    \end{pmatrix} \qquad\text{and}\qquad y_{\psi}:=d\psi_2\begin{pmatrix}
        0&0\\1&0
    \end{pmatrix},
\end{equation}
where $\psi_1 : \SL_2(\BC) \to \dualgroup{G}$ and $\psi_2 : \SL_2(\CC) \to \dualgroup{G}$ are defined by 
\[
\psi(1,x,1) = (\psi_1(x),1)
\qquad
\textrm{and}
\qquad
\psi(1,1,y) = (\psi_2(y),1),
\]
which is to say,
\[\psi_1=\psi\vert_{\{1\}\times \SL_2(\CC)\times \{1\}}: \SL_2(\BC) \to \dualgroup{G} \hspace{0.3cm}\text{and}\hspace{0.3cm} \psi_2=\psi\vert_{\{1\}\times \{1\}\times \SL_2(\BC)}: \SL_2(\BC) \to \Lgroup{G}.
\]
To see that the $H_\lambda$-orbit of $(x_\psi,y_\psi)$ is $\Lambda^\mic_\lambda$, argue as in the proof of \cite{CFMMX22}*{Proposition 6.9} and note that the proof of that proposition applies to $G$ equally as it applies to $G^0$, as it largely relies on dimension counting.
\end{proof}

For any local Arthur parameter $\psi$ for $G$, define
\[
\NEvs_{\psi} : \Per_{H_\lambda}(V_\lambda) \to \Rep(A_\psi)
\]
by the composition of $\NEvs_{C_{\phi_\psi}} : \Per_{H_\lambda}(V_\lambda) \to \Loc_{H_\lambda}(\Lambda^\reg_{C_{\phi_\psi}})$ with the equivalence 
\[
\Loc_{H_\lambda}(\Lambda_{C_{\phi_\psi}}^\reg) \simeq \Rep(A_\psi)
\]
determined by $y_\psi \in \Lambda_{x_{\phi_\psi}}^\reg$.
This is one of our main tools for understanding Arthur packets and the distributions attached to them.
Note that $\NEvs_\psi$ is the restriction of 
\[
\NEvs_{z_\psi}
: 
D^b_{c,H_\lambda}(V_\lambda) \to D \Rep(A_\psi)
\]
to $\Per_{H_\lambda}(V_\lambda)$.
In particular, for any pure inner form $\delta$ for $G$, Definition~\ref{definition: ABV-packet} may be written in the form
\[
\Pi^\ABV_{\phi_\psi}(G_\delta)
=
\{
\pi \in \Pi_\lambda(G_\delta) 
\tq
\NEvs_\psi\, \IC(\pi) \ne 0.
\}
\]

If $\psi$ a local Arthur parameter and $s\in Z_{\dualgroup{G}}(\psi)$ then
${}^{s}\phi_\psi = \phi_\psi$ and $y_\psi \in (\Lambda^\reg_{x_\phi})$ so the condition appearing in the definition of $\etaABV{\phi_\psi,s}{,\delta}$ is satisfied and
\[
\langle s,\pi\rangle_{\phi_\psi}^\ABV 
=
\trace(s, \NEvs_{z_\psi}\IC(\pi)).
\]

\begin{corollary}\label{corollary: transfer in Arthur type}
Let $G$ be a quasisplit group as in Definition~\ref{assu: groups we consider}. 
Assume Desiderata~\ref{desideratum: LLC} and Hypothesis~\ref{hypothesis: KLH} for $G$.
Let $(G',s,\xi)$ be an endoscopic triple for $G$.
Let $\psi'$ be a local Arthur parameter for $G'$ and set $\psi = \xi\circ\psi'$. 
Then, for every pure inner form $\delta$ for $G$,
    \[
    \Transfer{G',s}{G_\delta}( \etaABV{\phi_{\psi'}}{})
    =
    e(G_\delta)\,
    \etaABV{\phi_{\psi},s}{,\delta}.   
    \]
\end{corollary}

\begin{proof}
    This is a consequence of Theorem~\ref{theorem: transfer of ABV-distributions} using Lemma~
    \ref{lemma: Arthur is regular}.
\end{proof}

\begin{lemma}\label{lemma: Arthur is regular}
    Let $G$ be a quasisplit group as in Definition~\ref{assu: groups we consider}.
    Let $(G',s,\xi)$ be an endoscopic triple for $G$ as in Section~\ref{subsection: endoscopic triple}.
    Let $\psi'$ be an Arthur parameter for $G'$.
    Then $\phi_{\psi'}$ is microlocally generic for $G$ in the sense of Definition~\ref{definition: regular parameter}. 
\end{lemma}

\begin{proof}
The argument for $G^0$ given in \cite{CFMMX22}*{Proposition 6.9} adapts directly to $G$.
\end{proof}

Now, let $G$ a quasisplit classical group as reviewed in Section~\ref{subsection: classical groups}, so $G$ is $\U^E_{n}$, $\SO_{2n+1}$, $\Sp_{2n}$ or $\O^{\alpha}_{2n}$.
Recall the definition of $G^+$ from Section~\ref{subsection: classical groups}.
In Section~\ref{sec verification for hypotehsis on LLC for classical groups} we verified Desiderata~\ref{desideratum: LLC} for the groups $G$ and $G^+$.

\begin{remark}
In Theorem~\ref{theorem: Vogan's conjecture} we assume Hypothesis~\ref{hypothesis: KLH} for $G$ and $G^+$. 
See \cite{CHAI} for the proof of this hypothesis for unipotent representations of $\GL_N^+(\mathbb{Q}_p) = \GL_N(\mathbb{Q}_p)\rtimes\langle\theta_N\rangle$;
for all representations of classical groups see \cite{Sol25};
for all representations of general linear groups $\GL_N(F)$, see \cite{bal25};
for unipotent representations of connected reductive groups, see \cites{Lus02, CG10, Sol25}.
\end{remark}

\begin{theorem}[Vogan's conjecture on local Arthur packets]\label{theorem: Vogan's conjecture}
Let $G$ be a quasisplit classical group and $\delta$ a pure inner form of $G$.
Assume Hypothesis~\ref{hypothesis: KLH} for these groups and also for $G^+$. Also assume that Conjecture \ref{conj: Arthur} holds for $G_{\delta}$. 
Let $\psi$ be a local Arthur parameter for $G$.
\begin{enumerate}
    \item[(i)]
    Arthur packets are ABV-packets and the stable distributions associated to Arthur packets are ABV-distributions:
\[
\Pi^\Art_{\psi}(G_{\delta}) 
=
\Pi^\ABV_{\phi_{\psi}}(G_{\delta}),
\qquad\text{and}\qquad
\etaArt{\psi}{}
=
\etaABV{\psi}{}.
\]
\item[(ii)]
More generally, for any semisimple $s\in Z_{\dualgroup{G}}(\psi)$,
\[
\etaArt{\psi,s}{}
=
\etaABV{\phi_{\psi},s}{}.
\]
\item[(iii)]
The characters of $A_{\psi}$ defined by endoscopic character relations are given by traces of representations of $A_{\psi}$ given by $\NEvs_{z_{\psi}}$: 
    \[
    \langle s,\pi\rangle_{\psi}^{\Art}
    =
    \trace\left(s, \NEvs_{z_{\psi}}\IC(\pi)\right)
    =: \langle s,\pi\rangle_{\phi_{\psi}}^{\ABV}.
\]
\end{enumerate}
\end{theorem}
We explain the equality $\Pi^\Art_{\psi}(G_{\delta}) 
=
\Pi^\ABV_{\phi_{\psi}}(G_{\delta})$ in the following remark.

\begin{remark}\label{remark: multiset}
By definition, the ABV-packet 
\[\Pi_{\phi_{\psi}}^{\ABV}(G_{\delta}):=\{\pi \in \Pi_{\lambda}(G_{\delta})\ | \ \NEvs_{\psi}\IC(\pi)\neq 0\}\]
is a set, but a priori the Arthur packet $\Pi^\Art_{\psi}(G_{\delta}) $ can be a multi-set. To compare them, we turn ABV packets into a multi-set together with a function with value in $\Irr{A_{\phi_{\psi}}^{\mic}}$ as follows. If $\NEvs_{\psi}(\IC(\pi))$ is nonzero, we decompose it into a direct sum of irreducible representations of $A_{\phi_{\psi}}^{\mic}$:
\[ \NEvs_{\psi} \IC(\pi) = \varepsilon_1\oplus \cdots \oplus \varepsilon_{r}.\]
Then the multi-set contains $r$ copies of $\pi$, and the function take the $i$-th copies to the irreducible representation $\varepsilon_i$. We will show that the multi-set defined in this way agree with the multi-set $\Pi_{\psi}(G_{\delta})$.

We remark that it is proved in \cites{Moe11, Moe11b} that Arthur packets for pure inner forms of symplectic and orthogonal groups are indeed sets. However, we do not need this fact in our proof.
\end{remark}

The argument proceeds in 4 steps.

\subsubsection{\bf Step 1.}

Let $\psi$ be a local Arthur parameter for $G_n$.
Set $\psi^{+} = \xi\circ \psi$ where $(G_n,{\hat\theta},\xi)$ is the endoscopic triple for $G_n^+ = \GL_N^+$; set $\phi^+ = \phi_{\psi^+}$.
We begin by showing
\[
\etaArt{\psi^+,{\hat\theta}}{}
=
\etaABV{\phi^+,{\hat\theta}}{}.
\]

To simplify notation 
we write $\psi_{\GL}$ for the associated local Arthur parameter for $(G_n^+)^0 = \GL_N$.
Recall that
\begin{equation}\label{equation: step 1, a}
\etaArt{\psi^+,{\hat\theta}}{}
=
{\pi_{\psi_\GL}^+} - {\pi_{\psi_\GL}^-},
\end{equation}
where $\pi_{\psi_\GL}$ is the product of generalized Speh representations of $\GL_N(F)$ and $\pi_{\psi_\GL}^+$ is its Whittaker normalized extension to $\GL_N^+(F)$ and $\pi_{\psi_\GL}^-$ is the other representation appearing in 
\[
\Ind_{\GL_N}^{\GL_N^+}
\ \pi_{\psi_\GL} 
= \pi_{\psi_\GL}^+ \oplus \pi_{\psi_\GL}^-.
\]

\begin{lemma}
With notation as above, set $\phi = \phi_{\psi}$, $\phi_{\GL} = \phi_{\psi_{\GL}}$ and $\phi^+ = \phi_{\psi^+}$.
Then
\[
\etaABV{\phi^+,{\hat\theta}}{}
=
{\pi_{\psi_\GL}^+} -{\pi_{\psi_\GL}^-}.
\]
\end{lemma}

\begin{proof}
By Theorem~\ref{theorem: tranfer from the identity}, 
\[
\Pi^\ABV_{\phi^{+}}(\GL_N^+)
=
\Ind_{\GL_N}^{\GL_N^+} \ \Pi^\ABV_{\phi_{\GL}}(\GL_N).
\]
By \cites{CR2, Rid23, Lo24},
\[
\Pi^\ABV_{\phi_{\GL}}(\GL_N) = \{ \pi_{\psi_\GL} \}= \Pi_{\phi_{\GL}}(\GL_N),
\]
so 
\[
\Pi^\ABV_{\phi^{+}}(\GL_N^+) = \{ \pi_{\psi_\GL}^+, \pi_{\psi_\GL}^- \}= \Pi_{\phi^{+}}(\GL_N^+).
\]
As a consequence, Proposition \ref{prop L-packet inside ABV packet} gives
\[ \eta_{\phi^+,\hat{\theta}}^{\ABV}=\langle {\hat\theta}, \pi_\psi^+ \rangle_{\phi^+}^\ABV
 \ [\pi_\psi^+]
 +
 \langle {\hat\theta},
 \pi_\psi^- \rangle_{\phi^+}^\ABV 
 \ [\pi_\psi^-]=\langle {\hat\theta}, \pi_\psi^+ \rangle_{\phi^+}
 \ [\pi_\psi^+]
 +
 \langle {\hat\theta},
 \pi_\psi^- \rangle_{\phi^+} 
 \ [\pi_\psi^-]= \pi_{\psi_{\GL}}^+-\pi_{\psi_{\GL}}^-. \]
This completes the proof of the lemma.
\end{proof}

This shows \eqref{equation: step 1, a} and therefore 
completes Step 1.

\subsubsection{\bf Step 2.}\label{sec step 2}
    The second step in the proof of this theorem is to show 
\begin{equation}\label{equation: type A, step 2, a}
\etaArt{\psi}{} = \etaABV{\phi_{\psi}}{},
\end{equation}
for the classical groups $G$ we considered. This follows from the following proposition.
\begin{proposition}\label{prop step 2}
Let $G$ be a classical group we considered and recall the definition of $G^+$ from Section~\ref{subsection: classical groups}. Fix a local Arthur parameter $\psi$ of $G$ and let $\lambda$ be its infinitesimal parameter.
\begin{enumerate}
    \item  The vector space $K_{\BC}\Pi_{\lambda}(G)^{\std}$ is contained in the image of $\Ind_{G^0}^G: K_{\BC}\Pi(G^0)^{\st} \to K_{\BC}\Pi(G).$
    \item On the vector space $K_{\BC}\Pi_{\lambda}(G)^{\std}$, the two linear transformations coincide
    \[ \Lift{G,\hat{\theta}}{G^+}= \frac{1}{[G:G^0]}\Transfer{G^0,\hat{\theta}}{G^+} \circ \Res_{G^0}\]
\end{enumerate}
\end{proposition}

Indeed, Theorem \ref{theorem: Conjecture 2 qs} implies that $\etaABV{\phi_{\psi}}{} \in K_{\BC}\Pi_{\lambda}(G)^{\std}$. Then Proposition \ref{prop: Lift of ABV disconnected} and step 1 implies that
\[ \Lift{G,\hat{\theta}}{G^+}(\etaABV{\phi_{\psi}}{}) = \etaABV{\phi_{\psi^+,\hat{\theta}}}{}= \etaArt{\psi^+,\hat{\theta}}{}. \]
Now Proposition \ref{prop step 2} above implies that $\etaABV{\phi_{\psi}}{}$ satisfies the characterization of the unique distribution $\eta_{\psi}$ of $G$ in Conjecture \ref{conj: Arthur}(1), which verifies \eqref{equation: type A, step 2, a}.

Now we prove the proposition.
\begin{proof}

Recall that $K_{\BC}\Pi_{\lambda}(G)^{\std}$ has a basis $\{\mu_{\phi}\tq\phi \in \Phi_{\lambda}(G)\}$. Fixing a $\phi \in \Phi_{\lambda}(G)$ and write $\phi \leftrightarrow (P=MN, \phi_t, \nu)$ under the Langlands classification. Let $\phi^+:= \xi \circ \phi$ be the corresponding $L$-parameter of $G^+$.

Suppose first that $G$ is connected. Then 
\[ M \simeq \GL_{n_1}(E) \times \cdots \times \GL_{n_f}(E) \times G^-,\]
where $G^-$ is a connected classical group of the same type as $G$. Under this isomorphism, write 
\[ \phi_t= \phi_1 \boxtimes \cdots \boxtimes \phi_f \boxtimes \phi_{G^-},\ \ \nu( g_1,\ldots,g_f, g^-)= \prod_{i=1}^f |g_i|^{x_i}.\]
 Then by the definition of the $L$-packet $\Pi_{\phi}$, we have  $\mu_{\phi}= \Ind_P^G (\mu_{\phi_t}\otimes \nu)$, where
\[ \mu_{\phi_t}= \pi_{\phi_1} \otimes \cdots \otimes \pi_{\phi_f} \otimes  \eta_{\phi_{G^-}},\]
where $\Pi_{\phi_i}(\GL_{n_i}(E))= \{\pi_{\phi_i}\}$ and $\eta_{\phi_{G^-}}$ is the stable distribution on $G^-$ associate to the tempered local Arthur parameter $\phi_{G^-}$ as in Conjecture \ref{conj: Arthur}. In particular, $\mu_{\phi}$ is also stable, which proves Part (1) in this case. Moreover, by the compatibility of transfer and induction (Proposition \ref{proposition: induction commutes with transfer}), we have 
\[ \Transfer{G,\hat{\theta}}{ G^+} (\mu_{\phi})=\Transfer{G,\hat{\theta}}{ G^+}  \circ \Ind_{P}^G \mu_{\phi_t}\otimes \nu= \Ind_{P^+}^{G^+}  \circ \Trans_{M,\hat{\theta}}^{M^+} \left( \mu_{\phi_t} \otimes \nu \right), \]
where $P^+$ is the standard parabolic of $G^+$ with standard Levi 
\[ M^+ \simeq \big{(}\GL_{n_1}(E) \times \cdots \times \GL_{n_f}(E) \times \GL_{N^-}(E) \times \GL_{n_f}(E) \times \cdots \times \GL_{n_1}(E)  \big{)} \rtimes \langle \hat{\theta} \rangle.\]
Now by Theorem \ref{thm Arthur's theory for quasisplit and tempered} for the tempered local Arthur parameter $\phi^-$ of $G^-$ and Proposition \ref{prop twisted transfer from Siegel GL to disconencted GL}, we have 
\begin{equation}\label{eqn: Siegelish}
\Trans_{M,\hat{\theta}}^{M^+} (\mu_{\phi_t} \otimes \nu )= (\pi_{\xi\circ \phi_t\otimes \nu^+})^+ -(\pi_{\xi\circ \phi_t\otimes \nu^+})^-,
\end{equation}
where 
 \[\pi_{\xi\circ \phi_t \otimes \nu^+} =\pi_{\phi_1}\lvert\cdot \rvert^{x_1} \otimes \cdots \otimes \pi_{\phi_f}\lvert\cdot \rvert^{x_f} \otimes \pi_{\phi_{G^-}^+}  \otimes \pi_{\phi_f}^{\vee}\lvert\cdot \rvert^{-x_f} \otimes \cdots \otimes \pi_{\phi_1}^{\vee}\lvert\cdot \rvert^{-x_f}.\]
Note that $\phi^+ \leftrightarrow (P^+, \xi\circ \phi_t, \nu^+)$ under the Langlands classification. Thus, we conclude that 
 \[\Transfer{G,\hat{\theta}}{ G^+}(\mu_{\phi})= \Ind_{P^+}^{G^+}  \circ \Trans_{M,\hat{\theta}}^{M^+} \left( \mu_{\phi_t} \otimes \nu \right)= \mu_{\phi^+,\hat{\theta}}. \]
 On the other hand, recall that we have assumed Kazhdan-Lusztig hypothesis for $G$ and $G^+$. Then Theorem \ref{theorem: lift = transfer} gives 
 \[\Lift{G,\hat{\theta}}{ G^+}(\mu_{\phi}) =  \mu_{\phi^+,\hat{\theta}} \]
 as well. We conclude that the linear transformations $\Transfer{G,\hat{\theta}}{ G^+}$ and  $\Lift{G,\hat{\theta}}{ G^+}$ agrees on a basis of $K_{\BC}\Pi_{\lambda}(G)^{\std}$. This proves Part (2) in this case.

Next, we consider the case that $G$ is not connected, {\it i.e.}, $G= \O_{2n}^{\alpha}$. If $G^-$ is not trivial, then Lemma \ref{lem restriction of std mod even orthogonal}(1) implies that 
\[ \Res_{G^0}(\mu_{\phi})= \Ind_{P^0}^{G^0} \left( \pi_{\phi_1} \otimes \cdots \otimes \pi_{\phi_f} \otimes   \Res_{(G^-)^0}\eta_{\phi_{G^-}} \right),\]
and the rest of the argument follows with straightforward modifications, which we omit. If $G^-$ is trivial, then  $\mu_{\phi}= \Ind_P^G (\mu_{\phi_t}\otimes \nu)$, where
\[ M \simeq \GL_{n_1}(E) \times \cdots \times \GL_{n_f}(E),\]
is connected and
\[\mu_{\phi_t}= \pi_{\phi_1} \otimes \cdots \otimes \pi_{\phi_f}.\]
With the notation in Lemma \ref{lem restriction of std mod even orthogonal}(2), we have 
\[ \Res_{G^0}(\mu_{\phi})= \Ind_{P_1}^{G^0}(\mu_{\phi_t} \otimes \nu) + \Ind_{P_2}^{G^0}(\mu_{\phi_t} \otimes \nu).\]
Similarly, by Proposition \ref{proposition: induction commutes with transfer}, we have 
\[ \half{1}\Transfer{G^0,\hat{\theta}}{G^+} \circ \Res^{G}_{G^0}(\mu_{\phi}) = \half{1}\sum_{i=1}^2 \Transfer{G^0,\hat{\theta}}{G^+} \circ \Ind_{P_i}^{G^0} (\mu_{\phi_t} \otimes \nu)=   \half{1}\sum_{i=1}^2 \Ind_{P^+}^{G^+} \circ \Transfer{ M,\hat{\theta}}{M^+}  (\mu_{\phi_t} \otimes \nu), \]
where 
\[ M^+ \simeq \big{(}\GL_{n_1}(E) \times \cdots \times \GL_{n_f}(E) \times \GL_{n_f}(E) \times \cdots \times \GL_{n_1}(E)  \big{)} \rtimes \langle \hat{\theta} \rangle.\]
The rest of the argument follows with straightforward modifications, which we omit. This completes the proof of the proposition.
\end{proof}

\subsubsection{\bf Step 3.}\label{sec step 3}
The third step in the proof of the theorem is to show
\begin{equation}\label{equation: step 3, a}
    \etaArt{\psi,s}{,\delta}
    =
    \etaABV{\phi_\psi,s}{,\delta},
\end{equation}
for all semisimple $s\in Z_{\dualgroup{G}}(\psi)$. The equality follows from the following proposition.

\begin{proposition}\label{prop step two}
In the setting of Conjecture \ref{conj: Arthur}(2), let $\lambda$ (resp. $\lambda''$) be the infinitesimal parameter of $\psi$ of $G$ (resp. $\psi''$ of $G''$). Then the linear transformations $\Lift{G'',s}{G_{\delta}}$ and $\frac{1}{[G'':G']} \Transfer{G',s}{G_{\delta}} \circ \Res_{G'}$ agree on the vector space $ K_{\BC}\Pi_{\lambda''}(G'')^{\std}$. In particular, we have 
\begin{align}\label{eq lift of ABV step 3}
    \frac{1}{[G'':G']} \Transfer{G',s}{G_{\delta}} \circ \Res_{G'} (\etaABV{\phi_{\psi''}}{}) = \etaABV{\phi_{\psi},s}{,\delta}
\end{align}
which only depends on the image of $s$ in $A_{\phi_{\psi}}^{\mic}\simeq A_{\psi}$.
\end{proposition}
\begin{proof}
We have evaluated $\Lift{G'',s}{G}$ on a basis of $ K_{\BC}\Pi_{\lambda''}(G'')^{\std}$ in  Proposition \ref{prop lift=trans disconnected} that
\[\Lift{G'',s}{G_{\delta}} (\mu_{\phi''})= \mu_{\phi,s}^{\delta}. \]
Thus, to show the two linear transformations agree on  $ K_{\BC}\Pi_{\lambda''}(G'')^{\std}$,
it suffices to show that 
\[\frac{1}{[G'':G']} \Transfer{G',s}{G_{\delta}} \circ \Res_{G'} (\mu_{\phi''})= \mu_{\phi,s}^{\delta}\]
as well. This is a consequence of the compatibility of transfer and induction (Proposition \ref{proposition: induction commutes with transfer}) and Arthur's theory for tempered packets (Theorem \ref{thm Arthur's theory for quasisplit and tempered}) for Levi subgroups of $G$, similar to the proof of Proposition \ref{prop step 2}. We omit the straightforward modification. Finally, \eqref{eq lift of ABV step 3} follows from Proposition \ref{prop: Lift of ABV disconnected}. This completes the proof of the proposition.
\end{proof}

\subsubsection{\bf Step 4.}
In the last step, we show that
\[ \Pi_{\psi}(G_{\delta})= \Pi_{\phi_{\psi}}^{\ABV}(G_{\delta})\]
as multi-sets, and the functions from these multi-sets to $\Irr{A_{\psi}} = \Irr{A_{\phi_{\psi}}^{\mic}}$ agree.
See Remark \ref{remark: multiset}.

Let $\psi$ be any local Arthur parameter of $G$. Note that $s_{\psi}$ lies in the image of $\psi$. Therefore, $s_{\psi}$ is central in $Z_{\dualgroup{G}}(\psi)$; hence its image in $A_{\psi}$ is also central. We define and compute $\pi^{\ABV}(\psi,\varepsilon)$, the ABV-analogue of $\pi(\psi, \varepsilon)$. {The computation is general, not limited to classical groups.}
\begin{proposition}\label{prop pi psi varepsilon ABV}
     For any irreducible representation $\varepsilon\in 
    \widehat{A}_{\psi}[\delta]$, define
    \[ \pi^{\ABV}(\psi,\varepsilon):= \frac{1}{|A_{\psi}| } \sum_{x \in A_{\psi}} \trace(xs_{\psi}, \varepsilon) \etaABV{\phi_{\psi},x}{,\delta}.\]
    Then 
    \[\pi^{\ABV}(\psi,\varepsilon)= \sum_{\pi \in \Pi_{\phi_{\psi}}^{\ABV}(G_{\delta})}  m(\varepsilon, \NEvs_{\phi_{\psi}}\IC(\pi)) \chi_{\varepsilon}(s_{\psi}) (-1)^{d(\pi)-d(\psi)} \pi ,\]
    where $\chi_{\varepsilon}$ is the central character of $\varepsilon$ and where $m(\varepsilon, \NEvs_{\phi_{\psi}}\IC(\pi))$ is the multiplicity of $\varepsilon$ in $\NEvs_{\phi_{\psi}}(\IC(\pi))$ as a representation of $A_{\phi_{\psi}}^{\mic}$.
\end{proposition}
\begin{proof}
    This follows from a direct computation:
    \begin{align*}
         &\frac{1}{|A_{\psi}| } \sum_{x \in A_{\psi}} \trace(xs_{\psi}, \varepsilon) \etaABV{\phi_{\psi},x}{,\delta_{\varepsilon}}\\
         =&\frac{1}{|A_{\psi}| } \sum_{x \in A_{\psi}} \trace(xs_{\psi}, \varepsilon) \sum_{\pi \in \Pi_{\phi_{\psi}}^{\ABV}(G_{\delta})} (-1)^{d(\pi)-d(\psi)}\trace(x, \NEvs_{\phi_{\psi}}\IC(\pi)) \pi\\
         =& \sum_{\pi \in \Pi_{\phi_{\psi}}^{\ABV}(G_{\delta})}  \left(\frac{1}{|A_{\psi}| } \sum_{x \in A_{\psi}} \trace(x, \varepsilon)\trace(x, \NEvs_{\phi_{\psi}}\IC(\pi)) \right) \chi_{\varepsilon}(s_{\psi})(-1)^{d(\pi)-d(\psi)} \pi\\
         =& \sum_{\pi \in \Pi_{\phi_{\psi}}^{\ABV}(G_{\delta})}  m(\varepsilon, \NEvs_{\phi_{\psi}}\IC(\pi)) \chi_{\varepsilon} (s_{\psi})(-1)^{d(\pi)-d(\psi)} \pi.
    \end{align*}
    Here the third equality follows from the inner product for  characters of representations of the finite group $A_{\psi}$. This completes the proof of the proposition.
\end{proof}
\begin{remark}\label{rmk non-negative}
    Since $m(\varepsilon, \NEvs_{\phi_{\psi}}\IC(\pi))$ is non-negative, we see that Conjecture \ref{conj: Arthur}(3) holds for the pure inner form $G_{\delta}$ for $\psi$ if and only if 
    \begin{align}\label{eq non-negativity}
        (-1)^{d(\pi)-d(\psi)}
    \chi_\varepsilon(s_{\psi})=1
    \end{align} 
    for all irreducible subquotient $\varepsilon$ of $ \NEvs_{\phi_{\psi}}  (\IC(\pi))$. 
    It is a work in progress of the authors to show \eqref{eq non-negativity} in the general case.
\end{remark}

By the previous steps, for any  $\varepsilon \in \widehat{A}_{\psi}[\delta]$, we have 
\[ \pi(\psi, \varepsilon)= \pi^{\ABV}(\psi, \varepsilon).\]
Since we have assumed Conjecture \ref{conj: Arthur}(3) holds for $G_{\delta}$, Proposition \ref{prop pi psi varepsilon ABV} becomes 
\[\pi(\psi, \varepsilon)= \pi^{\ABV}(\psi, \varepsilon)=\sum_{\pi \in \Pi_{\phi_{\psi}}^{\ABV}(G_{\delta})}  m(\varepsilon, \NEvs_{\phi_{\psi}}\IC(\pi))  \pi. \]
Then it is clear from the construction in Remark \ref{remark: multiset} that 
\[ \Pi_{\psi}(G_{\delta})= \Pi_{\phi_{\psi}}^{\ABV}(G_{\delta})\]
as multi-sets, and the functions from these multi-sets to $\Irr{A_{\psi}} \iso \Irr{A_{\phi_{\psi}}^{\mic}}$ agree.

Finally, suppose that 
\[ \NEvs_{\psi}(\IC(\pi))= \varepsilon_1 \oplus \cdots \oplus \varepsilon_r,\]
where $\varepsilon_i \in \Irr{A_{\phi_{\psi}}^{\mic}}.$ 
Then the multiplicity of $\pi$ in $\pi(\psi,\varepsilon)=\pi^{\ABV}(\psi, \varepsilon)$ is equal to
\[m(\varepsilon, \NEvs_{\phi_{\psi}}\IC(\pi))= | \{1 \leq i \leq r \ | \ \varepsilon_i= \varepsilon\}|. \]
Therefore, for any $s \in A_{\psi}$, it follows from definition that
\[ \langle s, \pi \rangle_{\psi}^{\Art}= \trace(s, \bigoplus_{1\leq i\leq r}\varepsilon_i)= \trace(s,\NEvs_{\psi}\IC(\pi))=   \langle s, \pi \rangle_{\phi_\psi}^{\ABV}. \]
This concludes the proof of Theorem \ref{theorem: Vogan's conjecture}.\qed

\section{Arthur packets for even special orthogonal groups}

The proof of the Vogan's conjecture in previous sections is not straightforward when it involves even special orthogonal groups because of the absence of the local Langlands correspondence for them.
Throughout this section, we let $G$ be an even orthogonal group and let $G^0$ be its identity component. We modify the definition of Arthur packets for $G^0$ and prove the Vogan's conjecture, assuming the existence of the local Langlands correspondence for $G^0$ satisfying certain desiderata. The proof is much simpler based on the results in Part \ref{part: II Transfer}.

First, we explain the status for the local Langlands correspondence for even special orthogonal groups. The endoscopic theory of Arthur gives a {\it weak local Langlands correspondence}\IN{$\textrm{WLLC}$, weak local Langlands correspondence} for $G^0$:
\begin{align*}
    \textrm{WLLC}: \Pi(G^0)^\pure_{/\sim} 
    \ \longrightarrow\  \Phi^{\mathfrak{e}}(G^0)_{/ \sim}.
\end{align*}
Here, $\pi \sim \pi'$ (resp. $(\phi,\varepsilon) \sim (\phi', \varepsilon')$ ) if they are equivalent up to a conjugation by some $c \in G$ (resp. $\widehat{c}\in \dualgroup{G}$). See the discussion in \cite{AG17}*{\S 3.6}. A recent work of Peng (\cite{Pen25}) establish an unambiguous local Langlands correspondence extending Arthur's weak local Langlands correspondence
\begin{align}\label{eq LLC SO even}
     \textrm{LLC}: \Pi(G^0)^\pure \to \Phi^{\mathfrak{e}}(G^0),
\end{align}
assuming that $p \neq 2$ and $F/\Q_p$ is unramified. This assumptions may be dropped according to Remark(4) below \cite{Pen25}*{Theorem A}. In the rest of the section, we will assume the following Hypothesis.
\begin{hypothesis}\label{hyp LLC even SO}
    We assume that \eqref{eq LLC SO even} is well-defined in all cases with the following properties.
    \begin{enumerate}
        \item The Kazhdan-Lusztig hypothesis holds.
        \item The assumptions of Theorem \ref{thm: stable standards form a basis - quasisplit} holds. Note that this includes Part (1) and (2) of Corollary \ref{cor ABV SO even ECR} below holds when $\psi$ is tempered.
    \end{enumerate}
\end{hypothesis}

Under this Hypothesis, the distributions $\etaABV{\phi}{}$ is well-defined for all $L$-parameters of $\phi$. Moreover, we have the following endoscopic character relations for $L$-parameters of Arthur type.

\begin{corollary}\label{cor ABV SO even ECR}
Assume Hypothesis \ref{hyp LLC even SO}.
    Let $\psi$ be a local Arthur parameter of a quasisplit even special orthogonal group $G^0$, and let $\delta$ be a pure inner form of $G^0$.
    \begin{enumerate}
        \item [(1)] The stable distribution $\etaABV{\phi_{\psi}}{}$ satisfies the twisted endoscopic character relations
        \[ \Transfer{G^0, \theta}{G^+} (\etaABV{\phi_{\psi}}{}) = \pi_{\psi}^+- \pi_{\psi}^-,\]
        where $(G^0, \theta, \xi)$ is an endoscopic triple of the disconnected general linear group
        $G^+= \GL_N \rtimes \langle \theta \rangle$.
        \item [(2)] Suppose that $(G',s,\xi)$ is an endoscopic triple of $G^0$ so that $\psi= \psi' \circ \psi$ for some local Arthur parameter $\psi'$ of $G'$. Then we have the following endoscopic identity
        \[ \Transfer{G',s}{G_{\delta}^0} (\etaABV{\phi_{\psi}}{})= e(G_{\delta})\etaABV{\phi_{\psi},s}{}.\]
        \item [(3)] For $\pi= \pi(\phi_{\psi}, \varepsilon) \in \Pi_{\phi_{\psi}}(G^0)^\pure$ and $s$ as in Part (2), the coefficient of $\pi$ in $\etaABV{\phi_{\psi},s}{}$ is given by $\trace(s_{\psi}s, \varepsilon)$, where we regard $\varepsilon $ as a representation of $A_{\psi}$ via the restriction along $A_{\phi_{\psi}} \to A_{\psi}$.
    \end{enumerate}
\end{corollary}
\begin{proof}
    Parts (1) and (2) are special case of Theorem \ref{theorem: transfer of ABV-distributions}, where the assumptions are contained in Hypothesis \ref{hyp LLC even SO}. Part (3) is exactly Proposition \ref{prop L-packet inside ABV packet}.
\end{proof}

\appendix

\part{Appendices on geometry and functors, by Jörg Schürmann}\label{part: VI Schurmann}

\fancypagestyle{appendixstyle}{
    \fancyhf{}
    \fancyhead[CO]{\scriptsize APPENDICES ON GEOMETRY AND FUNCTORS}
    \fancyhead[CE]{\scriptsize SCHÜRMANN}
    \fancyhead[R]{\scriptsize\thepage}
    \renewcommand{\headrulewidth}{0pt}
}
\pagestyle{appendixstyle}
\thispagestyle{appendixstyle}

The work starting with \cite{CFMMX22} uses  Deligne's vanishing cycles functor in an equivariant derived category with $\ell$-adic coefficients to define ABV-packets for $p$-adic groups.
In this appendix we switch to the analytic site; for background on nearby and vanishing cycles for the analytic site, we refer to
\cite{MS22}*{Section 10.4} or \cite{Ach21}*{Chapter 4}.
For the corresponding comparison theorem between the $\ell$-adic and analytic nearby and vanishing cycles, see
 \cite{SGA7II}*{Expos\'{e}s XIII et XIV} for the case of torsion coefficients and see
\cite{BBD82}*{Section 6.1} for the general switch from torsion to 
$\bar{\mathbb{Q}}_\ell$-coefficients
(see also, e.g. \cite{Fu11}*{Chapter 9 and 10}, for the foundations on $\ell$-adic
  nearby and vanishing cycles).
We remark that, in fact, $\ev{C}$ comes in the $\ell$-adic context
equipped with an ``action of inertia'', which is to say, an action of $\Gal({\bar \eta}/\eta)$. Its topological counterpart is the monodromy action on the nearby and vanishing cycles. These will not be used in this Appendix.\\

The present paper makes crucial use of a trace formula for these vanishing cycles and also of the Kashiwara-Dubson type formula for comparison of local and microlocal Euler characteristics. 
In this appendix we explain how the inverse of the Kashiwara-Dubson index formula ports to the context of the present paper, provide a proof of the trace formula used in this paper and further show that the functor $\Evs_C$ extends to a sheaf complex on a larger target, namely, off of the regular part of the conormal to an orbit $C$.

This requires a movement from Deligne's equivariant vanishing cycles to their complex analytic counterparts which replaces sometimes the action of a complex algebraic group $G$ with the action of a maximal compact (real) subgroup $K$.
For this, the equivariant derived category of Bernstein-Lunts \cite{BL94} is well-adapted and so we use it here, together with the corresponding restriction and change of group functors, in particular, the forgetful functor $\operatorname{For}: D^b_{c,G}(M;k) \to D^b_{c}(M;k)$.
Note that the nearby and vanishing cycle functors can be lifted up to this equivariant derived context, commuting with the restriction and change of group structures, because they commute with smooth based change
(see, e.g. \cite{Sch03}*{Remark 4.3.7 on p. 265}), and because we only use invariant functions with respect to the group action. 

Consequently, this appendix explains how to move between three languages: local vanishing cycles (in complex algebraic and analytic flavours), stratified Morse theory and $D$-modules, as presaged in \cite{ABV92}*{Chapter 24}. 
In particular, all these approaches to characteristics cycles give the same microlocal Euler characteristics.
More broadly, we will also discuss the merits and demerits of these three different perspectives on microlocal geometry, especially as it applies to the definition of ABV-packets for real and $p$-adic groups.

This leads to the conclusion, among others, that ABV-packets, as defined in \cite{CFMMX22}, coincide with micropackets for $p$-adic groups as conceived in \cite{Vog93}*{Definition 8.5}. At the same time, this compliments the equivariant part of the discussion in \cite{ABV92} regarding the proof of \cite{ABV92}*{Theorem 24.8},
together with a non-equivariant comparison of these three different micro-local view points.
In fact the comparison is via some key results from stratified Morse theory for constructible sheaves, which at the moment are only  available
non-equivariantly.

\section{Applications of (equivariant) vanishing cycles}\label{section: equivariant vanishing cycles}

In the first part of this section we explain  characteristic cycles of (equivariant) constructible sheaf complexes and the inverse Kashiwara-Dubson formula
in the language of local vanishing cycles. In the second part we prove the key fixed-point formula for stalks of equivariant vanishing cycles.

\subsection{Remarks on the inverse Kashiwara-Dubson formula and local vanishing cycles}\label{subsection: KDF}

Although the following constructions apply to a more general context, for simplicity we restrict to the case of a complex algebraic group $G$, e.g. a (finite cover group of a) complex linear algebraic group, acting on a smooth algebraic variety $M$ with finitely many orbits $C$, as this is the only case needed for our applications. 
Then also the identity component $G^0$ is acting on $M$ with finitely-many connected orbits $C'$, so that each $G$-orbit $C$ is a finite disjoint union of $G^0$-orbits $C'$, which are all of the same dimension. The {\it closure ordering} on orbits is defined by
\[ 
C_1 \geq C_2 \Longleftrightarrow \overline{C_1} \supseteq C_2.
\]
Let $k$ be any field of characteristic $0$.

The following {\it commutative} diagram, including the right upper corner only in the case that $k = \BC$, presents all the key players.
\begin{equation}\label{equation: players}
\begin{tikzcd}
   K_0\Per_{G}(M;k) 
   \arrow[->>]{dd}[swap]{\chi^{\operatorname{loc}}}
   \arrow{ddrr}{\operatorname{CC}}
   && 
   \arrow[->>]{dd}{\operatorname{{CC}}}
   \arrow[bend left=5]{ll}{\operatorname{DR}}[swap]{\iso}
   \arrow[bend right=5]{ll}[swap]{\operatorname{Sol}} K_0\operatorname{Mod}_{c}(D_M,G) \\
   && \\
\operatorname{Fun}_{c}(M)^{G} \arrow{rr}{\operatorname{CC}}[swap]{\simeq}
&& 
\operatorname{Lan}(\Lambda)^{G}\\
\arrow[equal]{u}\ 
{\mathbb{Z}[\{\1_{C}\,\vert\, C\}]}
\arrow{rr}{\operatorname{CC}}[swap]{\simeq}
&& 
\arrow[equal]{u}\ 
{\mathbb{Z}[[\{\overline{T^*_C M}]\,\vert\, C\}]}
\end{tikzcd}
\end{equation}
Let us define these key players and explain how they fit into the diagram above.  
\begin{itemize}
    \item 
    In Diagram~\eqref{equation: players}, 
    $\Per_{G}(M;k)$ is the abelian category of $G$-equivariant perverse sheaves on $M$ with coefficients in a field $k$, such as $\bar{\mathbb{Q}}$ or $\BC$ or $\bar{\mathbb{Q}}_\ell$.
    Recall the isomorphisms of Grothendieck groups 
    \[
    K_0 \Per_{G}(M;k)
    \simeq
    K_0 D^b_{c,G}(M;k)
    \simeq 
    K_0 \operatorname{Sh}_{c,G}(M;k),
    \]
    coming from the (bounded and non-degenerate) perverse and standard t-structures for $D^b_{c,G}(M;k)$, where $\operatorname{Sh}_{c,G}(M;k)$ is the abelian category of $G$-equivariant constructible sheaves on $M$ (see, for example, \cite{Sch03}*{Lemma 3.3.1, p. 177}).
    Note that these Grothendieck groups are free abelian of finite rank in case $k$ is algebraically closed and of characteristic $0$, with two distinguished bases produced from the indexing set
    \[
    \Xi(M,G) \ceq \left\{ (C,\mathcal{L}) \tq \textrm{\ $H$-orbit } C , \mathcal{L} \in \Irrep(A_C) \right\},
    \]
    where $A_C \iso Z_H(x)/Z_{H}(x)^0$, for $x\in C$, is the $H$-equivariant fundamental group of $C$.
    From a pair $(C,\mathcal{L})$ the class of the perverse extensions $\IC(\mathcal{L}_C)$ gives a basis for $K_0\Per_{H}(M;k)$ while class of the extension by zero complexes $(j_C)_!\mathcal{L}_C$ gives a basis for $K_0\operatorname{{Sh}_{c,G}(M,k)}$, as recalled in \cite{ABV92}*{Definition 1.22, p.11}, for example.
    Observe that a base extension of such fields induces an isomorphism of these Grothendieck groups and bases, and consequently these groups and bases are independent of the chosen field $k$. This allows us to identify these groups for $k=\bar{\mathbb{Q}}_\ell$ and $k=\BC$ using field extensions $\bar{\mathbb{Q}}_\ell \leftarrow \bar{\mathbb{Q}} \rightarrow \BC$.
    In this way one can compare the topological context with the $\ell$-adic and $D$-module approach.
    \item 
In Diagram~\eqref{equation: players}, $\operatorname{Mod}_{c}(D_M,G)$ denotes the abelian category of $G$-equivariant coherent algebraic $D$-modules, which in our context are automatically regular holonomic by \cite{HTT08}*{Theorem 11.6.1}.
\item 
For $k=\BC$, the $\Z$-module homomorphisms  $\operatorname{DR}$ and $\operatorname{Sol}$ are induced from the the corresponding DeRham and shifted solution functors (to fit with perverse sheaves) of the associated analytic left D-module $(-)^{an}$
(and $\mathcal{O}_M$ and $\mathcal{D}_M$ the sheaf of holomorphic functions and differential-operators on $M$):
\begin{equation*} \begin{split}
DR(-)&=\Omega_M^\bullet\otimes_{\mathcal{O}_M}(-)^{an}[\dim M]\quad \text{and}\\ 
Sol(-)[\dim M]&=Rhom_{\mathcal{D}_M}((-)^{an},\mathcal{O}_M)[\dim M]
\end{split}\end{equation*}
as they appear in the (equivariant) Riemann-Hilbert Correspondence \cite{HTT08}*{Section 7.2} and 
\cite{BL94}*{Section 4.2}, giving the isomorphism across the top of the diagram.
Since the DeRham functor $DR(-)$
and shifted solution functor $Sol(-)[\dim M]$ are exchanged by duality of holonomic $D$-modules fitting with Verdier duality $\mathcal{D}$ \cite{HTT08}*{Proposition 4.2.1 and Corollary 4.6.5}, they agree after composing with $\chi^\mathrm{loc}$ or $\operatorname{CC}$.
Also note that $\chi^\mathrm{loc}$ and $\operatorname{CC}$ both factorize over the functor $\operatorname{For}$ that forgets the equivariant structure of the sheaf, while remembering the Whitney stratification $\mathscr{S}$ (or $\mathscr{S}'$)
of $M$ given by the $G$-orbit (or $G^0$-orbit) structure of $M$ (so that all cohomology sheaves of the underlying constructible sheaf complex $\operatorname{For}\mathcal{F}$ are locally constant on the $G$-orbits $C$ (or $G^0$-orbits $C'$).
\item 
In Diagram~\eqref{equation: players}, $\operatorname{Fun}_{c}(M)^{G}$ denotes the free abelian group of $G$-invariant $\mathbb{Z}$-valued functions on $M$ (which are then automatically constructible), with basis consisting of the indicator functions $1_C$ of the $G$-orbits $C$.
\item 
In Diagram~\eqref{equation: players}, $\chi^\mathrm{loc}$ is the surjective group homomorphism induced from the stalk-wise Euler characteristic; with respect to the basis $1_C$ above, this is given by the collections of group homomorphisms $\chi^\mathrm{loc}_C$.
Observe that 
\begin{equation}\label{equation: chiloc}
\chi^{\textrm{loc}}_{C}((j_{C_i})_!\mathcal{L}_{C_i}) 
= \begin{cases} 
\rank(\mathcal{L}_{C_i})\cdot \ 1_C & C_i = C \\
0 & C_i \ne C.
\end{cases}
\end{equation}
\item 
For $k=\BC$, the vertical $\operatorname{CC}$
is the characteristic cycle homomorphism on regular holonomic $D$-modules as in \cite{Ginsburg:vanishing} using good filtrations for holonomic $D$-modules and their associated graded cycles, fitting with the microlocal multiplicites of Kashiwara \cite{Kas83}*{Theorem 3.2.42, p. 93}. Note that multiplicities of a characteristic cycle for a holonomic $D$-module are non-negative.
Characteristic cycles appear in \cite{ABV92}*{Proposition 19.12} and \cite{Vog93}*{Definition 8.5}; that is, 
for $\mathcal{M}$ in $\operatorname{Mod}_{c}(D_M,G)$, we have 
\begin{equation}\label{equation: CCM}
\operatorname{CC}(\mathcal{M})
=
\sum_{C} \chi^\mic_C(\operatorname{For}\mathcal{M}) [\overline{T^*_C M}].
\end{equation}
For general $k$, the diagonal $\operatorname{CC}$
is the characteristic cycle homomorphism
for constructible sheaf complexes is given by \cite{MS22}*{Section 3.3} and \cite{Sch03}*{Section 5.2}
\begin{equation}
\operatorname{CC}(\mathcal{F}) = \sum_{C} (-1)^{\dim C} \chi(\operatorname{NMD}(\operatorname{For}\mathcal{F},df_x)) [\overline{T^*_C M}],
\end{equation}
where the sheaf theoretical normal Morse data of $\mathcal{F}$ in $D^b_{c,G}(M;k)$ is given by
\begin{equation}\label{equation: A 5}
\operatorname{NMD}(\operatorname{For}\mathcal{F},df_x)
=
\left(R\Phi_{f\vert_{N}}[-1](\operatorname{For}\mathcal{F}\vert_{N})\right)_x,
\end{equation}
as in \cite{MS22}*{(10.30), p. 714} and \cite{Sch03}*{(5.12), p. 283}.
Here, $N$ is normal slice to $C$ at  a chosen $x\in C$ and $f : (M,x) \to (\BC,0)$ is a holomorphic function germ with regular covector $df_x \in (T^*_C M)^\reg$, defined as the compliment in $T^*_C M$ of  $\overline{T^*_{C_i} M}$ for all $C_i \ne C$. If $x\in C'\subset C$ belongs to the $G^0$-orbit $C'$, then $(T^*_C M)^\reg$ restricts to $(T^*_{C'} M)^\reg$ over $C'$.
By stratified Morse theory, this vanishing cycles sheaf depends only on $df_x\in (T^*_{C'} M)^\reg$, for the connected component $C'$ of $C$ containing $x$ (given by the corresponding $G^0$-orbit), but not on the choice of $x$, $f$ and $N$ above (i.e. its isomorphism class is locally constant along $(T^*_{C'} M)^\reg$).
The $G$-equivariant structure of $\mathcal{F}$ implies that, in fact, it only depends on $C$, and not only on the component $C'$ as above.
The fact that $D$-module multiplicities fit with those coming from stratified Morse theory is explained in Remark \ref{micro-solution}
of this Appendix (later on,
see also \cite{Sch03}*{Example 5.3.4, p. 359-360}), thus implying the commutativity of the diagram. Note that the sign $(-1)^{\dim C}$ is coming from the fact that only the shifted restriction $\cdot \vert_{N}[-\dim C]$ is exact with respect to the perverse t-structure, so that also here the multiplicities of a characteristic cycle for a perverse sheaf are non-negative.
This motivates the definition
\begin{equation}\label{equation: stalks 5}
\chi^\mathrm{mic}_C(\mathcal{F}) \ceq (-1)^{\dim C} \chi(\operatorname{NMD}(\operatorname{For}\mathcal{F},df_x)) 
=
(-1)^{\dim C} \chi\left( \left(R\Phi_{f\vert_{N}}[-1](\operatorname{For}\mathcal{F}\vert_{N})\right)_x\right),
\end{equation}
with $df_x \in (T^*_C M)^\reg$; with this definition, we have
\begin{equation}
\operatorname{CC}(\mathcal{F})
=
\sum_{C} \chi^\mic_C(\operatorname{For}\mathcal{F}) [\overline{T^*_C M}].
\end{equation}
The diagonal $\operatorname{CC}$ factorizes over the group of $G$-invariant constructible functions since it is determined by the Euler characteristics of the normal Morse data as explained in \cite{MS22}*{Section 3.3} and \cite{Sch03}*{Section 5.2}; this defines the horizontal $\operatorname{CC}$.
Also note that Verdier duality $\mathcal{D}$ commutes with normal Morse data up to a shift
\cite{Sch03}*{(5.11), p. 281}
$$\mathcal{D}(\operatorname{NMD}(\operatorname{For}\mathcal{F},df_x))\simeq
\operatorname{NMD}(\mathcal{D}(\operatorname{For}\mathcal{F}),-df_x)[-2\dim C]$$
for $x\in C$, together with switching from $df_x$ to $-df_x$. In particular
$\chi^\mathrm{mic}_C$ and $\operatorname{CC}$
are invariant under Verdier duality $\mathcal{D}$.
\item 
Finally, in Diagram~\eqref{equation: players}, $\operatorname{Lan}(\Lambda)$ is the free abelian group generated by the irreducible components $\overline{T^*_{C'} M}$ of the conormal variety 
$$\Lambda:= \bigsqcup_{C'} T^*_{C'} M \subset T^*M$$ with respect to the $G^0$-orbits 
$C'$ in $M$. Note that $\Lambda$ is closed in $T^*M$ by the Whitney a-regularity of the stratification $\mathscr{S}'$.
These conormal varieties $\overline{T^*_{C'} M}$ are $\BC^*$-conic and $G^0$-invariant irreducible Lagrangian subvarieties of the cotangent bundle $T^*M$. We view these cycles as either top dimensional cycles  in Chow groups or as top dimensional  Borel-Moore homology classes of $\Lambda$.
Then $\operatorname{Lan}(\Lambda)$ gets an induced permutation action of $G/G^0$ (and therefore also of $G$),
with $\operatorname{Lan}(\Lambda)^{G}$ the free  abelian group generated by the cycle classes $[-]$ of conormal closures
$\overline{T^*_{C} M}$ of the $G$-orbits $C$ in $M$.
\end{itemize}

Having introduced all the key players, we now turn to a discussion of the index formula and explain how it can be used in the present paper.
The linear functionals $$\chi^{\textrm{loc}}_C : K_0 \Per_{H}(M;k) \to \Z$$ are linearly independent in $K_0 \Per_{H}(M;k)^\vee$  as can be seen by using the constant sheaves $\mathcal{L}_C = \1_{C}$ for all $G$-orbits $C$ in \eqref{equation: chiloc}.
Using the description of $\chi^\mic_C$ in terms of Euler characteristics of vanishing cycles, we next explain a geometric reformulation of a Kashiwara-Dubson type formula of the following shape 
\begin{equation}\label{equation: stalks 7}
\chi^\textrm{mic}_C = \sum_{C_i \geq C} c(C_i,C)\, \chi^\textrm{loc}_{C_i},
\end{equation}
for all $G$-orbits $C$, with the property $c(C,C) = (-1)^{\dim C}$; see, for example, \cite{MS22}*{(10.39), (10.40), p. 724}, \cite{Sch03}*{(5.36), p. 293, (5.38), p. 294}, which treat a more general complex analytic or algebraic context without group actions. 
This implies that the linear functionals $\chi^\textrm{mic}_C$ are linearly independent, as $C$ ranges over $G$-orbits, and generates the same sublattice spanned by $\chi^\textrm{loc}_C$, as $C$ ranges over $G$-orbits.
The following explanation of \ref{equation: stalks 7} will also pave the way for our proof of the Lefschetz trace formula for vanishing cycles, below, which will use traces in place of Euler characteristics.\\

We will now give a geometric description of stalks of nearby cycles $(R\Psi_f \mathcal{F})_x$ and vanishing cycles $(R\Phi_f \mathcal{F})_x$ for complex analytic constructible $\mathcal{F}\in D^b_{\mathscr{S}-c}(M;k)$, for $f: (M,x) \to (\CC,0)$ the germ of a holomorphic function on a complex manifold $M$, following \cite{MS22}*{Section 10.4}.
At this moment we do not ask for a group action on $M$ but 
we fix a complex analytic Whitney stratification $\mathscr{S} = \{ S \}$ of $M$, which in the equivariant context before can be taken to be (the connected components of) either the $G$-orbits or, alternatively, of the intersection of the $G$-orbits in $M$ with a normal slice $N$ at the given point $x$. 
Here, $\mathscr{S}$-constructibility for $\mathcal{F}$ means that the cohomology sheaves of $\mathcal{F}$ restricted to all strata $S$ are locally constant with all stalks finite-dimensional. 
Since we are working only locally, we can assume that we have finitely-many strata.

\begin{proposition}[\cite{MS22}*{Proposition 10.4.11 and Corollary 10.4.14}, \cite{Sch03}*{Lemma 1.1.1 and Example 1.1.3}]
\label{proposition: stalks}
    With the notation as above, fix a complex analytic Whitney stratification $\mathscr{S} = \{ S \}$ of $M$. 
    Choose a real analytic distance germ $r: (M,x) \to ([0,\infty),0)$ with the property $r^{-1}(0) = \{ x \}$.
    For small enough positive $\delta > 0$ (which we write as $0<\delta <<1$), for even smaller positive $\epsilon << \delta$, and $w\in \CC$ with $0< |w| < \epsilon$, one has a canonical isomorphism of distinguished triangles in $D^b_c(\{x\};k)$, for every $\mathcal{F} \in D^b_{\mathscr{S}-c}(M;k)$:
    \begin{equation}\label{diagram: stalks}
    \begin{tikzcd}
    & (R\Phi_f[-1]\mathcal{F})_x \arrow{r} & \mathcal{F}_x \arrow{r} & (R\Psi_f\mathcal{F})_x \arrow {r}{[1]} & \ \\
    & 
    \arrow{u}{\simeq} R\Gamma(T_{f,x}, M_{f,x};\mathcal{F}) \arrow{r} & \arrow{u}{\simeq}
    R\Gamma(T_{f,x};\mathcal{F}) \arrow{r}& \arrow{u}{\simeq}
    R\Gamma(M_{f,x};\mathcal{F}) \arrow{r}{[1]} & \
    \end{tikzcd}
    \end{equation}
where $T_{f,x} \ceq \{ r\leq \delta, |f| \leq \epsilon \}$, a local compact tube around the special fibre $\{ f=0\}$ at $x$, and $M_{f,x} \ceq \{ r\leq \delta, f = w \}$, a local compact nearby fibre called a Milnor fibre. 

Tensoring with a local system $\mathcal{L}$ of finite rank, the projection formula for the lower triangle implies that the canonical morphism
$(R\Psi_f\mathcal{F})\otimes \mathcal{L} 
\stackrel{\sim}{\to} R\Psi_f(\mathcal{F}\otimes \mathcal{L})$ is an isomorphism so that also
\begin{equation}\label{diagram: tensor-phi}
(R\Phi_f[-1]\mathcal{F})\otimes \mathcal{L} \simeq R\Phi_f[-1](\mathcal{F}\otimes \mathcal{L})\:.
\end{equation}
\end{proposition}

Taking the Euler characteristics above after noting that all spaces are finite-dimensional,  Proposition~\ref{proposition: stalks} gives, by additivity of Euler characteristics in the compact subanalytic context (e.g. by subanalytic triangulations or cell decompositions as in \cite{Sch03}*{Section 2.3}),
\begin{equation}\label{equation: stalks 8}
\begin{array}{rcl}
\chi((R\Phi_f[-1]\mathcal{F})_x)  &=& 
\chi(\mathcal{F}_x) 
-
\chi((R\Psi_f\mathcal{F})_x) \\
&=& 
\chi(\mathcal{F}_x) 
-
\chi(R\Gamma(M_{f,x};\mathcal{F})) \\
&=& 
\chi(\mathcal{F}_x) 
-
\mathop{\sum}\limits_{S_i} 
\chi_c(M_{f,x}\cap S_i) \cdot 
\chi(\mathcal{F}_{x_i}), 
\end{array}
\end{equation}
where the summation above is taken over all strata $S_i$ for which $S_i\cap M_{f,x} \ne \emptyset$ and for each such stratum $S_i$ we have chosen an element $x_i$ in this intersection. Here $\chi_c(-):=\chi(H^*_c(-;k))$ is the (additive) Euler characteristic of the compactly supported cohomology of a relatively compact subanalytic subset (like $M_{f,x}\cap S_i$).\\

We return now to the global context when $M$ is equipped with a $G$-action with finitely-many orbits $C$ and apply the machinery above to $\operatorname{For}\mathcal{F}$ for $\mathcal{F} \in D^b_{c,G}(M;k)$, and note that $\operatorname{For}\mathcal{F}$ is constructible with respect to the Whitney stratification $\mathscr{S}$
of $M$ give by (the connected components of) the $G$-orbits.
Given $x \in C \subset M$, we apply this to the normal slice $N$ to $C$ at $x$ with the induced stratification $\mathcal{S}$ given by (the connected components of) the intersections $C_i\cap N$. 
Since $\mathcal{F}$ is $G$-equivariant,  $\chi(\operatorname{For}\mathcal{F}_{x_i})$ is constant on the $G$-orbits.
Now, combining \eqref{equation: stalks 5} with \eqref{equation: stalks 8} now gives
\[
\begin{array}{rcl}
(-1)^{\dim C}\chi^\textrm{mic}_{C}(\mathcal{F})
&=&
\chi^\textrm{loc}_{C}(\mathcal{F})
-
\mathop{\sum}\limits_{C_i > C}\chi_c(M_{f,x}\cap C_i\cap N)\cdot
\chi^\textrm{loc}_{C_i}(\mathcal{F})
\end{array}
\]
Since
\[
(-1)^{\dim C}\chi_c(M_{f,x}\cap C_i\cap N)
=
\chi(\operatorname{NMD}((j_{C_i})_!\1_{C_i},df_x)),
\]
this expression only depends on $C$ and $C_i$, as discussed before, and now this proves \eqref{equation: stalks 7} with 
\[
c(C_i,C) \ceq  
\begin{cases}
(-1)^{\dim C +1} \chi_c(M_{f,x}\cap C_i\cap N) & C_i > C,\\
(-1)^{\dim C} & C_i = C.
\end{cases}
\]

\begin{remark}
All the discussion so far also works in the complex analytic context for a complex Lie group $G$ with finitely many connected components,
so that the equivariant derived categories of 
Bernstein-Lunts \cite{BL94}*{Part I, Section 3}
and their functors are available. 
For the Riemann-Hilbert correspondence one uses then of course analytic D-modules.
In the complex analytic context we explicitly assume that $G$ acts holomorphically on a complex manifold $M$ with finitely many locally closed orbits $C$, whose closure $\overline{C}$ is complex analytic with boundary
$\partial{C}= \overline{C}\backslash C$ a union of lower dimensional orbits.
\end{remark}

\subsection{A fixed-point formula for equivariant vanishing cycles}\label{subsection: FPF for Ev}

The purpose of this appendix is to prove the following result, which is specialized to our framework in Proposition~\ref{theorem: FPF} and used in the proof of Theorem~\ref{theorem: Lift of ABV}.

\begin{theorem}\label{theorem: Schurmann}
Let $M$ be a smooth complex algebraic variety equipped with the action of a diagonalizable group $G$. 
For $\sigma\in G$ (automatically semisimple), let $f: M\to \BC$ be $\sigma$-invariant morphism of complex varieties: $f =f\circ \sigma$.
Then, for every $\mathcal{F}\in D^b_{c,G}(M;k)$,
\begin{equation}\label{trace-van-alg}
\trace\left(\sigma, (\Rphi{f}[-1]{ \mathcal{F}})_x \right)
=
\trace\left(\sigma, (\Rphi{f\vert_{M^\sigma}}[-1]{ \mathcal{F}\vert_{M^\sigma}})_x \right),
\end{equation}
for all $x\in M^\sigma$.
\end{theorem}

\begin{proof}
Consider the Zariski closure $D$ of the subgroup $\langle \sigma \rangle$ generated by $\sigma$ in $G$.
Then $D$ is also a diagonalizable group and has cyclic component group generated by the image of $\sigma$ in the component group $D/D^0$.
In fact, the projection $\langle \sigma \rangle \to D/D^0$ is surjective since $\langle \sigma \rangle \cap D^i \ne \emptyset$ for each connected component $D^i$ of $D$, otherwise $D$ would not be the Zariski closure of $\langle \sigma \rangle$ in $G$.
Note that $f$ is invariant under the group $D$, since $f$ is invariant under $\sigma$.
Using a forgetful (change of group) functor, we now assume $D=G$. 

Then the stalks appearing above lie in $D^b_{c,G}(\{x\};k)$ and thus they have finite-dimensional cohomology.
Accordingly, trace means the alternating sum of the traces of $\sigma$ acting on cohomology groups. 
Since the $G$-action on $\{ x\}$ is trivial, this trace depends only on the image of $\sigma$ in the component group $G/G^0$.

Our method of proof uses a lifting of Diagram~\eqref{diagram: stalks} from $D^b_{c}(M;k)$ to $D^b_{c,G}(M;k)$. 
In order to do this, we only need to find the germ of a real analytic distance function $r$ as before, which in addition is now $G$-invariant. 
In that case the subspaces $T_{f,x}$ and $M_{f,x}$ are $G$-invariant and the arrows lift canonically to the $G$-equivariant derived category, since in the proof of Proposition~\ref{proposition: stalks} in \cite{Sch03}*{Lemma 1.1.1, p. 27} only uses the standard six functor formalism.
That even the equivariant vertical arrows in Diagram~\eqref{diagram: stalks} are isomorphisms in the equivariant derived category can be checked after forgetting from $D^b_{c,G}(\{x\};k)$ to $D^b_{c}(\{x\};k)$. 
In fact, such an equivariant morphism is an isomorphism if and only if its cone is $0 \in D^b_{c,G}(\{x\};k)$, but since we have a bounded non-degenerate t-structure commuting with forget, either the perverse of standard one, this can be checked after forgetting the $G$-action.

 Suppose first that $G$ is finite.
Consider the quotient map $p: (M,x)\to (M/G,x)$ with $M/G$ an algebraic variety (since $G$ is finite). Taking a closed embedding $(M/G,x)\hookrightarrow (\CC^n,0)$, and using the Euclidean distance function $d: (\CC^n,0) \to ([0,\infty),0)$, we choose $r = d\circ p$. Here we use $p^{-1}(\{x\}) = \{ x\}$, so $r$ is real analytic and $r^{-1}(0) = \{ x\}$, hence $r$ is a $G$-invariant distance function.
We can now conclude 
\[
\trace(\sigma, (R\Phi_f[-1]\mathcal{F})_x)
=
\trace(\sigma, R\Gamma(T_{f,x}, M_{f,x};\mathcal{F})),
\]
where now $T_{f,x}$ and $M_{f,x}$ are compact subanalytic subsets of $M$.
Applying the Lefschetz formula for global cohomology \cite{Sch03}*{(3.24), p.196}, we get
\[
\trace(\sigma, R\Gamma(T_{f,x}, M_{f,x};\mathcal{F}))
=
\trace(\sigma, R\Gamma(T_{f,x}^\sigma, M_{f,x}^\sigma;\mathcal{F}\vert_{M^\sigma})).
\]
Finally, 
\[
\trace(\sigma, R\Gamma(T_{f,x}^\sigma, M_{f,x}^\sigma;\mathcal{F}\vert_{M^\sigma}))
=
\trace\left(\sigma, (\Rphi{f\vert_{M^\sigma}}[-1]{ \mathcal{F}\vert_{M^\sigma}})_x \right),
\]
since we can use the germ of restriction of $r$ to the submanifold $M^\sigma$ with the same $\epsilon$ and $\delta$ as in the ambient context for $M$ (since they can be chosen arbitrarily small).
This proves the theorem in the case that $G$ is finite cyclic.

Finding a $G$-invariant distance function $r$ in this simple way does not work unless $G$ is finite. 
Now we prove the general case of the theorem, by a different method, by switching to a maximal compact real analytic Lie subgroup $K$ of the reductive group $G$. 

Then $K$ is Zariski dense in $G$ so that $M^\sigma = M^G = M^K$ and the component groups coincide: $A = G/G^0 = K/K^0$.
In our case, the abelian reductive group $G=D$ is a product $G = T\times A$, for cyclic $A$ and torus $T$, so that $K$ is a product $K = T_c\times A$, where $T_c$ compact torus. 

Since we cannot assume $\sigma$ lies in the compact subgroup $K$, we switch to an element $\sigma'\in K$ such that $\sigma/G^0 = \sigma'/G^0$ and $\sigma'$ is topologically dense in $K$, {\it i.e.}, the topological closure of $\langle \sigma'\rangle$ is $K$. 
Then $M^{\sigma'} = M^K = M^{\sigma}$.
Recalling that the traces in the theorem depend only on the component of $\sigma$, we can now replace $\sigma$ with $\sigma'$ and pass $\mathcal{F}$ from $D^b_{c,G}(M;k)$ to $D^b_{c,K}(M;k)$, where now we work with subanalytic constructible sheaves.
Now we need to produce the germ of a real analytic $K$-invariant function $r$ with $r^{-1}(0)= \{ x\}$.
For this we use a local linearization theorem of Cartan and Bochner (see \cite{Huc90}*{Section 4(a), pp. 154-155} or \cite{Akh95}*{Theorem, p.36}), which gives us a local $K$-equivariant holomorphic isomorphism $p: (M,x)\to (T_xM,0)$,  identifying $p: (M^K,x) \simeq (T_xM^K,0)$. Since $M^K = M^{\sigma'}=M^{\sigma}$, we can conclude that $M^{\sigma}$ is a complex submanifold of $M$
with $(TM)^\sigma = T(M^\sigma)$ for the induced action of $\sigma$ on the tangent bundle $TM$.
Arguing as above, we choose a $K$-invariant scalar product on $T_xM$ by averaging (using the normalized Haar measure on $K$) giving us a real analytic $K$-invariant Euclidean distance function $d$ on $T_xM$. 
Then $r = d\circ p$ is the germ of a $K$-invariant real analytic function as we want. 

As before, we can now lift Diagram~\ref{diagram: stalks} from $D^b_{c}(M;k)$ to $D^b_{c,K}(M;k)$. 
Now we may argue as above, but only using the standard $t$-structure, to reduce to the statement to the corresponding fixed-point formula
\[
\trace(\sigma', R\Gamma(T_{f,x}, M_{f,x};\mathcal{F}))
=
\trace(\sigma', R\Gamma(T_{f,x}^{\sigma'}, M_{f,x}^{\sigma'};\mathcal{F}\vert_{M^{\sigma'}})),
\]
where $T_{f,x}$ and $M_{f,x}$ are compact subanalytic subsets of $M$, for $\mathcal{F}\in D^b_{c,K}(M;k)$ $K$-equivariant subanalytic constructible. 
The theorem now follows from \cite{Sch03}*{Corollary 3.3.2 and (3.24), p. 196}.
\end{proof}

\begin{remark}
Theorem~\ref{theorem: Schurmann} is an analogue of the trace formula appearing in \cite{ABV92}*{Theorem 25.8} in the language of vanishing cycles, where \cite{ABV92}*{Theorem 25.8} use a description of microlocal traces in the language of stratified Morse theory with similar adapted pairs of spaces (just using the real part $re(f)$ of $f$ as in \cite{Sch03}*{Proposition 5.0.3, p. 281}),
to go to a trace formula for global cohomology. 
Note, however, that the automorphism $\sigma$ in \cite{ABV92}*{Theorem 25.8} is of finite order, thus semisimple, but Theorem~\ref{theorem: Schurmann} also applies in the more general case that $\sigma$ is a semisimple element in a linear complex algebraic group $G$ (which is always contained in a diagonalizable subgroup of $G$).
\end{remark}

\begin{remark}\label{smooth-fixed-points}
 Let $M$ be a smooth complex algebraic variety equipped with the action of a diagonalizable group $G$. Then we have shown in the proof above, that for
$\sigma\in G$ the fixed point set 
$M^{\sigma}$ is a complex submanifold of $M$
with $(TM)^\sigma = T(M^\sigma)$ for the induced action of $\sigma$ on the tangent bundle $TM$. 
This extends \cite{ABV92}*{Lemma 25.3} from a finite order automorphism to a semisimple element $\sigma$  in a linear complex algebraic group acting on $M$ (which is always contained in a diagonalizable subgroup).
Similarly, \cite{ABV92}*{Lemma 25.5 and Theorem 25.6} then also apply to this more general context.
\end{remark}

As a corollary of the method of proof of Theorem~\ref{theorem: Schurmann}, we get the following result.

\begin{theorem}\label{theorem: Schurmann analytic}
Let $M$ be a complex analytic manifold equipped with a holomorphic action of a compact real abelian Lie group $K = T_c\times A$, with $T_c$ a  compact real torus and $A$ a finite abelian component group. 
For $\sigma\in K$, let $f: M\to \BC$ be a $\sigma$-invariant holomorphic function.
Then, for every $\mathcal{F}\in D^b_{c,K}(M;k)$,
\begin{equation}\label{trace-van-an}
\trace\left(\sigma, (\Rphi{f}[-1]{ \mathcal{F}})_x \right)
=
\trace\left(\sigma, (\Rphi{f\vert_{M^\sigma}}[-1]{ \mathcal{F}\vert_{M^\sigma}})_x \right),
\end{equation}
for all $x\in M^\sigma$.
\end{theorem}

\section{Global vanishing cycles}\label{section: global vanishing cycles}

In this appendix we make some improvements to the functors $\NEvs_C$ and $\Ev_C$ introduced in \cite{CFMMX22}*{Chapter 7}
(in a more restricted context), largely by passing to the analytic site.

\subsection{Pyasetskii orbit duality and global vanishing cycles}\label{subsection: PDGV}

We now move to the geometric situation showing up in this project, assuming that $H$ is a complex (linear) algebraic group 
acting linearly on a (finite dimensional) complex vector space $V$ with finitely many orbits $C$. Then there is a bijection between the finite set $\{C'\}$ of $H^0$-orbits in $V$ and the finite set $\{\overline{T^*_{C'} M}\}$  of $\CC^*$-conic and $H^0$-invariant irreducible Lagrangian subvarieties of the cotangent bundle $T^*V=V\times V^*$
(see, e.g. \cite{Tev05}*{Theorem 1.10}).

We also set 
\[
\Lambda := \bigsqcup_{C'} T^*_{C'}V = \bigsqcup_{C} T^*_{C}V \subset T^*V = V\times V^*,
\]
and refer to this conormal variety in this context also as the {\it commuting variety}.

Note that these constructions are completely symmetric in $V$ and $V^*$ under the natural identification
$$T^*V = V\times V^* \simeq V^*\times V=T^*V^* \:,$$
preserving the canonical symplectic structure of the cotangent bundles up to a sign.

Now we privilege $V$ as a base parameter and $V^*$ as the cotangent direction. 
For each $H$-orbit $C\subseteq V$, the
conormal bundle  $\,_C\Lambda \ceq T^*_{C}V$ of $C$ is  the inverse image of the projection $\Lambda \to V$ above $C$.

The regular part of the conormal to $C$ is given by
${}_C\Lambda^\reg \ceq \,_C\Lambda \setminus {}_C\Lambda^\textrm{sing}$ for 
\[
{}_C\Lambda^\textrm{sing} \ceq \,_C\Lambda \cap \left( \operatorname{\bigcup}\limits_{C_i\ne C} (\,_{C}\Lambda \cap \overline{\,_{C_i}\Lambda})\right).
\]
Note that  $\,_C\Lambda^\textrm{sing}$ is closed in $\,_C\Lambda$ and $\,_C\Lambda^\reg$ is open and dense in $\,_C\Lambda$.
Similar definitions apply to the conormal bundle $\,_{C'}\Lambda \ceq T^*_{C'}V$ of a $H^0$-orbit $C'$.
Then the connected components of ${}_C\Lambda^\reg$ are given by the ${}_{C'}\Lambda^\reg$ for all $H^0$-orbits $C'\subset C$.

Instead of priviledging $V$ and the projection $p: V\times V^* \to V$ we could switch the role of $V$ with $V^*$ and use the projection $q : V\times V^* \to V^*$ and  the isomorphism $$a : V\times V^* \to V^*\times V$$ given by $(x,y) \mapsto (y,x)$ 
as above so that we may view $V^*$ as the base.
Here we identify $(V^*)^* =V$.
Accordingly, for an $H$-orbit $D$ (or $H^0$-orbit $D'$) in $V^*$, we set 
$\Lambda_{D}= \Lambda\cap q^{-1}(D)= T^*_DV^*$ (or $\Lambda_{D'}= \Lambda\cap q^{-1}(D')= T^*_{D'}V^*$)
for the projection $q : V\times V^*\simeq T^*V^* \to V^*$.
The Pyasetskii orbit duality \cite{Pja75} is defined by picking an irreducible component $\Lambda_i$ of the commuting variety $\Lambda$, which then can be written uniquely as
\[
\overline{\,_{C'}\Lambda}
=
\Lambda_i
= 
\overline{\Lambda_{D'}}
\]
for $H^0$-orbits $C'\subseteq V$ and $D'\subseteq V^*$, in which case we use the notation $C'\mapsto C'^* :=D'$.
For example $\{0\}=C'^*\subset V^*$ for the open and dense $H^0$-orbit $C'\subset V$.
Note that this orbit duality also implies that each $H^0$-orbit $C'$ in $V$ (or $D'$ in $V^*$) is
$\CC^*$-conic, which then is also true for the corresponding $H$-orbits (and then the Pyasetskii orbit duality can also be proved and formulated in terms of projective duality as in \cite{Tev05}*{Theorem 2.1 and Corollary 2.4}).
In particular $\Lambda\subset f^{-1}(0)$ for the canonical pairing $f: V\times V^* \to \CC$.
In fact, $f(v,w)=0$ for $v\in T_xC'$ and $(x,w)\in T^*_{C'}V$, with $x\in T_xC'\subset V$ for $x\in C'$, since this 
$H^0$-orbit is $\CC^*$-conic.

Moreover, this orbit duality commutes with the corresponding $G/G^0$-permutation actions inducing also a 
Pyasetskii orbit duality $C\mapsto C^* :=D$ for the possible disconnected (linear) algebraic group $H$, via 
\[
\overline{\,_{C}\Lambda}
=
G/G^0(\Lambda_i)
= 
\overline{\Lambda_{D}}
\]
for a $G/G^0$-orbit of an irreducible component $\Lambda_i$ of the commuting variety $\Lambda$.

Note that $\,_C\Lambda^\reg$ is identified under the identification $a: V\times V^* \to V^* \times V$ with the regular part $\Lambda_{C^*}^\reg$ of the conormal to $C^*$; in particular,  $\,_{C}\Lambda^\reg = \,\Lambda_{C^*}^\reg \subset C\times C^*$. Finally, the {\it commuting variety} $\Lambda$ is also the reduced zero fibre $\Lambda=\mu^{-1}(0)^{red}$
of the corresponding {\it moment map}\IN{$\mu$, moment map}
\[
\mu  : T^*V = V\times V^* \to \mathfrak{h}^*
=\textrm{Lie}(H)^*=\textrm{Lie}(H^0)^*
\]
as in \cite{Shu25} and \cite{Tev05}*{Remark 2.5}).\\

Fitting with the perspective of this paper, the preferred base from now on for $\Lambda$ is $V$, as reflected in the following definition.
\begin{definition}\label{definition: Appendix ev}  
For an $H$-orbit $C$ in $V$, consider the functor $$\ev{C} : D^b_{c,H}(V;k) \to D^b_{c,H}(f_{C^*}^{-1}(0);k)$$ defined by \IN{$\ev{C}$, vanishing cycle functor}
\[
\ev{C}\mathcal{F} \ceq  R\Phi_{f_{C^*}}[-1]p_{C^*}^*\mathcal{F},
\]
where $p_{C^*} : V\times C^* \to V$ is the projection and $f_{C^*}$ is the restriction of the canonical 
$H$-invariant pairing $f: V\times V^* \to \CC$ to $V\times C^*$.
\end{definition}

The following result will show that $\ev{C}$ may be thought of as a parametrized family of vanishing cycles $R\Phi_{f_y}$ for $y\in C^*$.
As the proof will show, here it is important that the parameter space $C^*$ is an $H$-orbit.

\begin{theorem}\label{theorem: family ev}
    For every $H$-orbit $C\subseteq V$, the equivariant base change arrow
    \[
    j_{(x,y)}^* R\Phi_{f_{C^*}} p_{C^*}^*\mathcal{F} \to j_{x}^* R\Phi_{f_y}\mathcal{F}
    \]
    is an $Z_{H}(x,y)$-equivariant isomorphism, for all $(x,y)\in f_{C^*}^{-1}(0)$, for all $\mathcal{F}\in D^b_{c,H}(V;k)$.
\end{theorem}

\begin{proof}
Note that $\operatorname{For} p_{C^*}^*\mathcal{F}$, as above, is automatically constructible with respect to the finite $H\times H$-orbit stratification of $V\times C^*$.
By $H\times H$-invariance, this is automatically a Whitney stratification.
But this is not adapted to $f_{C^*}^{-1}(0)$.
For this reason, we introduce new $H$-invariant strata in $f_{C^*}^{-1}(0)$ so that we get an $H$-invariant Whitney stratification $\mathscr{S}$ of $V\times C^*$ adapted to $f_{C^*}^{-1}(0)$, so that $f_{C^*}^{-1}(0)$ is a union of strata from $\mathscr{S}$, giving a Whitney stratification $\mathscr{S}_0$ of $f_{C^*}^{-1}(0)$. 
This implies, by \cite{MS22}*{Corollary 10.4.9, p. 753 and Corollary 10.4.12, p. 756} and \cite{Sch03}*{Lemma 4.2.1, p. 247-248}, the following refinement of Definition~\ref{definition: Appendix ev}
\begin{equation}\label{equation: A S-c}
\ev{C} = R\Phi_{f_{C^*}}[-1]p_{C^*}^* :  D^b_{c,H}(V;k) \to D^b_{\mathscr{S}_0-c,H}(f_{C^*}^{-1}(0);k),
\end{equation}
where the index on the target denotes constructability (after applying forget) with respect to the given $H$-invariant Whitney stratification $\mathscr{S}_0$ of $f_{C^*}^{-1}(0)$. 
Recall that this means that all cohomology sheaves of these vanishing cycles are locally constant on these strata.
The $H$-invariance is this stratification gives us the crucial property that the projection $q: S \to C^*$ is a submersion, for all $S\in \mathscr{S}$; {\it i.e.}, that $dq: T_z S \to T_{q(z)} C^*$ is surjective, for all $z\in S$.
For this, for each $z\in S$, consider the $H$-orbit $H\cdot z\subset S$ and note that $dq: T_{z}(H\cdot z) \to T_{q(z)} C^*$ is already surjective since $C^*$ is an $H$-orbit and $Z_{H}(z) \subseteq Z_{H}(q(z)) \subseteq H$.
This implies that for a given parameter $y\in C^*$, the  submanifold $q^{-1}(y) = V\times \{ y\}$ is transversal to all strata $S\in \mathscr{S}$ in the ambient manifold $V\times C^*$.
Now apply the  base change isomorphism for vanishing cycles 
\cite{MS22}*{Proposition 10.4.19 (3.), p. 758} and \cite{Sch03}*{Lemma 4.3.4, p. 265} in the context of intersecting with the transversal submanifold.
This finishes the proof because, as in the context of local vanishing cycles, it is enough to prove that the base change arrow is an isomorphism non-equivariantly, {\it i.e.}, after applying the forgetfull functor $D^b_{c,H}(V;k) \to D^b_{c}(V;k)$. 
\end{proof}

\begin{remark}
Compared to the equivariant induction argument from \cite{CFMMX22}*{Proposition 7.8},
the advantage of the argument above is that it applies non-equivariantly more generally to $H$-monodromic  instead of $H$-equivariant sheaf complexes, i.e. to sheaf complexes whose cohomology sheaves are locally constant along the  $H$-orbits.
\end{remark}

\begin{remark}\label{remark: bad}
    Let us give some hints for how we can find an $H$-invariant Whitney stratification adapted to $f_{C^*}^{-1}(0)$, {\it i.e.}, so that $f_{C^*}^{-1}(0)$ is a union of $H$-invariant strata (as in \cite{Tes82}*{Section 3 p. 478}.
    Such a stratification is constructed by induction, by first removing a possible singular locus and starting  in top dimension, comparing a pair of strata $S,S'$ of pure dimensional complex algebraic strata, with $S\subset \partial S' = \overline{S'}\setminus S'$.
    Then the bad subset $B(S,S')\subset S$ where the Whitney condition fails, is a constructible subset of lower dimension by \cite{Tes82}*{Théorème 1.2, p. 455, Proposition 6.1.2, p. 452}
    whose closure has to be removed in the induction step.
    In case $S$ and $S'$ are both $H$-invariant, by the characterization given in {\it loc. sit.}, also the bad set $B(S,S')$ is $H$-invariant, since the construction is canonical.
\end{remark}

Note that Theorem~\ref{theorem: family ev} can be expressed in this form:
\[
j_{(x,y)}^*\, \ev{C}\mathcal{F} = j_{x}^*\, R\Phi_{f_{y}}[-1]\mathcal{F},
\]
for all $(x,y)\in f_{C^*}^{-1}(0)$.

\begin{corollary}\label{corollary: conormal support}
    For every $\mathcal{F}\in D^b_{c,H}(V;k)$, the support of $\ev{C}\mathcal{F}$ is contained in 
    $$\Lambda \cap (V\times C^*) =\Lambda_{C^*} \subset f_{C^*}^{-1}(0)\:, $$
    with $T^*_{C^*}V \simeq \Lambda_{C^*}$.
\end{corollary}

\begin{proof}
    By Theorem~\ref{theorem: family ev}, $(\ev{C}\mathcal{F})_{(x,y)} = (R\Phi_{f_y}[-1]\mathcal{F})_x$, for all $(x,y)\in f_{C^*}^{-1}(0)$.
    Now take a point $(x,y)\in f_{C^*}^{-1}(0)$ which is not in the closed subset  $f_{C^*}^{-1}(0)\cap \Lambda$.
    The vanishing of this last stalk, working non-equivariantly, follows from \cite{MS22}*{Proposition 10.4.16, p. 757} and \cite{Sch03}*{Remark 4.2.4, p. 248}, since $x$ is not a stratified critical point of $f_y$, with respect to the orbit stratification {\it i.e.}, $df_y(x)$ is not contained in $\Lambda$.   
\end{proof}

\begin{theorem}\label{theorem: refinement}
    We now make an important refinement of the construction above about a $H$-invariant Whitney stratification of $V\times C^*$ adapted to $f_{C^*}^{-1}(0)$, having in addition the following properties: it is also adapted to $\Lambda_{C^*}$ and further $\,_{C}\Lambda^\reg =\Lambda_{C^*}^\reg$ is a stratum of the stratification of $\Lambda_{C^*}$.
    In particular, 
\begin{equation} \label{evc-refined}
j^*_{\Lambda_{C^*}} \ev{C} = R\Phi_{f_{C^*}}[-1]p_{C^*}^* :  D^b_{c,H}(V;k) \to D^b_{\mathscr{S}_0\cap \Lambda_{C^*}-c,H}(\Lambda_{C^*};k)
\end{equation}
and
\begin{equation} \label{evc-refined-perv}
j^*_{\Lambda_{C^*}} \ev{C}[\dim C^*] : \Per_{H}(V;k) \to \Per_{\mathscr{S}_0\cap \Lambda_{C^*}-c,H}(\Lambda_{C^*};k).
\end{equation}
Composing with the open restriction $j_{\Lambda^\reg_{C^*}} : \,_{C}\Lambda^\reg ={\Lambda^\reg_{C^*}} \hookrightarrow \Lambda_{C^*}$ now gives
\begin{equation} \label{evcs-refined}
\evs{C} \ceq j^*_{\Lambda^\reg_{C^*}} \ev{C}[-\codim C^*] = R\Phi_{f_{C^*}}[-1]p_{C^*}^* :  D^b_{c,H}(V;k) \to D^b_{\Lambda^\reg_{C^*}-c,H}(\Lambda^\reg_{C^*};k)
\end{equation}
and 
\begin{equation} \label{evcs-refined-perv}
\evs{C} = j^*_{\Lambda_{C^*}^\reg} \ev{C}[-\codim C^*] : \Per_{H}(V;k) \to \Loc_{H}(\Lambda^\reg_{C^*};k)\:,
\end{equation}
where $\Loc_{H}(\Lambda^\reg_{C^*};k)$ is the abelian category of $H$-equivariant local systems (of $k$-vector spaces) on 
$\,_{C}\Lambda^\reg = \,\Lambda^\reg_{C^*}$.
Here the induced stratification of $\Lambda^\reg_{C^*}$ just consists of one smooth stratum so that the indicated 
$\Lambda^\reg_{C^*}$-constructibility just means that all underlying cohomology sheaves are locally constant (with finite dimensional stalks) on $\,_{C}\Lambda^\reg = \,\Lambda^\reg_{C^*}$.
\end{theorem}

\begin{proof}
The proof that such a stratification exists, which is quite technical, will be given in Lemma~\ref{lemma: Whitney}. 
The next claim is a consequence of Equation~\eqref{equation: A S-c} together with Corollary~\ref{corollary: conormal support}.
We have seen $\ev{C} = R\Phi_{f_{C^*}}[-1]p_{C^*}^* : D^b_{c,H}(V;k) \to D^b_{c,H}(\Lambda_{C^*};k)$.
    Smooth pullback by $p_{C^*}^*[\dim C^*]$ preserves perversity, by \cite{MS22}*{Corollary 10.2.25, p. 696} and \cite{Sch03}*{Lemma 6.0.3 (3), p. 386}.
    Vanishing cycles $R\Phi_{f_{C^*}}[-1]$ also preserves perversity, by \cite{MS22}*{Theorem 10.4.22, p. 759} and \cite{Sch03}*{Theorem 6.0.2, p. 402}. Taken together, this proves the second claim.
Since $\Lambda_{C^*}^\reg \hookrightarrow \Lambda_{C^*}$ is the inclusion of an open stratum, the third claim now follows immediately. 
Then constructibility with respect to this open stratum $\Lambda_{C^*}^\reg$ just means that all cohomology sheaves are locally constant on this open stratum of dimension $\dim V$.
Then the last claim follows immediately.
\end{proof}

\begin{remark}
    The last part of Theorem~\ref{theorem: refinement} has the important consequence that $T^*_C(V_\lambda)_\gen = \,_{C}\Lambda^\reg = \Lambda^\reg_{C^*}$, where $T^*_C(V_\lambda)_\gen$ was defined in \cite{CFMMX22}*{Section 7.9}.
\end{remark}

\begin{lemma}\label{lemma: Whitney}
There is an $H$-invariant Whitney stratification of $V\times C^*$ that is adapted to $f_{C^*}^{-1}(0)$, adapted to $\Lambda_{C^*}$ and for which $\,_{C}\Lambda^\reg =\Lambda_{C^*}^\reg$ is a stratum.
\end{lemma}

\begin{proof}
    Here we are referring to the Whitney b-regularity condition, as in \cite{Sch03}*{Definition 4.1.2, p. 224}, for a complex algebraic stratification $\mathscr{S}$ of $V\times C^*\subset V\times V^*$, with locally closed, pure dimensional complex algebraic submanifolds $S\subset V\times C^*$ called the strata $S$ of $\mathscr{S}$.
    A pair of strata $S$, $S'$ with $S\subseteq \partial S' = \overline{S'}\setminus S'$ (so that $\dim S< \dim S'$) satisfies the Whitney-b condition over $S$ if the following property holds.
    \begin{itemize}
        \item Let $v_n\in S$ and $w_n\in S'$ be two sequences converging to the same limit point $v\in S$.
        \item Assume the tangent spaces $T_{w_n} = T_{w_n}(S')$ converge to some limiting vector space $\tau \subset V\times V^*$ and the secant lines $\ell_n = \CC(v_n-w_n) \subset V\times V^*$ convergent to some limiting line $\ell\subset V\times V^*$.
        \item These limiting lines and vector spaces are taken in the corresponding compact Grassmanian for $V\times V^*$.
    \end{itemize}
Under these assumptions, $S'$ is Whitney b-regular over $S$ if 
\begin{equation}
\lim_{n\to \infty} \ell_n = \ell \subseteq \tau = \lim_{n\to \infty} T_{w_n}
\end{equation}
    for all such sequences.

    Let's start with some conditions where we get the Whitney stratification for free.
\begin{enumerate}
    \item 
    $S$ is an algebraic submanifold of an ambient algebraic manifold $M$ (both of pure dimension) and we take $S' \ceq M\setminus S$, then $S'$ is Whitney b-regular over $S$. This will be applied later on for $S= T^*_{C}(V) \cap (V\times C^*) = \Lambda\cap (C\times C^*)$ in $M = C\times C^*$.
    \item The stratification of $V\times C^*$ given by $H\times H$-orbits $D\times C^*$ is a Whitney stratification, because the Whitney b-condition holds generically as mentioned in Remark~\ref{remark: bad}. By $H\times H$-invariance, it holds therefore everywhere.
    \item The Whitney property is preserved if we restrict to the open subset $\{ f_{C^*} \neq 0\}$ in $V\times C^*$. Also note that all critical points of $f_{C^*}\vert_{D\times C^*}$ are contained in $\{ f_{C^*} = 0\}$ since $D\times C^*$ is $\CC^*$-conic and $f_{C^*}$ is homogeneous.
    This critical locus of $f_{C^*}\vert_{D\times C^*}$ is $\Lambda \cap (D\times C^*) = T^*_DV \cap (V\times C^*)$
    (as can be easily shown similarly to \cite{CFMMX22}*{Proof of Lemma 7.15}).
\item    Now we get an induced Whitney stratification of $V\times C^*\setminus (\Lambda\cap (V\times C^*))$ with strata in $\{ f_{C^*} \neq 0\}$ as before, and adding the new strata $(D\times C^* )\cap \left(\{ f_{C^*} = 0 \} \setminus (\Lambda\cap (V\times C^*))\right)$ in $\{ f_{C^*} = 0 \} \setminus (\Lambda\cap (V\times C^*))$. Since we avoided all critical points, this is Whitney b-regular by \cite{Sch03}*{Example 4.3.3, p. 257}. All these strata are $H$-invariant by construction.
\end{enumerate}
What remains is to stratify $\Lambda \cap (V\times C^*)$ in a suitable way.
Note that $\Lambda \cap (V\times C^*) $ is of dimension $\dim V$ with the union of the top-dimensional irreducible components given by $\overline{\Lambda^\reg_{C^*}} \cap (V\times C^*)$.
Here we introduce as a new stratum $N= \Lambda^\reg_{C^*}\subset C\times C^*$.
Note that this is an open subset of the critical locus of $f_{C^*}\vert_{C\times C^*}$ with complement of lower dimension. 
We need to show that all other strata introduced before are Whitney b-regular over $N$.
If this is done, we can inductively refine $\Lambda_{C^*} \setminus N$ as in the proof of Theorem~\ref{theorem: family ev}, so it remains to show that all the strata $S'$ as above are Whitney b-regular over $N$. 
\begin{enumerate}
    \item 
    The stratum $(C\times C^*)\setminus (\Lambda\cap (V\times C^*))$ is Whitney b-regular over $N$ as in Item (1) above.
    \item 
    A stratum $(D\times C^*) \cap \{f_{C^*} \neq 0 \}$ as in Item (3) above is Whitney b-regular over $N$ since $D\times C^*$
    already was Whitney b-regular over $C\times C^*$.
\end{enumerate}
 
    The final case to show is that a stratum
    $S' = (D\times C^* )\cap \left(\{ f_{C^*} = 0 \} \setminus (\Lambda\cap (V\times C^*))\right)$ as in Item (4) above is Whitney b-regular over $S= N$.
Now suppose \begin{itemize}
        \item Let $v_n\in S$ and $w_n\in S'$ are two sequences converging to the same limit point $v\in S$.
        \item Assume the tangent spaces $T_{w_n} = T_{w_n}(S')$ converge to some limiting vector space $\tau \subset V\times V^*$ and the secant lines $\ell_n = \CC(v_n-w_n) \subset V\times V^*$ converget to some limiting line $\ell\subset V\times V^*$.
    \end{itemize}
Going over to subsequences (since a Grassmanian is compact) we can assume $T'_{w_n} \ceq T_{w_n}(D\times C^*)$ is also convergent to a limiting vector space $\tau'$.
By our construction $T_{w_m} = T'_{w_n}\cap \ker d_{w_n}f_{C^*}$.
In the limit we get $\tau \subseteq \tau'\cap \ker d_{v}f_{C^*}$.
This inclusion is indeed an equality since $d_vf_{C^*}$ is not vanishing on $\tau'$ because it is a regular covector.
We have chosen sequences $t_n\in \CC^*$ so that $t_n(v_n-w_n)$ converges to a vector $0\ne z \in V\times V^*$.
Again we can also assume, passing to a subsequence, that $t_n v_n$ converges to $0\ne v'$ and $t_n w_n$ converges to $0\ne w'$ with $z = v'-w'$.
By the Whitney b-regularity of the original product stratification with strata $D\times C^*$, we know that the limiting line $\ell$ generated by $z$ is contained in $\tau'$, {\it i.e.}, $\ell\subseteq \tau'$.
This Whitney-b property also implies the Whitney-a property; see \cite{Sch03}*{p. 226} for example.
So $T_{v} \subseteq \tau'$.
But $C\times C^*$ is $\CC^*$-conic so that $t_n v_n \in T_{v_n}$.
In the limit we get $v'\in T_{v}$ so that $d_vf_{C^*}(v') =0$ since $v$ is a critical point of $f_{C^*}\vert_{C\times C^*}$.
Since $f_{C^*}$ is homogeneous, $(D\times C^*)\cap \{ f_{C^*} =0 \}$ is $\CC^*$-conic so that the sequence $t_n w_n\in T_{w_n} = T'_{w_n} \cap \ker d_{w_n}f_{C^*}$. 
In the limit we get $w'\in \tau' \cap \ker (df_{C^*})_{v} = \tau$.
All together we get $z= v'-w'\in \tau' \cap \ker d_{v}f_{C^*} = \tau$.
This means that the limiting line $\ell$ is contained in $\tau$, which finishes the proof of the Whitney-b property.
\end{proof}

Let us finish this subsection with the following {\it multiplicativity} alias {\it Thom-Sebastiani property} of 
the micro-local vanishing cycle functors  $\ev{C}$ and $\evs{C}$.
Assume $H_i$  is a complex linear algebraic group 
acting linearly on a (finite dimensional) complex vector space $V_i$ with finitely many orbits $C_i$ for $i=1,2$.
Then the same is true for $H=H_1\times H_2$ acting on
$V=V_1 \times V_2$ with finitely many orbits $C=C_1\times C_2$, with
$$(C_1\times C_2)^*=C_1^*\times C_2^*$$
for the corresponding dual orbits. Similarly 
$$T^*V\simeq TV^*_1\times T^*V_2 = (V_1\times V_1^*)\times (V_2\times V_2^*)\:,$$
with $\,_C\Lambda \ceq T^*_{C}V\simeq \,_{C_1}\Lambda\times \,_{C_2}\Lambda$ and
$\,_C\Lambda^\reg \ceq (T^*_{C}V)^\reg \simeq \,_{C_1}\Lambda^\reg\times \,_{C_2}\Lambda^\reg$
(and similarly for $\Lambda_{C^*}$ and $\Lambda_{C^*}^\reg$).
Finally the canonical evaluation pairing $f:V\times V^*\to \CC$ is given by the sum function 
$$f=f_1\circ pr_1 + f_2\circ pr_2 \:,$$
with $pr_i: V\times V^*\to V_i\times V_i^*$ and $f_i: V_i\times V_i^*\to \CC$ the corresponding projection and evaluation pairing for $i=1,2$.
Then the (equivariant version of the) {\it Thom-Sebastiani-Theorem} (see \cite{Ma01} and 
\cite{Sch03}*{Corollary 1.3.4 on p.72}) implies for $\mathcal{F}_i\in D^b_{c,H_i}(V_i;k)$ for $i=1,2$:

\begin{equation}\label{TS-theorem}
  \ev{C}(\mathcal{F}_1\boxtimes \mathcal{F}_2)
  \simeq \ev{C_1}(\mathcal{F}_1) \boxtimes
  \ev{C_2}(\mathcal{F}_2) \in  D^b_{c,H}(\Lambda_{C^*};k) \:.
\end{equation}
Similarly
\begin{equation}\label{TS-theorem2}
  \evs{C}(\mathcal{F}_1\boxtimes \mathcal{F}_2)
  \simeq \evs{C_1}(\mathcal{F}_1) \boxtimes
  \evs{C_2}(\mathcal{F}_2) \in  D^b_{\Lambda^\reg_{C^*}-c,H}(\Lambda^\reg_{C^*};k) \:,
\end{equation}
with $\evs{C} \ceq j^*_{\Lambda_{C^*}^\reg} \ev{C}[-\codim C^*]$.

\subsection{Remarks on slicing and normalization} 
\label{subsection: normalization}
In this section we restrict our context further to  a linear representation $V$ of a complex linear algebraic group $H$, whose underlying representation of $H^0$ is coming in the following way from a graded Lie algebra.

Let $G$ be a connected and reductive complex linear algebraic group. 
Consider the $\Z$-grading $\mathfrak{g} = \mathop{\oplus}_{n\in \Z} \mathfrak{g}(n)$ on $\mathfrak{g} = \operatorname{Lie} G$ given by a cocharacter $\lambda : \CC^* \to G$. 
Let our $H^0$ also be the centralizer
$H^0 = G^\lambda$  of $\lambda$, which then also is a connected and reductive complex linear algebraic group.
Observe that $\operatorname{Lie} H^0 = \mathfrak{h} = \mathfrak{g}(0)$.
For $n$ non-zero, $\mathfrak{g}(n)$ is an instance of a prehomogeneous vector space  with finitely many orbits for the linear action of $H^0$ (see, e.g. \cite{Vin76}*{} or 
\cite{Lus95}*{Section 3}). We fix now a nondegenerate and $G$-invariant symmetric bilinear form on $\mathfrak{g}$ (e.g. a suitable extension of the
Killing form for $[\mathfrak{g},\mathfrak{g}]$ 
as in \cite{CFMMX22}*{Proof of Proposition 6.2}):
$$\langle \cdot , \cdot\rangle : \mathfrak{g}\times \mathfrak{g} \to \CC\:.$$
Then also $\langle \cdot , \cdot\rangle : \mathfrak{g}(n)\times \mathfrak{g}(-n) \to \CC$ is non-degenerate, allowing us to identify  $\mathfrak{g}(-n) \simeq \mathfrak{g}(n)^*$.
We assume now $V=\mathfrak{g}(n)$ for a fixed $n\ne 0$, with a given $H$-action  extending this $H^0$-action (with $H^0$ the identity component of $H$)
so that $V$ also has finitely many $H$-orbits as in the section before.
Moreover, the identification $V^*\simeq \mathfrak{g}(-n)$ induces then a dual $H$-action on $V^*\simeq \mathfrak{g}(-n)$ so that the $H$-invariant evaluation pairing $f: V\times V^*\to \CC$
is identified with $\langle \cdot , \cdot\rangle : \mathfrak{g}(n)\times \mathfrak{g}(-n) \to \CC$.
Using the pairing $\langle \cdot , \cdot\rangle$, the moment map $\mu: T^*V=V\times V^* \to \mathfrak{h}^*$ 
is given by the Lie bracket for $\mathfrak{g}$ (see \cite{Tev05}*{Remark 2.7}):
\[
\mu = [\cdot , \cdot ] : \mathfrak{g}(n)\times \mathfrak{g}(-n) \to \mathfrak{g}(0)\simeq \mathfrak{g}(0)^*.
\]
Then the {\it commuting variety} is given by
\[
\Lambda = \mu^{-1}(0)^{red} = \{ (x,y) \in V\times V^* \tq [x,y] =0 \},
\]
together with $\,_{C}\Lambda= T^*_CV = \{ (x,y) \in C\times V^* \tq [x,y] =0 \}$ and
$T^*_{C^*}V \simeq \Lambda_{C^*} = \{ (x,y) \in V\times C^* \tq [x,y] =0 \}$ for an $H$-orbit $C$ in $V$.\\

While $\mathfrak{g}(n)\times \mathfrak{g}(-n)$ carries a natural $H\times H$-action, $f=\langle \cdot , \cdot\rangle$ is only $H$-invariant for the diagonal $H$-action on $\mathfrak{g}(n)\times \mathfrak{g}(-n)=T^*V$.
This will allow us to use $H$-equivariant vanishing cycles in the context of Bernstein-Lunts equivariant derived categories \cite{BL94} (as before, since nearby and vanishing cycles commute with smooth pull-back). 
We now comment in more detail on the expectations of \cite{CFMMX22}*{Section 7.11}.\\

As explained in the section before, the functor $\ev{C}$ produces a parametrized family of vanishing cycles $j_x^* R\Phi_{f_y}[-1]$ for $(x,y)\in \Lambda_{C^*}$.
In the language of stratified Morse theory for constructible sheaves, as in \cite{Sch03}*{Chapter 5}, these are the local Morse data of $f_y$ \cite{MS22}*{(10.24), p. 710} (non-equivariantly, since stratified Morse theory is at the moment not available equivariantly): 
\begin{equation}\label{LMD-family}
\operatorname{For} \left((\ev{C}\mathcal{F})_{(x,y)}\right)
=
\operatorname{For} (\Phi_{f_y}[-1]\mathcal{F})_x )
=
\operatorname{LMD}(\operatorname{For}\mathcal{F},f_y,x)
\in D^b_{c}(V;k),
\end{equation}
for all $\mathcal{F}\in D^b_{c,H}(V;k)$ and $(x,y)\in \Lambda_{C^*}$.
This situation is very specific to the context of  a linear $H$-representation with finitely many orbits.
Here, the vanishing cycles functor $\Phi_{f_y}[-1]$ is of global nature, depending on the function $f_y$, for $y\in C^*$.
In microlocal geometry, a similar statement is only available for the so-called normal Morse data as in \eqref{equation: A 5}
\begin{equation}\label{NMD-family}
\operatorname{NMD}(\operatorname{For}\mathcal{F},(df)_x)
=
\left(R\Phi_{f\vert_{N}}[-1](\operatorname{For}\mathcal{F}\vert_{N})\right)_x,
\end{equation}
where the role of $f$ is played by $f_y$.
Here, $N$ is a local normal slice to the $H$-orbit $C$ at $x\in C$ and the power of microlocal geometry, here in the language of stratified Morse theory for constructible sheaves, is that this depends {\it only} on the non-degenerate covector $(df)_x \in \,_{C}\Lambda^\reg$ as in \cite{Bry86}*{Remarque 1.13, p. 17}, motivating the introduction of the functor $\Ev{}$, as explained in \cite{CFMMX22}*{Section 7.3}.
The comparison between these two versions, of local and global views, asks for a comparison between local and normal Morse data, respectively local vanishing cycles with or without slicing, for regular covectors 
$(df)_x \in \,_{C}\Lambda^\reg=\Lambda_{C^*}^\reg$.
This hints already at an important difference: the microlocal view restricts to regular covectors whereas the global vanishing cycles approach works more generally, in the sense that it can be studied outside of regular covectors, namely for $(df_y)_x \in \Lambda_{C^*}$. \\

Such a comparison is best possible in the more restricted context of this section
(with the $H^0$-action on $V= \mathfrak{g}(n)$ for $n \ne 0$ coming from a graded Lie algebra, fitting with Vogan varieties as they appear in this paper and in \cite{CFMMX22}).
Let's first explain the comparison of local and normal Morse data in the context of global vanishing cycles on $V$ via twisting by a rank-$1$ equivariant local system
\begin{equation}\label{twisting-loc}
\mathcal{T}_{C}
=
\evs{C}\IC(\1_{C}) \in \Loc_{H}(\Lambda_{C^*}^\reg;k),
\end{equation}
which works equivariantly and avoids the language of normal Morse data.
\begin{remark}\label{rem-rank-one}
Note that the property of   $\mathcal{T}_{C} = \evs{C}\IC(\1_{C}) $ being a rank-$1$ equivariant local system
is a local statement, not depending on the choice of the nondegenerate and $G$-invariant symmetric bilinear form $\langle \cdot , \cdot\rangle$ on $\mathfrak{g}$, which can be checked non-equivariantly. Similarly it also applies to an algebraic group $G'$ acting via an algebraic group homomorphism $G'\to H$ on $V=\mathfrak{g}(n)$ for a fixed $n\ne 0$.
\end{remark}

It comes from a local calculation of a description of $f_{C^*}$ along the smooth critical locus $\Lambda_{C^*}^\reg$, in \cite{CFMMX22}*{Section 7.7}, showing that $f_{C^*}\simeq \langle \cdot , \cdot\rangle: C\times C^*\to \CC $
is a holomorphic Morse-Bott function.
In other words, in suitable local coordinates in $C\times C^*$,  $f_{C^*}$ is given by a sum of squares in these variables $g \ceq \sum_{i=1}^e u_i z_i^2$ \cite{CFMMX22}*{Lemma 7.18}, for units $u_i$ in suitable local coordinates on an open $U\subseteq C\times C^*$, where $z_1,\ldots, z_e$ are coordinates in the normal direction, so $e =e_C = \dim C\times C^* - \dim \Lambda$, called the eccentricity of $C$. 
Because of this, the vanishing cycles
\[
R\Phi_{g}[-1] \1_{U}
= \mathcal{L}_{Z}[-e_C],
\]
where $Z = U \cap g^{-1}(0)$, 
with $\mathcal{L}_{Z}$ is a quadratic local system on $U$ with monodromy sign character coming from a $2$-group of order $2e_C$. 

\begin{remark}\label{local-sys}
 Let $\mathcal{L}$ be an $H$-equivariant local system on $C$. Then the equivariant counterpart of equation (\ref{diagram: tensor-phi}) implies:
\begin{equation}\label{twisting1}
 \evs{C}\IC(\mathcal{L}) =(\evs{C}\IC(\1_{C}))\otimes j^*_{\Lambda_{C^*}^\reg}(\mathcal{L}\boxtimes \1_{C^*})
 \in \Loc_{H}(\Lambda_{C^*}^\reg;k)
\end{equation}
 as in \cite{CFMMX22}*{Proposition 7.13(b)}.
 In particular
 \begin{equation}\label{twisting2}
 \mathcal{T}_{C}^\vee \otimes \evs{C}\IC(\mathcal{L}) = j^*_{\Lambda_{C^*}^\reg}(\mathcal{L}\boxtimes \1_{C^*})
 \in \Loc_{H}(\Lambda_{C^*}^\reg;k)
 \end{equation}
 as in \cite{CFMMX22}*{Theorem 7.22(d)}.
\end{remark}
We now  define the {\it normalized version}
\begin{equation}\label{equation: appendix NEvs}
\nevs{C} \ceq 
\mathcal{T}_{C}^\vee \otimes j^*_{\Lambda^\reg_{C^*}} R\Phi_{f_{C^*}}[-1]p_{C^*}^*[-\codim C^*] :  D^b_{c,H}(V;k) \to D^b_{\Lambda^\reg_{C^*}-c,H}(\Lambda^\reg_{C^*};k).
\end{equation}
Here the induced stratification of $\Lambda^\reg_{C^*}$ just consists of one smooth stratum so that the indicated 
$\Lambda^\reg_{C^*}$-constructibility just means that all underlying cohomology sheaves are locally constant (with finite dimensional stalks) on $\,_{C}\Lambda^\reg = \,\Lambda^\reg_{C^*}$.
This definition is designed so that
\[
\nevs{C} : \Per_{H}(V;k) \to \Loc_{H}(\Lambda^\reg_{C^*};k) = \Rep(A^\reg_{C^*})
\]
is an exact functor of abelian categories.
As before, $\Loc_{H}(\Lambda^\reg_{C^*};k)$ is the abelian category of $H$-equivariant local systems (of $k$-vector spaces) on $\Lambda^\reg_{C^*}$, and where we write $A^\reg_{C^*}$ for the equivariant fundamental group of $\Lambda^\reg_{C^*;k}$.
Note that the twisting by $\mathcal{T}_{C}^\vee \in \Loc_{H}(\Lambda^\reg_{C^*;k})$ only makes sense on $\Lambda^\reg_{C^*}$, but not on $\Lambda_{C^*}$.\\

For all applications in this paper, it is not the equivariant local system which comes into play but its restriction to a chosen $H$-orbit in $\Lambda_{C^*}^\reg$. For this reason, we introduce
\begin{equation}\label{nevs-mic}
\nevsmic{(x,y)} \ceq j^*_{\mathcal{O}_{H}(x,y)}\, \nevs{C}: \Per_{H;k}(V) \to \Loc_{H}(\mathcal{O}_{H}(x,y);k) = \Rep(A^\mic_{(x,y)})
\end{equation}
where $\mathcal{O}_H(x,y)$ is the $H$-orbit of $(x,y)\in \Lambda_{C^*}^\reg$ and where $j_{\mathcal{O}_{H}(x,y)}$ is the inclusion of this orbit in $\Lambda_{C^*}^\reg$.
Here, $A^\mic_{(x,y)}$ is the microlocal fundamental group of $\mathcal{O}_{H}(x,y)$ at the base point $(x,y)$, so $A^\mic_{(x,y)} = Z_{H}(x,y)/Z_{H}(x,y)^0$.
For every $(x,y)\in \Lambda_{C^*}^\reg$,  restriction along the inclusion $j_{\mathcal{O}_{H}(x,y)}$
defines 
\[
\Rep(A_{C^*}^\reg) = \Loc_{H}(\Lambda^\reg_{C^*};k) \to \Loc_{H}(\mathcal{O}_H(x,y);k) = \Rep(A^\mic_{(x,y)}),
\]
and thus defines $A^\mic_{(x,y)} \to A^\reg_{C^*}$.\\

Since $A_{(x,y)}^\mic$ can depend dramatically on $(x,y)\in \Lambda_{C^*}^\reg$ (except up to the $H$-action), this motivates the introduction of an non-empty open $H$-invariant subvariety $U\subseteq \Lambda_{C^*}^\reg$, as in \cite{ABV92}*{Lemma 24.3, p.254} such that $A^\mic_{(x,y)}$ is a local system of finite groups on $U$, in which case the image of 
$A^\mic_{(x,y)} \to \pi_1^{H}(U,(x,y))$ is constant, as in \cite{ABV92}*{24.6, p. 257}.
In the important special case, which applies to all local Arthur parameters, when there is an open $H$-orbit $\mathcal{O}_{H}(x,y) \subset\Lambda_{C^*}^\reg$, this is the preferred choice of $U$, above, and this coincides with the variety $\Lambda_{C^*}^\sreg \subseteq \Lambda_{C^*}^\reg$.

The compatibility with $\chi^\mic_{C}$ is given by the following commutative diagram, which does not depend on the choice of base point $(x,y)\in \Lambda_{C^*}^\reg$:

\begin{equation}
\begin{tikzcd}
K_0\Per_{H}(V;k)  
\arrow{d}[swap]{\chi^\mic_C} 
\arrow{rr}{\nevs{C}} 
&&
K_0\Loc_{H}(\Lambda_{C^*}^\reg;k) \arrow{d} \\
\Z &&  \arrow{ll}{\chi} K_0\Loc_{H}(\mathcal{O}_{H}(x,y);k) 
\end{tikzcd}
\end{equation}
We can  use this as a definition of $\chi^\mic_{C}$ via global vanishing cycles or, as we explain now, it fits with the old definition of $\chi^\mic_{C}$ given above.
In fact, this is only the numerical shadow of the comparison between $\nevs{C^*}$ and normal Morse data, as we now explain.

By construction, for a chosen base point $(x,y)\in U$,
\[
\operatorname{For}(j_{(x,y)}^*\mathcal{T}_{C}) 
= 
\operatorname{For}\left( j_{(x,y)}^* R\Phi_{f}[-1]\1_{U}[e_{C}] \right)
=
\operatorname{TMD}(\1_{C\times C^*},f,(x,y))[e_{C}],
\]
with the tangential Morse data as in \cite{Sch03}*{(5.3), p. 272} written there in terms of the real part $re(f)$ of $f$, fitting with the vanishing cycles description in the complex context by \cite{Sch03}*{Lemma 5.4.1, p. 362}:
\[
\operatorname{TMD}(\1_{C\times C^*},f,(x,y)) = 
\operatorname{LMD}(\1_{C\times C^*},f,(x,y))
= 
j_{(x,y)}^* \, R\Phi_{f}[-1]\1_{U}.
\]
Equivalently, by the base change of Theorem~\ref{theorem: family ev},
\[
\operatorname{TMD}(\1_{C},f_y,x) = 
\operatorname{LMD}(\1_{C},f_y,x)
= 
j_{x}^* \, R\Phi_{f_y}[-1]\1_{C}.
\]
So, tangential Morse data is just local Morse data on the stratum $C\times C^*$ for the constant sheaf $\1_{C\times C^*}$.
In this way we forget the equivariant structure as well as the monodromy action on the local system living on $\Lambda^\reg_{C^*}$.

Now we give important special cases of tangential Morse data $\operatorname{TMD}(\1_{S},g,z)$ on a given stratum $S$ in the point $z\in S$ for a function $g : S\to \CC$ given as $\sum_{i=1}^s z_i^2$ for $0\leq s\leq \dim S$:
\[
\operatorname{TMD}(\1_{S},g,z)= \1_{z}[-s].
\]
Especially, for $s=0$ we take $g=0$, which applies for example to the case $C$ or $C^*$ an open stratum with dual orbit $\{0\}$ and $s=e_C=0$.
In the second case we take $s= \dim S$, corresponding to a holomorphic Morse critical point of $g$, as in \cite{Sch03}*{p. 368} or, the intermediate case, $s = e_C$, for the case discussed above, as in \cite{CFK22}*{Theorem 7.19, p. 66}.\\

This gives light to the expectation of \cite{CFMMX22}*{(7.33), p. 72}, which is a special case of \cite{Sch03}*{Theorem 5.4.1, p. 367}.
Specializing to the context of this paper with stratum $C\times C^*\subset V\times C^*$, $f = f_{C^*}$ and $(x,y)\in \Lambda_{C^*}^\reg$, 
we get for all $\mathcal{F}\in D^b_{c,H}(V;k)$:
\[
\operatorname{LMD}( \operatorname{For}(\mathcal{F}\boxtimes \1_{C^*}),f,(x,y))
=
\operatorname{TMD}(\operatorname{For}(\1_{C\times C^*}),f,(x,y))
\otimes
\operatorname{NMD}(\operatorname{For}(\mathcal{F}\boxtimes \1_{C^*}),(df)_{(x,y)}),
\]
or similarly
\[
\operatorname{LMD}( \operatorname{For}(\mathcal{F}),f_y,x)
=
\operatorname{TMD}(\operatorname{For}(\1_{C}),f_y,x)
\otimes
\operatorname{NMD}(\operatorname{For}(\mathcal{F}),(df_y)_x).
\]
This is precisely the non-equivariant confirmation of the expectations 
\cite{CFMMX22}*{(7.34) and (7.35), p. 72}, with
\[
\operatorname{TMD}(\operatorname{For}(\1_{C}),f_y,x) = \1_{x}[-e_C].
\]
In particular, we get for all $\mathcal{F}\in D^b_{c,H}(V;k)$:
\begin{equation}\label{equation: NMD}
\begin{array}{rcl}
\operatorname{For}(j_{(x,y)}^*\, \nevs{C}\mathcal{F})
&=& 
\operatorname{NMD}(\operatorname{For}(\mathcal{F}\boxtimes\1_{C^*}),(df)_{(x,y)})
[-e_C-\codim C^*]\\
&=&
\operatorname{NMD}(\operatorname{For}(\mathcal{F}\boxtimes\1_{C^*}[\dim C^*]),(df)_{(x,y)})[-\dim (C\times C^*)]\\
&\simeq&
\operatorname{NMD}(\operatorname{For}(\mathcal{F}),(df_y)_x)[-\dim C]\\
&=& 
\left(R\Phi_{f_y\vert_{N}}[-1](\operatorname{For}\mathcal{F}\vert_{N})\right)_x
[-\dim C] \in D^b_{c}(\{x\};k),
\end{array}
\end{equation}
for any local normal slice $N$ at $x\in C$, with regular covector $(df_y)_x \in \Lambda_C^\reg$.
The isomorphism in the third line follows from the fact that normal Morse data $NMD$ commute by definition with the smooth pullback $p: C\times C^*\to C$, with $p^*(df_y)_x=(df)_{(x,y)}$ also regular.
On the right hand side the shift by $[-\dim C]$ fits with the sign $(-1)^{\dim C}$ showing up in the definition of characteristic cycles for perverse sheaves.
In our global context here,  there is even a preferred  Slodowy slice (restricted from $\mathfrak{g}$ to $V$),
which is affine linear and $Z_{H}(x)$-invariant (see \cite{CFMMX22}*{p. 72}).\\

The identification above shows that, after applying the forgetful functor (since stratified Morse theory is not developed equivariantly so far),  the twisting by $\mathcal{T}_C$ in the approach by global vanishing cycles gives the translation to the microlocal viewpoint with a dependence only on the regular differential $(df)_{(x,y)}\in \Lambda^\reg_{C^*}$. In particular, the functor
$\nevs{C}$ satisfies in the context of our Vogan varieties all the properties asked for the functor $Q^{\mic}$ in 
\cite{ABV92}*{Theorem 24.8, p. 259} (with $Q^{\mic}$ being the subject of the next section in this Appendix).\\

\section{\texorpdfstring{The functor $Q^\mic$ via equivariant microlocalization}{}}\label{section: Qmic}
We are coming back now to the more general context of
Subsection \ref{subsection: KDF} from the Appendix, and explain the construction of the functor $Q^\mic=Q^\mic_C$
from \cite{ABV92}*{Theorem 24.8} via the equivariant version of the sheaf theoretic 
microlocalization functor $\mu_C$ from \cite{KaS94}*{Section 4.3} (fitting with the corresponding D-module functors as in \cite{MF96}). We also point out the similarities and differences to the (normalized) vanishing cycle functor from Section \ref{section: global vanishing cycles} of the Appendix.
In particular, non-equivariantly also the stalk of 
$Q^\mic_C$ at a regular covector $df_x\in \,_{C}\Lambda^\reg$ fits with the corresponding 
shifted normal Morse datum $\operatorname{NMD}(\operatorname{For}(\mathcal{F}),df_x)[-\dim C]$
(as in equation (\ref{equation: NMD}) for the corresponding  normalized vanishing cycles $\nevs{C}$).\\

For simplicity we restrict again to the case of a complex algebraic group $G$, e.g. a (finite cover group of a) complex linear algebraic group, acting on a smooth algebraic variety $M$ with finitely-many orbits $C$, as this is the only case needed for our applications. These $G$-orbits give us a Whitney stratification $\mathscr{S}$ of $M$, with
 the {\it closure ordering} on these orbits defined by
\[ 
C_1 \geq C_2 \Longleftrightarrow \overline{C_1} \supseteq C_2.
\]
Recall that
$$\Lambda:=  \bigsqcup_{C} T^*_{C} M \subset T^*M$$ 
denotes the conormal variety 
with respect to the  $G$-orbits 
$C$ in $M$. Note that $\Lambda$ is closed in $T^*M$ by the Whitney a-regularity of the stratification $\mathscr{S}$.
For each $G$-orbit $C\subseteq V$, the
conormal bundle  $\,_C\Lambda \ceq T^*_{C}V$ of $C$ is  the inverse image of the projection $\Lambda \to V$ above $C$.
The regular part of the conormal to $C$ is given by
${}_C\Lambda^\reg \ceq \,_C\Lambda \setminus {}_C\Lambda^\textrm{sing}$ for 
\[
{}_C\Lambda^\textrm{sing} \ceq \,_C\Lambda \cap \left( \operatorname{\bigcup}\limits_{C_i\ne C} (\,_{C}\Lambda \cap \overline{\,_{C_i}\Lambda})\right).
\]
Note that  $\,_C\Lambda^\textrm{sing}$ is closed in $\,_C\Lambda$ and $\,_C\Lambda^\reg$ is open and dense in $\,_C\Lambda$.\\

We now recall the real and complex deformation to the normal bundle $T_CM$ for a given $G$-orbit $C$ in $M$,
with the real deformation living on a real-analytic manifold $\tilde{M}'$ used in \cite{KaS94}*{Section 4.2} for the definition of the corresponding real specialization functor 
$\nu_S$ of Kashiwara-Schapira, and the complex deformation living on a complex-algebraic (or complex-analytic) manifold $\tilde{M}$ used in \cite{Ver83}
by Verdier for the definition of his complex specialization functor $Sp_C$. For both definitions we have to assume that $C$ is closed in $M$,
which we can do after restriction from $M$ to the open subset $M\backslash \partial C$, with $\partial C= \overline{C}\backslash C$ the boundary of $C$ (given by a union of other $G$-orbits of smaller dimension.
Also note that the complex deformation to the normal cone of Verdier works more generally for a possible singular complex algebraic or analytic subset, with $\tilde{M}$ possible singular, but this is not needed here). All involved spaces can be best visualized by the following commutative diagram, whose right part is cartesian for the corresponding deformation parameters, and whose left part indicates the corresponding projections to $C$ or $M$ (and all vertical arrows are given by natural inclusions):
\begin{equation}\begin{CD}
    C @<<< T_CM @>>> \{0\} \\
    @VVV @V i VV @VVV \\
    M @< p' << \tilde{M}' @> t' >> \R \\
    @| @VVV  @VVV \\
    M @< p << \tilde{M} @> t >> \CC \\
    @| @AAA @AAA \\
    M @< \pr_1 << M\times \CC^* @> \pr_2 >> \CC^*\:.
\end{CD}\end{equation}
By naturality, all arrows are $G$-equivariant for the induced $G$-action on these spaces (with
$G$ acting trivially on $\CC$). Let $j: \{t'>0\} \hookrightarrow  \tilde{M}'$ be the corresponding open inclusion, and $k$ be any field of characteristic $0$. Then the real specialization functor of Kashiwara-Schapira is defined as
$$\nu_S= i^*Rj_*j^*p'^*: D^b_c(M;k)\to D^b_{\R^+-mon}(T_CM;k) \:,$$
whereas the complex specialization functor of Verdier is defined as
$$ Sp_S= R\Psi_t \,p^*: D^b_c(M;k)\to D^b_{c, mon}(T_CM;k) \:,$$
with $R\Psi_t$ the nearby cycle functor with respect to the deformation parameter $t$.
Here $mon$ (or $\R^+-mon$) stands for {\it monodromic} with respect to the $\CC^*$-action (or
$\R^+$-action) on the fibres of the complex normal bundle $T_CM\to C$ (i.e. all cohomology sheaves are locally constant on the orbits of these actions). Since $\R^+$ is contractible, the inclusion $\R^+\hookrightarrow \CC^*$ can be lifted up to a section
$s: \R^+\to \tilde{\CC}^*$ of the universal cover $\pi: \tilde{\CC}^*\to \CC^*$, which we fix from now on. Then one gets for a constructible sheaf complex 
$\mathcal{F}\in D^b_c(M;k)$ a natural comparison isomorphism
\begin{equation} \label{iso-sp-nu}
\begin{CD}
  Sp_C  \mathcal{F} @> \sim >> \nu_C \mathcal{F}
  \in D^b(T_CM;k) 
\end{CD}\end{equation}
(with the adjunction morphism constructed as in \cite{Sch03}*{p. 69}, which then is stalkwise an isomophism by \cite{Sch03}*{Lemma 1.1.1 on p. 27}, based on the local Milnor fibre description as in Proposition \ref{proposition: stalks},
together with its real counterpart for the functor $i^*Rj_*j^*$).
The {\it microlocalization} functor $\mu_C$ of \cite{KaS94}*{Section 4.3} is then defined as $\mu_C\mathcal{F} = \left(\nu_C\mathcal{F}\right)^\wedge$, with $(-)^\wedge$ the Fourier-Sato transformation
$$(-)^\wedge: D^b_{\R^+-mon}(T_CM;k)\to D^b_{\R^+-mon}(T^*_CM;k)\:.$$
Using the isomorphism \eqref{iso-sp-nu}
and the open inclusion $j_C: M\backslash \partial C \hookrightarrow M$,
we finally get our complex version of the  {\it microlocalization} functor $\mu_C:=(-)^\wedge \circ \Sp_C \circ  j_C^*$ defined as the composition
\begin{equation}\label{micloc}
\begin{CD}
  D^b_c(M;k) @> j_C^* >> D^b_c(M\backslash \partial C;k) @> Sp_C >> D^b_{c, mon}(T_CM;k)
 @> (-)^\wedge >>   D^b_{c, mon}(T^*_CM;k) \:. 
\end{CD}\end{equation}

Note that $\nu_C$ and $\mu_C[\codim C]$ are t-exact for the perverse t-structure (\cite{KaS94}*{Proposition 10.3.19} or
\cite{Sch03}*{Remark 6.0.8 on p. 411}), fitting for $k=\CC$ with corresponding D-module functors in terms of the solution functor $Sol(-)=Rhom_\mathcal{D}((-)^{an}, \mathcal{O})$ as in 
\cite{MF96}*{Theorem 1.1.2}.
Here $(-)^{an}$ denotes the corresponding analytic left D-module
(with $\mathcal{O}$ and $\mathcal{D}$ the sheaf of holomorphic functions and differential-operators).
\\

These nearby and vanishing cycle functors, as well as the specialization, Fourier-Sato 
and microlocalization functors $\nu_S, Sp_C, (-)^\wedge$ and $\mu_C$ can be lifted up to the equivariant derived context
of Bernstein-Lunts \cite{BL94}, commuting with the restriction and change of group structures, like the forgetful functor $\operatorname{For}: D^b_{c,G}(M;k) \to D^b_{c}(M;k)$.
This is the case
because they commute with smooth based change
(see, e.g. \cite{Sch03}*{Remark 4.3.7 on p. 265} for $R\Psi_t, R\Phi_t$, \cite{KaS94}*{Proposition 4.2.5 and Proposition 4.3.5} for $\nu_C, \mu_C$, and \cite{KaS94}*{Proposition 3.7.13} for $(-)^\wedge$ commuting with any pullback)
and we only use invariant functions (like the deformation parameter) with respect to the group actions. 
Similarly one gets for an equivariant constructible sheaf complex 
$\mathcal{F}\in D^b_{c,G}(M;k)$ a natural comparison isomorphism
\begin{equation} \label{G-iso-sp-nu}
\begin{CD}
  Sp_C  \mathcal{F} @> \sim >> \nu_C \mathcal{F}
  \in D^b_G(T_CM;k) \:,
\end{CD}\end{equation}
coming from the fact that the non-equivariant version \eqref{iso-sp-nu} commutes with smooth based change
(and gives this isomorphism after applying $\operatorname{For}$). Altogether we get in this way from 
\eqref{micloc} the
{\it equivariant microlocalization functor} $\mu_C$ given as the composition
\begin{equation}\label{G-micloc}
\begin{CD}
  D^b_{c,G}(M;k) @> j_C^* >> D^b_{c,G}(M\backslash \partial C;k) @> Sp_C >> D^b_{c, mon, G}(T_CM;k)
 @> (-)^\wedge >>   D^b_{c, mon, G}(T^*_CM;k) \:, 
\end{CD}\end{equation}
with the indices $c$ and $mon$ referring to the corresponding properties of the underlying non-equivariant sheaf complexes.
Again $\mu_C[\codim C]$ is t-exact for the perverse t-structure:
\begin{equation}\label{G-micloc-perv}
\mu_C[\codim C]: \Per_{G}(M;k)\to \Per_{mon,G}(T^*_CM;k)\:.
\end{equation}

This equivariant microlocalization functor $\mu_C$ should be seen as a micro-local counterpart of the vanishing cycle functor
$\ev{C}$ studied in Corollary \ref{corollary: conormal support} of this Appendix. But the target of the vanishing cycle functor
$\ev{C}$ are equivariant constructible sheaf complexes on $\Lambda \cap (V\times C^*) =\Lambda_{C^*}$ (the conormal bundle to the dual orbit $C^*$
in $V^*$),
whereas the target of our equivariant microlocalization functor $\mu_C$ are equivariant constructible sheaf complexes on the conormal bundle $T^*_CM=\,_C\Lambda$.
Similarly to Theorem \ref{theorem: refinement} of this Appendix,
we get the derived version of our functor $Q^\mic_C$
by composing with $j^*_{\,_C\Lambda^\reg}[-\dim C]$ for the open inclusion of the regular part 
$$j^*_{\,_C\Lambda^\reg}: \,_C\Lambda^\reg=(T^*_CM)^\reg \hookrightarrow T^*_CM= \,_C\Lambda \:.$$
We now  define 
\begin{equation}\label{equation: derived Qmic}
Q^\mic_C := j^*_{\,_C\Lambda^\reg}\,\mu_C[-\dim C] :  
D^b_{c,G}(M;k) \to 
D^b_{\,_C\Lambda^\reg-c,G}(\,_C\Lambda^\reg;k).
\end{equation}

Here the used stratification of $\,_C\Lambda^\reg$  consists by definition just of one smooth stratum so that the indicated 
$\,_C\Lambda^\reg$-constructibility just means that all underlying cohomology sheaves are locally constant (with finite dimensional stalks) on $\,_C\Lambda^\reg=(T^*_CM)^\reg$ (as will be explained in a moment).
This definition is designed so that
\begin{equation}\label{equation: def-Qmic}
Q^\mic_C  : \Per_{G}(M;k) \to \Loc_{G}(\,_C\Lambda^\reg;k) = \Rep(A^\reg_{C})
\end{equation}
is an exact functor of abelian categories
(as asked for in \cite{ABV92}*{Theorem 24.8(c) on p. 259}). As before, $\Loc_{G}(\,_C\Lambda^\reg;k)$ is the abelian category of $G$-equivariant local systems (of $k$-vector spaces) on $\,_C\Lambda^\reg$, and where we write $A^\reg_{C}$ for the equivariant fundamental group of $\,_C\Lambda^\reg$.\\

Let us now explain why all cohomology sheaves of $j^*_{\,_C\Lambda^\reg}\,\mu_C\operatorname{For}(\mathcal{F})$ are locally constant on $\,_C\Lambda^\reg$ for a given
 $\mathcal{F}\in D^b_{c,G}(M;k)$ (see also \cite{Sch03}*{(5.78) on p. 359}).
 This follows from a suitable estimate of 
 $$\mu\supp(\operatorname{For}(\mathcal{F}))\subset T^*M\:,$$
with $\mu supp$ the {\it micro-support} in the sense of \cite{KaS94}*{Definition 5.1.2}.
Here, $\mu\supp(\mathcal{G})\subset T^*M$ for $\mathcal{G}\in D^b(M,k)$ is by definition the complement of the largest open subset $U\subset T^*M$ such that each stalk $(R\Gamma_{\{g\geq g(x)}\mathcal{G})_x\simeq 0$
 for any differentiable function germ $g : (M,x) \to (\R,g(x))$ 
 with $dg_x \in U$. For example $\mu\supp(\mathcal{G})\subset T^*_MM$ is contained in the zero-section $T^*_MM$ of $T^*M$, if and only if all
 cohomology sheaves of $\mathcal{G}$ are locally constant on $M$ (see \cite{Sch03}*{Proposition 4.1.1 on p.217}).
 Since $\operatorname{For}(\mathcal{F})$ is constructible with respect to the Whitney stratification 
$\mathscr{S}$ of $M$ given by the $G$-orbits $C$, one has by \cite{Sch03}*{Corollary 4.0.2 on p.212} the estimate
\begin{equation}\label{musupp1}
   \mu\supp(\operatorname{For}(\mathcal{F}))\subset \Lambda=  \bigsqcup_{C} T^*_{C}M \subset  T^*M\:. 
\end{equation}
In fact this estimate is even equivalent to the  property that $\operatorname{For}(\mathcal{F})$ is $\mathscr{S}$-constructible. Similarly \cite{Sch03}*{Proposition 5.0.1 on p.273} or \cite{MS22}*{Proposition 3.13} shows that in this case 
$\mu\supp(\operatorname{For}(\mathcal{F}))$ is the union over all conormal closures $\overline{T^*_CM}=\overline{\,_C\Lambda}$ with
non-trivial normal Morse data $\operatorname{NMD}(\operatorname{For}(\mathcal{F}),df_x)\neq 0$ for $df_x\in \,_C\Lambda^\reg$. In particular
\begin{equation}\label{CC-essimate}
 \supp(CC(\operatorname{For}(\mathcal{F})) \begin{cases} \subset \mu\supp(\operatorname{For}(\mathcal{F}))
  \quad \text{for $\mathcal{F}\in D^b_{c,G}(M;k)$,}\\
  = \mu\supp(\operatorname{For}(\mathcal{F}))
  \quad \text{for $\mathcal{F}\in \Per_{G}(M;k)$.}
\end{cases}\end{equation}

After restriction to the open subset $M\backslash \partial C \hookrightarrow M$ for a fixed $G$-orbit $C$, we can assume that $C\subset M$ is closed.
Then \cite{KaS94}*{Theorem 6.4.1} gives for  $\mathcal{G}\in D^b(M,k)$ the estimate
$$\mu \supp(\mu_C\, \mathcal{G}) \subset C_{T^*_CM}(\mu \supp(\mathcal{G}))\subset T_{T^*_CM}(T^*M) \simeq T^*(T^*_CM)\:, $$
with $C_{T^*_CM}(-)$ the corresponding normal cone (in the sense of \cite{KaS94}*{Definition 4.1.1}) along the closed submanifold $T^*_CM$ of $T^*M$.
Here the isomorphism $T_{T^*_CM}(T^*M) \simeq T^*(T^*_CM)$ from \cite{KaS94}*{(6.2.2) on p. 259} is used, which also identifies  the corresponding zero-sections. Since $\,_C\Lambda^\reg$ is open in $T^*_CM=\,_C\Lambda$ 
as well as in $\Lambda$, we can conclude from the estimate
\eqref{musupp1} that $\mu \supp(j^*_{\,_C\Lambda^\reg}\,\mu_C\,\operatorname{For}(\mathcal{F}) )$ is contained in the zero-section of 
$T^*(\,_C\Lambda^\reg)$ for a given
 $\mathcal{F}\in D^b_{c,G}(M;k)$.
But as already remarked before, this implies (by \cite{Sch03}*{Proposition 4.1.1 on p.217}) that
all cohomology sheaves of $j^*_{\,_C\Lambda^\reg}\,\mu_C\operatorname{For}(\mathcal{F})$ are locally constant on $\,_C\Lambda^\reg$.

\begin{remark}\label{rem:lc}
 Note that  \eqref{CC-essimate} holds more generally for an algebraic group $G'$ acting on a complex algebraic manifold $M'$, without assuming that it only acts with finitely many orbits. Then one gets for  $\mathcal{F}\in \Per_{G'}(M';k)$ the implication: assume the stalk Euler characteristic 
 $\chi(\operatorname{For}(\mathcal{F})_{x'}) $
 is a constant function independent of $x'\in M'$. Then
 $$\mu\supp(\operatorname{For}(\mathcal{F})) =
 \supp(CC(\operatorname{For}(\mathcal{F})) \subset T^*_{M'}M'$$
 is contained in the zero-section of $T^*M'$ so that
 $\mathcal{F}[-\dim M']$  is an $G'$-equivariant local system on $M'$.
\end{remark}

Coming back to the context of this section, using \cite{Sch03}*{(5.52) on p. 311 and (5.76) on p. 356}, we get the identification of the non-equivariant stalk of the functor 
$Q^\mic_C$ with the corresponding shifted {\it normal Morse data}:
\begin{equation}\label{Qmic=NMD}
\operatorname{For}(Q^\mic_C (\mathcal{F})_{df_x})\simeq   \operatorname{NMD}(\operatorname{For}(\mathcal{F}),df_x)[-\dim C]\in D^b_{c}(\{x\};k)
\end{equation}
for $df_x\in \,_C\Lambda^\reg$ and $\mathcal{F}\in D^b_{c,G}(M;k)$.

\begin{remark} \label{micro-solution}
Assume $\mathcal{F}=Sol(\mathcal{M})=Rhom_{\mathcal{D}_M}(\mathcal{M}^{an},\mathcal{O}_M)\in D^b_{c,G}(M;\CC)$ 
is the {\it solution complex} of a $G$-equivariant coherent D-module $\mathcal{M}\in \operatorname{Mod}_{c}(D_M,G)$,
with $(-)^{an}$ the corresponding analytic left D-module (and $\mathcal{O}_M$ and $\mathcal{D}_M$ the sheaf of holomorphic functions and differential-operators on $M$).
Then by \cite{KaS94}*{Theorem 11.3.3} the micro-support
$$\mu \supp(\operatorname{For}\mathcal{F}) = \mathrm{char}(\operatorname{For}\mathcal{M}^{an})$$
is just the {\it characteristic variety} of $\operatorname{For}\mathcal{M}^{an}$, and
the stalk $\,\operatorname{For}(\mu_C (\mathcal{F})_{df_x})$ for $df_x\in \,_C\Lambda^\reg$
corresponds to the {\it micro-solution} of $\operatorname{For}\mathcal{M}^{an}$ in $df_x$ (see \cite{MF95}*{Corollary 2.4 and Remark 4}), with
$$\mu_C(\operatorname{For}(\mathcal{F})) \simeq Rhom_{\pi^{-1}\mathcal{D}_M}(\pi^{-1}\operatorname{For}\mathcal{M}^{an},\mu_C(\mathcal{O}_M))\:.$$
Here $\pi: T^*M\to M$ is the projection. Moreover, the complex $\mu_C(\mathcal{O}_M)$
is concentrated in degree $c=\codim C$, with
$$\mathscr{C}^\R_{C|M}=\mathcal{H}^c(\mu_C(\mathcal{O}_M))$$
the sheaf of holomorphic microfunctions along $C$ (compare  \cite{KaS94}*{Proposition 11.4.1 and
Definition 11.4.2}).
By \cite{Kas83}*{Theorem 3.2.42 on p. 93} one has for $df_x\in \,_C\Lambda^\reg$:
\begin{equation*}\begin{split}\mathcal{E}xt^i_{\pi^{-1}\mathcal{D}_M}(\pi^{-1}\operatorname{For}\mathcal{M}^{an},\mathscr{C}^\R_{C|M})_{df_x}
&\simeq \mathcal{E}xt^i_{\mathcal{E}_M}(\mathcal{E}_M \otimes_{\pi^{-1}\mathcal{D}_M} \pi^{-1}\operatorname{For}\mathcal{M}^{an},\mathscr{C}^\R_{C|M})_{df_x}\\
&=\begin{cases}
0 \quad &i\neq 0,\\
\CC^m \quad &i=0.
\end{cases}\end{split}\end{equation*}
Here $m$ is the {\it multiplicity} of $\operatorname{For}\mathcal{M}$ along  $T^*_CM$, with $\mathcal{E}_M$ the ring of finite-order microdifferential operators (which is flat over $\pi^{-1}\mathcal{D}_M$). Together with \eqref{Qmic=NMD} this implies directly the equality of the different notions of characteristic cycles for the equivariant D-module $\mathcal{M}\in \operatorname{Mod}_{c}(D_M,G)$ and the equivariant perverse sheaf
$$\mathcal{F}[\dim M]=Sol(\mathcal{M})[\dim M]\in \Per_{G}(M;\CC)$$ 
as stated in the discussion after the commutative diagram \eqref{equation: players} in Subsection \ref{subsection: KDF}  of this Appendix.
\end{remark}

In particular our functor $Q^\mic_C$ satisfies the  property from \cite{ABV92}*{Theorem 24.8(a) on p. 259}, that the rank (or Euler characteristik) of $\operatorname{For}(Q^\mic_C (\mathcal{F})_{df_x})$ for $df_x\in  \,_C\Lambda^\reg$ and an equivariant perverse sheaf 
$\mathcal{F}\in \Per_{G}(M;k)$
(or equivariant constructible sheaf complex $\mathcal{F}\in D^b_{c,G}(M;k)$) is equal to the corresponding multiplicity $\chi^\mic_C(\mathcal{F})$
of $\mathcal{F}$ in its characteristic cycle.

\begin{remark} \label{rem-comparison}
Let us compare \eqref{Qmic=NMD} with \eqref{equation: NMD} in the more restricted context of Subsection \ref{subsection: normalization}
from this Appendix, with $G=H$ acting linearly on a (finite dimensional) complex vector space $M=V$ with finitely many orbits $C$, whose underlying representation of $H^0$ is coming  from a $\Z$-graded Lie algebra, e.g. as in the context of a {\it Vogan variety} as in the main part of this paper (or in \cite{Vin76}*{} and
\cite{Lus95}*{Section 3}).
Then we get for 
$$\mathcal{F}\in D^b_{c,H}(V;k) \quad \text{and} \quad (x,y)\in \,_C\Lambda^\reg = \Lambda^\reg_{C^*}\subset T^*V=V\times V^*$$
the non-equivariant identification between the stalks of the {\it normalized vanishing cycles} $ \nevs{C}$ and of the {\it micro-localization functor} $Q^\mic_C$, which are both isomorphic to the corresponding shifted {\it normal Morse data}:
\begin{equation} \label{stalk-news-qmic}
  \operatorname{For}(j_{(x,y)}^*\, \nevs{C}\mathcal{F}) \simeq   \operatorname{NMD}(\operatorname{For}(\mathcal{F}),(df_y)_x)[-\dim C]
  \simeq \operatorname{For}(Q^\mic_C (\mathcal{F})_{(df_y)_x})\:.
\end{equation} 
Note that all cohomology sheaves of $\nevs{C}\mathcal{F}$ and also of $Q^\mic_C (\mathcal{F})$ are $H$-equivariant local systems on $\,_C\Lambda^\reg = \Lambda^\reg_{C^*}$.
In future work \eqref{stalk-news-qmic} will be extended to the {\it equivariant} context for the {\it stalk along an $H$-orbit} in $\,_C\Lambda^\reg = \Lambda^\reg_{C^*}$,
once the basic results of stratified Morse theory for constructible sheaves are extended to such an equivariant context. In particular, this equivariant version of 
\eqref{stalk-news-qmic} will imply for $(x,y)\in \,_C\Lambda^\reg = \Lambda^\reg_{C^*}$ fixed by a semi-simple element $s\in H$ the following important equality of micro-local traces:
\begin{equation}\label{mic-traces}
  \trace({s},j_{(x,y)}^*\, \nevs{C}\mathcal{F})
= \trace({s},Q^\mic_C (\mathcal{F})_{(df_y)_x}) \:.  
\end{equation}
Similarly,  the trace formula appearing in \cite{ABV92}*{Theorem 25.8} will also be a special case of 
Theorem~\ref{theorem: Schurmann} from this Appendix.
\end{remark}

Next we show 
that our functor $Q^\mic_C$ satisfies the normalization property from \cite{ABV92}*{Theorem 24.8(b) on p. 259},
which is the counterpart for $Q^\mic_C$ of the equation \eqref{twisting2} from Remark \ref{local-sys} in this Appendix for the normalized vanishing cycle functor $\nevs{C}$. Fix a $G$-orbit $C$ in $M$ and an equivariant perverse sheaf $\mathcal{F}\in \Per_{G}(M)$, which is supported on the closure
$\overline{C}$ of $C$, with $\mathcal{F}|C \simeq \mathcal{L}[\dim C]$ for an $H$-equivariant local system $\mathcal{L}$ on $C$. After restriction to the open subset $M\backslash \partial C \hookrightarrow M$, we can assume that $i: C\hookrightarrow M$ is closed 
with $\mathcal{F}\simeq i_*\mathcal{L}[\dim C]$. Then the specialization $Sp_{id}\mathcal{L}$ of $\mathcal{L}$ with respect to the identity inclusion
$id: C\to C$ is given by $\mathcal{L} \simeq Sp_{id}\mathcal{L}$ under the zero-section identification $C\simeq T_CC$ (see \cite{Ver83}*{Property (SP5) on p. 353}. The counterpart for the real specialization functor $\nu_{id}$ follows from \cite{KaS94}*{Theorem 4.2.3(iv)}). 
Moreover the Verdier (and also the real) specialization functor is functorial for the closed inclusion $i: C\hookrightarrow M$
 (see \cite{Ver83}*{Property (SP2) on p. 352} and \cite{KaS94}*{Proposition 4.2.4}), so that 
$$Sp_C\mathcal{F}\simeq Ti_*(\mathcal{L}[\dim C]) \simeq s_*\mathcal{L}[\dim C] = s_!\mathcal{L}[\dim C]\:,$$
with $s: C\simeq T_CC\to T_CM$ the zero-section embedding identified with the map $Ti: T_CC\to T_CM$ of the normal bundles induced by the closed embedding $i: C\hookrightarrow M$. The transpose map of $s$ is the projection map $\,^ts=\pi: T^*_CM\to C$  of the conormal bundle,
so that \cite{KaS94}*{Proposition 3.7.14(i)} implies 
$$\mu_C \mathcal{F}  \simeq ( s_!\mathcal{L}[\dim C])^\wedge  \simeq \pi^*\mathcal{L}[\dim C] \in D^b_{c,mon,G}(T^*_CM;k)\:.$$
In terms of our functor $Q^\mic_C$ this translates to the {\it normalization} property of $Q^\mic_C$:
\begin{equation} \label{Qmic-normalization}
 Q^\mic_C  \mathcal{F} = j^*_{\,_C\Lambda^\reg}\,\mu_C  \mathcal{F} [-\dim C] \simeq j^*_{\,_C\Lambda^\reg}\pi^*\mathcal{L}\in \Loc_{G}(\,_C\Lambda^\reg;k) \:.
\end{equation}\\

Let us finish this subsection with the following {\it multiplicativity property} of the micro-localization functors $\mu_C$ and $Q^\mic_C$
(similarly to \eqref{TS-theorem} and \eqref{TS-theorem} for 
the micro-local vanishing cycle functors  $\ev{C}$ and $\evs{C}$ in Subsection \ref{subsection: normalization} of this Appendix).
Assume $G_i$  is a complex  algebraic group 
acting linearly on a complex algebraic manifold $M_i$ with finitely many orbits $C_i$ for $i=1,2$.
Then the same is true for $G=G_1\times G_2$ acting on
$M=M_1 \times M_2$ with finitely many orbits $C=C_1\times C_2$. Moreover
$$T^*M\simeq TM^*_1\times T^*M_2 \:,$$
with $\,_C\Lambda \ceq T^*_{C}M\simeq \,_{C_1}\Lambda\times \,_{C_2}\Lambda$ and
$\,_C\Lambda^\reg \ceq (T^*_{C}M)^\reg \simeq \,_{C_1}\Lambda^\reg\times \,_{C_2}\Lambda^\reg$.

By the (equivariant version of) \cite{KaS94}*{Proposition 4.2.6} one has  for $\mathcal{F}_i\in D^b_{c,G_i}(M_i;k)$ for $i=1,2$
a natural (iso)morphism
$$Sp_{C}( j_C^*(\mathcal{F}_1\boxtimes \mathcal{F}_2))
  \simeq Sp_{C_1}( j_{C_1}^* \mathcal{F}_1) \boxtimes
  Sp_{C_2}( j_{C_2}^* \mathcal{F}_2) \in  D^b_{c,G}(T_CM;k) \:.$$
Here it follows from
\cite{Sch03}*{Theorem 1.0.4 on p.23},
that the general multiplicativity morphism from  \cite{KaS94}*{Proposition 4.2.6} is in our case an isomorphism.
Then the multiplicativity of the Fourier-Sato transformation $(-)^{\wedge}$ (as in \cite{KaS94}*{Proposition 3.7.15},
together with the pullback functoriality of $(-)^{\wedge}$ as in \cite{KaS94}*{Proposition 3.7.13(i)}) implies
the {\it multiplicativity} of the micro-localization functor $\mu_C$:
\begin{equation} \label{micro-multiplicativity}
  \mu_{C}(\mathcal{F}_1\boxtimes \mathcal{F}_2)
  \simeq \mu_{C_1}(\mathcal{F}_1) \boxtimes
  \mu_{C_2}(\mathcal{F}_2) \in  D^b_{c,G}(T^*_CM;k) \:,
\end{equation}
and similarly for the functor  $Q^\mic_C$:
\begin{equation} \label{micro-multiplicativity2}
  Q^\mic_{C}(\mathcal{F}_1\boxtimes \mathcal{F}_2)
  \simeq Q^\mic_{C_1}(\mathcal{F}_1) \boxtimes
  Q^\mic_{C_2}(\mathcal{F}_2) \in  D^b_{\,_C\Lambda^\reg-c,G}(\,_C\Lambda^\reg;k) \:.
\end{equation}
The (non-equivariant) counterpart in the context of D-modules is discussed in \cite{MF96}*{Proposition 1.2.5}.

\begin{remark}
All the results of this last Section \ref{section: Qmic} of this Appendix also hold in the complex analytic context for a complex Lie group $G$ with finitely many connected components,
so that the equivariant derived categories of 
Bernstein-Lunts \cite{BL94}*{Part I, Section 3}
and their functors are available. 
For the Riemann-Hilbert correspondence one uses then of course analytic D-modules.
In this complex analytic context we explicitly assume that $G$ acts holomorphically on a complex manifold $M$ with finitely many locally closed orbits $C$, whose closure $\overline{C}$ is complex analytic with boundary
$\partial{C}= \overline{C}\backslash C$ a union of lower dimensional orbits.
\end{remark}

\fancypagestyle{referencesstyle}{
    \fancyhf{}
    \fancyhead[C]{\scriptsize REFERENCES}
    \fancyhead[R]{\scriptsize\thepage}
    \renewcommand{\headrulewidth}{0pt}
}
\pagestyle{referencesstyle}
\thispagestyle{referencesstyle}
\bibliographystyle{alpha}
\bibliography{ArXivreferences-v1}

\fancypagestyle{indexstyle}{
    \fancyhf{}
    \fancyhead[C]{\scriptsize INDEX}
    \fancyhead[R]{\scriptsize\thepage}
    \renewcommand{\headrulewidth}{0pt}
}
\fancypagestyle{indexfirststyle}{
    \fancyhf{}
    \fancyhead[R]{\scriptsize\thepage}
    \renewcommand{\headrulewidth}{0pt}
}
\clearpage
\pagestyle{indexstyle}
\begingroup
\makeatletter
\let\ps@plain\ps@indexfirststyle
\makeatother
\printindex
\clearpage
\endgroup

\fancypagestyle{authorstyle}{
    \fancyhf{}
    \fancyhead[C]{\scriptsize AUTHOR INFORMATION}
    \fancyhead[R]{\scriptsize\thepage}
    \renewcommand{\headrulewidth}{0pt}
}
\pagestyle{authorstyle}
\thispagestyle{authorstyle}

\end{document}